\documentclass[reqno,psamsfonts]{memo-l}
\usepackage{amscd,amssymb,eucal,latexsym,upref}

\def\vv<#1>{\langle#1\rangle}

\newcommand{\cm}{com\-mu\-ta\-tive mo\-no\-id}
\newcommand{\poag}{partially ordered abelian group}

\newcommand{\crm}{co\-ni\-cal re\-fi\-ne\-ment mo\-no\-id}
\newcommand{\wpr}{\to_{\mathrm{w}}}
\newcommand{\wpru}{\nearrow_{\mathrm{w}}}
\newcommand{\wprd}{\searrow_{\mathrm{w}}}

\newcommand{\Vemb}{V-em\-bed\-ding}
\newcommand{\Vhom}{V-ho\-mo\-mor\-phism}
\newcommand{\Vmod}{V-mod\-u\-lar}
\newcommand{\Vmeas}{V-meas\-ure}
\newcommand{\mcont}{meet-con\-tin\-u\-ous}
\newcommand{\jcont}{join-con\-tin\-u\-ous}
\newcommand{\res}{\mathbin{\restriction}}

\newcommand{\zero}{\ensuremath{\mathbf{0}}}
\newcommand{\one}{\ensuremath{\mathbf{1}}}
\newcommand{\two}{\ensuremath{\mathbf{2}}}
\def\dnw{\mathop{\downarrow}}

\def\rps{refined partial semigroup}

\def\S{\mathbf{S}}
\def\C{\mathbf{C}}
\def\sd{\smallsetminus}
\def\lo{\left]}
\def\ro{\right[}
\def\open#1{\lo#1\ro}

\newcommand{\tr}{\vartriangleleft}
\newcommand{\ntr}{\ntriangleleft}
\newcommand{\utr}{\trianglelefteq}
\newcommand{\nutr}{\ntrianglelefteq}
\newcommand{\dtr}{\mathbin{\vartriangleleft\kern-10pt
{\lower 3pt\hbox{$\scriptscriptstyle\ne$}}\kern3pt}}

\newcommand{\Dim}{{\Delta}}

\DeclareMathOperator{\nor}{Nor}
\DeclareMathOperator{\con}{Con}
\DeclareMathOperator{\ccon}{Con_c}
\DeclareMathOperator{\Id}{Id}
\DeclareMathOperator{\Idc}{Id_c}
\DeclareMathOperator{\FP}{FP}
\DeclareMathOperator{\rect}{Rect}
\DeclareMathOperator{\NEq}{NEq}
\DeclareMathOperator{\diag}{diag}
\DeclareMathOperator{\FL}{FL}

\DeclareMathOperator{\DD}{Dim}
\DeclareMathOperator{\MM}{M}

\def\tvi{\vrule height 12pt depth 6pt width 0pt}

\newcommand{\diam}[2]{\diamondsuit_{#1}({#2})}
\newcommand{\NN}{\mathbb{N}}
\newcommand{\PP}{\mathbb{P}}

\newcommand{\RR}{\mathbb{R}}
\newcommand{\ZZ}{\mathbb{Z}}
\newcommand{\ZZb}{{\overline{\ZZ}^+}}
\newcommand{\Ft}{\widetilde{\mathbf{F}}}
\newcommand{\ee}{\varepsilon_{QP}}
\newcommand{\ff}{\eta_{QP}}
\newcommand{\fft}{\tilde{\eta}_{QP}}
\newcommand{\ind}{\mathrm{Ind}}
\newcommand{\fin}{\mathrm{Fin}}

\newcommand{\ii}[1]{\emph{#1}}

\numberwithin{equation}{chapter}
\numberwithin{figure}{chapter}

\theoremstyle{plain}
\newtheorem{theorem}{Theorem}[chapter]
\newtheorem{proposition}[theorem]{Proposition}
\newtheorem{lemma}[theorem]{Lemma}
\newtheorem{corollary}[theorem]{Corollary}
\newtheorem{examplepf}[theorem]{Example}
\newtheorem{claim}{Claim}

\theoremstyle{definition}
\newtheorem{definition}[theorem]{Definition}
\newtheorem{example}[theorem]{Example}
\newtheorem{notation}[theorem]{Notation}
\newtheorem{problem}{Problem}

\theoremstyle{remark}
\newtheorem*{remark}{Remark}
\newtheorem*{note}{Note}

\newcommand{\qedc}{{\qed}~{\rm Claim~{\theclaim}.}}

\newenvironment{cproof}
{\begin{proof}[Proof of Claim.]}
{\qedc\renewcommand{\qed}{}\end{proof}}

\makeindex

\begin{document}
\frontmatter
\title[Dimension monoid]{The Dimension Monoid\\
                         of a Lattice}
\author{Friedrich Wehrung}
\address{C.N.R.S.\\
         Universit\'e de Caen, Campus II\\
         D\'epartement de Math\'ematiques\\
         B.P. 5186\\
         14032 CAEN cedex\\
         FRANCE}

\email{wehrung@math.unicaen.fr}
\urladdr{http://www.math.unicaen.fr/\~{}wehrung}

\begin{abstract}
We introduce the \emph{dimension monoid} of a lattice $L$,
denoted by \(\DD L\).  The monoid
\(\DD L\) is commutative and conical, the latter meaning that the
sum of any two nonzero elements is nonzero. Furthermore, \(\DD L\) is
given along with the \emph{dimension map}, $\Dim$, from \(L\times L\)
to
\(\DD L\), which has the intuitive meaning of a \emph{distance}
function.
The maximal semilattice quotient of $\DD L$ is
isomorphic to the semilattice $\ccon L$ of compact congruences of
$L$; hence \(\DD L\) is a precursor of the congruence lattice of
$L$. Here are some additional features of this construction:

\begin{enumerate}
\item Our dimension theory provides a
generalization to \emph{all} lattices of the von Neumann dimension
theory of continuous geometries. In particular, if $L$ is an
irreducible continuous geometry, then \(\DD L\) is either
isomorphic to \(\ZZ^+\) or to \(\RR^+\).

\item If $L$ has no infinite bounded chains, then \(\DD L\) embeds
(as an ordered monoid) into a power of \(\ZZ^+\cup\{\infty\}\).

\item If $L$ is modular or if $L$ has no infinite bounded chains,
then \(\DD L\) is a refinement monoid.

\item If $L$ is a simple geometric lattice, then \(\DD L\) is
isomorphic to \(\ZZ^+\), if $L$ is modular, and to the
two-element semilattice, otherwise.

\item If $L$ is an $\aleph_0$-\mcont\ complemented modular lattice,
then both \(\DD L\) and the dimension function $\Dim$ satisfy
(countable) completeness properties.
\end{enumerate}

If $R$ is a von Neumann regular ring and if
$L$ is the lattice of principal right ideals of the matrix ring
\(M_2(R)\), then \(\DD L\) is isomorphic to the
monoid \(V(R)\) of isomorphism classes of finitely generated
projective right $R$-modules. Hence the dimension theory of
lattices provides a wide lattice-theoretical generalization of
nonstable K-theory of regular rings.

\end{abstract}

\maketitle

\tableofcontents

\chapter*{Introduction}

In this work, we associate with any lattice its
\emph{dimension monoid}\index{dimension monoid!of a lattice},
\(\DD L\).
\footnote{The corresponding concept for \emph{rings}
originates in B. Blackadar's paper \cite{Blac90} and it is
denoted there by $V(R)$, for a ring $R$. However, for $L$ a
lattice, the notation $V(L)$ (or $\mathbf{V}(L)$) is used very
widely in lattice theory to denote the \emph{variety}
generated by $L$, so the author of the present work gave up
this notation for the dimension monoid of a lattice.}
The \cm\ \(\DD L\)
is generated by elements \(\Dim(a,b)\), where
$a$, $b\in L$, subjected to a few very simple relations. The
element \(\Dim(a,b)\) may be viewed as a monoid-valued ``distance''
between $a$ and $b$.
It is important to note that the relations
defining the $\Dim$ function are, in particular, satisfied by the
mapping $\Theta$, that with any elements $a$ and
$b$ of $L$, associates the principal congruence
$\Theta(a,b)$ generated by the pair $\vv<a,b>$.
In particular, the dimension monoid%
\index{dimension monoid!of a lattice} \(\DD L\) is a
\emph{precursor} of the congruence
lattice\index{congruence!lattice} \(\con L\) of
$L$; more precisely, the semilattice $\ccon L$ of all finitely
generated congruences of $L$ is isomorphic to the maximal
semilattice quotient of \(\DD L\) (see Corollary~\ref{C:congquotV}).
The dimension monoid gives much more information about the
lattice than the congruence lattice does.

\section{Dimension monoids of special classes of lattices}
\label{S:SpClass}

\subsection*{Modular lattices}
\index{lattice!modular (not necessarily complemented) ---}
We shall study dimension monoids of modular lattices in
Chapter~\ref{ModLatt}. The main observation is that if $L$ is a
modular lattice, then $\DD L$ is a \emph{refinement monoid}%
\index{monoid!refinement ---},
that is, if $a_0$, $a_1$, $b_0$, $b_1$ are elements of $\DD L$ such
that $a_0+a_1=b_0+b_1$, then there are elements $c_{ij}$ (for $i$,
$j<2$) such that $a_i=c_{i0}+c_{i1}$ and
$b_i=c_{0i}+c_{1i}$, for all $i<2$. Furthermore, the algebraic reason
for an equality
$\sum_{i<m}\Dim(a_i,b_i)=\sum_{j<n}\Dim(c_j,d_j)$ is
the existence of ``Schreier-like''
common refinements\index{refinement!Schreier ---} of the
corresponding sequences of intervals, $\vv<[a_i,\,b_i]\mid i<m>$
and $\vv<[c_j,\,d_j]\mid j<n>$.
The sizes of the potential refinements may no longer be predictable,
for example, checking the equality \(\Dim(a,b)=\Dim(a',b')\)
may require large
common refinements of the intervals \([a,\,b]\) and \([a',\,b']\).
This complication does not occur, for example, for complemented
modular lattices.

In general, our alternative construction of $\DD L$ for $L$
modular requires an abstract extension procedure from a partial
semigroup\index{semigroup!partial ---} into a total
semigroup\index{semigroup}. The details of this are worked out
in Chapter~\ref{MaxCommQuot}. Once these monoid-theoretical details
are settled, an alternative presentation of the dimension monoid of
any modular lattice is given as a particular case.

\subsection*{Finite lattices; BCF lattices}
By definition, a partially ordered set is
BCF\index{lattice!BCF ---|see {BCF}}, if it does not have
any infinite bounded chain. Hence every finite partially ordered
set is, of course, BCF. As in the case of modular lattices, we find,
for any BCF lattice $L$, an alternative presentation of $\DD L$. It
should probably not be a surprise that \(\DD L\) is defined in
terms of certain finite configurations of the lattice $L$, which we
shall call
\emph{caustic pairs}\index{caustic!pair}. The form of this set
of relations shows that the corresponding dimension monoids are
always refinement\index{monoid!refinement ---} monoids.
Furthermore, these are not arbitrary refinement monoids, but
monoids of a very special kind: they are the
\emph{primitive}\index{monoid!primitive ---} monoids studied
by R. S. Pierce \cite{Pier89}. It follows, in particular, that
$\DD L$ enjoys many other properties, such as antisymmetry,
unperforation\index{unperforation}
\cite{Good86}, separativity\index{separativity}
\cite{Wehr94}, pseudo-cancellation\index{pseudo-cancellation}
\cite{ShWe94,Wehr92a,Wehr92b}, or the interval
axiom\index{interval axiom} \cite{Wehr96a}. Moreover,
\(\DD L\), endowed with its
algebraic preordering (which turns out to be an
\emph{ordering}), embeds into a power of
\(\ZZ^+\cup\{\infty\}\). Thus, dimensionality is given by a
\emph{family} of numerical
(\(\ZZ^+\cup\{\infty\}\)-valued) dimension functions,
as opposed to only one such function. Finally,
the results about dimension monoids of BCF lattices
(especially the fact that they satisfy the \emph{interval
axiom}\index{interval axiom}, see Corollary~\ref{C:pptiesV(L)})
implies that the natural map from a finite lattice
into its rectangular extension has the ``dimension extension
property'' (see Corollary~\ref{C:RectDEP}).\medskip

These two sections already illustrate how much more
information the dimension monoid carries than the
congruence lattice. This feature can be crystallized by the
following result (see Corollary~\ref{C:SimpleGeom}):

\begin{quote}
\em If $L$ is a simple%
\index{lattice!simple ---}
geometric%
\index{geometric lattice|see {lattice}}%
\index{lattice!geometric ---} lattice, then \(\DD L\) is
isomorphic to \(\ZZ^+\), if $L$ is modular, and to the
two-element semilattice $\two$\index{tzzwo@$\two$|ii},
otherwise.
\end{quote}

Each half of this result requires a quite
careful look into the dimension theory of the corresponding
class of lattices (modular%
\index{lattice!modular (not necessarily complemented) ---},
BCF\index{lattice!BCF ---}, resp.).\medskip

\subsection*{Continuous geometries}
A continuous geometry\index{continuous!geometry|ii} $L$ is defined
as a complete\index{lattice!complete ---},
\mcont\index{lattice!meet-continuous ---},
\jcont\index{lattice!join-continuous ---},
complemented modular lattice. If $L$ is, in addition,
\emph{indecomposable}%
\index{continuous!geometry!indecomposable ---}, then one can
define a \([0,\,1]\)-valued \emph{dimension function} $\Dim$
on $L$ satisfying the following properties:

\begin{enumerate}
\item \(\Dim(1_L)=1\).

\item The map $\Dim$ is \emph{additive}, that is, if $a$ and
$b$ are independent elements of $L$, then
\(\Dim(a\vee b)=\Dim(a)+\Dim(b)\).

\item For all $a$, $b\in L$, \(\Dim(a)=\Dim(b)\) if and only
if $a$ and $b$ are \emph{perspective}\index{perspectivity}.
\end{enumerate}

This \emph{unary} dimension function can also be viewed as a
\emph{binary} distance function on $L$, still denoted by $\Dim$,
defined by
 \[
 \Dim(a,b)=\Dim(x),\text{ for every sectional complement }x
 \text{ of }a\wedge b\text{ in }a\vee b.
 \]
In particular, $L$ is \emph{complete} with respect to this
metric. When $L$ is no longer assumed to be indecomposable%
\index{continuous!geometry!indecomposable ---}, then
one can still define a dimension function $\Dim$ satisfying
2 and 3 above. However, $\Dim$ is no longer numerical, but
takes its values in an interval (of the form \([0,\,u]\)) of a
Dedekind-complete lattice-ordered group; this result is
contained (although this is not explicitly stated in this way) in
\cite{Neum60,Iwam44}.

The dimension theory of continuous geometries originates in
J. von~Neumann's ground-breaking work \cite{Neum60}; it was
developed in further, considerable, work of I. Amemiya, I. Halperin,
T. Iwamura. References for this are, for example,
\cite{AmHa59,Halp38,Halp39a,Halp39b,Halp61,HaNe40,Iwam44}.
Another exposition of von~Neumann's work can be found in F.~Maeda's
book \cite{FMae58}. Further work on the dimension theory of more
general lattices equipped with an orthocomplementation can be found
in \cite{Fill65,Kalm83,Loom55,SMae55}.
The methods introduced by von~Neumann were used to a
large extent in important subsequent work, especially on modular
lattices. References for this are, for example,
\cite{Free76,Free79,Free80, Free87,Herr84,HeHu74,JiRo92,Jons60}.

As it will turn out in this work, the von~Neumann dimension and
ours are, up to a canonical isomorphism, \emph{identical}.
Of course, although our definition of dimensionality is much
shorter than von~Neumann's definition, there is no magic
behind this, since one cannot avoid von~Neumann's difficult theory
of continuous geometries anyway.

Once again, the dimension monoid gives more information about the
structure of the lattice than the congruence lattice: in the
case of an indecomposable
\index{continuous!geometry!indecomposable ---} continuous
geometry\index{continuous!geometry} $L$, the congruence lattice of
$L$ is trivial (that is, $L$ is simple), while the dimension
monoid of $L$ is isomorphic either to $\ZZ^+$ or to $\RR^+$.

\section[Relatively complemented modular lattices]%
{Relatively complemented modular lattices and\\
von~Neumann regular rings}

\subsection*{Fundamental open problems}
We shall state in this section two well-known hard open problems.
The first problem is stated in lattice theory, the second in ring
theory. However, after a closer look, these problems seem to live
simultaneously in different worlds.

\begin{itemize}
\item \emph{Lattices.}
The Congruence Lattice Problem%
\index{congruence!Congruence Lattice Problem|ii}, raised by
R. P. Dilworth in the early forties (see \cite{GrSc0}):
if $A$ is a distributive algebraic lattice, does there exist a
lattice $L$ whose congruence%
\index{congruence!lattice}
lattice is isomorphic to $A$?

\item \emph{Rings.}
The Representation Problem of
Refinement Monoids\index{monoid!refinement ---}
by (von~Neumann) regular%
\index{ring!(von Neumann) regular ---} rings, raised by
K.~R. Goodearl, see \cite{Good95}, in 1995,  but implicit in several
still unsolved problems of \cite{Good91}: if $M$ is a \crm%
\index{monoid!conical refinement ---} with order-unit, does
there exist a regular ring $R$ such that \(V(R)\)%
\index{dimension monoid!of a ring} is
isomorphic to $M$?
\end{itemize}

In the second problem above, $V(R)$, the \emph{dimension monoid} of
$R$, denotes the monoid of isomorphism
types of finitely generated projective right $R$-modules.

At the present time, the status of the two problems are related:
the second problem has counterexamples in size $\aleph_2$ and above
(this starts with \cite{WehrB} and continues in a more
lattice-theoretical fashion in
\cite{WehrC,PTWe,TuWe}) but the answer below $\aleph_2$ is not
known. The first problem is solved positively in the case
where $A$ has at most $\aleph_1$ compact elements,
see A.~P. Huhn \cite{Huhn89a,Huhn89b}, but it is
still open for higher cardinalities.

Many other ring-theoretical
problems which can be formulated lattice-the\-o\-ret\-i\-cal\-ly can
be found in
\cite{Good91}. Note that the appropriate
interface between (complemented, modular) lattice theory and
(regular)%
\index{ring!(von Neumann) regular ---} ring theory is provided
by von~Neumann's Coordinatization%
\index{von Neumann Coordinatization Theorem!original form}
Theorem \cite{Neum60,FMae58}. A positive solution to the
representation problem of \crm s%
\index{monoid!conical refinement ---} by regular%
\index{ring!(von Neumann) regular ---} rings would require the
construction of complicated regular%
\index{ring!(von Neumann) regular ---} rings;
however, such methods are not known at present.
On the other hand, the field of lattice theory is
rich with many construction methods which do not seem to have
any analogues in ring theory, most of them based on
\emph{partial lattices} (see \cite{Grat}).

\subsection*{Dimension theory of lattices and rings}

The dimension theory of regular rings is, really, the dimension
theory of complemented modular lattices. Therefore, from
Chapter~\ref{BasicSCL} on, all our lattices will be
relatively complemented and modular as a rule.

For the convenience of the reader,
Chapter~\ref{BasicSCL} recalls some basic properties of relatively
complemented modular lattices that will be used in the sequel,
mostly continuity conditions, independency%
\index{independent} and perspectivity\index{perspectivity}. In
Chapter~\ref{RelCompl}, it is shown that the dimension map
$\Dim$ is a \Vmeas%
\index{Vmeas@\Vmeas}\ as defined in H. Dobbertin
\cite{Dobb83}. In particular, the dimension
range\index{dimension range} is a lower subset of the
dimension monoid. The preparatory
Corollary~\ref{C:V(I)inV(L)} is the first of a whole series
that will be discussed further.

In Chapter~\ref{RegRings}, we establish, in particular, that in
many cases, the (nonstable) K-theory of a regular%
\index{ring!(von Neumann) regular ---} ring (that is,
essentially the study of \(V(R)\)
for $R$ regular\index{ring!(von Neumann) regular ---}) is
coherent with the dimension theory of the associated lattices.
The basic idea is the following. Let
$R$ be a regular\index{ring!(von Neumann) regular ---} ring.
Consider the (complemented, modular) lattice
\(L=\mathcal{L}(R_R)\)
\index{LzzofR@$\mathcal{L}(R_R)$, $\mathcal{L}(M)$|ii} of all
principal right ideals\index{ideal!of a ring} of
$R$ and the relation of isomorphy $\cong$ on $L$, and let
$I$, $J\in\mathcal{L}(M)$. Are $I$ and $J$ isomorphic if
and only if \(\Dim(\zero,I)=\Dim(\zero,J)\) holds? The answer
is negative, in general
(see Corollary~\ref{C:RegRingCounterex}), but positive in many
cases. The reason for this is the following: the relation
of isomorphy $\cong$ is what we call a \emph{normal}%
\index{normal!equivalence} equivalence on $L$,
that is, it is additive, refining, and it contains
perspectivity\index{perspectivity}, and any two independent
isomorphic ideals\index{ideal!of a ring} are
perspective\index{perspectivity}. Furthermore, in quite a
number of cases, there is
\emph{at most} one normal\index{normal!equivalence} equivalence
on a given lattice. This is sufficient to ensure, for example,
that \(V(R)\cong\DD L\) if
$L$ satisfies a certain condition weaker than the existence of
a homogeneous basis%
\index{homogeneous!basis} with at least two elements
(see Corollary~\ref{C:uniqnor}). Furthermore, this makes it possible
to give a precise characterization of isomorphy in
lattice-theoretical terms---more precisely, in terms of the
square \(\sim_2\)%
\index{Pzzersp2@$\sim_2$|ii} of the relation of
perspectivity\index{perspectivity} in $L$.

To illustrate this, we now present two applications of this theory.
The context for our first example is Section~\ref{S:LatMonInd}, where
we relate A.~P. Huhn's definition of
$n$-distributivity\index{distributive!$n$- ---} with a certain,
simpler, monoid-theoretical concept of index
(see Definition~\ref{D:MonInd}), so that for a relatively
complemented modular lattice $L$, every closed interval of $L$ is
$n$-distributive, for some $n$, if and only if every element of
$\DD L$ has finite index. In particular, by using this, together
with some machinery in the theory of partially ordered abelian
groups (that will appear elsewhere), one can prove the following
statement:

 \begin{quote}
 \em Let $L$ be a relatively complemented lattice. If the congruence
 lattice of $L$ is isomorphic to the three-element chain, then there
 exists a closed interval $I$ of $L$ such that $I$ is not
 $n$-distributive, for all $n\in\omega$.
 \end{quote}

Another application, given in Corollary~\ref{C:NoReprRing}, is the
following. In \cite{WehrB}, the author has shown that there exists
a dimension group\index{group!dimension ---} $G$ of size $\aleph_2$
with order-unit that is not isomorphic to $K_0(R)$ for any
von~Neumann regular ring $R$. The example given there is a
complicated free construction. However, by using the results of
\cite{WehrB,WehrC,PTWe} and of
Section~\ref{S:RegRing}, one can obtain such
dimension groups $G$ that are, in addition, easy to define:
namely, take the Grothendieck group of the dimension monoid of any
free lattice with at least $\aleph_2$ generators in any locally
finite, modular, non-distributive variety.

\subsection*{Generalizations of continuous geometries}

Most of the dimension theory of regular rings can be carried out
for arbitrary relatively complemented modular lattices, without any
assumption that the lattice is coordinatizable.

Chapters~\ref{FinDistr} to \ref{CtbleMeetCo} are devoted to
the dimension theory of relatively complemented modular
lattices that satisfy additional \emph{completeness}
assumptions. A classical example of such a study is of course
von~Neumann's theory of continuous\index{continuous!geometry}
geometries, but the latter does not include, for example,
lattices of subspaces of a vector space (with trivial dimension
theory). Neither does the more general framework of
complemented modular lattices which are both
$\aleph_0$-\mcont%
\index{lattice!meet-continuous ---!$\aleph_0$-\mcont\ ---} and
$\aleph_0$-\jcont%
\index{lattice!join-continuous ---!$\aleph_0$-\jcont\ ---}
studied by I. Halperin
\cite{Halp38,Halp39a,Halp39b}. On the other hand, the
ring-theoretical analogue of this is the dimension theory of
right-continuous%
\index{continuous!right- ---|see {ring}}%
\index{ring!right-continuous ---}
regular \index{ring!(von Neumann) regular ---}
rings and it has been
quite extensively studied,
see, for example, the works by P. Ara, C. Busqu\'e and K.~R.
Goodearl in \cite{Ara87,Busc90,Good82,Good91}. In the latter
context, the countable additivity of the dimension function for
$\aleph_0$-right-con\-tin\-u\-ous%
\index{continuous!$\aleph_0$-right- ---|see {ring}}%
\index{ring!$\aleph_0$-right-continuous ---}
regular\index{ring!(von Neumann) regular ---} rings is obtained
by Goodearl in the case of directly finite%
\index{directly finite!ring} regular\index{ring!(von Neumann) regular ---}
rings \cite[Theorem 2.2]{Good82}, and then extended by Ara to
arbitrary $\aleph_0$-right-con\-tin\-u\-ous%
\index{ring!$\aleph_0$-right-continuous ---} regular rings in
\cite[Theorem 2.12]{Ara87}.

One of the main ideas of the
corresponding arguments is, roughly speaking, to reduce first the
problem to $\aleph_0$-right self-injective regular rings;
in this context, we can use the uniqueness of injective hulls. As
there do not seem to be corresponding direct lattice-theoretical
translations of the latter concepts, the treatment of the dimension
theory of non-coordinatizable lattices cannot be done in this
framework. The lack of a unified theory for all these different
situations may be highlighted by Halperin's following remark in his
1961 survey paper \cite{Halp61}:

\begin{quote}
One suspects that the $\aleph_0$-descending
continuity axiom should be dropped or replaced by some weaker
axiom, just as in measure theory we drop the restriction that
the whole space be of finite measure. We would then presumably
get a wider theory of lattice-structure and dimension theory
with the dimension no longer finite; the lattice of all linear
subspaces of an infinite dimensional vector space would then
be included, as it should be.
\end{quote}

This we do in Chapters~\ref{FinDistr} to \ref{CtbleMeetCo}: in
fact, the $\aleph_0$-descending continuity
(= $\aleph_0$-join-con\-ti\-nu\-i\-ty) can be, as hoped by
Halperin, \emph{completely} dropped, so that in particular, we get a
unified lattice-theoretical setting for all the contexts mentioned in
the previous paragraph. In particular, the dimension function
is countably additive (see Proposition~\ref{P:V(L)GCA}).

In order to prove that the corresponding dimension monoids
also satisfy some completeness conditions, the core of the
problem turns out to be the following \emph{normality
problem}: prove that all lattices under consideration satisfy
the sentence
\[
(\forall x,y)\bigl((x\approx y\text{ and }x\wedge y=0)
\Rightarrow x\sim y\bigr)
\]
where $\sim$ (resp., $\approx$) denotes the relation of
perspectivity\index{perspectivity}
(resp., projectivity\index{projectivity!of elements}).
Not every complemented modular lattice is normal%
\index{normal!lattice}%
\index{lattice!normal ---|see {normal}} (see
Section~\ref{S:NonNorm} for a
counterexample, pointed to the author by Christian Herrmann),
but we prove here normality for a class of lattices large
enough for our purpose, that is, for $\aleph_0$-\mcont%
\index{lattice!meet-continuous ---!$\aleph_0$-\mcont\ ---}
complemented modular lattices (see Theorem~\ref{T:MeetContNor}).
Our proof uses an improvement of the von~Neumann
Coordinatization%
\index{von Neumann Coordinatization Theorem!
J\'onsson improvement}
Theorem established by B. J\'onsson in his 1960 paper
\cite{Jons60}, which makes it possible to reduce the normality problem
to \emph{$3$-distributive}%
\index{distributive!$n$- ---} lattices (see
Corollary~\ref{C:L/Nor(L)FinInd}); and for the latter,
perspectivity\index{perspectivity} is transitive (see
Corollary~\ref{C:nDistrTrans}).\smallskip

Once this normality property is established, the dimension
theory of $\aleph_0$-\mcont%
\index{lattice!meet-continuous ---!$\aleph_0$-\mcont\ ---}
complemented modular (or, more generally, conditionally
$\aleph_0$-\mcont%
\index{lattice!meet-continuous ---!$\aleph_0$-\mcont\ ---}
relatively complemented modular) lattices can be started. This
we do in Chapter~\ref{CtbleMeetCo}, where we also extend the
results to $\aleph_0$-\emph{join}-continuous%
\index{continuous!$\aleph_0$-join- ---|see {lattice}}%
\index{lattice!join-continuous ---!$\aleph_0$-\jcont\ ---}
complemented modular lattices. The idea of this last step is
to note that if $L$ is a $\aleph_0$-\mcont%
\index{lattice!meet-continuous ---!$\aleph_0$-\mcont\ ---}
complemented modular lattice, then its dimension monoid
\(\DD L\) is a special kind of commutative monoid called a
\emph{generalized cardinal algebra} (GCA)%
\index{generalized cardinal algebra (GCA)}, as defined by A.~Tarski
in \cite{Tars49}. Therefore, \(\DD L\) satisfies a certain
monoid-theoretical axiom, given in Theorem~\ref{T:AxImplNorm}, that
implies in turn normality\index{normal!lattice} of the lattice. On
the other hand, we shall see that in order to establish the fact that
\(\DD L\) is a GCA%
\index{generalized cardinal algebra (GCA)}
if $L$ is $\aleph_0$-\mcont%
\index{lattice!meet-continuous ---!$\aleph_0$-\mcont\ ---},
one needs to prove \emph{first} the normality%
\index{normal!lattice} of $L$. Note that this is not the first
time that GCA's appear in dimension theoretical considerations
about lattices, see, for example, P.~A. Fillmore \cite{Fill65}, which
deals with not necessarily modular lattices equipped with an
operation of relative orthocomplementation, \(\vv<x,y>\mapsto y-x\)
for \(x\leq y\), along with an equivalence relation $\sim$ subjected to
certain axioms.

At this point, an important difference with the study of
continuous geometries appears:

\begin{quote}
\em In the context of $\aleph_0$-\mcont\ complemented modular
lattices, dimensionality is no longer equivalent to
perspectivity.
\end{quote}

We prove instead that two elements $x$ and $y$ are
equidimensional (in symbol, \(\Dim(x)=\Dim(y)\)) if and only
if there are decompositions \(x=x_0\oplus x_1\) and
\(y=y_0\oplus y_1\) such that \(x_0\sim x_1\) and
\(y_0\sim y_1\), see Theorem~\ref{T:normeqform}.

This theory will make it possible, in particular, to
extend some of the results of \cite{Ara87,Good82}, to
$\aleph_0$-\emph{left}-continuous%
\index{ring!$\aleph_0$-left-continuous ---}
regular\index{ring!(von Neumann) regular ---} rings. These
results are applied in Section~\ref{S:ContGeom}
of Chapter~\ref{CtbleMeetCo}, where we apply
the previous results of Chapter~\ref{CtbleMeetCo} to lattices
satisfying both $\kappa$-completeness assumptions (for $\kappa$
an infinite cardinal) \emph{and} a finiteness condition (no
nontrivial homogeneous\index{homogeneous!sequence} sequence),
as the situation without this finiteness condition seems to be
much more complicated. The results of this short section
generalize classical results about
continuous\index{continuous!geometry} geometries.

\subsection*{Acknowledgments}
The author wishes to thank warmly Christian Herrmann, who
communicated to him quite a number of very useful comments, such
as, in particular, the proof of the main result of
Section~\ref{S:NonNorm}---namely, the presentation of a non-normal
modular ortholattice. Warm thanks are also due to George
Gr\"atzer, who flooded the manuscript with
red ink, tirelessly pointing thousands of various style errors.

\chapter*{Notation and terminology}

We shall denote by $\omega$\index{ozzmega@$\omega$|ii} the set
of all natural numbers.

If \(\vv<P,\leq>\) is a partially preordered set, a
\emph{lower subset}\index{lower subset|ii}
of $P$ is a subset $X$ such that if
\(x\in X\) and \(p\leq x\) in $P$, then \(p\in X\). If $x$ and
$y$ are elements of $P$, we put
\[
[x,\,y]=\{z\in P\mid x\leq z\leq y\},
\]
and write \(x\prec y\), if \(x<y\) and \([x,\,y]=\{x,y\}\).

We say that $P$ satisfies the \emph{interpolation property}%
\index{interpolation!property|ii}, if for all positive
integers $m$ and $n$ and all elements $a_i$ (for $i<m$) and $b_j$
(for $j<n$) of $P$ such that $a_i\leq b_j$ (for all $i$, $j$), there
exists an $x\in P$ such that \(a_i\leq x\leq b_j\) (for all
$i$, $j$); note that it suffices to verify this for $m=n=2$.
\smallskip

A \emph{semigroup}\index{semigroup|ii} is a nonempty set
endowed with an associative binary operation, while a
\emph{monoid}\index{monoid|ii} is a
semigroup with a
zero element. All our semilattices will be join-semilattices.
We shall need in Chapter~\ref{MaxCommQuot} to introduce
\emph{partial semigroups}%
\index{semigroup!partial ---} which may fail to be
commutative. The definition of a partial semigroup will be
similar to the definition of a category, except for the
existence of identities (such ``categories without identities''
are sometimes called \emph{allegories}). If $S$ is a
semigroup, we shall denote by
\(S^\circ\)\index{szzzero@$S^\circ$|ii} the monoid
consisting of $S$ with a new zero element $\mathsf{O}$
adjoined. A monoid $M$ is
\emph{conical}\index{monoid!conical ---|ii}, if it satisfies
the sentence
\[
(\forall x,y)(x+y=0\Longrightarrow x=y=0).
\]
Every \cm\ $M$ can be endowed with its
\emph{algebraic}\index{algebraic!preordering|ii} preordering
\(\leq_{\mathrm{alg}}\)%
\index{azzlgsym@$\leq_{\mathrm{alg}}$|ii}, defined by the
formula
\[
x\leq_{\mathrm{alg}}y\Longleftrightarrow (\exists z)(x+z=y).
\]
The \emph{Grothendieck group} of $M$ is the abelian group of
all elements of the form \([x]-[y]\), $x$, $y\in M$, where
\([x]-[y]=[x']-[y']\) if and only if there exists $z$ such
that \(z+x+y'=z+x'+y\). It is a partially preordered abelian
group, with positive cone the submonoid
\(\{[x]\mid x\in M\}\). The monoid $M$ can also be endowed
with two binary relations,
$\propto$%
\index{pzzropto@$\propto$|ii} and
$\asymp$\index{azzsymp@$\asymp$|ii}, defined by
\begin{gather*}
x\propto y\Leftrightarrow
(\exists n\in\NN)(x\leq_{\mathrm{alg}}ny);\\
x\asymp y\Leftrightarrow
(x\propto y\text{ and }y\propto x).
\end{gather*}
In a given semigroup\index{semigroup}, a
\emph{refinement matrix}\index{refinement!matrix|ii} is an
array of the form
\[
\begin{tabular}{|c|c|c|c|c|}
\cline{2-5}
\multicolumn{1}{l|}{} & $b_0$ & $b_1$ & $\ldots$ &
$b_{n-1}$\tvi\\
\hline
$a_0$\tvi & $c_{00}$ & $c_{01}$ & $\ldots$ &
$c_{0,n-1}$\\
\hline
$a_1$\tvi & $c_{10}$ & $c_{11}$ & $\ldots$ &
$c_{1,n-1}$\\
\hline
$\vdots$ & $\vdots$ & $\vdots$ & $\ddots$ & $\vdots$\\
\hline
$a_{m-1}$\tvi & $c_{m-1,0}$ & $c_{m-1,1}$ &
$\ldots$ & $c_{m-1,n-1}$\\
\hline
\end{tabular}
\]
where \(a_i=\sum_{j<n}c_{ij}\) for all \(i<m\) and
\(b_j=\sum_{i<m}c_{ij}\) for all \(j<n\).

The \emph{refinement property}\index{refinement!property|ii} is
the semigroup-theoretical\index{semigroup} axiom that states
that every equation \(a_0+a_1=b_0+b_1\) admits a refinement,
that is, elements \(c_{ij}\) (for $i$, $j<2$) such that
\(a_i=c_{i0}+c_{i1}\) and \(b_i=c_{0i}+c_{1i}\), for all
\(i<2\). Equivalently,
for any finite sequences \(\vv<a_i\mid i<m>\) and
\(\vv<b_j\mid j<n>\) of elements of the underlying
semigroup\index{semigroup} such that $m$, $n>0$ and
\(\sum_{i<m}a_i=\sum_{j<n}b_j\), there exists a refinement
matrix as above. A \emph{refinement
monoid}\index{monoid!refinement ---|ii} is a \cm\ satisfying
the refinement property.

A submonoid $I$ of a \cm\ $M$ is an
\emph{ideal}\index{ideal!of a monoid|ii}, if it is, in
addition, a lower subset of $M$ for the algebraic preordering.
One associates with
$I$ a monoid congruence\index{congruence!monoid ---}
\(\equiv_I\)\index{EzzquivI@$\equiv_I$|ii} defined by
\[
x\equiv_Iy\Leftrightarrow(\exists u,v\in I)(x+u=y+v),
\]
and then \(M/I=M/\equiv_I\)\index{MzzoverI@$M/I$|ii} is a
conical\index{monoid!conical ---}
\cm. In addition, if
$M$ satisfies the refinement property, then so does
\(M/I\). Denote by \(\Id M\)
\index{IzzdofM@$\Id M$|ii} the
(algebraic) lattice of ideals of $M$, and by \(\Idc M\)
\index{IzzdofMc@$\Idc M$|ii} the semilattice of all compact
elements of \(\Id M\).
\smallskip

For every partially preordered group $G$, we shall denote by
$G^+$ the positive cone of $G$. In particular,
\(\ZZ^+\)\index{zzzeeplus@$\ZZ^+$|ii} is the monoid of all
natural numbers; we shall put \(\NN=\ZZ^+\setminus\{0\}\)%
\index{nzzatural@$\NN$|ii}. It is well-known
\cite[Proposition 2.1]{Good86} that $G$ satisfies the
interpolation\index{interpolation!property} property if and
only if $G^+$ satisfies the refinement property; we say then
that $G$ is an
\emph{interpolation group}\index{group!interpolation ---|ii}.
Furthermore, we say that
$G$ is \emph{unperforated}\index{unperforation|ii}, if $G$
satisfies the statement \((\forall x)(mx\geq0\Rightarrow x\geq0)\),
for every positive integer $m$. A
\emph{dimension group}\index{group!dimension ---|ii} is a
directed, unperforated interpolation group.
\smallskip

Our lattices will not necessarily be bounded. For every lattice
$L$, we define a subset
$\diag L$\index{dzziag@$\diag L$|ii} of $L\times L$ by the
following formula:
\[
\diag L=\{\vv<x,y>\mid x\leq y\ \text{in}\ L\}.
\]
Moreover, we shall put
\(\S(L)=\{[x,\,y]\mid x\leq y\mbox{ in }L\}\)%
\index{SzzofL@$\S(L)$|ii}. An \emph{ideal}
\index{ideal!of a semilattice|ii} of a semilattice is a
nonempty lower subset closed under the join operation. A
lattice $L$ is
\emph{complete}\index{lattice!complete ---|ii}, if every
subset of $L$ admits a supremum (thus every subset of $L$ also
admits an infimum). An element $a$ of a lattice is
\emph{compact}\index{compact (element of a lattice)|ii},
if for every subset $X$ of $L$ such that the supremum
\(\bigvee X\) exists and \(a\leq\bigvee X\), then there exists
a finite subset $Y$ of $X$ such that \(a\leq\bigvee Y\). A
lattice is
\emph{algebraic}\index{lattice!algebraic ---|ii}, if it is
complete and every element is a supremum of compact elements. If
$L$ is a lattice, we shall denote by
\(\con L\)\index{czzonL@$\con L$|ii} the lattice of all
congruences\index{congruence!lattice} of $L$ (it is a
distributive algebraic lattice) and by
\(\ccon L\)%
\index{czzconL@$\ccon L$|ii} the
semilattice\index{congruence!semilattice} of all compact
( = finitely generated) congruences of $L$
(it is a distributive join-semilattice with zero). A lattice
$L$ is \emph{complemented}\index{lattice!complemented ---|ii}, if it
is bounded (that is, it has a least and a largest
element) and every element of $L$ has a complement. A
lattice $L$ is \emph{relatively
complemented}\index{lattice!relatively complemented ---|ii},
if every closed interval \([a,\,b]\) of $L$ is complemented;
\emph{sectionally complemented}
\index{lattice!sectionally complemented ---|ii}, if it has a
least element, denoted by $0$, and every interval of the form
\([0,\,a]\) is complemented. Note that for a \emph{modular} lattice,
the following chain of implications is valid:
\[
\text{(complemented)}\Rightarrow\text{(sectionally complemented)}
\Rightarrow\text{(relatively complemented)}.
\]
Of course, none of these implications remains valid for arbitrary
lattices.\smallskip

All our rings will be \emph{associative} and \emph{unital}
(but not necessarily commutative). If $R$ is a
regular\index{ring!(von Neumann) regular ---} ring, we shall
denote it as \(R_R\) when viewed as a right
$R$-module over itself. Furthermore, we shall denote by
$\FP(R)$%
\index{FzzPR@$\FP(R)$|ii} the class of finitely generated
projective right $R$-modules. For every
$M\in\FP(R)$, we shall denote by
$\mathcal{L}(M)$
\index{LzzofR@$\mathcal{L}(R_R)$, $\mathcal{L}(M)$|ii} the
(complemented, modular) lattice of all finitely generated
right submodules of $M$, by $nM$ the direct sum
$M\oplus\cdots\oplus M$ ($n$ times), for all $n\in\NN$, and by
$[M]$ the isomorphism class of $M$. The set
$V(R)$%
\index{dimension monoid!of a ring}%
\index{VzzR@$V(R)$, $R$ ring|ii}
of all isomorphism classes of elements of
$\FP(R)$ (\emph{dimension monoid} of $R$)
can be endowed with an addition defined by
\([A]+[B]=[A\oplus B]\) (see \cite{Good95}); the structure
$V(R)$ is a \crm%
\index{monoid!conical refinement ---}\
\cite[Theorem 2.8]{Good91}, with order-unit $[R]$. Then
\(K_0(R)\) is defined as the Grothendieck
\index{group!Grothendieck ---}
group of \(V(R)\). In particular, \(K_0(R)\) is a partially
preordered abelian group with order-unit.
\smallskip

Our proofs will sometimes be divided into claims.
For the sake of readability, proofs of claims will be
written within the proof environment beginning with
\medskip

\noindent\textsc{Proof of Claim}
\medskip

\noindent and ending with
\medskip

\begin{equation*}
\tag*{\qedsymbol\,Claim\,$n$}
\end{equation*}
\medskip

\noindent where $n$ is the number of the claim.

\mainmatter

\chapter{The dimension monoid of a lattice}\label{DimMon}

We start out by defining, for any lattice, its dimension monoid,
$\DD L$. The map $L\mapsto\DD L$ can be extended to a functor on
the category of lattices and lattice homomorphisms. This functor
satisfies many properties satisfied by the functor
$L\mapsto\ccon L$\index{czzconL@$\ccon L$}, although they are harder to
establish as a rule.

\section{Basic categorical properties}

The congruence semilattice functor, $L\mapsto\ccon L$, from
lattices and lattice homomorphisms to $\{\vee,0\}$-semilattices and
$\{\vee,0\}$-homomorphisms, preserves direct products of finitely
many factors, as well as direct limits. In this section, we shall
define the dimension functor and establish that it also satisfies
these properties, see Proposition~\ref{P:Vfunctor}.

\begin{definition}\label{D:V(L)}
Let $L$ be a lattice. Then
the \emph{dimension monoid}
\index{dimension monoid!of a lattice|ii} of
$L$, denoted by \(\DD L\)\index{DzzL@$\DD L$, $L$ lattice|ii},
is the commutative monoid\index{monoid} defined by
generators \(\Dim(a,b)\)%
\index{Dzzimab@$\Dim(a,b)$, $\Dim_L(a,b)$|ii},
\(a\leq b\) in $L$, and the following relations:

\begin{itemize}
\item[(D0)]\index{D012@(D0), (D1), (D2)|ii}
\(\Dim(a,a)=0\), for all \(a\in L\).
\item[(D1)] \(\Dim(a,c)=\Dim(a,b)+\Dim(b,c)\), for all
$a\leq b\leq c$ in $L$.
\item[(D2)] \(\Dim(a,a\vee b)=\Dim(a\wedge b,b)\),
for all $a$, $b$ in $L$.
\end{itemize}
\end{definition}

We shall write $\Dim_L(a,b)$%
\index{Dzzimab@$\Dim(a,b)$, $\Dim_L(a,b)$|ii} for
$\Dim(a,b)$ in case $L$ is not understood.
If $L$ has a zero element $0$, then we shall write
$\Dim(x)$\index{Dzzimabx@$\Dim(x)$|ii} for
$\Dim(0,x)$.

The \emph{dimension range}\index{dimension range|ii} of
$L$ will be, by definition, the set
\(\{\Dim(a,b)\mid a\leq b\ \mathrm{in}\ L\}\).

As usual, say that two intervals are
\emph{transposes of each other},
if they are of the form $[a\wedge b,\,b]$ and
$[a,\,a\vee b]$---in symbol,
\([a\wedge b,\,b]\nearrow[a,\,a\vee b]\)%
\index{tzzransppup@$\nearrow$|ii} and
\([a,\,a\vee b]\searrow[a\wedge b,\,b]\)%
\index{tzzransppdn@$\searrow$|ii}. The rule (D2) can then be
read as ``if two intervals are transposes of each other, then they
have the same length''.

We denote by $\sim$%
\index{PzzerspInt@$\sim$ (for intervals)|ii}
(resp., $\approx$)%
\index{PzzrojInt@$\approx$ (for intervals)|ii} the symmetric
(resp., equivalence) relation generated by
$\nearrow$). The relation $\approx$ is called
\emph{projectivity}%
\index{projectivity!of intervals|ii}.

\begin{notation}\label{N:K0(L)}
Let $L$ be a lattice. Then we shall denote by
$K_0(L)$\index{KzzzeroL@$K_0(L)$|ii} the Grothendieck group%
\index{group!Grothendieck ---} of
$\DD L$.
\end{notation}

Therefore, $K_0(L)$ is, in fact, a
\emph{directed partially preordered abelian group}, of which
the positive cone is the maximal cancellative quotient of
$\DD L$,
that is, the quotient of $\DD L$ by the monoid
congruence\index{congruence!monoid ---}
$\equiv$ defined by
$x\equiv y$ if and only if there exists $z$ such that
\(x+z=y+z\). As we shall see in Chapter~\ref{RegRings}, this
notation is in some sense consistent with the corresponding
notation for von~Neumann regular%
\index{ring!(von Neumann) regular ---} rings; in particular,
there are examples where the natural preordering of
$K_0(L)$ is not a partial order.

\begin{proposition}\label{P:V(L)cone}
For every lattice $L$, $\DD L$ is a conical \cm
\index{monoid!conical ---}; furthermore, $\Dim(a,b)=0$ if and only $a=b$,
for all $a\leq b$ in $L$.
\end{proposition}

\begin{proof}
Let $\nu$ be the \two-valued function defined on \(\diag L\)
by the rule
\[
\nu(\vv<a,b>)=\left\{
\begin{array}{ll}
0,&\mbox{(if $a=b$)}\\
\infty,&\mbox{(if $a<b$)}.
\end{array}
\right.
\]
Then $\nu$ obviously satisfies (D0), (D1) and (D2) above, thus there
exists a unique monoid\index{monoid} homomorphism
$\phi\colon \DD L\to\two$ such that $\phi(\Dim(a,b))=\nu(\vv<a,b>)$,
for all $a\leq b$ in $L$; in particular,
${\phi}^{-1}\{0\}=\{0\}$. The conclusion follows
immediately.
\end{proof}

If $f\colon A\to B$ is a lattice homomorphism, then it is easy
to see that there exists a unique monoid\index{monoid}
homomorphism
$\DD f\colon\DD A\to\DD B$
such that
\[
\DD(f)(\Dim_A(a,a'))=\Dim_B(f(a),f(a')),
\]
 for all $a\leq a'$ in $A$.

\begin{proposition}\label{P:Vfunctor}
The operation $\DD$
defines a (covariant) functor from the category of lattices
(with lattice homomorphisms) to the category of
conical\index{monoid!conical ---} \cm s (with monoid
homomorphisms). Furthermore, $\DD$ preserves direct limits and
finite direct products.
\end{proposition}

\begin{proof}
The verification of the fact that $\DD$ is a
direct limit preserving functor is tedious but
straightforward. As to direct products, it suffices
to prove that, if $A$, $B$, and $C$ are lattices such that
$C=A\times B$ and if $p\colon C\twoheadrightarrow A$ and
$q\colon C\twoheadrightarrow B$ are the natural
projections, then
\[
\epsilon\colon\DD C\to\DD A\times\DD B,\quad
\gamma\mapsto(\DD(p)(\gamma),\DD(q)(\gamma))
\]
is an isomorphism.

To prove this, first let $a_0\leq a_1$ in $A$ and
$b_0\geq b_1$ in $B$. Then applying (D2) to
$c_0=\vv<a_0,b_0>$ and $c_1=\vv<a_1,b_1>$ yields
\[
\Dim_C(\vv<a_0,b_0>,\vv<a_1,b_0>)=
\Dim_C(\vv<a_0,b_1>,\vv<a_1,b_1>),
\]
whence it follows easily that the value of
\(\Dim_C(\vv<a_0,b>,\vv<a_1,b>)\) is independent of the
element \(b\in B\). Then, the function which with
\(\vv<a_0,a_1>\) (for \(a_0\leq a_1\) in $A$) associates
\(\Dim_C(\vv<a_0,b>,\vv<a_1,b>)\) (for some \(b\in B\))
is easily seen to satisfy (D0), (D1) and (D2), thus
there exists a unique monoid homomorphism
\(\varphi\colon \DD A\to\DD C\)
such that the equality
\[
\varphi(\Dim_A(a_0,a_1))=\Dim_C(\vv<a_0,b>,\vv<a_1,b>)
\]
holds for all \(a_0\leq a_1\) in $A$ and all $b\in B$.
Similarly, there exists a unique monoid homomorphism
$\psi\colon\DD B\to\DD C$ such
that the equality
\[
\psi(\Dim_B(b_0,b_1))=\Dim_C(\vv<a,b_0>,\vv<a,b_1>)
\]
holds for all $a\in A$ and all $b_0\leq b_1$ in $B$.
Thus let $\eta\colon\DD A\times\DD B\to\DD C$ be defined by
$\eta(\vv<\alpha,\beta>)=\varphi(\alpha)+\psi(\beta)$. It is
easy to verify that $\epsilon$ and $\eta$ are mutually
inverse.
\end{proof}

Note also that a proof similar to the one of the
functoriality of $\DD$ yields the following result:

\begin{proposition}
\label{P:V(Lop)} Let
\(L^\mathrm{op}\)\index{Lzzop@$L^\mathrm{op}$|ii}
be the dual lattice of $L$. Then there exists an
isomorphism from \(\DD L\) onto
\(\DD(L^\mathrm{op})\) that with
every \(\Dim_L(a,b)\) (for \(a\leq b\) in $L$)
associates \(\Dim_{L^\mathrm{op}}(b,a)\).
\qed\end{proposition}

\section{Basic arithmetical properties}

Let us first note that in any lattice one can extend the
definition of $\Dim(a,b)$ to all pairs $\vv<a,b>$, and define a
new operation $\Dim^+$\index{Dzzimabplus@$\Dim^+(a,b)$|ii}, by
putting
\[
\Dim(a,b)=\Dim(a\wedge b,a\vee b)\ \text{ and }\
\Dim^+(a,b)=\Dim(a\wedge b,a)
\]
(note that the new $\Dim$ notation is consistent with the
old one). Throughout this section, we shall fix a lattice $L$. Our
main purpose is to provide some computational ease with the $\Dim$
symbol. In particular, we shall see that $\Dim$ is a monoid-valued
``distance'' on $L$, see Proposition~\ref{P:triangineq}, and that
the meet and the join operations are \emph{continuous} with respect
to this operation, see Proposition~\ref{P:modularlaw}.

\begin{lemma}\label{L:dimfromdimplus}
The equality
\[
\Dim(a,b)=\Dim^+(a,b)+\Dim^+(b,a)
\]
holds for all $a$, $b\in L$.
\end{lemma}

\begin{proof}
A simple calculation:
\begin{align*}
\Dim(a,b)=\Dim(a\wedge b,a\vee b)&=
\Dim(a\wedge b,a)+\Dim(a,a\vee b)\\
&=\Dim(a\wedge
b,a)+\Dim(a\wedge b,b)\\
&=\Dim^+(a,b)+\Dim^+(b,a).\tag*{\qed}
\end{align*}
\renewcommand{\qed}{}
\end{proof}

In the following two lemmas, for all elements $a$,
$b$, and $c$ of
$L$, define three elements \(\mathrm{Mod}(a,b,c)\)%
\index{Mzzodabc@$\mathrm{Mod}(a,b,c)$|ii},
$\mathrm{Distr}(a,b,c)$%
\index{Dzzistrabc@$\mathrm{Distr}(a,b,c)$|ii} and
$\mathrm{Distr}^*(a,b,c)$%
\index{Dzzistrsabc@$\mathrm{Distr}^*(a,b,c)$|ii} of
$\DD L$ by putting
\begin{align*}
\mathrm{Mod}(a,b,c)&=\Dim(b\vee(a\wedge c),a\wedge(b\vee c)),\\
\mathrm{Distr}(a,b,c)&=
\Dim((a\wedge c)\vee(b\wedge c),(a\vee b)\wedge c),\\
\mathrm{Distr}^*(a,b,c)&=
\Dim((a\wedge b)\vee c,(a\vee c)\wedge(b\vee c)).\\
\end{align*}

\begin{lemma}
\label{L:modularlaw0}
The following equality
\[
\Dim(a,b) =\Dim(a\wedge c,b\wedge c)+\Dim(a\vee c,b\vee c)
+\mathrm{Mod}(a,b,c)
\]
holds for all $a$, $b$, $c\in L$ such that $a\geq b$.
\end{lemma}

Note that in particular, if $L$ is
modular\index{lattice!modular (not necessarily complemented)
---}, then the remainder term\goodbreak $\mathrm{Mod}(a,b,c)$
is always equal to $0$ (this follows from $a\geq b$).\smallskip

\begin{proof}
Since
$a\geq a\wedge(b\vee c)\geq b\vee(a\wedge c)\geq b$,
\begin{align*}
\Dim(a,b) &=\Dim(b,b\vee(a\wedge c))+
\Dim(b\vee(a\wedge c),a\wedge(b\vee c))\\
&\quad+\Dim(a\wedge(b\vee c),a)\\
&=\Dim(a\wedge b\wedge
c,a\wedge c)+
\Dim(b\vee(a\wedge c),a\wedge(b\vee c))\\
&\quad+\Dim(b\vee
c,a\vee b\vee c)\\
&=\Dim(a\wedge c,b\wedge c)+\mathrm{Mod}(a,b,c)+
\Dim(a\vee c,b\vee c).\tag*{\qed}
\end{align*}
\renewcommand{\qed}{}\end{proof}

\begin{proposition}\label{P:modularlaw}
The following equality
\begin{align*}
\Dim(a,b) &=\Dim(a\wedge c,b\wedge c)+\Dim(a\vee c,b\vee c)\\
&\quad+\bigl[\mathrm{Mod}(a\vee b,a\wedge b,c)+
\mathrm{Distr}(a,b,c)+\mathrm{Distr}^*(a,b,c)\bigr]
\end{align*}
holds for all $a$, $b$, $c\in L$.
\end{proposition}

\begin{proof}
Again a simple sequence of calculations:
\begin{align}
\Dim(a,b) &=\Dim(a\wedge b,a\vee b)
\tag{\text{by definition}}\\
&=\Dim(a\wedge b\wedge c,(a\vee b)\wedge c)+
\Dim((a\wedge b)\vee c,a\vee b\vee c)
\notag\\
&\quad+\mathrm{Mod}(a\vee b,a\wedge b,c)
\tag{\text{by Lemma~\ref{L:modularlaw0}}}\\
&=\Dim(a\wedge b\wedge c,(a\wedge c)\vee(b\wedge c))+
\mathrm{Distr}(a,b,c)+\mathrm{Distr}^*(a,b,c)
\notag\\
&\quad+\Dim((a\vee c)\wedge(b\vee c),a\vee b\vee c)+
\mathrm{Mod}(a\vee b,a\wedge b,c),\notag
\end{align}

while we also have, by the definition of the extended $\Dim$,
\begin{align*}
\Dim(a\wedge c,b\wedge c)+\Dim(a\vee c,b\vee c)&=
\Dim(a\wedge b\wedge c,(a\wedge c)\vee(b\wedge c))\\
&\quad+\Dim((a\vee c)\wedge(b\vee c),a\vee b\vee c);
\end{align*} the conclusion follows immediately.\end{proof}

In particular, we see that we always have the inequality
\[
\Dim(a,b)
\geq\Dim(a\wedge c,b\wedge c)+
\Dim(a\vee c,b\vee c)\]
for the algebraic\index{algebraic!preordering} preordering of
$\DD L$.

\begin{proposition}\label{P:triangineq}
The following equality
\begin{itemize}
\item[\rm(i)] $\Dim^+(a,c)\leq\Dim^+(a,b)+\Dim^+(b,c)$;
\item[\rm(ii)] $\Dim(a,c)\leq\Dim(a,b)+\Dim(b,c)$
\end{itemize}
holds for all $a$, $b$, $c\in L$.
\end{proposition}

\begin{proof} Because of Lemma~\ref{L:dimfromdimplus}, it
suffices to prove (i). Put $c'=a\wedge c$. Then $a\geq c'$ and
$\Dim^+(a,c')=\Dim(c',a)$. We shall now perform a few
computations which can be followed on Figure~\ref{Fig:Free21}
(especially for applications of the rule (D2)).

\begin{figure}[hbt]
\begin{picture}(100,200)(-20,-20)
\thicklines
\put(27.88,2.12){\line(-1,1){25.76}}
\put(2.12,32.12){\line(1,1){25.76}}
\put(2.12,132.12){\line(1,1){25.76}}
\put(32.12,157.88){\line(1,-1){25.76}}
\put(57.88,127.88){\line(-1,-1){25.76}}
\put(62.12,127.88){\line(1,-1){45.76}}
\put(30,97){\line(0,-1){34}}
\put(27.88,102.12){\line(-1,1){25.76}}
\put(32.12,57.88){\line(1,-1){25.76}}
\put(32.12,2.12){\line(1,1){25.76}}
\put(62.12,32.12){\line(1,1){45.76}}

\put(30,-8){\makebox(0,0){$b\wedge c'$}}
\put(-5,30){\makebox(0,0)[r]{$c'$}}
\put(65,30){\makebox(0,0)[tl]{$a\wedge b$}}
\put(-5,130){\makebox(0,0)[r]{$a$}}
\put(35,65){$(a\wedge b)\vee c'$}
\put(35,90){$a\wedge(b\vee c')$}
\put(115,80){\makebox(0,0)[l]{$b$}}
\put(65,130){\makebox(0,0)[bl]{$b\vee c'$}}
\put(30,168){\makebox(0,0){$a\vee b$}}

\put(30,0){\circle{6}}
\put(0,30){\circle{6}}
\put(60,30){\circle{6}}
\put(30,60){\circle{6}}
\put(30,100){\circle{6}}
\put(110,80){\circle{6}}
\put(0,130){\circle{6}}
\put(60,130){\circle{6}}
\put(30,160){\circle{6}}

\end{picture}
\caption{Lattice generated by one chain and one element}
\label{Fig:Free21}
\end{figure}
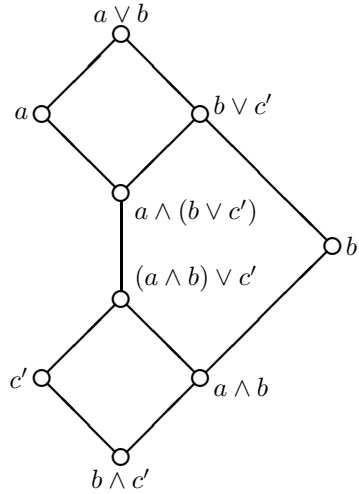

To start with, \(a\geq a\wedge(b\vee c')\geq a\wedge b\), thus
\begin{align*}
\Dim^+(a,b)
&=\Dim(a\wedge b,a\wedge(b\vee c'))+
\Dim(a\wedge(b\vee c'),a)\\
&\geq\Dim((a\wedge b)\vee c',a\wedge(b\vee c'))+
\Dim(a\wedge(b\vee c'),a),\\
\end{align*}
while
\begin{align*}
\Dim^+(b,c)=\Dim(b\wedge c,b)&\geq
\Dim(a\wedge b\wedge c,a\wedge b)
&\text{(by Proposition~\ref{P:modularlaw})}\\
&=\Dim(b\wedge c',a\wedge b)&\\
&=\Dim(c',(a\wedge b)\vee c')
&\text{(by (D2))},
\end{align*}
whence we obtain that
\begin{align*}
\Dim^+(a,b)+\Dim^+(b,c)&\geq
\Dim(c',(a\wedge b)\vee c')+
\Dim((a\wedge b)\vee c',a\wedge(b\vee c'))\\
&\quad+\Dim(a\wedge(b\vee c'),a)\\
&=\Dim(c',a)\\
&=\Dim^+(a,c).
\end{align*}
So (i) holds, thus so does (ii).\end{proof}

\chapter{Dimension monoids and
congruences\index{congruence!lattice ---}}\label{DimCong}

As mentioned in the Introduction, if $L$ is any lattice, then the
maximal semilattice quotient of $\DD L$ is isomorphic to
$\ccon L$\index{czzconL@$\ccon L$}.
We provide a proof of this in Section~\ref{S:IdDimMon}, see
Corollary~\ref{C:congquotV}. Furthermore, in
Section~\ref{S:VmodLatt}, we give a dimension analogue of the axiom
of ``weak modularity'' used to study congruences, see
\cite{Grat}. In Section~\ref{S:DistrLatt}, we settle the rather
easy question of characterization of dimension monoids of
\emph{distributive} lattices, which, of course, generalizes the
characterization result of congruence semilattices of distributive
lattices.

\section{Ideals of the dimension monoid}\label{S:IdDimMon}

In this section, we shall introduce a natural correspondence
between congruences\index{congruence!lattice ---} of a
lattice $L$ and ideals\index{ideal!of a monoid} of
$\DD L$
\index{VzzL@$\DD L$, $L$ lattice}. This correspondence is
given as follows:

\begin{itemize}
\item For every congruence\index{congruence!lattice ---}
$\theta$ of $L$, let
$\DD(\theta)$%
\index{Dzztheta@$\DD(\theta)$ ($\theta$ congruence)|ii} be the
ideal\index{ideal!of a monoid} of
$\DD L$ generated by all
elements $\Dim(x,y)$, where
$\vv<x,y>\in\theta$.

\item For every ideal\index{ideal!of a monoid} $I$ of
$\DD L$, let $\C(I)$%
\index{CzzofI@$\C(I)$|ii} be the binary relation on $L$ defined
by
\[
\vv<x,y>\in\C(I)\Longleftrightarrow\Dim(x,y)\in I.
\]
By Propositions \ref{P:modularlaw} and \ref{P:triangineq}, it
is immediate that $\C(I)$ is, in fact, a
congruence\index{congruence!lattice ---} of $L$.
\end{itemize}

In \ref{P:congsonV}--\ref{C:congquotV}, let $L$ be a lattice.

\begin{proposition}\label{P:congsonV}
The two maps $\DD$ and
$\C$ are mutually inverse lattice isomorphisms between
$\Id(\DD L)$ and $\con L$.
\end{proposition}

\begin{proof}
Note first that obviously, the map $\Theta$ satisfies
(D0), (D1) and (D2): this means that \(\Theta(a,a)=0\) (for all
$a\in L$), $\Theta(a,c)=\Theta(a,b)\vee\Theta(b,c)$ (for all
$a\leq b\leq c$ in $L$), and
$\Theta(a,a\vee b)=\Theta(a\wedge b,b)$ (for all $a$, $b\in L$).
Therefore, there exists a unique monoid homomorphism
$\rho\colon\DD L\to\ccon L$ such that
$\rho(\Dim(a,b))=\Theta(a,b)$ (for all $a\leq b$ in
$L$); hence, this also holds without any order restriction on
$a$ and $b$ (because $\Theta(a,b)=\Theta(a\vee b,a\wedge b)$).
Thus let
$\theta$ be a congruence\index{congruence!lattice ---} of $L$.
For all $\vv<x,y>\in\C(\DD(\theta))$,
that is, $\Dim(x,y)\in\DD(\theta)$, there are
$n\in\NN$ and elements $\vv<x_i,y_i>\in\theta$ (for all $i<n$) such
that $\Dim(x,y)\leq\sum_{i<n}\Dim(x_i,y_i)$; it follows, by
applying the homomorphism $\rho$, that
$\Theta(x,y)\subseteq\bigvee_{i<n}\Theta(x_i,y_i)
\subseteq\theta$, so that $\vv<x,y>\in\theta$. So we have
proved that $\C(\DD(\theta))\subseteq\theta$. The
converse being trivial, we have
$\theta=\C(\DD(\theta))$.

Conversely, let $I$ be an ideal\index{ideal!of a monoid} of
$\DD L$. By definition,
$\DD(\C(I))$ is the ideal\index{ideal!of a monoid} of
$\DD L$ generated by all
$\Dim(x,y)$, where
$\vv<x,y>\in\C(I)$, that is,
$\Dim(x,y)\in I$; since $\DD L$
is generated as a monoid by all the elements of the form
$\Dim(x,y)$, we obtain immediately that $I=\DD(\C(I))$.

Finally, it is obvious that $\DD$ and $\C$ are homomorphisms of
partially ordered sets (with the inclusion on each side): the
conclusion follows.
\end{proof}

This can be used, for example, to characterize the natural
monoid homomorphism
$\rho\colon \DD L\to\ccon L$ of the proof above (defined by
$\rho(\Dim(a,b))=\Theta(a,b)$). We first prove a lemma:

\begin{lemma}\label{L:|x-y|,|a-b|} Let $a$, $b$, $x$, and $y$ be
elements of $L$. Then $\vv<x,y>\in\Theta(a,b)$ if and only if
there exists $n\in\NN$ such that
$\Dim(x,y)\leq n\cdot\Dim(a,b)$, that is,
$\Dim(x,y)\propto\Dim(a,b)$.
\end{lemma}

\begin{proof}
Without loss of generality, we can assume that $a\leq b$. If
$\Dim(x,y)\leq n\cdot\Dim(a,b)$ for some $n\in\NN$, then,
applying $\rho$, we obtain $\Theta(x,y)\subseteq\Theta(a,b)$,
that is, $\vv<x,y>\in\Theta(a,b)$.

To prove the converse, let us first recall the definition of
\emph{weak projectivity}\index{weak projectivity|ii} between
closed intervals used (see, for example,
\cite[Chapter III]{Grat}) to describe
congruences\index{congruence!lattice ---} of lattices. If
\(u\leq v\) and
\(u'\leq v'\) in $L$, we say that
\begin{itemize}
\item \([u,\,v]\wpru[u',\,v']\)%
\index{wzzpr@$\wpr$, $\wprd$, $\wpru$|ii}, if \(v\leq v'\) and
\(u=u'\wedge v\).

\item \([u,\,v]\wprd[u',\,v']\)
\index{wzzpr@$\wpr$, $\wprd$, $\wpru$|ii}, if \(u'\leq u\) and
\(v=u\wedge v'\).

\item \([u,\,v]\to[u',\,v']\), if
\([u,\,v]\wpru[u',\,v']\)
\index{wzzpr@$\wpr$, $\wprd$, $\wpru$|ii} or
\([u,\,v]\wprd[u',\,v']\)%
\index{wzzpr@$\wpr$, $\wprd$, $\wpru$|ii}.
\end{itemize}

The first two cases can be visualized on Figure~\ref{Fig:wpr}.

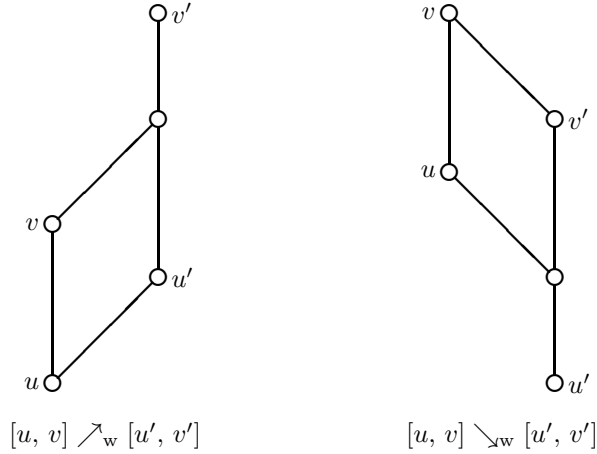
\begin{figure}[hbt]
\begin{picture}(200,200)(-50,-120)
\thicklines
\put(-30,-80){\circle{6}}
\put(-30,-20){\circle{6}}
\put(10,-40){\circle{6}}
\put(10,20){\circle{6}}
\put(10,60){\circle{6}}

\put(-30,-77){\line(0,1){54}}
\put(10,-37){\line(0,1){54}}
\put(10,23){\line(0,1){34}}
\put(-27.88,-77.88){\line(1,1){35.76}}
\put(-27.88,-17.88){\line(1,1){35.76}}

\put(-35,-80){\makebox(0,0)[r]{$u$}}
\put(-35,-20){\makebox(0,0)[r]{$v$}}
\put(15,-40){\makebox(0,0)[l]{$u'$}}
\put(15,60){\makebox(0,0)[l]{$v'$}}

\put(-10,-100){\makebox(0,0){$[u,\,v]\wpru[u',\,v']$}}

\put(120,0){\circle{6}}
\put(120,60){\circle{6}}
\put(160,20){\circle{6}}
\put(160,-40){\circle{6}}
\put(160,-80){\circle{6}}

\put(120,3){\line(0,1){54}}
\put(160,-37){\line(0,1){54}}
\put(160,-77){\line(0,1){34}}
\put(122.12,-2.12){\line(1,-1){35.76}}
\put(122.12,57.88){\line(1,-1){35.76}}

\put(115,0){\makebox(0,0)[r]{$u$}}
\put(115,60){\makebox(0,0)[r]{$v$}}
\put(165,-80){\makebox(0,0)[l]{$u'$}}
\put(165,20){\makebox(0,0)[l]{$v'$}}

\put(140,-100){\makebox(0,0){$[u,\,v]\wprd[u',\,v']$}}

\end{picture}
\caption{Weak projectivity of intervals}\label{Fig:wpr}
\end{figure}

Then, let $\wpr$\index{wzzpr@$\wpr$, $\wprd$, $\wpru$|ii}
(\emph{weak projectivity}) be the
transitive closure of $\to$. We first need the
following claim:

\setcounter{claim}{0}
\begin{claim}
If $[u,\,v]\wpr[u',\,v']$, then
$\Dim(u,v)\leq\Dim(u',v')$.
\end{claim}

\begin{cproof}
It is sufficient to prove the
conclusion for $[u,\,v]\to[u',\,v']$. If
\([u,\,v]\wpru[u',\,v']\),
that is, \(v\leq v'\) and
\(u=u'\wedge v\), then we have
\[\Dim(u,v)=\Dim(u',u'\vee v)\leq\Dim(u',v').\]
Similarly, if
\([u,\,v]\wprd[u',\,v']\), that is,
\(u'\leq u\) and \(v=u\wedge v'\), then we have
\[\Dim(u,v)=\Dim(u\wedge v',v')\leq\Dim(u',v'),\]
thus concluding the proof.
\end{cproof}

Now, in the context of the proof of our Lemma, suppose that
$\vv<x,y>\in\Theta(a,b)$. By
\cite[Theorem III.1.2]{Grat}, there exist
$n\in\NN$ and elements $z_i$ (for $i\leq n$) of $L$ such that
$x\wedge y=z_0\leq z_1\leq\cdots\leq z_n=x\vee y$ and
$[z_i,\,z_{i+1}]\wpr[a,\,b]$,
for all $i<n$. It follows, using Claim~1, that
\begin{equation}
\Dim(x,y)=\Dim(x\wedge y,x\vee y)
=\sum_{i<n}\Dim(z_i,z_{i+1})\leq n\cdot\Dim(a,b).\tag*{\qed}
\end{equation}
\renewcommand{\qed}{}\end{proof}

Now, the promised characterization of $\rho$:

\begin{corollary}\label{C:congquotV}
The natural homomorphism
\(\rho\colon \DD L\to\ccon L\) induces an isomorphism from
the maximal semilattice quotient
\(\DD L/\asymp\) onto \(\ccon L\).
\end{corollary}

\begin{proof}
It is obvious that $\rho$ is surjective; thus, it
suffices to characterize properly the kernel of
\(\rho\). Let $\pi$ be the natural surjective homomorphism from
\(\DD L\) to the semilattice
\(\Idc(\DD L)\) of all
compact%
\index{compact (element of a lattice)} (= finitely
generated) elements of
\(\Id(\DD L)\), so that it sends every
\(\alpha\in\DD L\) to the ideal\index{ideal!of a monoid}
generated by $\alpha$. In particular, for all \(a\leq b\) in
$L$, \(\C\circ\pi(\Dim(a,b))\) is the set of all
\(\vv<x,y>\in L\times L\) such that
\(\Dim(x,y)\in \pi(\Dim(a,b))\), that is,
\(\Dim(x,y)\propto\Dim(a,b)\), hence it is just
\(\Theta(a,b)\) by Lemma~\ref{L:|x-y|,|a-b|}. Hence,
\(\rho=\C\circ\pi\). Since \(\C\) is a lattice isomorphism, we
obtain that
\(\rho(\alpha)\leq\rho(\beta)\) if and only if
\(\pi(\alpha)\leq\pi(\beta)\), for all
\(\alpha\) and \(\beta\) in \(\DD L\); this means that
\(\alpha\propto\beta\). The conclusion follows
immediately.\end{proof}

Hence, in particular, \(\DD L\)
cannot be an arbitrary \cm! Indeed, it is well known (see for
example
\cite[Theorem II.3.11]{Grat}) that for every
lattice $L$, the congruence\index{congruence!lattice} lattice
\(\con L\) of $L$ is distributive. This
is in turn equivalent to saying that
\(\ccon L\) is a distributive
semilattice. Hence, \emph{for every lattice $L$, the maximal
semilattice quotient of
\(\DD L\) is a distributive semilattice}. Thus a natural
question which arises from this is whether
\(\DD L\) is itself a refinement
monoid for every lattice $L$. In the following chapters, we
shall give some partial positive answers to that question, in
two very different cases (in the case of \emph{modular}%
\index{lattice!modular (not necessarily complemented) ---}
lattices and in the case of lattices \emph{without infinite bounded
chains}).\smallskip

Now let us go back to the isomorphisms $\DD$ and $\C$ of
Proposition~\ref{P:congsonV}. They can be used to compute
dimension monoids of quotients:

\begin{proposition}\label{P:latt<>mon} Let $\theta$ be a
congruence\index{congruence!lattice ---} of a lattice $L$. Then
$\DD(L/\theta)\cong\DD L/\DD(\theta)$.
\end{proposition}

\begin{proof} Put $I=\DD(\theta)$. Recall that the monoid
congruence\index{congruence!monoid ---}
$\equiv\pmod I$ on $\DD L$ is defined
by
$\alpha\equiv\beta\pmod I$ if and only if there are
$\alpha'$, $\beta'\in I$ such that
$\alpha+\alpha'=\beta+\beta'$, and that then, by definition,
$\DD L/I=\DD L/\equiv_I$.

\setcounter{claim}{0}
\begin{claim}
For all $x$, $x'$, $y$, $y'\in L$,
if $x\equiv x'\pmod\theta$ and $y\equiv y'\pmod\theta$, then
$\Dim(x,y)\equiv\Dim(x',y')\pmod I$.
\end{claim}

\begin{cproof}
We prove it first for
$x\leq x'$ and $y=y'$. Indeed, in this case, we have
\begin{align*}
\Dim(x',y)+\Dim(x\wedge y,x'\wedge y)&=
\Dim(x\wedge y,x'\wedge y)+\Dim(x'\wedge y,x'\vee y)\\
&=\Dim(x\wedge y,x'\vee y)\\
&(\mathrm{because}\ x'\vee y\geq
x'\wedge y\geq x\wedge y)\\
&=\Dim(x\wedge y,x\vee
y)+\Dim(x\vee y,x'\vee y)\\
&=\Dim(x,y)+\Dim(x\vee y,x'\vee y),
\end{align*} with
\(\Dim(x\wedge y,x'\wedge y)\leq\Dim(x,x')\in I\) and
\(\Dim(x\vee y,x'\vee y)\leq\Dim(x,x')\in I\), so that
\(\Dim(x,y)\equiv\Dim(x',y)\pmod I\).

Now, in case $y=y'$ (but there are no
additional relations on $x$ and $x'$), we have
\begin{align*}
\Dim(x,y) &\equiv\Dim(x\vee x',y)\pmod I\\
&\equiv\Dim(x',y)\pmod I.\\
\end{align*}
Hence, in the general case, we have, using the fact that
$\Dim$ is symmetric,
\begin{align*}
\Dim(x,y)&\equiv\Dim(x',y)\pmod I\\
&\equiv\Dim(x',y')\pmod I,
\end{align*}
which concludes the proof of the claim.
\end{cproof}

Claim~1 above makes it possible to define a map $\varphi_0$ on
$\diag(L/\theta)$ by sending
\(\vv<[x]_\theta,[y]_\theta>\) to
\(\Dim_L(x,y)\ \mathrm{mod}\,I\). It is obvious that
\(\varphi_0\) satisfies
(D0), (D1) and (D2), hence
there exists a unique monoid homomorphism
\(\varphi\colon\DD(L/\theta)\to\DD L/I\) such that the equality
\(\varphi(\Dim_{L/\theta}([x]_\theta,[y]_\theta))=
\Dim_L(x,y)\ \mathrm{mod}\,I\) holds for all
\([x]_\theta\leq [y]_\theta\) in \(L/\theta\). Conversely, let
\(\pi\colon L\twoheadrightarrow L/\theta\) be the natural
homomorphism. Then \(\DD(\pi)\) vanishes on all elements of the
form \(\Dim_L(x,y)\), where
\(x\equiv y\pmod\theta\), thus on all elements of \(I\), hence
it induces a monoid homomorphism
\(\psi\colon \DD L/I\twoheadrightarrow\DD(L/\theta)\), which
sends every
\(\Dim_L(x,y)\ \mathrm{mod}\,I\) (for \(x\leq y\) in $L$) to
\(\Dim_{L/\theta}([x]_\theta,[y]_\theta)\). Then it is trivial
that $\varphi$ and $\psi$ are mutually inverse.
\end{proof}

It follows immediately that dimension monoids of
\emph{countable} lattices are exactly all quotients
$\DD(\FL(\omega))/I$, where $\FL(\omega)$ is the free lattice on
$\omega$ generators and $I$ is an ideal\index{ideal!of a
monoid} of $\DD(\FL(\omega))$.

\begin{lemma}\label{L:genV(I)} Let $L$ be a lattice and let
$X\subseteq L\times L$. Let $\theta=\Theta(X)$ be the
congruence\index{congruence!lattice ---} of $L$ generated by
$X$. Then $\DD(\theta)$ is the
ideal\index{ideal!of a monoid} of
$\DD L$ generated by all
elements $\Dim(x,y)$, where $\vv<x,y>\in X$.
\end{lemma}

\begin{proof} Let $I$ be the ideal\index{ideal!of a monoid} of
$\DD L$ generated by all
elements
$\Dim(x,y)$, $\vv<x,y>\in X$. By definition, we have
$X\subseteq\C(I)$, thus, since $\C(I)$ is a
congruence\index{congruence!lattice ---} of $L$,
$\theta\subseteq\C(I)$, whence
$\DD(\theta)\subseteq\DD(\C(I))=I$. On the
other hand, it is trivial that
$I\subseteq\DD(\theta)$, so that
equality holds.\end{proof}

For example, for any lattice $L$, if $L_\mathrm{mod}$ is the
maximal modular\index{lattice!modular (not necessarily
complemented) ---} quotient of $L$, then
$\DD(L_\mathrm{mod})$ is
isomorphic to $\DD L/I$, where $I$ is the ideal\index{ideal!of
a monoid} of $\DD L$ generated
by all elements
$\mathrm{Mod}(a,b,c)=\Dim(b\vee(a\wedge c),a\wedge(b\vee c))$,
where $a\geq b$ and $c$ are elements of $L$.

\begin{corollary}\label{C:presdiag}
The $\DD$ functor preserves
all pushout diagrams of lattices and lattice homomorphisms of
the form
\[
\begin{CD}
A @>f>> C\\ @VeVV @VV{\bar e}V\\ B @>{\bar f}>> D
\end{CD}
\]
such that $e$ is \emph{surjective}.
\end{corollary}

\begin{proof} It is easy to verify that every pushout diagram
of the form above is isomorphic to one of the form
\[
\begin{CD}
A @>f>> C\\ @VeVV @VV{\bar e}V\\ A/\alpha @>
{\bar f}>> C/\gamma
\end{CD}
\]
where $\alpha$ is a congruence
\index{congruence!lattice ---} of
$A$, $\gamma$ is the congruence\index{congruence!lattice ---}
of $C$ generated by the image of $\alpha$ under $f$, and
$e\colon A\to A/\alpha$ and $\bar e\colon C\to C/\gamma$ are
the natural projections. Hence, by
Proposition~\ref{P:latt<>mon}, the image of this diagram under
$\DD$ is (isomorphic to) a diagram of the form
\[
\begin{CD}
\DD A @>{\DD f}>> \DD C\\ @VVV @VVV\\
\DD A/\DD\alpha @>>> \DD C/\DD\gamma
\end{CD}
\]
where vertical arrows are projections from a monoid to a
quotient under an ideal\index{ideal!of a monoid}
(respectively, $\DD\alpha$ and $\DD\gamma$).
Now, by Lemma~\ref{L:genV(I)}, $\DD\gamma$ is
the ideal\index{ideal!of a monoid} of
$\DD C$ generated by all elements
$\Dim(f(x),f(y))$, where
$\vv<x,y>\in\alpha$, thus it is the ideal\index{ideal!of a
monoid} of $\DD C$ generated by the
image of the ideal\index{ideal!of a monoid} $\DD\alpha$ under
the map
$\DD f$. Thus we are reduced to prove the following
monoid-the\-o\-ret\-i\-cal fact: if $M$ and
$N$ are conical\index{monoid!conical ---} \cm s,
$f\colon M\to N$ is a monoid homomorphism and $I$ is an
ideal\index{ideal!of a monoid} of $M$ and $J$ is the
ideal\index{ideal!of a monoid} of
$N$ generated by $f[I]$, then the following diagram
\[
\begin{CD}
M @>f>> N\\ @VVV @VVV\\ M/I
@>{\bar f}>> N/J
\end{CD}
\]
where the vertical arrows are the natural ones, is a pushout
in the category of conical\index{monoid!conical ---} \cm s with
monoid homomorphisms. Again, this is an easy
verification.\end{proof}

\begin{note} The $\DD$ functor is very far from preserving
arbitrary pushouts.
\end{note}

\section{V-modular lat\-tices}
\label{S:VmodLatt}

The following definition of V-modularity is a slight strengthening
of weak modularity.

\begin{definition}\label{D:Pmod}
A lattice $L$ is
\emph{\Vmod}\index{lattice!V-modular ---|ii}, if
\(\Dim(a,b)\) belongs to the canonical image of
\(\DD [c,\,d]\) in \(\DD L\), for all
elements \(a\leq b\) and \(c\leq d\) of $L$ such that
\([a,\,b]\wpr[c,\,d]\).
\end{definition}

The latter statement means, of course, that there are
\(n\in\omega\) and elements \(c_i\leq d_i\) (for all $i<n$) of
\([c,\,d]\) such that \(\Dim(a,b)=\sum_{i<n}\Dim(c_i,d_i)\).

It is easy to see, using \cite[Corollary III.1.4]{Grat}, that
the conclusion that \(\Dim(a,b)\) belongs to the canonical
image of \(\DD [c,\,d]\) in \(\DD L\) holds then under the
weaker assumption that
\(\vv<a,b>\in\Theta(c,d)\).

Two well-known classes of \Vmod\ lattices are provided by the
following essential observation:

\begin{proposition}\label{P:ExVmod} Every lattice which is
either modular%
\index{lattice!modular (not necessarily complemented) ---} or
relatively complemented is \Vmod%
\index{lattice!V-modular ---}.
\end{proposition}

\begin{proof}
In fact, one proves that for all elements
\(a\leq b\) and \(c\leq d\) of $L$, if
\([a,\,b]\wpr[c,\,d]\),
then there exists a subinterval \([c',\,d']\) of \([c,\,d]\)
such that \([a,\,b]\approx[c',\,d']\). This is classical (see,
for example, \cite[proof of Lemma III.1.7]{Grat} and also
\cite[Exercise III.1.3]{Grat}).
\end{proof}

\begin{proposition}\label{P:VmodVtheta} Let $L$ be a
\Vmod\index{lattice!V-modular ---} lattice and let $X$ be a
lower subset of \(\vv<\S(L),\subseteq>\). Let
$\theta$ be the congruence\index{congruence!lattice ---} of
$L$ generated by $X$. Then
\(\DD\theta\) is generated,
\emph{ as a monoid}, by all elements of the form \(\Dim(x,y)\),
where \(\vv<x,y>\in X\).
\end{proposition}

\begin{proof}
By Lemma~\ref{L:genV(I)}, \(\DD\theta\) is the
\emph{ideal}\index{ideal!of a monoid} of
\(\DD L\) generated by all
elements of the form
\(\Dim(x,y)\), where \([x,\,y]\in X\). Thus, if $M$ denotes the
submonoid of \(\DD L\)
generated by all elements of the form
\(\Dim(x,y)\), where \([x,\,y]\in X\), it suffices to prove
that for all \(a\leq b\) in $L$, all \(m\in\omega\) and all
\([a_i,\,b_i]\) (for \(i<m\)) in $X$, if
\(\Dim(a,b)\leq\sum_{i<n}\Dim(a_i,b_i)\), then
\(\Dim(a,b)\in M\). First, by applying
Corollary~\ref{C:congquotV} to the previous inequality, one
obtains that
\(\vv<a,b>\in\bigvee_{i<n}\Theta(a_i,b_i)\subseteq\theta\),
thus, by \cite[Corollary III.1.4]{Grat}, there exists a
decomposition \(a=c_0\leq c_1\leq\cdots\leq c_n=b\)
satisfying that for all \(j<n\), there exists
\([x,\,y]\in X\) such that \([c_j,\,c_{j+1}]\wpr[x,\,y]\).
But by \Vmod ity\index{lattice!V-modular ---}, it follows then
that \(\Dim(c_j,c_{j+1})\) is a finite sum of elements of the
form \(\Dim(x',y')\), where
\([x',\,y']\subseteq[x,\,y]\). Since $X$ is a lower subset of
\(\vv<\S(L),\subseteq>\), it follows that
\([x',\,y']\in X\); thus \(\Dim(c_j,c_{j+1})\in M\).
Hence $\Dim(a,b)=\sum_{j<n}\Dim(c_j,c_{j+1})$ belongs to $M$.
\end{proof}

\begin{corollary}\label{C:VKidVL}
Let $L$ be a \Vmod\index{lattice!V-modular ---} lattice and
let $K$ be a convex sublattice of $L$. Then the range of the
canonical map from \(\DD K\) to \(\DD L\) is an
ideal\index{ideal!of a monoid} of~\(\DD L\).
\end{corollary}

\begin{proof}
Put \(X=\diag K\)
and let $\theta$ be  the congruence\index{congruence!lattice
---} of $L$ generated by $X$. Then, by
Proposition~\ref{P:VmodVtheta},
\(\DD\theta\) is the set of
all finite sums of elements of the form
\(\Dim(x,y)\), where
\(\vv<x,y>\in X\); thus it is exactly the range of the natural
map from \(\DD K\) to \(\DD L\).
However,
\(\DD\theta\) is, by
definition, an ideal\index{ideal!of a monoid} of \(\DD L\).
\end{proof}

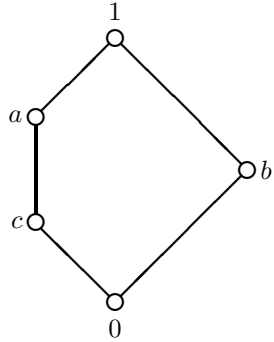
\begin{figure}[hbt]
\begin{picture}(100,140)(0,-70)
\thicklines
\put(27.88,-47.88){\line(-1,1){25.76}}
\put(0,-17){\line(0,1){34}}
\put(2.12,22.12){\line(1,1){25.76}}
\put(32.12,47.88){\line(1,-1){45.76}}
\put(77.88,-2.12){\line(-1,-1){45.76}}

\put(30,-50){\circle{6}}
\put(0,-20){\circle{6}}
\put(0,20){\circle{6}}
\put(30,50){\circle{6}}
\put(80,0){\circle{6}}

\put(-5,-20){\makebox(0,0)[r]{$c$}}
\put(30,-60){\makebox(0,0){$0$}}
\put(-5,20){\makebox(0,0)[r]{$a$}}
\put(30,60){\makebox(0,0){$1$}}
\put(85,0){\makebox(0,0)[l]{$b$}}
\end{picture}

\caption{Dimension without V-modularity}\label{Fig:noVmod}
\end{figure}

If $L$ is not \Vmod\index{lattice!V-modular ---}, then the
conclusion of Corollary~\ref{C:VKidVL} may not hold. For
example, consider the case where $L$ is the pentagon
($L=\{0,a,b,c,1\}$ with $a>c$),
$K$ is the convex sublattice defined by $K=\{0,b\}$, and
\(e:K\hookrightarrow L\) is the inclusion map, see
Figure~\ref{Fig:noVmod}.

Then the range of $\DD e$ is not convex in $\DD L$, as
$\Dim_L(c,a)$ belongs to the lower subset generated
by the range of $\DD e$ but not to the range of $\DD e$.

\begin{corollary}\label{C:simpleVL}
Let $L$ be a simple\index{lattice!simple ---}
\Vmod\index{lattice!V-modular ---}
lattice and let $a$ and $b$ be elements of $L$ such that
\(a\prec b\). Then
\(\DD L=\ZZ^+\cdot\Dim(a,b)\).
\qed
\end{corollary}

Further consequences of this fact can be found in
Corollary~\ref{C:succVhom}, Corollary~\ref{C:SimpleVmodShort}
and Corollary~\ref{C:SimpleGeom}.

\section{Dimension monoids of distributive lat\-tices}
\label{S:DistrLatt}

In this section, we shall compute dimension monoids of
\emph{distributive} lattices. Recall that a \emph{generalized
Boolean lattice} is a sectionally complemented distributive
lattice. It is well-known that every distributive lattice $L$
embeds into a unique (up to isomorphism) generalized Boolean
lattice $B=\mathbf{B}(L)$%
\index{BzzofL@$\mathbf{B}(L)$|ii} which it generates (with
respect to the lattice operations and the operation which with
$\vv<a,b>$ associates the unique relative complement $a\sd b$%
\index{szzmallsm@$a\sd b$|ii} of
$a\wedge b$ in $a$, see, for example, \cite[II.4]{Grat}). Furthermore,
the operation $L\mapsto\mathbf{B}(L)$ is given by a direct limit
preserving \emph{functor}. Then, with every generalized Boolean
lattice $B$, associate the commutative monoid $\ZZ^+[B]$%
\index{zzzeeplusB@$\ZZ^+[B]$|ii} defined by generators $\one_u$
(for $u\in B$) and relations $\one_0=0$,
$\one_{u\vee v}=\one_u+\one_v$, for all $u$, $v\in B$ such that
$u\wedge v=0$. Alternatively, $\ZZ^+[B]$ is isomorphic to the
monoid of all non-negative continuous functions with compact
support from the Stone space of $B$ to $\ZZ^+$. Again, the
operation $B\mapsto\ZZ^+[B]$ is given by a direct limit
preserving functor.

\begin{lemma}\label{L:projint}
Let $L$ be a distributive
lattice and let \([a_0,\,b_0]\), \([a_1,\,b_1]\) in
\(\S(L)\). Then \([a_0,\,b_0]\approx[a_1,\,b_1]\) if and only
if \(b_0\smallsetminus a_0=b_1\smallsetminus a_1\) in
\(\mathbf{B}(L)\).
\end{lemma}

\begin{proof} If
\([a_0,\,b_0]\nearrow[a_1,\,b_1]\), then
\(b_0\smallsetminus a_0=b_1\smallsetminus a_1\)
(this is the identity
\((x\vee y)\smallsetminus x=y\smallsetminus(x\wedge y)\)).
Conversely, suppose that
\(b_0\smallsetminus a_0=b_1\smallsetminus a_1\) in
\(\mathbf{B}(L)\). Put
\(a=a_0\wedge a_1\) and
\(b=b_0\wedge b_1\); it is immediate to verify that
\([a,\,b]\nearrow[a_i,\,b_i]\),
for all \(i<2\), whence
\([a_0,\,b_0]\approx[a_1,\,b_1]\).
\end{proof}

\begin{theorem}\label{T:V(L)} Let $L$ be a distributive
lattice, let \(B=\mathbf{B}(L)\). Then there exists a unique
monoid homomorphism
\(\varphi\colon \DD L\to\ZZ^+[B]\) such that
\[
\varphi(\Dim(a,b))=
\one_{b\smallsetminus a},
\]
for all \(a\leq b\) in $L$; furthermore, $\varphi$ is an
isomorphism.
\end{theorem}

\begin{proof}
The existence and uniqueness of \(\varphi\) are
easy. Moreover, since $B$ is generated by $L$ as a generalized
Boolean algebra, every element of $B$ is a disjoint union of
elements of the form \(b\smallsetminus a\) (for \(a\leq b\) in
$L$), see \cite[II.4, Lemma 3]{Grat}, thus $\varphi$ is surjective.
Now we prove that $\varphi$ is one-to-one. Since the $\DD$ functor
preserves direct limits and since every finitely generated
distributive lattice is finite, it suffices to prove
it in the case where $L$ is finite.
However, in this case, up to
isomorphism, $L$ is the distributive lattice of all lower
subsets of a finite partially ordered set $P$, and
\(\varphi(\Dim(a,b))=\one_{b\setminus a}\in(\ZZ^+)^P\), for all
\(a\subseteq b\) in $L$. Now let
\(\psi\colon (\ZZ^+)^P\to\DD L\) be the unique monoid
homomorphism defined by
\(\psi(\one_{\{p\}})=\Dim(a_p,b_p)\), where we put
\(b_p=\{x\in P\mid x\leq p\}\) and
\(a_p=\{x\in P\mid x<p\}\). To verify that
\(\psi\varphi=\mathrm{id}_{\DD L}\), it suffices to verify that
\(\psi\varphi(\Dim(a,b))=\Dim(a,b)\), for every prime interval
\([a,\,b]\) of $L$. However, the latter condition means that
\(b=a\cup\{p\}\), for some element
\(p\) of \(P\setminus a\), so that we obtain
\(\psi\varphi(\Dim(a,b))=\psi(\one_{\{p\}})=\Dim(a_p,b_p)\).
But
\(b\setminus a=b_p\setminus a_p=\{p\}\), thus, by
Lemma~\ref{L:projint},
\(\Dim(a,b)=\Dim(a_p,b_p)\), whence we obtain that
\(\psi\varphi=\mathrm{id}_{\DD L}\). Therefore, \(\varphi\) is
an isomorphism.
\end{proof}

In particular, by using Corollary~\ref{C:congquotV}, we recover
the classical fact that the
congruence\index{congruence!semilattice} semilattice of any
distributive lattice is a generalized Boolean lattice, see
\cite{CrDi73}.

\begin{corollary}
Let $L$ be a distributive lattice, let $M$
be a \cm, let \(\mu\colon\diag L\to M\) be a map such that
the following conditions
\begin{itemize}
\item[\rm (i)] \(\mu(\vv<a,a>)=0\);
\item[\rm (ii)] \(\mu(\vv<a,c>)=\mu(\vv<a,b>)+\mu(\vv<b,c>)\),
if \(a\leq b\leq c\);
\item[\rm (iii)]
\(\mu(\vv<a\wedge b,b>)=\mu(\vv<a,a\vee b>)\)
\end{itemize}
hold, for all $a$, $b$, and $c$ in $L$.
Then there exists a unique monoid homomorphism
\(\bar\mu\colon\ZZ^+[\mathbf{B}(L)]\to M\) such that
\(\bar\mu(b\smallsetminus a)=\mu(\vv<a,b>)\) holds, for all
\(a\leq b\) in $L$.\qed
\end{corollary}

\chapter{Basic properties of
refinement\index{monoid!refinement ---} monoids}\label{RefMon}

Refinement monoids have been formally introduced, independently,
in H.~Dobbertin \cite{Dobb82} and P.~A. Grillet \cite{Gril70}, and probably
in other places as well.
Nevertheless, their origin can probably be traced back to
A. Tarski \cite{Tars49}. We shall review in this chapter
what we need about refinement monoids.

\section[Systems of inequalities]%
{Solutions of basic systems
of equations and inequalities}

We start with a classical lemma; see, for example,
\cite[Lemma 1.9]{Wehr92a}.

\begin{lemma}\label{L:a+b=nc}
Let $M$ be a
refinement\index{monoid!refinement ---} monoid, let $a$, $b$,
and $c$ be elements of $M$, let $n$ be a positive integer. Then
\(a+b=nc\) if and only if there are elements $x_k$ (for \(0\leq
k\leq n\)) of $M$ such that
\begin{align}
a&=\sum_{0\leq k\leq n}kx_k,\notag\\
b&=\sum_{0\leq k\leq n}(n-k)x_k,\notag\\
c&=\sum_{0\leq k\leq n}x_k.\tag*{\qed}
\end{align}
\end{lemma}

The following consequence will be of interest for dimension monoids
of sectionally complemented modular lattices:

\begin{corollary}\label{C:charasymp}
Let $M$ be a refinement\index{monoid!refinement ---} monoid,
let $a$ and $b$ be elements of $M$. Then the following are equivalent:

\begin{itemize}
\item[\rm (i)] There exists $n\in\NN$ such that $a\leq nb$ and
$b\leq na$.

\item[\rm (ii)] There are finite sets $X$ and $Y$, a subset
$\Gamma\subseteq X\times Y$ with domain $X$ and image
$Y$, and elements $a_i$ (for $i\in X$) and $b_j$ (for $j\in Y$) such
that
\begin{gather*}
a=\sum_{i\in X}a_i,\\
b=\sum_{j\in Y}b_j,\\
(\forall\vv<i,j>\in\Gamma)(a_i=b_j).
\end{gather*}
\end{itemize}
\end{corollary}

\begin{proof}
(ii)$\Rightarrow$(i) is trivial. Conversely,
suppose that (i) holds. By Lemma~\ref{L:a+b=nc}, there are
decompositions of the form
\begin{align*}
a&=\sum_{0\leq k\leq n}kx_k,\\
b&=\sum_{0\leq k\leq n}x_k,
\end{align*}
for $x_0$,\dots, $x_n\in M$. Therefore,
\(x_0\leq b\leq na\leq n^2\sum_{1\leq k\leq n}x_k\), thus,
applying again Lemma~\ref{L:a+b=nc}, there are decompositions
\(x_0=\sum_{0\leq l\leq n^2}ly_l\) and
\(\sum_{1\leq k\leq n}x_k=\sum_{0\leq l\leq n^2}y_l\). Applying
refinement\index{refinement!property} to the latter equation
yields a refinement matrix\index{refinement!matrix} of the form
\[
\begin{tabular}{|c|c|}
\cline{2-2}
\multicolumn{1}{l|}{} & $y_l\ (0\leq l\leq n^2)$\tvi\\
\hline
$x_k\ (1\leq k\leq n)$ & $z_{kl}$\tvi\\
\hline
\end{tabular}
\]
Therefore, we compute
\[
a=\sum_{1\leq k\leq n}kx_k=
\sum_{1\leq k\leq n;\ 0\leq l\leq n^2}kz_{kl},
\]
while
\[
b=x_0+\sum_{1\leq k\leq n}x_k=
\sum_{0\leq l\leq n^2}ly_l+\sum_{1\leq k\leq n}x_k=
\sum_{1\leq k\leq n;\ 0\leq l\leq n^2}(l+1)z_{kl}.
\]
Thus we have obtained decompositions of the form
\(a=\sum_{0\leq k<m}p_kc_k\) and \(b=\sum_{0\leq k<m}q_kc_k\),
with $p_k$, $q_k\in\NN$, and $c_k\in M$ (for all $k<m$). Now put
\begin{align*}
X&=\bigcup_{0\leq k<m}(\{k\}\times p_k),\\
Y&=\bigcup_{0\leq k<m}(\{k\}\times q_k),\\
x_{\vv<k,s>}&=c_k\qquad
(\text{for all}\ \vv<k,s>\in X),\\
y_{\vv<k,s>}&=c_k\qquad
(\text{for all}\ \vv<k,s>\in Y),\\
\Gamma&=\{\vv<\vv<k,s>,\vv<l,t>>\in X\times Y\mid k=l\}.
\end{align*}
Since the $p_k$'s and the $q_k$'s are nonzero,
$\Gamma$ has domain $X$ and image $Y$.
The rest is trivial.
\end{proof}

We refer to \cite[Lemma 1.11]{Wehr92a} for the following lemma.

\begin{lemma}\label{L:ApproxCanc}
Let $M$ be a
refinement\index{monoid!refinement ---} monoid, endowed with
its algebraic\index{algebraic!preordering} preordering $\leq$.
Let $a$, $b$, and
$c$ be elements of $M$, let $n$ be a positive integer.

\begin{itemize}

\item[\rm (i)] If \(a+c=b+c\), then there are elements $d$, $u$,
and $v$ of $M$ such that $a=d+u$, $b=d+v$, and $nu$,
$nv\leq c$.

\item[\rm (ii)] If \(a+c\leq b+c\), then there exists $d\in M$
such that $nd\leq c$ and \(a\leq b+d\).\qed

\end{itemize}
\end{lemma}

The following lemma is one of the basic results of
\cite{Wehr92b}, see \cite[Theorem~1]{Wehr92b}.
Its origin goes back to A. Tarski \cite{Tars49}.

\begin{lemma}
\label{L:ArchPOM} Let $\vv<M,+,0,\leq>$ be a \cm\ endowed with
a preordering $\leq$ satisfying the following two axioms:
\begin{enumerate}

\item \((\forall x)(x\geq 0)\);

\item \((\forall x,y,z)(x\leq y\Rightarrow x+z\leq y+z)\).

\end{enumerate}
Then $M$ embeds into a
power of $\vv<[0,\,\infty],+,0,\leq>$ if and only if it is
antisymmetric and satisfies the following axiom:
\begin{equation}
(\forall x,y)\bigl( (\forall
n\in\NN)(nx\leq(n+1)y)\Longrightarrow x\leq y\bigr).
\tag*{\qed}
\end{equation}
\end{lemma}

In \cite{Wehr92a,Wehr92b}, preordered monoids satisfying
1 and 2 above are called \emph{positively ordered monoids}
(P.O.M.s).

If the condition of Lemma~\ref{L:ArchPOM} is satisfied, we
then say that $\vv<M,+,0,\leq>$ is \emph{Archimedean}%
\index{Archimedean|ii}. This terminology is consistent with
the one used for positive cones of \poag s (namely, if $G$ is a
\poag, then $G$ is Archimedean if and only if $G^+$ satisfies
our Archimedean condition presented here),
thus we prefer it to the somewhat uninspired ``regular'' used
in \cite{Wehr92a,Wehr92b}.

The following lemma has been proved in
\cite[Theorem 2.16]{Wehr92b}.

\begin{lemma}\label{L:ArchRefPOM}
Let $M$ be a refinement
monoid\index{monoid!refinement ---} in which the
algebraic\index{algebraic!preordering} ordering
$\leq$ is \emph{antisymmetric}. Then $M$ is
Archimedean if and only if for all elements
$a$ and $b$ of $M$, and for every sequence
$\vv<c_n\mid n\in\omega>$ of elements of $M$, if
\[
a\leq b+c_n\ \text{ and }\ \sum_{i<n}c_i\leq b
\]
holds for all $n\in\omega$, then $a\leq b$.\qed
\end{lemma}

\section[From interval to monoid]%
{Transfer from a generating interval to the whole monoid}

We shall formulate in this section a series of lemmas whose
meaning is always that the whole structure of a
refinement\index{monoid!refinement ---} monoid can be somehow
reflected by the additive structure of a generating interval.
Lemma~\ref{L:cancref} is a prototype of this kind of situation.

\begin{lemma}[folklore]
\label{L:cancref}
Let $M$ be a
refinement\index{monoid!refinement ---} monoid and let $U$ be a
lower subset of $M$ generating $M$ as a monoid. Suppose that
$U$ is cancellative, that is, \(a+c=b+c\in U\)
implies that $a=b$, for all $a$, $b$, $c\in U$. Then $M$ is cancellative.
\end{lemma}

\begin{proof}
Put
\[
nU=\left\{\sum_{i<n}x_i\mid (\forall i<n)
(x_i\in U)\right\},
\]
for all $n\in\NN$. Since $U$ generates $M$ as a monoid, it suffices to
prove that
$2U$ is cancellative (because then, by induction, $2^nU$ is
cancellative, for all $n$). We first prove that \(a+c=b+c\) implies
that $a=b$, for all $a$, $b\in 3U$, and $c\in U$. Indeed,
by the refinement\index{refinement!property} property, there
exists a refinement\index{refinement!matrix} matrix
\[
\begin{tabular}{|c|c|c|}
\cline{2-3}
\multicolumn{1}{l|}{} & $b$ & $c$\tvi\\
\hline
$a$ & $d$ & $a'$\tvi\\
\hline
$c$ & $b'$ & $c'$\tvi\\
\hline
\end{tabular}
\]
so that \(a'+c'=b'+c'=c\in U\). By assumption, $a'=b'$,
whence
\(a=d+a'=d+b'=b\).

Next, we prove that if $a$, $b$, $c\in 2U$ such that $a+c=b+c$,
then
$a=b$. Indeed, since $c\in 2U$, there are $c_i$ (for $i<2$) in $U$
such that $c=c_0+c_1$. Then \((a+c_0)+c_1=(b+c_0)+c_1\) with
$a+c_0$, $b+c_0\in 3U$, and $c_1\in U$, whence, by the previous
result, \(a+c_0=b+c_0\). Again by previous result,
$a=b$. Hence $2U$ is cancellative.\end{proof}

\begin{note}
A similar result holds for
\emph{separativity}\index{separativity}, see
\cite[Proposition 3.3]{AGPO}. However, this may not hold for other
properties than separativity.
For example, there are non stably%
\index{monoid!stably finite ---} finite \crm s%
\index{monoid!conical refinement ---}
with a directly finite%
\index{directly finite!element in a monoid} order-unit (even
in the simple\index{monoid!simple ---} case, see \cite{WehrA}).
\end{note}

We shall next recall a definition, due essentially to
Dobbertin \cite{Dobb83}:

\begin{definition}\label{D:V-hom} A homomorphism $f\colon M\to
N$ of
\cm s is a \emph{\Vhom}\index{Vhom@\Vhom|ii} at an element
$c\in M$, if for all $\xi$, $\eta\in N$ such that
\(f(c)=\xi+\eta\), there are
$x$ and $y$ in $M$ such that $c=x+y$, $\xi=f(x)$, and
$\eta=f(y)$. We say that $f$ is a \Vhom, if
it is a
\Vhom\ at every element of $M$.
We say that $f$ is a
\emph{\Vemb}\index{Vemb@\Vemb}, if it is a one-to-one
\Vhom.
\end{definition}

\begin{note}
If $f$ is a \Vhom, then the
range of $f$ is necessarily an ideal\index{ideal!of a monoid}
of $N$. Furthermore, if $f$ is one-to-one, then $f$ is also
an embedding for the
\emph{algebraic}\index{algebraic!preordering} preorderings of
$M$ and $N$.
\end{note}

The following lemma is often useful in order to prove that a given
homomorphism is a \Vhom.

\begin{lemma}\label{L:V-homat} Let $M$ be a \cm\ and let $N$
be a refinement\index{monoid!refinement ---} monoid; let
$f\colon M\to N$ be a monoid homomorphism. Then the set of
elements of $M$ at which $f$ is a \Vhom\ is
closed under addition.\qed
\end{lemma}

In a similar spirit, the following lemma is often useful in order to
prove that a given \Vhom\ is one-to-one.

\begin{lemma}[folklore]\label{L:onetoone}
Let $M$ be a \cm\
and let $N$ be a refinement\index{monoid!refinement ---}
monoid. Let $f\colon M\to N$ be a \Vhom\ of
\cm s, let $U$ be a lower subset of $M$ generating $M$ as a
monoid. If $f\res_U$ is one-to-one, then $f$ is
one-to-one.
\end{lemma}

\begin{proof} Our first step is to prove that if
$x\in U$ and $y\in 2U$, then $f(x)=f(y)$ implies that $x=y$.
Indeed, let $y_0$ and $y_1$ in $U$ such that $y=y_0+y_1$. Then
\(f(x)=f(y_0)+f(y_1)\), thus, since $f$ is a
\Vhom, there are elements $x_0$ and $x_1$ of
$M$ such that
$x=x_0+x_1$ (thus $x_0$, $x_1\in U$) and $f(x_i)=f(y_i)$ (for all
$i<2$). By assumption, $x_i=y_i$, whence
\(x=x_0+x_1=y_0+y_1=y\).

Next, let $x$ and $y$ be elements of $2U$ such that
$f(x)=f(y)$. There are decompositions $x=x_0+x_1$ and
$y=y_0+y_1$ with $x_i$, $y_j\in U$. Since $N$ is a refinement
monoid\index{monoid!refinement ---}, there exists a refinement
matrix\index{refinement!matrix} of format
\[
\begin{tabular}{|c|c|c|}
\cline{2-3}
\multicolumn{1}{l|}{} & $f(y_0)$ & $f(y_1)$\tvi\\
\hline
$f(x_0)$\tvi & $\zeta_{00}$ & $\zeta_{01}$\\
\hline
$f(x_1)$\tvi & $\zeta_{10}$ & $\zeta_{11}$\\
\hline
\end{tabular}
\]
with $\zeta_{ij}\in N$. Fix $i<2$. Since $f$ is a
\Vhom, there are elements $z_{ij}$ (for $j<2$)
such that
\(x_i=z_{i0}+z_{i1}\) and \(f(z_{ij})=\zeta_{ij}\) (for all
$j<2$). Thus, for all $j<2$, \(f(y_j)=f(z_{0j}+z_{1j})\) with
$y_j\in U$ and \(z_{0j}+z_{1j}\in 2U\). By the first step
above, \(y_j=z_{0j}+z_{1j}\). It follows that
\(x=x_0+x_1=y_0+y_1=y\), so that $f\res_{2U}$ is
one-to-one. Hence, by induction on
$n$, \(f\res_{2^nU}\) is one-to-one for every $n$. But
$U$ generates $M$ as a monoid, whence $f$ is
one-to-one.
\end{proof}

\section{Index of an element in a monoid}\label{S:LMInd}

\begin{definition}\label{D:MonInd}
Let $M$ be a
conical\index{monoid!conical ---} \cm. For all $x\in M$,
we define the \emph{index} of $x$ in $M$, and we denote it by
$\ind_M(x)$%
\index{index!of an element in a monoid}%
\index{izzndm@$\ind_M(x)$|ii}, as the largest
$n\in\ZZ^+$ such that there exists
$y\in M\setminus\{0\}$ such that $ny\leq x$, if it exists;
otherwise, we put $\ind_M(x)=\infty$.
\end{definition}

Let us first establish a few elementary properties of the monoid
index:

\begin{lemma}\label{L:addInd} Let $M$ be a \crm%
\index{monoid!conical refinement ---}. Then the following
holds:

\begin{itemize}

\item[\rm (i)] The following inequality
\[
\max\{\ind_M(x),\ind_M(y)\}\leq\ind_M(x+y)
\leq\ind_M(x)+\ind_M(y)
\]
holds for all $x$, $y\in M$.

\item[\rm (ii)] The equality
\(\ind_M(nx)=n\cdot\ind_M(x)\)%
\index{index!of an element in a monoid} holds for all $x\in M$ and all
$n\in\NN$.

\item[\rm (iii)] The equality \(\ind_M(x)=0\) holds if and only
if $x=0$, for all $x\in M$.

\end{itemize}
\end{lemma}

\begin{proof} (i) It is obvious that
\(\max\{\ind_M(x),\ind_M(y)\}\leq\ind_M(x+y)\). Next,
let $n\in\NN$ such that \(n>\ind_M(x)+\ind_M(y)\). We prove
that \(n>\ind_M(x+y)\). Thus let $z\in M$ such that
$nz\leq x+y$, we prove $z=0$. Since $M$ is a
refinement\index{monoid!refinement ---} monoid, there are
$x'\leq x$ and $y'\leq y$ such that $nz=x'+y'$. By
Lemma~\ref{L:a+b=nc}, there are decompositions
$z=\sum_{k\leq n}z_k$, $x'=\sum_{k\leq n}kz_k$, and
$y'=\sum_{k\leq n}(n-k)z_k$. For all $k\leq n$, either
$k>\ind_M(x)$ or $n-k>\ind_M(y)$; in the first case,
$kz_k\leq x$ implies that $z_k=0$; in the second case,
$(n-k)z_k\leq y$ implies again that $z_k=0$. Thus $z=0$, which
completes the proof of (i).

(ii) is trivial for \(\ind_M(x)=\infty\); thus suppose that
\(\ind_M(x)<\infty\). Let \(m\in\NN\) such that
\(m\leq n\cdot\ind_M(x)\). Thus there exists \(k\in\NN\) such
that
\(k\leq\ind_M(x)\) and
\(m\leq nk\). Since \(k\leq\ind_M(x)\), there exists \(y>0\) in
$M$ such that \(ky\leq x\). Therefore, \(my\leq nky\leq nx\),
so that \(m\leq\ind_M(nx)\). This proves that
\(n\cdot\ind_M(x)\leq\ind_M(nx)\)%
\index{index!of an element in a monoid}; the converse follows
from (i) above.

(iii) is trivial.
\end{proof}

\begin{corollary}\label{C:addInd} Let $M$ be a \crm%
\index{monoid!conical refinement ---}. Then the set $\fin(M)$%
\index{fzzinM@$\fin(M)$|ii} of all elements of
$M$ with finite index is an ideal\index{ideal!of a monoid} of
$M$.\qed
\end{corollary}

In fact, $\fin(M)$ is a very special sort of \crm%
\index{monoid!conical refinement ---}:

\begin{proposition}\label{P:finIndDim} Let $M$ be a \crm%
\index{monoid!conical refinement ---}. Then $\fin(M)$ is the
positive cone of an Archimedean\index{Archimedean} dimension
group%
\index{group!dimension ---}.
\end{proposition}

\begin{proof} We first prove that \(\fin(M)\) is
cancellative. It suffices to prove that if $a$, $b$,
and $c$ are elements of $M$ such that
\(a+c=b+c\) and \(c\in\fin(M)\), then \(a=b\). Put
\(n=\ind_M(c)+1\)%
\index{index!of an element in a monoid}. By
Lemma~\ref{L:ApproxCanc}, there are decompositions
\(a=d+u\), \(b=d+v\) such that \(nu\leq c\) and \(nv\leq c\).
Since \(\ind_M(c)<n\), this implies that \(u=v=0\), so that
\(a=b\).

But $\fin(M)$ is an ideal\index{ideal!of a monoid} of $M$ and
$M$ is a \crm%
\index{monoid!conical refinement ---}, thus \(\fin(M)\) is also
a \crm; hence, it is the positive cone of an
interpolation\index{group!interpolation ---} group.
To obtain that
\(\fin(M)\) is Archimedean\index{Archimedean}, it suffices, by
Lemma~\ref{L:ArchRefPOM}, to prove that if $a$,
$b\in\fin(M)$, and
\(\vv<c_n\mid n\in\omega>\) is a sequence of elements of
\(\fin(M)\) such that \(a\leq b+c_n\)
and \(\sum_{i<n}c_i\leq b\) for all $n$, then \(a\leq b\). Take
\(n>\ind_M(b)\)%
\index{index!of an element in a monoid}. Since \(\fin(M)\) is
the positive cone of an
interpolation\index{group!interpolation ---} group, there
exists
\(c\in\fin(M)\) such that
\(a\leq b+c\) and \((\forall i<n)(c\leq c_i)\). It follows that
\(nc\leq b\), thus, since \(n>\ind_M(b)\)%
\index{index!of an element in a monoid}, \(c=0\), so that
\(a\leq b\).
\end{proof}

Toward Proposition~\ref{P:Ind1}, let us first establish a lemma.

\begin{lemma}\label{L:genLgroup}
Let $G$ be a \poag, let $U$ be a directed lower subset of
$G^+$ generating $G$ as a group such that the
interval $[0,\,u]$ is lattice-ordered for all $u\in U$. Then $G$ is
lattice-ordered.
\end{lemma}

\begin{proof}
Let $\mathcal{E}$ be the set of all lower subsets
$X$ of $G^+$ such that any two elements of $X$ have an infimum
in $X$. Fix first an element $X$ of $\mathcal{E}$. Then for all
$x$, $y\in X$, the pair $\{x,y\}$ admits an infimum in $X$, say
$z$. Then $z$ is also the infimum of $\{x,y\}$ in $G$: indeed,
for all $s\in G$ such that $s\leq x$, $y$, there exists, by
interpolation\index{interpolation!property}, $t\in G$ such that
$0,s\leq t\leq x$, $y$. Since
$X$ is a lower subset of $G^+$, $t\in X$. But $t\leq z$ by the definition
of $z$, which proves our claim. Hence, for all
$X\in\mathcal{E}$, any two elements of $X$ have an infimum in
$G$.

Next, for any elements $A$ and $B$ of $\mathcal{E}$, we write
$A\mid B$ if $\{x,y\}$ admits
an infimum in $G$, for all $x\in A$ and $y\in B$. Then we prove the
following claim:

\setcounter{claim}{0}
\begin{claim}
Let $A$, $B$, and $C$ be elements of $\mathcal{E}$. If
$A\mid B$ and $A\mid C$, then $A\mid B+C$.
\end{claim}

\begin{cproof}
Let $a\in A$, $b\in B$, and
$c\in C$. Since $A\mid B$, $a\wedge b$ exists; then put
$u=a-a\wedge b$. Note that $u\in A$, thus, since $A\mid C$,
$u\wedge c$ exists. Now let \(d=a\wedge b+u\wedge c\), we
prove that $d=a\wedge(b+c)$. First, it is easy to see that
$d\leq a$ and $d\leq b+c$. Next, let $z\in G$ such that
$z\leq a$ and $z\leq b+c$. Thus $z\leq a+c,b+c$, thus, since
the translation of vector $a$ is an order-isomorphism of $G$,
\(z\leq a\wedge b+c\). But \(z\leq a=a\wedge b+u\), whence
$z\leq a\wedge b+u\wedge c=d\), so that \(d=a\wedge(b+c)\).
\end{cproof}

Now, since $U$ is directed, any two elements of $U$ have an
infimum, that is, $U\mid U$. Thus, by using Claim~1,
one deduces easily that $U\mid nU$, for all $n\in\NN$
(where $nU$ is the set of all sums of $n$ elements of $U$);
then, by using again Claim~1,
$nU\mid nU$ for all $n\in\NN$.
But by assumption, $G^+=\bigcup_{n\in\NN}nU$, so
that $G^+\mid G^+$. Since $U$ generates $G$ as a group, $G$ is
directed: thus $G$ is lattice-ordered.
\end{proof}

\begin{proposition}\label{P:Ind1} Let $M$ be a \crm%
\index{monoid!conical refinement ---}. If $M$ is generated
by a directed lower subset $U$ of elements of index at most $1$,
then $M$ is the positive cone of a lattice-ordered group%
\index{group!lattice-ordered ---}.
\end{proposition}

See also \cite[Corollary 3.4]{Good94} for a similar statement.

\begin{proof} First, by Corollary~\ref{C:addInd} and
Proposition~\ref{P:finIndDim}, $M$ is the positive cone of an
Archimedean\index{Archimedean} dimension group%
\index{group!dimension ---} $G$. Therefore, by
Lemma~\ref{L:genLgroup}, it suffices to prove that $U$ is
lattice-ordered. Thus it suffices to prove that for any
elements $a$ and $b$ of $U$, $a\vee b$ exists (because in any
partially ordered group, $a\wedge b$ exists if and only if
$a\vee b$ exists); so we proceed: since
$U$ is directed, there exists $u\in U$ such that
$a$, $b\leq u$. Since $G$ has
interpolation\index{interpolation!property}, there exists
$c\in G$ such that $a$, $b\leq c\leq a+b$, $u$. Since
$c\leq a+b$ and $a$, $b$, $c\in G^+$, there are, by Riesz
decomposition, elements
$a'$ and $b'$ such that $0\leq a'\leq a$, $0\leq b'\leq b$,
and $c=a'+b'$. Since $c\leq u$ and $u$ has index at most $1$,
we have
$a'\wedge b'=0$. We prove now that $c=a\vee b$. Indeed, by
definition, $a$, $b\leq c$. Let $x\in M$ such that $a$,
$b\leq x$. It follows that $a'$, $b'\leq x$,
thus there are $u$, $v\in M$ such that
\(x=a'+u=b'+v\). Applying Riesz decomposition to the inequality
$a'\leq b'+v$ and noting that \(a'\wedge b'=0\) yields that
$a'\leq v$, we obtain \(u=x-a'\geq x-v=b'\), whence
\(x\geq a'+b'=c\), thus completing the proof that $c=a\vee b$.
\end{proof}

By contrast, in the case of index $2$, we have the following
counterexample (it appears in the final comments of
\cite[Section 3]{Good94}):

\begin{examplepf}\label{E:Index2} There exists a non
lattice-ordered dimension group%
\index{group!dimension ---}%
\index{group!lattice-ordered ---} with or\-der-unit $(G,u)$
such that \(\ind_{G^+}(u)=2\)%
\index{index!of an element in a monoid}.
\end{examplepf}

\begin{proof}
Let $G$ be the set of all sequences
\(\vv<x_n\mid n\in\omega>\) of integers such that
$x_n=x_0+x_1$ for all
large enough $n$, ordered componentwise.
Put $G_n=\ZZ\times\ZZ\times\ZZ^n$ ordered componentwise,
for all $n\in\omega$, and
let $f_n\colon G_n\to G$ such that $f(x)$ is defined by
concatenation of $x$ with the constant sequence with value
$x_0+x_1$. Define similarly $f_{mn}\colon G_m\to G_n$ for
$m\leq n$. Then it is easy to verify that $G$ is the direct
limit of the direct system $\vv<G_m,f_{mn}>_{m\leq n}$ with
limiting maps $f_n$ (for $n\in\omega$), thus $G$ is a dimension
group%
\index{group!dimension ---} (it is even
\emph{ultrasimplicial}, that is, a direct limit of
simplicial groups\index{group!simplicial ---} with embeddings
of \poag s). Now, define elements
\(a=\vv<a_n\mid n\in\omega>\),
\(b=\vv<b_n\mid n\in\omega>\) of \(G^+\) by
\(a_0=b_1=1\), \(a_1=b_0=0\), and
\(a_n=b_n=1\) for all \(n\geq 2\).
Then \(u=a+b\) is an order-unit of $G$ of index
$2$, but \(\{a,b\}\) has no infimum in $G$.
\end{proof}

\section{Generalized cardinal algebras}\label{S:GCAs}

The results of this section will be used mainly in
Chapter~\ref{CtbleMeetCo}. Generalized cardinal algebras have been
introduced in Tarski's monograph \cite{Tars49} as certain
(partial) refinement monoids equipped with a (partial)
infinite addition defined on countable sequences satisfying
quite natural properties. Let us be more specific:

\begin{definition}\label{D:GCA}
A \emph{generalized cardinal algebra}%
\index{generalized cardinal algebra (GCA)|ii} (GCA, for short) is a
quadruple
\(\vv<A,+,0,\sum>\) such that $0$ is an element of $A$, $+$ is a
\emph{partial} binary operation on $A$, $\sum$ is a
\emph{partial} function from
\(A^\omega\) to $A$ (we write \(\sum_na_n\) instead of
\(\sum\vv<a_n\mid n\in\omega>\)) satisfying the following
rules:

\begin{itemize}
\item[\rm (GCA1)] If \(a\in A\) and
\(\vv<a_n\mid n\in\omega>\in A^\omega\), then
\(a=\sum_na_n\) implies that \(a=a_0+\sum_na_{n+1}\).

\item[\rm (GCA2)] If \(c\in A\) and
\(\vv<a_n\mid n\in\omega>,\vv<b_n\mid n\in\omega>\in
A^\omega\), then
\(c=\sum_n(a_n+b_n)\) implies that
\(c=\sum_na_n+\sum_nb_n\).

\item[\rm (GCA3)] \(a+0=0+a=a\), for all \(a\in A\).

\item[\rm (GCA4)] (Refinement
postulate)\index{refinement!postulate|ii} If $a$, $b\in A$, and
\(\vv<c_n\mid n\in\omega>\in A^\omega\) such that
\(a+b=\sum_nc_n\), then there exists a
refinement\index{refinement!matrix} matrix (that is, the
obvious infinite generalization) of the form
\[
\begin{tabular}{|c|c|}
\cline{2-2}
\multicolumn{1}{l|}{} & \(c_n\,(n\in\omega)\)\tvi\\
\hline
$a$\tvi & \(a_n\)\\
\hline
$b$\tvi & \(b_n\)\\
\hline
\end{tabular}
\]
\item[\rm (GCA5)] (Remainder postulate) Let \(\vv<a_n\mid
n\in\omega>\in A^\omega\) and
\(\vv<b_n\mid n\in\omega>\in A^\omega\) be sequences such
that \(a_n=b_n+a_{n+1}\) for all
\(n\in\omega\). Then all infinite sums
\(\sum_kb_{n+k}\) (for \(n\in\omega\)) are defined and there exists
\(a\in A\) such that \(a_n=a+\sum_kb_{n+k}\) for all $n$.
\end{itemize}

For all \(a\in A\), we put \(\infty a=\sum_na\)%
\index{izznftya@$\infty a$|ii}, if the sum is defined.\smallskip

\noindent A \emph{cardinal algebra} (CA)%
\index{cardinal algebra (CA)|ii} is a GCA where both
$+$ and $\sum$ are \emph{total} operations.
\end{definition}

Tarski's monograph \cite{Tars49} contains many often non
trivial consequences of this system of axioms. Here are some of
them:

\begin{proposition}\label{P:basicGCA1}
Let \(\vv<A,+,0,\sum>\) be a GCA%
\index{generalized cardinal algebra (GCA)}. Then the following
properties hold in $A$:

\begin{itemize}
\item[\rm (a)] \(\infty0=0\).

\item[\rm (b)] \(\vv<A,+,0>\) is a commutative, associative
partial monoid.

\item[\rm (c)] If \(\vv<a_k\mid k\in\omega>\) is a
sequence of elements of $A$ and if \(n\in\omega\) such that
\(a_k=0\) for all \(k\geq n\), then \(\sum_ka_k\) is defined if and
only if \(a_0+\cdots+a_{n-1}\) is defined, and then the two
are equal.

\item[\rm (d)] For all $a$, $b\in A$, \(a+b=b\) if and only if
\(\infty a\) is defined and \(\infty a\leq b\).

\item[\rm (e)] The binary relation $\leq$ defined on $A$ by
``\(a\leq b\) if and only if there exists $x$ such that
\(a+x=b\)'' is a partial ordering on $A$.

\item[\rm (f)] If \(\vv<a_n\mid n\in\omega>\) and
\(\vv<b_n\mid n\in\omega>\) are two sequences of elements
of
$A$ such that \((\forall n\in\omega)(a_n\leq b_n)\) and
\(\sum_nb_n\) is defined, then \(\sum_na_n\) is defined and
\(\sum_na_n\leq\sum_nb_n\).

\item[\rm (g)] Let $\sigma$ be a permutation of $\omega$ and
let \(\vv<a_n\mid n\in\omega>\in A^\omega\). Then
\(\sum_na_n\) is defined if and only if \(\sum_na_{\sigma(n)}\)
is defined, and then the two are equal.\smallskip

\noindent{\rm This makes it possible to
define unambiguously the notation
\(a=\sum_{i\in I}a_i\) for every (at most) countable set $I$,
every \(a\in A\), and every
\(\vv<a_i\mid i\in I>\in A^I\).}
\smallskip

\item[\rm (h)] Let $I$ and $J$ be (at most) countable
sets and let \(\vv<a_{ij}\mid \vv<i,j>\in I\times J>\) be
a family of elements of $A$. Then
\(\sum_{\vv<i,j>\in I\times J}a_{ij}\) is defined if and only
if \(\sum_{i\in I}\sum_{j\in J}a_{ij}\) is defined, and then
the two are equal.

\item[\rm (i)]
If \(\vv<a_k\mid k\in\omega>\in A^\omega\) and
\(b\in A\), then \(\sum_ka_k\) is defined and \(\leq b\) if and
only if all partial sums \(\sum_{k<n}a_k\) (for \(n\in\omega\)) are
defined and \(\leq b\). In particular, \(\sum_ka_k\), if
defined, is the supremum of all partial sums \(\sum_{k<n}a_k\)
(for \(n\in\omega\)).

\item[\rm (j)] Let \(\vv<a_m\mid m\in\omega>\) and
\(\vv<b_n\mid n\in\omega>\) be sequences of elements
of
$A$ such that \(\sum_ma_m=\sum_nb_n\). Then there exists an
infinite refinement\index{refinement!matrix} matrix of the
following form:
\[
\begin{tabular}{|c|c|}
\cline{2-2}
\multicolumn{1}{l|}{} & \(b_n\,(n\in\omega)\)\tvi\\
\hline
\(a_m\,(m\in\omega)\)\tvi & \(c_{mn}\)\\
\hline
\end{tabular}
\]
\end{itemize}\qed
\end{proposition}

It is easy to see that every lower subset of a GCA%
\index{generalized cardinal algebra (GCA)} is a GCA. Quite
naturally, every GCA appears, in fact, as a lower subset of a CA.
As in \cite[Section 7]{Tars49}, we say that a
\emph{closure}%
\index{generalized cardinal algebra (GCA)!closure of a ---} of
a GCA $A$ is a CA%
\index{cardinal algebra (CA)} \(\overline A\) such that
$A$ is a subset of \(\overline A\), every element of
\(\overline A\) is a countable sum of elements of $A$, and
\(a=\sum_na_n\)
in $A$ if and only if \(a=\sum_na_n\) in \(\overline A\),
for all \(a\in A\) and
\(\vv<a_n\mid n\in\omega>\in A^\omega\). In
particular, it follows easily that $A$ is a lower subset of
\(\overline A\). We have the following result, see
\cite[Theorem 7.8]{Tars49}:

\begin{proposition}\label{P:basicGCA2}
Let $A$ be a GCA%
\index{generalized cardinal algebra (GCA)}. Then there exists
a unique closure of $A$, up to isomorphism.\qed
\end{proposition}

In fact, Tarski proves existence of the closure (an essential
tool is Proposition~\ref{P:basicGCA1}.(j) above, that is,
the infinite refinement\index{refinement!property} property),
but uniqueness is nearly trivial. Now, the formulation of
``arithmetical'' results is much more convenient in CA's%
\index{cardinal algebra (CA)} than in GCA's\index{generalized
cardinal algebra (GCA)}. Let us list a few more of these:

\begin{proposition}\label{P:basicCA}
Every CA%
\index{cardinal algebra (CA)} satisfies the following
properties:
\begin{itemize}
\item[\rm (a)] (A slight strengthening of the
pseudo-cancellation\index{pseudo-cancellation} property:)
\[
(\forall a,b,c)(a+c=b+c\Rightarrow
(\exists d) (d=2d\leq c\text{ and }a+d=b+d)).
\]

\item[\rm (b)]
\((\forall x,y)(mx\leq my\Rightarrow x\leq y)\), for all \(m\in\NN\).

\item[\rm (c)] For all sequences \(\vv<a_m\mid
m\in\omega>\) and \(\vv<b_n\mid n\in\omega>\) such that
\((\forall m,n\in\omega)(a_m\leq b_n)\), there exists $c$ such
that for all $m$, $n\in\omega$, we have
\(a_m\leq c\leq b_n\).\qed
\end{itemize}
\end{proposition}

In fact (see \cite{RaSh92} or \cite[Corollary 2.17]{Wehr92b})
every CA\index{cardinal algebra (CA)} admits an embedding (for
both the order and the addition but not necessarily for the
infinite addition) into a power of
\(\vv<[0,\,+\infty],+,0,\leq>\).

\section[Infinite interpolation]%
{Infinite $\kappa$-interpolation and monotone
$\kappa$-completeness}

The main reference for this section is \cite{Good86}, from which
most of the results of this section are easy generalizations.

If $\kappa$ is any ordinal, we say that a partially ordered set
$P$ is \emph{monotone $\kappa$-complete}%
\index{monotone@monotone $\kappa$-completeness|ii}, if for
every increasing (resp., decreasing)
$\kappa$-sequence of elements of $P$ which is bounded above
(resp., below) has a supremum (resp., an infimum). It is obvious
that monotone $\kappa$-completeness is equivalent
to monotone \(\mathrm{cf}(\kappa)\)-completeness, where
\(\mathrm{cf}(\kappa)\) is the \emph{cofinality} of $\kappa$
(thus a regular cardinal). In the case where $P$ is a
partially ordered group, it suffices to verify one of the two
dual conditions above. Next, we say that $P$ satisfies the
\emph{$\kappa$-interpolation property}%
\index{interpolation!$\kappa$- --- property|ii}, if
for all families \(\vv<x_\xi\mid \xi<\kappa>\) and
\(\vv<y_\eta\mid \eta<\kappa>\) of elements of $P$ such
that \(x_\xi\leq y_\eta\) (for all $\xi$, $\eta$), there exists
\(z\in P\) such that \(x_\xi\leq z\leq y_\eta\) (for all
$\xi$, $\eta$). The proof of the following lemma is similar to
\cite[Lemma 16.2]{Good86} (the possibility \(\kappa>\aleph_0\)
implies that an additional induction proof is required, but
this proof is easy):

\begin{lemma}\label{L:kInterp} Let $\kappa$ be a infinite
cardinal and let $P$ be a partially ordered set. Then $P$
satisfies the $\kappa$-interpolation
property if and
only if whenever $\alpha$, $\beta\leq\kappa$ are ordinals and
\(\vv<x_\xi\mid \xi<\alpha>\) and
\(\vv<y_\eta\mid \eta<\beta>\) are families of elements of
\(P\) such that \(x_\xi\leq y_\eta\) for all $\xi$, $\eta$,
and one of the following possibilities

\begin{itemize}
\item[\rm (a)] \(\alpha=\beta=2\),

\item[\rm (b)] \(\alpha=2\) and
\(\vv<y_\eta\mid \eta<\beta>\) is decreasing,

\item[\rm (c)] \(\beta=2\) and \(\vv<x_\xi\mid
\xi<\alpha>\) is increasing,

\item[\rm (d)] \(\vv<x_\xi\mid \xi<\alpha>\) is
increasing and
\(\vv<y_\eta\mid \eta<\beta>\) is decreasing
\end{itemize}

\noindent holds, then there exists \(z\in P\) such that
\(x_\xi\leq z\leq y_\eta\), for all $\xi$, $\eta$.\qed
\end{lemma}

It follows, in particular, that every monotone
$\kappa$-complete partially ordered set with interpolation
also satisfies $\kappa$-interpolation.

Now, we say that a partially ordered set $P$ is
\emph{$\kappa$-directed}, if every subset of $P$ of size at
most $\kappa$ is bounded above. The proof of the following
lemma is similar to the proof of
\cite[Proposition 16.3]{Good86} (the directedness assumption
allows to reduce the problem to the case of an order-unit):

\begin{lemma}\label{L:kInt,UtoG}
Let $G$ be an interpolation
group\index{group!interpolation ---}, let $U$ be a lower
subset of \(G^+\) that generates
\(G^+\) as a monoid. If $U$ is
$\kappa$-directed and satisfies the $\kappa$-interpolation
property, then
$G$ satisfies the $\kappa$-interpolation property.\qed
\end{lemma}

The proof of the following lemma is similar to the proof
of \cite[Proposition 16.9]{Good86}:

\begin{lemma}\label{L:MsC,UtoG} Let $G$ be an
interpolation\index{group!interpolation ---} group, let $U$ be
a lower subset of \(G^+\) which generates
\(G^+\) as a monoid. If $U$ is $\kappa$-directed and monotone%
\index{monotone@monotone $\kappa$-completeness}
$\kappa$-complete, then $G$ is monotone%
\index{monotone@monotone $\kappa$-completeness}
$\kappa$-complete.\qed
\end{lemma}

\chapter{Dimension theory of refined partial semigroups%
\index{semigroup!refined partial ---}}\label{MaxCommQuot}

The main purpose of this chapter is to construct, from a given
\emph{partial} semigroup $\vv<S,\oplus>$, endowed with an
equivalence relation, $\approx$, a commutative monoid $M$, together
with a ``universal, $\approx$-invariant, $M$-valued, finitely
additive measure'' from $S$ to $M$. Of course, the existence of
such a measure in general is a trivial, categorical statement, so
that we shall be more interested in the \emph{construction} of $M$,
for certain classes of structures $\vv<S,\oplus,\approx>$, which
we shall call \emph{\rps s}.

\section[Refined partial semigroups]%
{Refined partial semigroups; the
refinement\index{refinement!property} property}

Define a \emph{magma} to be a set endowed with a binary
operation. For any set $A$, let
$\mathcal{M}(A)$ be the free magma on $A$; we shall slightly
stray from the usual notational conventions by denoting by
$\vv<s,t>\mapsto s\oplus t$%
\index{ozzplus@$\oplus$!partial operation (general)|ii} the
natural binary operation on
$\mathcal{M}(A)$. The reason for this introduction lies in the
fact that we shall often be concerned in this paper about
\emph{partial} binary operations. If $\oplus^A$ is a partial
binary operation on a set $A$, then one defines inductively a
binary relation between elements of
$\mathcal{M}(A)$ and elements of $A$, reading ``the term $s$ is
defined with value $a$'', in notation $s\dnw=a$%
\index{dzzownrw@$s\dnw=a$|ii}, as follows: if
$s=b\in A$, then $s\dnw=a$ if and only if $b=a$; if
$s=s_0\oplus s_1$, then $s\dnw=a$ if and only if there exist
$a_0$ and $a_1$ such that $(\forall i<2)(s_i\dnw=a_i)$ and
$a=a_0\oplus^Aa_1$. Thus, for all $s$, there exists at most one
$a$ such that $s\dnw=a$, called naturally the \emph{value} of
$s$; if there exists such an $a$, we shall say that $s$ is
\emph{defined}, and write $s\dnw$; otherwise we shall say that
$s$ is \emph{undefined}. If $s$ and $t$ are terms on $A$,
$s\cong t$ will be the statement
\((\forall a\in A)(s\dnw=a\Leftrightarrow t\dnw=a)\). In the
sequel, we shall often omit the superscript $A$ in the notation
above, just writing $\oplus$ instead of
$\oplus^A$. We shall also often identify a defined term on $A$
with its value.

\begin{definition}\label{D:parsemi}
A \emph{partial semigroup}\index{semigroup!partial ---|ii}
is a nonempty set endowed with a partial binary operation
$\oplus$ which is \emph{associative}%
\index{associative!partial operation|ii}, that is, it
satisfies
\[
(\forall a,b,c) ((a\oplus b)\oplus c\cong a\oplus(b\oplus
c)).
\]
Similarly, we say that $S$ is \emph{commutative}%
\index{commutative!partial operation|ii}, if it satisfies
\[
(\forall a,b)(a\oplus b\cong b\oplus a).
\]
\end{definition}

This leads, for any finite sequence $\vv<a_i\mid i<n>$ of
elements of $A$ and for any element $a$ of $A$, to the
nonambiguous notations $\oplus_{i<n}a_i\dnw=a$ and
$\oplus_{i<n}a_i\dnw$ (we shall put
$\oplus_{i\leq n}a_i=(\oplus_{i<n}a_i)\oplus a_n$).
These notations are in general \emph{non-commutative}
($a\oplus b\dnw$ is not the same as
$b\oplus a\dnw$). There are many natural examples of partial
semigroups as defined above, but
also many examples of ``semigroup-like'' partial binary
operations that nevertheless fail to be partial
semigroups in the above sense. For
example, if $G$ is a group, define
$a\oplus b$ to be the group product
$ab$, if $a$ and $b$ \emph{commute}; or, if $L$ is a lattice
with $0$, define $a\oplus b=a\vee b$, if $a\wedge b=0$. For
both examples, associativity fails in general;
nevertheless, it holds in the second case if $L$ is
\emph{modular}%
\index{lattice!modular (not necessarily complemented) ---}
(this is well known; see also Proposition~\ref{P:oplus}).

\begin{definition}\label{D:refining}
A binary relation
$\sim$ on a partial semigroup\index{semigroup!partial ---}
$\vv<S,\oplus>$ is \emph{left refining}%
\index{refining (binary relation)!left ---|ii}, if it
satisfies the statement
\[
(\forall x_0,x_1,y)\bigl(y\sim x_0\oplus x_1\Rightarrow
(\exists y_0,y_1)(y_0\sim x_0\text{ and }y_1\sim x_1
\text{ and }y=y_0\oplus y_1)\bigr);
\]
\emph{right refining}%
\index{refining (binary relation)!right ---|ii}, if the dual
relation is left refining, and \emph{refining}%
\index{refining (binary relation)|ii}, if it is left
refining and right refining.

Furthermore, we say that $\sim$ is \emph{additive}%
\index{additive (binary relation)|ii}, if for all elements
$x$, $x'$, $y$, and $y'$ of $S$, if $x\sim x'$ and
$y\sim y'$ and both $x\oplus y$ and $x'\oplus y'$ are defined,
then $x\oplus y\sim x'\oplus y'$.

A \emph{\rps}\index{semigroup!refined partial ---|ii}
\(\vv<S,\oplus,\approx>\) is a
partial semigroup\index{semigroup!partial ---} $\vv<S,\oplus>$
endowed with a refining equivalence relation $\approx$ on $S$.
\end{definition}

If $\vv<S,\oplus>$ is a partial semigroup%
\index{semigroup!partial ---} and $\sim$ is a binary relation
on $S$, then a
\emph{$\sim$-refinement matrix}\index{refinement!matrix|ii} is
an array of the form
\begin{equation}\label{Eq:RefMatr}
\begin{tabular}{|c|c|c|c|c|}
\cline{2-5}
\multicolumn{1}{l|}{} & $b_0$ & $b_1$ & $\ldots$ &
$b_{n-1}$\tvi\\
\hline
$a_0$\tvi & $c_{00}/d_{00}$ & $c_{01}/d_{01}$ & $\ldots$ &
$c_{0,n-1}/d_{0,n-1}$\\
\hline
$a_1$\tvi & $c_{10}/d_{10}$ & $c_{11}/d_{11}$ & $\ldots$ &
$c_{1,n-1}/d_{1,n-1}$\\
\hline
$\vdots$ & $\vdots$ & $\vdots$ & $\ddots$ & $\vdots$\\
\hline
$a_{m-1}$\tvi & $c_{m-1,0}/d_{m-1,0}$ & $c_{m-1,1}/d_{m-1,1}$ &
$\ldots$ & $c_{m-1,n-1}/d_{m-1,n-1}$\\
\hline
\end{tabular}
\end{equation}
where $m$, $n\in\omega\setminus\{0\}$, the $a_i$'s, the
$b_j$'s, the $c_{ij}$'s and the $d_{ij}$'s are elements of $S$
such that the equalities $a_i=\oplus_{j<n}c_{ij}$,
$b_j=\oplus_{i<m}d_{ij}$, and $c_{ij}\sim d_{ij}$ hold for
all values of $i$ and $j$ in the corresponding domains.

If $m$ and $n$ are positive integers, if
$\vv<S,\oplus>$ is a partial semigroup%
\index{semigroup!partial ---} and $\sim$ is a binary relation
on $S$, then we shall say that
$\vv<S,\oplus,\sim>$ satisfies the
\emph{refinement\index{refinement!property|ii} property of
order $\vv<m,n>$}, if for all finite sequences
$\vv<a_i\mid i<m>$ and $\vv<b_j\mid j<n>$ of
elements of $S$ such that \(\oplus_{i<m}a_i=\oplus_{j<n}b_j\),
one can form a \(\sim_{(m-1)(n-1)}\)
refinement matrix of format
(\ref{Eq:RefMatr}) (we denote by $\sim_k$%
\index{Pzzersp2k@$\sim_k$|ii} the $k$th power of
$\sim$ as a binary relation, so that $\sim_0=\mathrm{id}_S$,
$\sim_1=\sim$ and $\sim_{k+1}=\sim\circ\sim_k$).

Furthermore, $\vv<S,\oplus,\sim>$ satisfies the
refinement property if and only
if it satisfies the refinement
property of every order $\vv<m,n>$, for
$m$, $n\in\omega\setminus\{0\}$. We say that a partial semigroup
\index{semigroup!partial ---}
$\vv<S,\oplus>$ satisfies the
refinement property if
$\vv<S,\oplus,=>$ satisfies the
refinement property (this definition
is compatible with the usual concept of
\emph{refinement\index{monoid!refinement ---} monoid},
that is, a \cm\ satisfying the
refinement property).

We shall encounter later in this chapter the following concept:

\begin{definition}\label{D:dimsgrp}
If $\vv<S,\oplus>$ is a
partial semigroup\index{semigroup!partial ---} and $\sim$ is a
binary relation on $S$, then the
\emph{dimension semigroup of} $\vv<S,\oplus,\sim>$, which we
shall denote by $\DD(S,\oplus,\sim)$%
\index{dimension semigroup of a \rps|ii}%
\index{DzzS@$\DD(S,\oplus,\sim)$ (for \rps s)|ii}, is the
commutative semigroup\index{semigroup} defined by generators
$|a|$ (for $a\in S$) and relations $|a\oplus b|=|a|+|b|$ (for all $a$,
$b$ in $S$ such that
$a\oplus b\dnw$) and $|a|=|b|$ (for all $a\sim b$ in $S$).
\end{definition}

The proof of the following lemma is about the same as the
classical one for commutative (total)
semigroups\index{semigroup}. We reproduce it here for convenience.

\begin{lemma} \label{L:refppty} Let $\vv<S,\oplus>$ be a
partial semigroup\index{semigroup!partial ---} and let $\sim$
be a refining binary relation on
$S$. Then $\vv<S,\oplus,\sim>$ satisfies the
refinement property if and only if
it satisfies the refinement property
of order $\vv<2,2>$.
\end{lemma}

\begin{proof} Let us prove the nontrivial direction. Suppose
that $\vv<S,\oplus,\sim>$ satisfies the
refinement\index{refinement!property} property of order
$\vv<2,2>$. We prove by induction on $m+n$ that
$\vv<S,\oplus,\sim>$ satisfies the
refinement\index{refinement!property} property of order
$\vv<m,n>$. If $m=1$ or $n=1$, this is trivial. If $m=n=2$,
then this is
the hypothesis. Thus let us, for example, prove the
refinement\index{refinement!property} property of order
$\vv<m+1,n>$, for $m\geq2$.
Thus let $a_i$ (for $i\leq m$) and $b_j$ (for $j<n$) in $S$
such that $\oplus_{i\leq m}a_i=\oplus_{j<n}b_j$. This can be
written as $a\oplus a_m=\oplus_{j<n}b_j$, where
$a=\oplus_{i<m}a_i$, so that by the induction hypothesis,
there exists a $\sim_{n-1}$-refinement%
\index{refinement!matrix} matrix of the form
\[
\begin{tabular}{|c|c|}
\cline{2-2}
\multicolumn{1}{l|}{} & $b_j\,(j<n)$\tvi\\
\hline
$a$ & $b'_j/b''_j$\tvi\\
\hline
$a_m$ & $c_{mj}/d_{mj}$\tvi\\
\hline
\end{tabular}
\]
Now, since $\oplus_{i<m}a_i=\oplus_{j<n}b'_j$ and by the
induction hypothesis, there exists a
$\sim_{(m-1)(n-1)}$-refinement\index{refinement!matrix} matrix
of the form
\[
\begin{tabular}{|c|c|}
\cline{2-2}
\multicolumn{1}{l|}{} & $b'_j\,(j<n)$\tvi\\
\hline
$a_i\,(i<m)$ & $c_{ij}/c'_{ij}$\tvi\\
\hline
\end{tabular}
\]
Thus, the equality $a_i=\oplus_{j<n}c_{ij}$ holds for all $i\leq m$. On the
other hand, the equality
$\oplus_{i<m}c'_{ij}=b'_j\sim_{n-1}b''_j$ holds for all $j<n$;
thus, since $\sim_{n-1}$ is refining, there are $d_{ij}$ (for $i<m$)
such that $b''_j=\oplus_{i<m}d_{ij}$ and $c'_{ij}\sim_{n-1}d_{ij}$, for
all $i$. Since $c_{ij}\sim_{(m-1)(n-1)}c'_{ij}$, for all $i$,
$j$, we also have $c_{ij}\sim_{m(n-1)}d_{ij}$. Moreover, the equality
$b_j=b''_j\oplus d_{mj}=\oplus_{i\leq m}d_{ij}$ holds for
all $j<n$, thus the following array
\[
\begin{tabular}{|c|c|}
\cline{2-2}
\multicolumn{1}{l|}{} & $b_j\,(j<n)$\tvi\\
\hline
$a_i\,(i\leq m)$ & $c_{ij}/d_{ij}$\tvi\\
\hline
\end{tabular}
\]
is a $\sim_{m(n-1)}$-refinement\index{refinement!matrix}
matrix.
\end{proof}

\section[Presentation using refinements]%
{An alternative
presentation of the dimension monoid}
\label{S:PresRef}

In this section, until \ref{C:chardim}, let
$\vv<S,\oplus,\approx>$ be a
\rps\index{semigroup!refined partial ---} satisfying the
refinement\index{refinement!property} property. Let
$[\omega]^{<\omega}$%
\index{ozzmegfin@$[\omega]^{<\omega}$|ii} denote the set of all
finite subsets of $\omega$ and let $\mathbb{S}$%
\index{Szzbig@$\mathbb{S}$|ii} be the set of all maps
\(p\colon a\to S\) for some \emph{nonempty}
\(a\in[\omega]^{<\omega}\) (we write
\(a=\mathrm{dom}(p)\)). For every such $p$, denote by
$\Sigma p$\index{Szzigmap@$\Sigma p$|ii} the set of all
elements $x$ of $S$ such that there exists a bijection
\(\sigma\colon n\to a\), where \(n=|a|\), and elements
\(x_i\) (for \(i<n\)) of $S$ such that
\(x_i\approx p(\sigma(i))\), for all $i$, and
\(x=\oplus_{i<n}x_i\). Note that since $\approx$ is refining,
\(\Sigma p\) is closed under $\approx$; it may also be empty.

For $p$, $q$ in $\mathbb{S}$, we write $p\perp q$%
\index{pzzerp@$\perp$ (for \rps s)|ii}, if
$\mathrm{dom}(p)\cap\mathrm{dom}(q)=\varnothing$, and then,
$p\oplus q=p\cup q$. Note that $\mathbb{S}$, thus equipped, is a
partial commutative semigroup\index{semigroup!partial ---}.
For all $p$, $q$ in $\mathbb{S}$, we write $p\to q$, if
there exists a surjective map
\(\alpha\colon
\mathrm{dom}(q)\twoheadrightarrow\mathrm{dom}(p)\)
such that $p\to^\alpha q$%
\index{tzzoalpha@$\to^\alpha$|ii},
where $p\to^\alpha q$ is the statement
\((\forall i\in\mathrm{dom}(p))
\bigl(p(i)\in
\Sigma(q\res_{{\alpha}^{-1}\{i\}})\bigr)\). Intuitively,
$p\to q$ means that $q$ ``refines'' $p$ up to $\approx$ and
up to the ordering of the summands.

\begin{lemma}\label{L:basicto} The following statements hold:

\begin{itemize}
\item[\rm (i)] The relation $\to$ is a preordering of
$\mathbb{S}$.

\item[\rm (ii)] The relation $\to$ is compatible with
$\oplus$, that is, if $p_0\to q_0$, $p_1\to q_1$,
$p_0\perp p_1$, and $q_0\perp q_1$, then
$p_0\oplus p_1\to q_0\oplus q_1$.

\item[\rm (iii)] The relation $\to$ is right refining (as in
Definition~\ref{D:refining}).

\item[\rm (iv)] The relation $\to$ is \emph{confluent},
that is, if $p\to q$ and $p\to r$, then there exists $s$
such that $q\to s$ and $r\to s$.
\end{itemize}
\end{lemma}

\begin{proof} (i) It is trivial that $\to$ is reflexive. Now
let $p\to^\alpha q$ and $q\to^\beta r$; we shall prove that
\(p\to^{\alpha\circ\beta}r\). First, it is clear that
\(\alpha\circ\beta\colon \mathrm{dom}(r)\twoheadrightarrow
\mathrm{dom}(p)\). Now let
$i\in\mathrm{dom}(p)$. If we put
$m=|{\alpha}^{-1}\{i\}|$, then there exist a bijection
$\sigma\colon m\to{\alpha}^{-1}\{i\}$ and elements
\(x_j\approx q(\sigma(j))\) (for all $j<m$) such that
$p(i)=\oplus_{j<m}x_j$. Similarly, for every $j<m$, put
\(n_j=|{\beta}^{-1}\{\sigma(j)\}|\); then there exist a
bijection
\(\tau_j\colon n_j\to{\beta}^{-1}\{\sigma(j)\}\) and
elements
\(y_{jk}\approx r(\tau_j(k))\) (for all $k<n_j$) such that
\(q(\sigma(j))=\oplus_{k<n_j}y_{jk}\). Since $\approx$ is
refining, for all $j<m$, there are elements
\(x_{jk}\approx y_{jk}\) (for all $k<n_j$) such that
\(x_j=\oplus_{k<n_j}x_{jk}\). Therefore,
\(p(i)=\oplus_{j<m}\oplus_{k<n_j}x_{jk}\) and the right-hand
side belongs to
\(\Sigma(r\res_{(\alpha\circ\beta)^{-1}\{i\}})\). Hence
\(p\to^{\alpha\circ\beta}r\), which proves the transitivity
of $\to$, thus completing the proof of (i).

Next, under hypothesis of (ii), if
\(p_i\to^{\alpha_i}q_i\) (for all $i<2$), then it is easy to see
that if we put \(\alpha=\alpha_0\cup\alpha_1\), then
\(p_0\oplus p_1\to^\alpha q_0\oplus q_1\). This proves (ii).

Now let \(p_0\oplus p_1\to^\alpha q\). Let $q_i$ (for $i<2$) be
the restrictions of $q$ of respective domains
\(\mathrm{dom}(q_i)={\alpha}^{-1}[\mathrm{dom}(p_i)]\)
(for $i<2$); put
\(\alpha_i=\alpha\res_{\mathrm{dom}(q_i)}\). Then
\(q=q_0\oplus q_1\) and \(p_i\to^{\alpha_i}q_i\): hence
$\to$ is right refining, thus proving (iii).

Finally we prove (iv). Let us do it first when
$|\mathrm{dom}(p)|=1$, say $p(i)=a$, where $i$ is the unique
element of
$\mathrm{dom}(p)$. Thus, putting $m=|\mathrm{dom}(q)|$ and
$n=|\mathrm{dom}(r)|$, there are bijections
\(\sigma\colon m\to\mathrm{dom}(q)\) and
\(\tau\colon n\to\mathrm{dom}(r)\) and elements \(x_j\approx
q(\sigma(j))\) (for all $j<m$) and
\(y_k\approx r(\tau(k))\) (for all $k<n$) such that
\(a=\oplus_{j<m}x_j=\oplus_{k<n}y_k\). By the
refinement\index{refinement!property} property, there exists a
$\approx$-refinement\index{refinement!matrix} matrix of the
form
\[
\begin{tabular}{|c|c|}
\cline{2-2}
\multicolumn{1}{l|}{} & $y_k\,(k<n)$\tvi\\
\hline
$x_j\,(j<m)$ & $z_{jk}/z'_{jk}$\tvi\\
\hline
\end{tabular}
\]

The relation \(q(\sigma(j))\approx\oplus_{k<n}z_{jk}\) holds
for all $j<m$, thus there are
elements $u_{jk}\approx z_{jk}$ such that
\(q(\sigma(j))=\oplus_{k<n}u_{jk}\); similarly,
the relation
\(r(\tau(k))\approx\oplus_{j<m}z'_{jk}\) holds for all $k<n$,
thus there are
elements \(v_{jk}\approx z'_{jk}\) such that
\(r(\tau(k))=\oplus_{j<m}v_{jk}\). Thus let $s\colon mn\to S$
be defined by the rule \(s(nj+k)=u_{jk}\) (for all $j<m$, $k<n$).
Let
\(\alpha\colon
\mathrm{dom}(s)\twoheadrightarrow\mathrm{dom}(q)\) be defined
by the rule \(\alpha(nj+k)=\sigma(j)\), and, similarly, let
\(\beta\colon
\mathrm{dom}(s)\twoheadrightarrow\mathrm{dom}(r)\) be defined
by the rule \(\beta(nj+k)=\tau(k)\). Then the equality
\(q(\sigma(j))=\oplus_{k<n}u_{jk}\in
\Sigma(s\res_{{\alpha}^{-1}\{\sigma(j)\}})\)
holds for all $j<m$,
while \(r(\tau(k))=\oplus_{j<m}v_{jk}\in
\Sigma(s\res_{{\beta}^{-1}\{\tau(k)\}})\) for all $k<n$, so
that \(q\to^\alpha s\) and \(r\to^\beta s\).

In the general case, let \(p\to q\) and \(p\to r\). Then there
exists a decomposition \(p=\oplus_{i<m}p_i\), where
\(m\in\omega\setminus\{0\}\) and each \(p_i\) has a one-element
domain; by (iii), there are corresponding decompositions
\(q=\oplus_{i<n}q_i\) and \(r=\oplus_{i<n}r_i\)
such that \(p_i\to q_i\) and \(p_i\to r_i\), for
all \(i<n\). Therefore, by the
previous result, there exists \(s_i\) such that \(q_i\to s_i\)
and \(r_i\to s_i\); without loss of generality, the domains of
the \(s_i\)'s are mutually disjoint. Hence, by (ii), if
\(s=\oplus_{i<m}s_i\), we have \(q\to s\) and \(r\to s\).
\end{proof}

\begin{note} The preordering $\to$ is not antisymmetric in
general: for example, for every \(p\in\mathbb{S}\) and every
bijection $\sigma$ from a set in \([\omega]^{<\omega}\)
onto the domain of $p$, then we have \(p\to p\circ\sigma\to p\).
\end{note}

Now, we define a binary relation $\equiv$ on $\mathbb{S}$ by
putting
\[
p\equiv q\Longleftrightarrow (\exists r\in\mathbb{S})
(p\to r\text{ and }q\to r).
\]

\begin{lemma}\label{L:basicequiv} The relation $\equiv$ is an
equivalence on $\mathbb{S}$, compatible with $\oplus$. The
quotient structure $\vv<\mathbb{S},\oplus>/\equiv$ is a (total)
commutative semigroup\index{semigroup}.
\end{lemma}

\begin{proof}
The fact that $\equiv$ is an equivalence follows
from the fact that $\to$ is a confluent preordering. If
$p_0\perp p_1$, $q_0\perp q_1$, and $p_i\equiv q_i$ (for all
$i<2$), there exist $r_0$, $r_1$ such that $p_i\to r_i$ and
$q_i\to r_i$ (for all $i<2$); by replacing the domain of $r_1$
with one with the same cardinality disjoint from
$\mathrm{dom}(r_0)$, one may assume that $r_0\perp r_1$. By
Lemma~\ref{L:basicto},
\(p_0\oplus p_1\to r_0\oplus r_1\) and
\(q_0\oplus q_1\to r_0\oplus r_1\), whence
\(p_0\oplus p_1\equiv q_0\oplus q_1\); thus $\equiv$ is
compatible with $\oplus$. Denote by
$[p]$ the $\equiv$-equivalence class of $p$,
for all $p\in\mathbb{S}$. For all
$p$, $q\in\mathbb{S}$, there exists $q'\in\mathbb{S}$ such that
$q\equiv q'$ (and $|\mathrm{dom}(q)|=|\mathrm{dom}(q')|$) and
$p\perp q'$, so that $[p]+[q]=[p\oplus q']$ is defined, whence
the totality assertion. Since $\oplus$ is commutative on
$\mathbb{S}$, the commutativity assertion is trivial.
\end{proof}

We shall denote by $\hat S$%
\index{Szzhat@$\hat S$|ii} the commutative
semigroup\index{semigroup}
$\vv<\mathbb{S},\oplus>/\equiv$. There is a natural map
$\epsilon\colon S\to\hat S$ defined by the rule
$\epsilon(x)=[\vv<x>]$. Then note that the equality
\([p]=\sum_{x\in\mathrm{rng}(p)}\epsilon(x)\) holds for all
$p\in\mathbb{S}$, so that the range of
$\epsilon$ generates $\hat S$ as a semigroup\index{semigroup}.
Note also that
$\epsilon$ is $\approx$-invariant, that is, $x\approx y$
implies that $\epsilon(x)=\epsilon(y)$. In fact, $\hat S$ is the
universal such semigroup\index{semigroup} (as defined in
Definition~\ref{D:dimsgrp}):

\begin{proposition}\label{P:charS}
$\hat S=\DD(S,\oplus,\approx)$%
\index{dimension semigroup of a \rps}.
\end{proposition}

\begin{proof}
It suffices to show that, for every commutative
semigroup\index{semigroup} $\vv<T,+>$, and for every
$\approx$-invariant homomorphism
$f\colon\vv<S,\oplus>\to\vv<T,+>$, there exists a unique
semigroup\index{semigroup} homomorphism
$g\colon \hat S\to T$ such that
$g\circ\epsilon=f$. Necessarily, the equality
\(g([p])=\sum_{x\in\mathrm{rng}(p)}g\circ\epsilon(x)
=\sum_{x\in\mathrm{rng}(p)}f(x)\) holds for all $p\in\mathbb{S}$,
whence $g$ is unique. As to
existence, we must prove first that the formula above defines
a map $g\colon \hat S\to T$, that is, that $p\equiv q$
implies that
\(\sum_{x\in\mathrm{rng}(p)}f(x)=
\sum_{x\in\mathrm{rng}(q)}f(x)\). By the definition of $\equiv$,
it suffices to prove the same conclusion for
$p\to q$. However, since $T$ is commutative and $f$ is
$\approx$-invariant, this is clear. Then it is trivial that
$g$ is a semigroup\index{semigroup} homomorphism and that
$g\circ\epsilon=f$.
\end{proof}

In accordance with Definition~\ref{D:dimsgrp}, we shall write
$|x|$ instead of $\epsilon(x)$ for $x\in S$.

\begin{corollary}\label{C:chardim} Suppose that, in addition,
$\vv<S,\oplus>$ is commutative and $\approx$ is additive.
Then $|x|=|y|$ if and only if $x\approx y$, for all $x$, $y\in S$.
\end{corollary}

\begin{proof} We prove the nontrivial direction. So suppose
that $|x|=|y|$. By definition, there exists
$p\in\mathbb{S}$ such that $\vv<x>\to p$ and
$\vv<y>\to p$. This means, by definition, that both $x$ and
$y$ belong to $\Sigma p$,
that is, if $n=|\mathrm{dom}(p)|$, there are bijections
$\sigma$, $\tau\colon n\to\mathrm{dom}(p)$ and elements
$x_i\approx p(\sigma(i))$ and $y_i\approx p(\tau(i))$ (for $i<n$)
such that $x=\oplus_{i<n}x_i$ and $y=\oplus_{i<n}y_i$. In
particular, $x_{\sigma^{-1}(i)}\approx y_{\tau^{-1}(i)}$
for all $i<n$. But by
assumption, $\vv<S,\oplus>$ is associative and commutative,
thus $\oplus_{i<n}x_{\sigma^{-1}(i)}$ (resp.,
$\oplus_{i<n}y_{\tau^{-1}(i)}$) is defined and equal to $x$
(resp., to $y$). Moreover, since $\approx$ is by assumption
additive, we obtain that $x\approx y$.
\end{proof}

\begin{note} The hypothesis of Corollary~\ref{C:chardim} could
have been weakened into a ``commutativity modulo $\approx$'' of
$\oplus$, but we shall not need that version.
\end{note}

\begin{theorem}\label{T:refmon} Let $\vv<S,\oplus,\approx>$ be
a \rps\index{semigroup!refined partial ---}\ satisfying the
refinement\index{refinement!property} property. Then
$\DD(S,\oplus,\approx)^\circ$ also satisfies the
refinement\index{refinement!property} property.
\end{theorem}

\begin{proof} By Proposition~\ref{P:charS}, it suffices to
prove that $\hat{S}^\circ$ satisfies the
refinement\index{refinement!property} property. Thus let us try
to refine in $\hat{S}^\circ$ an equation of the form
$\xi_0+\xi_1=\eta_0+\eta_1$, where the $\xi_i$'s and the
$\eta_j$'s belong to $\hat{S}^\circ$. If one of the $\xi_i$'s
or $\eta_j$'s is equal to $\mathsf{O}$, then it is trivial that
one can do so, thus let us suppose that they all belong to
$\hat S$.  Let $p_i\in\xi_i$ and $q_i\in\eta_i$, for all $i<2$;
one can take
$p_0\perp p_1$ and $q_0\perp q_1$, so that
\(p_0\oplus p_1\equiv q_0\oplus q_1\), that is, there
exists
$r$ such that $p_0\oplus p_1\to r$ and $q_0\oplus q_1\to r$.
By Lemma~\ref{L:basicto}.(iii) ($\to$ is right refining),
there are elements $p'_i$ and $q'_i$ (for $i<2$) such that
$p_i\to p'_i$,
$q_i\to q'_i$, and \(r=p'_0\oplus p'_1=q'_0\oplus q'_1\). Put
$\zeta_{ij}=[p'_i\cap q'_j]$, for all $i$, $j<2$ (\emph{we put
$[\varnothing]=\mathsf{O}$, showing the utility of the
extra zero!}). Then the following is a
refinement\index{refinement!matrix} matrix in $\hat{S}^\circ$:
\[
\begin{tabular}{|c|c|c|}
\cline{2-3}
\multicolumn{1}{l|}{} & $\eta_0$ & $\eta_1$\tvi\\
\hline
$\xi_0$\tvi & $\zeta_{00}$ & $\zeta_{01}$\\
\hline
$\xi_1$\tvi & $\zeta_{10}$ & $\zeta_{11}$\\
\hline
\end{tabular}
\]
which concludes the proof.
\end{proof}

However, if $S$ is a singleton and $\oplus$ is the nowhere
defined partial binary operation on $S$, then
$\vv<S,\oplus,=>$ is a \rps\index{semigroup!refined partial
---} satisfying the refinement\index{refinement!property}
property, while
$\DD(S,\oplus,=)\cong\NN$ does
\emph{not} satisfy the refinement\index{refinement!property}
property (one cannot refine in $\NN$ the equation
$1+1=1+1$). As we shall see in Corollary~\ref{C:commquot}, this
cannot happen for \emph{total} semigroups\index{semigroup}.

\section{Further properties of the dimension function}

Our next lemma will make it possible, in particular, to transfer the
refinement\index{refinement!property} property from
$\DD(S,\oplus,\approx)^\circ$ to
$\DD(S,\oplus,\approx)$.

\begin{lemma}\label{L:refwoO} Let $S$ be a commutative
semigroup\index{semigroup} satisfying the following axiom:
\[
(\forall a,b)(\exists u,v,c) (a=u+c\text{ and }b=c+v).
\]
If $S^\circ$ satisfies the
refinement\index{refinement!property} property, then so does
$S$.
\end{lemma}

For example, if $S$ has a zero element, then the refinement
property\index{refinement!property} of $S$ follows from the
refinement\index{refinement!property} property of
$S^\circ$, but a direct proof is much easier in this case.
\smallskip

\begin{proof}
Let $a_i$, $b_i\in S$ (for $i<2$) such that
$a_0+a_1=b_0+b_1$. If $S^\circ$ satisfies the
refinement\index{refinement!property} property, then there is a
refinement\index{refinement!matrix} matrix of the form
\[
\begin{tabular}{|c|c|c|}
\cline{2-3}
\multicolumn{1}{l|}{} & $b_0$ & $b_1$\tvi\\
\hline
$a_0$\tvi & $c_{00}$ & $c_{01}$\\
\hline
$a_1$\tvi & $c_{10}$ & $c_{11}$\\
\hline
\end{tabular}
\]
where the $c_{ij}$'s belong to $S^\circ$. If the $c_{ij}$'s
all belong to $S$, then we are done. Now suppose that
$c_{ij}=\mathsf{O}$, for some $\vv<i,j>$. Without loss of
generality,
$c_{01}=\mathsf{O}$, so that we have a
refinement\index{refinement!matrix} matrix
\[
\begin{tabular}{|c|c|c|}
\cline{2-3}
\multicolumn{1}{l|}{} & $b_0$ & $b_1$\tvi\\
\hline
$a_0$\tvi & $a_0$ & $\mathsf{O}$\\
\hline
$a_1$\tvi & $c_{10}$ & $b_1$\\
\hline
\end{tabular}
\]
where \(c_{10}\in S^\circ\). Now, by the assumption on $S$,
there are elements $u$, $v$, and $c$ of $S$ such that $a_0=u+c$
and $b_1=c+v$. Thus the following
 \[
 \begin{tabular}{|c|c|c|}
 \cline{2-3}
 \multicolumn{1}{l|}{} & $b_0$ & $b_1$\tvi\\
 \hline
 $a_0$\tvi & $u$ & $c$\\
 \hline
 $a_1$\tvi & $c+c_{10}$ & $v$\\
 \hline
 \end{tabular}
 \]
is a refinement\index{refinement!matrix} matrix in
$S$.\end{proof}

The following lemma will be used in Chapter~\ref{ModLatt}.

\begin{lemma}\label{L:Vmeas} Let $\vv<S,\oplus,\approx>$ be a
\rps\index{semigroup!refined partial ---} satisfying the
refinement\index{refinement!property} property. Let $c\in S$
and let $\xi$, $\eta\in\DD(S,\oplus,\approx)$ such that
$|c|=\xi+\eta$. Then there exist a decomposition
$c=\oplus_{i<k}c_i$ and subsets $I$ and $J$ of
$k$ such that $I\cap J=\varnothing$,
$I\cup J=k$, \(\xi=\sum_{i\in I}|c_i|\), and
\(\eta=\sum_{i\in J}|c_i|\).
\end{lemma}

\begin{proof} Let \(p=\vv<a_i\mid i<m>\) and
\(q=\vv<b_i\mid m\leq i<m+n>\) be such that $\xi=[p]$ and
$\eta=[q]$. Then we have \(\vv<c>\equiv p\oplus q\),
that is, there exist an element \(r\in\mathbb{S}\) and maps
$\alpha$, $\beta$ such that
\(\vv<c>\to^\alpha r\) and \(p\oplus q\to^\beta r\). Let
\(k=|\mathrm{dom}(r)|\). Then there are bijections
$\sigma$, $\tau\colon k\to\mathrm{dom}(r)$ and elements
\(c_j\approx r(\sigma(j))\) and
\(d_j\approx r(\tau(j))\) of $S$ (for all $j<k$) such that
\(c=\oplus_{j<k}c_j\), and
\((p\oplus q)(i)=\oplus_{j\in{\beta}^{-1}\{i\}}d_j\), for all
$i<m+n$. Note that
\(d_j\approx r(\tau(j))\approx c_{\sigma^{-1}\tau(j)}\) for all $j<k$,
whence \(|d_j|=|c_{\sigma^{-1}\tau(j)}|\). It follows that
\begin{align*}
\xi=\sum_{i<m}|p\oplus q(i)|&=
\sum_{i<m}\left|
\oplus_{j\in{\beta}^{-1}\{i\}}d_j\right|\\
&=\sum_{i<m}\sum_{j\in{\beta}^{-1}\{i\}}|d_j|\\
&=\sum_{j\in{\beta}^{-1}\{0,1,\ldots,m-1\}}
|c_{\sigma^{-1}\tau(j)}|,
\end{align*}

while

\begin{align*}
\eta=\sum_{m\leq i<m+n}|p\oplus q(i)|&=
\sum_{m\leq i<m+n}
\left|\oplus_{j\in{\beta}^{-1}\{i\}}d_j\right|\\
&=\sum_{m\leq i<m+n}\sum_{j\in{\beta}^{-1}\{i\}}|d_j|\\
&=\sum_{j\in{\beta}^{-1}\{m,m+1,\ldots,m+n-1\}}
|c_{\sigma^{-1}\tau(j)}|.
\end{align*}

Therefore, the conclusion holds for $I$ and $J$ defined by
\begin{align*}
I=\sigma^{-1}
\tau[{\beta}^{-1}\{0,1,\ldots,m-1\}],\\
J=\sigma^{-1}\tau[{\beta}^{-1}\{m,m+1,\ldots,m+n-1\}].
\tag*{\qed}
\end{align*}
\renewcommand{\qed}{}\end{proof}

\begin{example}\label{E:decomptypes}
Let a group $G$ act by
automorphisms on a generalized Boo\-le\-an algebra
$\mathcal{B}$ (see Section~\ref{S:DistrLatt}).
In $\mathcal{B}$, we write
\(a\oplus b=a\vee b\), if
\(a\wedge b=0\). Moreover, let $\approx$ be the orbital
equivalence associated with $G$. Then
\(\vv<\mathcal{B},\oplus,\approx>\) is a
\rps\index{semigroup!refined partial ---} and it is very easy
to see that it satisfies the
refinement\index{refinement!property} property. Thus, the
dimension semigroup\index{semigroup}
\(\DD(S,\oplus,\approx)\) is a
refinement\index{monoid!refinement ---} monoid. It is called, in
general, the \emph{monoid of equidecomposability types of
elements of $\mathcal{B}$ modulo $G$}, see
\cite{Wago85,Wehr90}. Note that the known proofs of refinement
in this monoid use either an extension of the Stone space
\(\Omega\) of $\mathcal{B}$ (as shown in
\cite[Chapter 8]{Wago85}) or the
refinement\index{refinement!property} property on the monoid
of bounded continuous \(\ZZ^+\)-valued functions on $\Omega$,
see \cite{Wehr90}, while on the other hand, the proof of the
fact that
\(\vv<\mathcal{B},\oplus>\) satisfies refinement is trivial.
Moreover, some generalizations of this proof to other
situations, as, for example, when $G$, instead of being a
group, is an
\emph{inverse semigroup}\index{semigroup!inverse ---} (for
example, partial isometries on a metric space, see
\cite{Tars49,Wago85}), are trivial as well.
\end{example}

\begin{example}\label{E:brookfield} Let $R$ be any ring.
Denote by \(R{\mbox{\bf-Noeth}}\)%
\index{RzzNoeth@$R{\mbox{\bf-Noeth}}$|ii} the category of all
Noetherian left $R$-modules. As in
\cite{Broo}, we say that finite submodule series of same
length, say \(A_0\subseteq A_1\subseteq\cdots\subseteq A_n\) and
\(B_0\subseteq B_1\subseteq\cdots\subseteq B_n\), of elements
of \(R{\mbox{\bf-Noeth}}\) are \emph{isomorphic}, if there
exists a permutation $\sigma$ of $n$ such that \(A_{i+1}/A_i\cong
B_{\sigma(i)+1}/B_{\sigma(i)}\) for all
\(i<n\).
If, in addition, \(A_0=\{0\}\) and \(B_0=\{0\}\), then say
that $A_n$ and $B_n$
\emph{have isomorphic submodule series} and write
\(A_n\sim B_n\). Then, in
\cite{Broo}, G. Brookfield endows the set
\(M(R{\mbox{\bf-Noeth}})\) of all
$\sim$-equivalence classes of elements of
\(R{\mbox{\bf-Noeth}}\) with an addition defined by
\([A]+[B]=[A\oplus B]\) (\([A]\) denotes the $\sim$-equivalence
class of $A$) and proves, among other things, that
\(M(R{\mbox{\bf-Noeth}})\) is a
refinement\index{monoid!refinement ---} monoid. In this
particular case, his approach is closely related to ours, in
the sense that it is, roughly speaking, equivalent to considering
the partial semigroup\index{semigroup!partial ---} $S$ of all
ordered pairs
\(\vv<A,B>\) of elements of \(R{\mbox{\bf-Noeth}}\) such that
\(A\subseteq B\), endowed with the (partial) addition $\oplus$
defined by
\[
\vv<A,B>\oplus\vv<C,D>=\vv<A,D>\mbox{ provided that }B=C
\]
and the equivalence relation $\approx$ defined by
\[
\vv<A,B>\approx\vv<C,D>\Longleftrightarrow B/A\cong D/C,
\]
so that
\(M(R{\mbox{\bf-Noeth}})=\DD(S,\oplus,\approx)\) is a
refinement\index{monoid!refinement ---} monoid. One difference
is that the elements of
\(\DD(S,\oplus,\approx)\) can already be
defined as quotients modulo $\sim$ of elements of
\(R{\mbox{\bf-Noeth}}\), essentially because if $A$ and $B$ are
any elements of \(R{\mbox{\bf-Noeth}}\), then there always exists an
element $C$ of \(R{\mbox{\bf-Noeth}}\) such that
\([C]=[A]+[B]\); this consideration simplifies the proof to
some extent, as, for example, the introduction of
$\mathbb{S}$ is no longer needed.
\end{example}

To conclude this chapter, let us mention another corollary of
Lemma \ref{L:refwoO}:

\begin{corollary}\label{C:commquot} Let $S$ be a (total)
semigroup\index{semigroup} satisfying the
refinement\index{refinement!property} property. Then the
maximal commutative quotient of $S$ also satisfies the
refinement\index{refinement!property} property.
\end{corollary}

\begin{proof} By Theorem~\ref{T:refmon} and
Lemma~\ref{L:refwoO}, it suffices to verify that the maximal
commutative quotient $T$ of $S$ satisfies the hypothesis of
Lemma~\ref{L:refwoO}. Since $T$ is a homomorphic image of
$S$, it suffices to see that $S$ satisfies this axiom. Thus
let $a$, $b\in S$. Since $S$ satisfies the
refinement\index{refinement!property} property, there is a
refinement\index{refinement!matrix} matrix in $S$ of the form

\[
\begin{tabular}{|c|c|c|}
\cline{2-3}
\multicolumn{1}{l|}{} & $a$ & $b$\tvi\\
\hline
$a$\tvi & $c_{00}$ & $c_{01}$\\
\hline
$b$\tvi & $c_{10}$ & $c_{11}$\\
\hline
\end{tabular}
\]
so that we obtain the conclusion for \(u=c_{00}\),
\(c=c_{01}\), \(v=c_{11}\).\end{proof}

\chapter{Dimension monoids of modular%
\index{lattice!modular (not necessarily complemented) ---}
lat\-tices}\label{ModLatt}

In this chapter, we shall see that the description of dimension
monoids of \emph{modular} lattices turns out to be, with the help
of the tools introduced in Chapter~\ref{MaxCommQuot}, a rather easy
task.

\section{The associated refined partial semigroup%
\index{semigroup!refined partial ---}}

For every lattice $L$, we shall endow the set
$\S(L)$ with the structure of a partial
semigroup\index{semigroup!partial ---}, defined as follows:
\index{ozzplus@$\oplus$!on intervals of a lattice|ii}
\[
[a,\,b]\oplus[c,\,d]=[a,\,d]\qquad
\text{in case}\ b=c.
\]

Note that $\oplus$ is very far from being commutative. Then
the following result holds trivially:

\begin{lemma}\label{L:charV(L)} The semigroup\index{semigroup}
$\DD L$ is isomorphic to the
dimension semigroup
$\DD(\S(L),\oplus,\approx)$.\qed
\end{lemma}

Note that $\DD L$ always has a
zero element, thus, by Lemma~\ref{L:refwoO}, to find sufficient
conditions for
$\DD L$ to satisfy the
refinement\index{refinement!property} property, it suffices to
find sufficient conditions for
$\vv<\S(L),\oplus,\approx>$ to satisfy the hypotheses of
Theorem~\ref{T:refmon}, that is, that
$\vv<\S(L),\oplus,\approx>$ is a \rps%
\index{semigroup!refined partial ---} satisfying the
refinement\index{refinement!property} property. Most of this
statement is, in fact, always true:

\begin{lemma}
Let $L$ be a lattice. Then
$\vv<\S(L),\oplus>$ is a partial
semigroup\index{semigroup!partial ---} satisfying the
$\approx$-refinement\index{refinement!property} property.
\end{lemma}

\begin{proof} It is obvious that $\vv<\S(L),\oplus>$ is a
partial semigroup\index{semigroup!partial ---}. Now let us
verify the refinement\index{refinement!property} property.
Thus let
$a\leq c_i\leq b$ (for all $i<2$); we try to refine the identity
\([a,\,c_0]\oplus[c_0,\,b]=[a,\,c_1]\oplus[c_1,\,b]\). But it
is obvious that the following is always a
$\sim$-refinement\index{refinement!matrix} matrix,
\[
\begin{tabular}{|c|c|c|}
\cline{2-3}
\multicolumn{1}{l|}{} & $[a,\,c_1]$ & $[c_1,\,b]$\tvi\\
\hline
$[a,\,c_0]$\tvi & $[a,\,c_0\wedge c_1]/[a,\,c_0\wedge c_1]$ &
$[c_0\wedge c_1,\,c_0]/[c_1,\,c_0\vee c_1]$\\
\hline
$[c_0,\,b]$\tvi & $[c_0,\,c_0\vee c_1]/[c_0\wedge c_1,\,c_1]$
&
$[c_0\vee c_1,\,b]/[c_0\vee c_1,\,b]$\\
\hline
\end{tabular}
\]
thus, \emph{a fortiori}, a
$\approx$-refinement\index{refinement!matrix} matrix.
\end{proof}

Note that this is the way the classical Schreier%
\index{refinement!Schreier ---} Refinement Theorem is proved.

Now, modularity%
\index{lattice!modular (not necessarily complemented) ---} is
used in the proof of the following lemma:

\begin{lemma}\label{L:projref} Let $L$ be a modular%
\index{lattice!modular (not necessarily complemented) ---}
lattice. Then $\approx$ is a
refining equivalence relation on $\S(L)$.
\end{lemma}

\begin{proof} It suffices to prove that
$\nearrow$ is both left and right
refining. Thus let
$[a,\,b]\nearrow[a',\,b']$ in
$\S(L)$. If $a\leq c\leq b$, then
$[a,\,c]\nearrow[a',\,c']$, where
$c'=a'\vee c$, and, by modularity%
\index{lattice!modular (not necessarily complemented) ---},
$c=c'\wedge b$ so that
$[c,\,b]\nearrow[c',\,b']$: thus
$\nearrow$ is right refining. Similarly, if $a'\leq c'\leq
b'$, then
$[c,\,b]\nearrow[c',\,b']$, where
$c=c'\wedge b$, and, by modularity%
\index{lattice!modular (not necessarily complemented) ---},
$c'=a\vee c$ so that
$[a,\,c]\nearrow[a',\,c']$; thus
$\nearrow$ is left refining.
\end{proof}

Putting together \ref{L:charV(L)}--\ref{L:projref},
Theorem \ref{T:refmon}, and Lemma \ref{L:refwoO}, we obtain
immediately the following result:

\begin{theorem}\label{T:V(L)ref} Let $L$ be a modular%
\index{lattice!modular (not necessarily complemented) ---}
lattice. Then $\DD L$ is a
\crm\index{monoid!conical refinement ---}.\qed
\end{theorem}

In fact, Chapter~\ref{MaxCommQuot} gives a much more precise
information about \(\DD L\). Indeed, it gives a test for the
equality of any ``dimension words''
\(\xi=\sum_{i<m}\Dim(a_i,b_i)\) (with \(a_i\leq b_i\)) and
\(\eta=\sum_{j<n}\Dim(c_j,d_j)\) (with \(c_j\leq d_j\)): the
criterion studied in Chapter~\ref{MaxCommQuot} gives that
\(\xi=\eta\) if and only if there exists $\zeta$ such that
\(\xi\to\zeta\) and \(\eta\to\zeta\) (the $\to$ relation was
introduced in Chapter~\ref{MaxCommQuot}).

Reformulated in lattice-the\-o\-ret\-i\-cal terms, this gives
the following: if\linebreak
\(p=\vv<[a_i,\,b_i]\mid i\in I>\) and
\(q=\vv<[c_j,\,d_j]\mid j\in J>\) are finite
sequences of elements of \(\S(L)\), we say that $p$ and $q$ are
\emph{isomorphic}, if there exists a bijection
\(\sigma\colon I\to J\) such that
\([a_i,\,b_i]\approx[c_{\sigma(i)},\,d_{\sigma(i)}]\) for all
\(i\in I\). Then, in the general case, we obtain that
\(\sum_{i<m}\Dim(a_i,b_i)=\sum_{j<n}\Dim(c_j,d_j)\) if and
only if the sequences \(\vv<[a_i,\,b_i]\mid i<m>\) and
\(\vv<[c_j,\,d_j]\mid j<n>\) have isomorphic refinements.

In fact, this alone makes it possible to prove
Theorem~\ref{T:V(L)ref}: one defines formally the equality of
two dimension words as above, and then one verifies that this
definition of equality is transitive (thus proving the analogue
of Lemma~\ref{L:basicequiv} and then
Proposition~\ref{P:charS}), by using a mild extension of the
Schreier Refinement%
\index{refinement!Schreier ---} Theorem to finite dimension
words instead of mere intervals. Then, the proof of refinement
goes as in the proof of Theorem~\ref{T:refmon} followed by
Lemma~\ref{L:refwoO}. All the technical details put together,
although not fundamentally difficult, take some space; thus
our choice, in this work, is to present the general argument of
Chapter~\ref{MaxCommQuot}, more likely to be immediately
generalizable to new situations without any further
calculation (see Examples
\ref{E:decomptypes} and \ref{E:brookfield} or the case of the
$\oplus$ operation in the sectionally complemented modular
case from Chapter~\ref{RelCompl} on, indeed used in this work).

\section[Modularity and dimension]%
{Links between modularity%
\index{lattice!modular (not necessarily complemented) ---} and
the dimension monoid}

As the results of this section will show, there are certain classes
of lattices, such as lattices with finite height, where modularity
can be recognized on the dimension monoid.
More elaborate results of this
kind can be found, for example, in Chapter~\ref{ShortLatt}.

\begin{proposition}\label{P:shortmod} Let $L$ be a modular%
\index{lattice!modular (not necessarily complemented) ---}
lattice without infinite bounded chains. Let $P$ be the set  of
all projectivity%
\index{projectivity!of intervals} classes of prime intervals of
$L$. Then $\DD L$ is isomorphic
to $(\ZZ^+)^{(P)}$, the free commutative monoid on $P$.
\end{proposition}

\begin{proof} If $\xi$ is a projectivity%
\index{projectivity!of intervals} class of prime intervals of
$L$ and if $a\leq b$ in
$L$, then the number of occurrences of an interval in $\xi$ in
any maximal chain from $a$ to $b$ (finite by assumption) is
independent on the chosen chain (this follows, for example, from
the Schreier\index{refinement!Schreier ---} Refinement
Theorem); denote it by $|a,b|_\xi$. Then it is easy to see that
$\vv<a,b>\mapsto|a,b|_\xi$ satisfies
(D0), (D1) and (D2), thus there exists a unique monoid homomorphism
$\pi\colon\DD L\to(\ZZ^+)^{(P)}$ such that
$\pi(\Dim(a,b))=\vv<|a,b|_\xi\mid\xi\in P>$ for all
$a\leq b$ in $L$. Note that
$\pi$ is obviously surjective. Conversely, note that for all
$\xi\in P$, the value of $\Dim(a,b)$, for $[a,\,b]\in\xi$, does
not depend on the choice of $\vv<a,b>$; denote it by
$\rho(\delta_\xi)$ (where $\vv<\delta_\xi\mid \xi\in P>$
is the canonical basis of $(\ZZ^+)^{(P)}$). Then $\rho$
extends to a unique monoid homomorphism from $(\ZZ^+)^{(P)}$ to
$\DD L$, still denoted by
$\rho$. Then for every prime interval
$[a,\,b]$ of $L$ with projectivity%
\index{projectivity!of intervals} class $\xi$, we have
$\rho\pi(\Dim(a,b))=\rho(\delta_\xi)=\Dim(a,b)$. Therefore,
$\rho\pi(\Dim(a,b))=\Dim(a,b)$, for all $a\leq b$ in
$L$, so that $\rho\pi=\mathrm{id}_{\DD L}$. Thus
$\pi$ is an embedding.\end{proof}

Moreover, we see that, in particular, if $L$ is modular%
\index{lattice!modular (not necessarily complemented) ---}
without infinite bounded chains, then $\DD L$ is
\emph{cancellative}. Conversely, we have the following result,
that generalizes the well-known fact, see
\cite[Theorem X.2, page 232]{Birk93}, that ($\RR$-valued)
metric lattices are modular%
\index{lattice!modular (not necessarily complemented) ---}.
We say that an element $u$ of a \cm\ $M$ is
\emph{directly finite}%
\index{directly finite!element in a monoid|ii}, if $M$
satisfies the statement
$(\forall x)(x+u=u\Rightarrow x=0)$. (This is weaker than
to say that $M$ is cancellative.)

\begin{proposition}\label{P:dirfin} Let $L$ be a lattice such
that every element of the dimension range\index{dimension
range} of
$L$ is directly finite. Then $L$ is
modular%
\index{lattice!modular (not necessarily complemented) ---}.
\end{proposition}

\begin{proof}
Let $a\geq c$ and $b$ be elements of
$L$. Then the intervals $[a\wedge b,\,b]$ and
$[a\wedge(b\vee c),\,b\vee c]$ on the one hand, and
$[a\wedge b,\,b]$ and $[(a\wedge b)\vee c,\,b\vee c]$ on the
other, are tranposes, thus, using (D2) and
then (D1), we obtain that
\begin{align*}
\Dim(a\wedge b,b) &=\Dim(a\wedge(b\vee c),b\vee c)\\
&=\Dim((a\wedge b)\vee c,b\vee c)\\
&=\Dim((a\wedge b)\vee
c,a\wedge(b\vee c))+
\Dim(a\wedge(b\vee c),b\vee c),
\end{align*} whence, since the element \(\Dim(a\wedge(b\vee
c),b\vee c)\) is directly finite,
\(\Dim((a\wedge b)\vee c,a\wedge(b\vee c))=0\), thus
\(a\wedge(b\vee c)=(a\wedge b)\vee c\) by Proposition
\ref{P:V(L)cone}.
\end{proof}

In particular, if $L$ is a lattice without infinite bounded
chains, $L$ is modular%
\index{lattice!modular (not necessarily complemented) ---} if
and only if
$\DD L$ is cancellative. Of
course, there are many modular%
\index{lattice!modular (not necessarily complemented) ---}
lattices with non-cancellative dimension lattices (for example
the lattice of all linear subspaces of any infinite
dimensional vector space over any field).

Let us now turn to \Vhom s\index{Vhom@\Vhom}
(Definition~\ref{D:V-hom}). One consequence of
Chapter~\ref{MaxCommQuot} is the following result:

\begin{proposition}\label{P:exV-hom}
Let \(e\colon K\to L\) be a lattice homomorphism
with convex range, with $L$ modular%
\index{lattice!modular (not necessarily complemented) ---}.
Then $\DD e$ is a \Vhom\index{Vhom@\Vhom}.
\end{proposition}

\begin{proof} Put $f=\DD e$. Since $f$ is a monoid
homomorphism and since $\DD L$
satisfies the refinement\index{refinement!property} property
(see Theorem~\ref{T:V(L)ref}), the set of all elements of
$\DD K$ at which $f$ is a
\Vhom\index{Vhom@\Vhom} is by Lemma~\ref{L:V-homat} closed
under addition. Therefore, it suffices to verify that
$f$ is a \Vhom\index{Vhom@\Vhom} at all elements of $\DD K$ of
the form
$\Dim_K(u,v)$, where $u\leq v$ in
$K$. Thus let
$\xi$, $\eta\in\DD L$ such that
$f(\Dim_K(u,v))=\xi+\eta$,
that is, $\Dim_L(e(u),e(v))=\xi+\eta$. By
Lemma~\ref{L:Vmeas}, there exist $k\in\omega$, subsets $I$ and
$J$ of $k$ such that
$I\cap J=\varnothing$, $I\cup J=k$ and elements $w'_i$
(for $i\leq k$) such that
$e(u)=w'_0\leq w'_1\leq\cdots\leq w'_k=e(v)$,
$\xi=\sum_{i\in I}\Dim_L(w'_i,w'_{i+1})$, and
$\eta=\sum_{i\in J}\Dim_L(w'_i,w'_{i+1})$. Since $e[K]$ is
convex in $L$, the $w'_i$'s belong to $e[K]$, thus one can
write $w'_i=e(w_i)$, $w_i\in K$. It is easy to see that one can
take $w_0=u$, $w_k=v$, and $w_i\leq w_{i+1}$, for all $i<k$.
Hence, if we put $\xi^*=\sum_{i\in I}\Dim_L(w_i,w_{i+1})$ and
$\eta^*=\sum_{i\in J}\Dim_L(w_i,w_{i+1})$, then
$\xi^*+\eta^*=\Dim_K(u,v)$ and $f(\xi^*)=\xi$,
$f(\eta^*)=\eta$.\end{proof}

\begin{corollary}\label{C:succVhom}
Let $L$ be a lattice. Then the following are equivalent:
\begin{itemize}
\item[\rm (i)] $L$ is simple\index{lattice!simple ---},
modular,
\index{lattice!modular (not necessarily complemented) ---}
and $L$ has a prime interval.

\item[\rm (ii)] \(\DD L\cong\ZZ^+\).
\end{itemize}
\end{corollary}

\begin{proof}
(ii)$\Rightarrow$(i) If \(\DD L\cong\ZZ^+\), then, by
Proposition~\ref{P:dirfin}, $L$ is modular%
\index{lattice!modular (not necessarily complemented) ---}.
Furthermore, by Corollary~\ref{C:congquotV}, $L$ is simple.
Let $\mu$ be the unique isomorphism from $\DD L$ onto
$\ZZ^+$. There are $a<b$ in $L$ such that $\mu\Dim(a,b)=1$.
If $c$ is an element of $L$ such that \(a\leq c\leq b\), then
\(1=\mu\Dim(a,b)=\mu\Dim(a,c)+\mu\Dim(c,b)\), thus $a=c$ or
$b=c$. Hence $a\prec b$.

(i)$\Rightarrow$(ii) Suppose the hypothesis of (i) satisfied,
and let $[a,\,b]$ be a prime interval in $L$. The
relation\(\vv<x,y>\in\Theta(a,b)\) holds for all
\(x\leq y\) in $L$, thus,
by \cite[Corollary III.1.4]{Grat}, there exists a finite
decomposition \(x=z_0<z_1<\cdots<z_n=y\) such that
\([z_i,\,z_{i+1}]\wpr[a,\,b]\) holds for all \(i<n\);
since $L$ is modular%
\index{lattice!modular (not necessarily complemented) ---}
and \(z_i<z_{i+1}\), we have necessarily
\([z_i,\,z_{i+1}]\approx[a,\,b]\), thus \(z_i\prec z_{i+1}\).
By the Jordan-H\"older property, it follows that $L$ has no
infinite bounded chains. By Proposition~\ref{P:shortmod}, we
obtain that \(\DD L\cong\ZZ^+\).
\end{proof}

\chapter{Primitive refinement monoids, revisited}
\label{PrimMon}

This monoid-theoretical chapter will serve as a transition towards
the study of dimension monoids of finite or, more generally, of
BCF, lattices, undertaken in Chapter~\ref{ShortLatt}. It
turns out that the dimension monoids of these lattices are very
special monoids, the \emph{primitive monoids} introduced by R.~S.
Pierce. These structures are, in some sense, the
monoid-theoretical analogues of the $\{\vee,0\}$-semilattices, where
every element is a finite join of join-irreducible elements.

This chapter will be devoted to recall some of Pierce's classical
results about primitive monoids, as well as proving a few new ones.
In particular, it will turn out that every primitive monoid embeds,
as a partially ordered monoid, into a power of
$\ZZb=\ZZ^+\cup\{\infty\}$, and satisfies many first-order axioms
that are not implied by the refinement property alone.

\section{The monoid constructed from a QO-system}
\label{S:MonFromQO}

A large part of the terminology that we shall use in this
chapter is borrowed from R.~S. Pierce \cite{Pier89}. In
particular, we shall say that an element $p$ of a \cm\ $M$ is
\emph{pseudo-indecomposable}%
\index{pseudo-indecomposable (PI)|ii}
(or \emph{PI}), if $M$ satisfies
\[
(\forall x,y)\bigl(x+y=p\Rightarrow(x=p\ \mathrm{or}\
y=p)\bigr).
\]
We shall say that $M$ is \emph{PI-generated}, if
the set of all nonzero PI elements of $M$ generates $M$ as a
monoid. A \emph{primitive\index{monoid!primitive ---|ii}
monoid} is a PI-generated refinement monoid for which
the algebraic\index{algebraic!preordering} preordering is
\emph{antisymmetric}.

Let us record some results from \cite{Pier89}. A
\emph{quasi-ordered system}, or \emph{QO-system}%
\index{quasi-ordered (QO) -system|ii}, is a set equipped with
a transitive binary relation. If
$\vv<P,\tr>$ is an antisymmetric QO-system, then one can
associate with it a natural partial ordering $\utr$ defined by
$x\utr y$ if and only if $x\tr y$ or $x=y$, along with a subset
$P^{(0)}=\{x\in P\mid x\tr x\}$. Then we put
$P^{(1)}=P\setminus P^{(0)}$. Conversely, if $\vv<P,\utr>$ is
a partially ordered set and if $P_0\subseteq P$, then one can
define an antisymmetric and transitive binary relation
$\tr$ on $P$ by putting $x\tr y$ if and only if $x\dtr y$ or
$x=y\in P_0$, where $x\dtr y$ is short for $x\utr y$ and
$x\ne y$.
Furthermore, these transformations are easily
seen to be mutually inverse.

If $\vv<P,\tr>$ is a QO-system, then one defines
$\mathbf{E}(P)=\mathbf{E}(P,\tr)$%
\index{EzzofP@$\mathbf{E}(P)$|ii} as the \cm\ defined by
generators $e_p$ (for $p\in P$) and relations $e_p\ll e_q$ for all
$p\tr q$ in $P$ (where $x\ll y$ is short for
$x+y=y$, ``$y$ \emph{absorbs} $x$''). In the case where there
may be any ambiguity on $P$, we shall write $e_p^P$ instead of
$e_p$.

For the following theorem, we refer to
R.~S. Pierce \cite[3.4--3.6]{Pier89}.

\begin{theorem}\label{T:Pierce}\hfill

\begin{itemize}
\item[\rm (i)] The primitive
monoids are exactly the monoids of the form
$\mathbf{E}(P)$, where $P$ is an
antisymmetric QO-system.
\item[\rm (ii)] If $\vv<P,\tr>$ is an antisymmetric QO-system%
\index{quasi-ordered (QO) -system}, then the nonzero PI%
\index{pseudo-indecomposable (PI)} elements of
$\mathbf{E}(P)$ are exactly the
$e_p$, $p\in P$, and $p\mapsto e_p$ is an embedding from
$\vv<P,\tr>$ into $(\mathbf{E}(P),\ll)$. Furthermore, every
element of $\mathbf{E}(P)$ admits a representation of the
form $\sum_{i<n}e_{p_i}$, where $n\in\omega$ such that
$p_i\ntr p_j$ holds for all
$i\ne j$. Moreover, this representation, the
\emph{reduced representation} of $x$, is unique up to
permutation.\qed
\end{itemize}
\end{theorem}

In particular, if $P$ is an antisymmetric QO-system%
\index{quasi-ordered (QO) -system}, then
$\mathbf{E}(P)$ is an antisymmetric
refinement\index{monoid!refinement ---} monoid.

\begin{example}\label{E:basicex}
For every $n\in\NN$ and for
every additive subgroup $G$ of $\ZZ^n$,
$G^+=G\cap(\ZZ^+)^n$ is generated by the set of its atoms
(minimal elements of $G^+\setminus\{0\}$), thus,
\emph{a fortiori}, it is PI-generated%
\index{pseudo-indecomposable (PI)}. It does not always satisfy
refinement---in fact, it does satisfy refinement if and only if
$G^+\cong(\ZZ^+)^k$, for some non-negative integer $k$.
\end{example}

The objective of the following lemma is to show that many
monoids are primitive. In particular, dimension monoids of BCF
lattices will turn out to admit a presentation as in the assumption
of Lemma~\ref{L:genrel}.

\begin{lemma}\label{L:genrel}
Let $I$ be any set, let $X$ and
$Y$ be subsets of $I\times I$. Then the \cm\ defined by
generators $e_i$ (for $i\in I$) and relations $e_i=e_j$ (for all
$\vv<i,j>\in X$) and $e_i\ll e_j$ (for all $\vv<i,j>\in Y$) is
primitive\index{monoid!primitive ---}.
\end{lemma}

\begin{proof}
Denote by $\mathbf{M}(I;X,Y)$%
\index{MzzIXY@$\mathbf{M}(I;X,Y)$|ii} the \cm\ above. First,
let
$\bar X$ be the equivalence relation generated by $X$. Then
$e_i=e_j$ holds in
$\mathbf{M}(I;X,Y)$ for all $\vv<i,j>\in\bar X$.
Thus, for all $i\in I$, if $\bar\imath$ denotes the
$\bar X$-equivalence class of $i$, then note
$e_{\bar\imath}=e_i$. Put $\bar I=I/\bar X$ and
$\bar Y=\{\vv<\bar\imath,\bar\jmath>\mid\vv<i,j>\in Y\}$;
then one can define a monoid homomorphism from
$\mathbf{M}(I;X,Y)$ to
$\mathbf{M}(\bar I;\varnothing,\bar Y)$ by
$e_i\mapsto e_{\bar\imath}$, and it is easy to verify that this
homomorphism is, in fact, an isomorphism. Thus
$\mathbf{M}(I;X,Y)\cong\mathbf{M}(\bar I;\varnothing,\bar Y)$.

Therefore, it suffices to prove that
$\mathbf{M}(I;\varnothing,X)$ is
primitive\index{monoid!primitive ---},
for all $X\subseteq I\times I$. First, note that
$\mathbf{M}(I;\varnothing,X)
\cong\mathbf{M}(I;\varnothing,\prec)$,
where $\prec$ is the transitive closure of $X$. Next, define
binary relations $\preceq$, $\sim$ on $I$ by
\begin{align*}
i\preceq j&\Leftrightarrow i=j\ \mathrm{or}\ i\prec j,\\
i\sim j&\Leftrightarrow i\preceq j\text{ and }j\preceq i.
\end{align*}
Thus $\vv<I,\preceq>$ is a preordered set with associated
equivalence
$\sim$; let $\vv<P,\leq>$ be the quotient partially ordered set
$\vv<I,\preceq>/\sim$. For every
$i\in I$, let $[i]$ be the $\sim$-equivalence class of $i$.
Note that $i\prec j\preceq k$ or $i\preceq j\prec k$ implies
that $i\prec k$. This makes it possible to define a binary
relation $\tr$ on $P$ by $[i]\tr[j]$ if and only if $i\prec j$.
Note also that
$\tr$ is antisymmetric and transitive. Then it is easy to
verify that one can define a monoid isomorphism from
$\mathbf{M}(I;\varnothing,\prec)$ onto $\mathbf{E}(P,\tr)$ by
$e_i\mapsto e_{[i]}$.
\end{proof}

\section[Numerical dimension functions]%
{Representation of
primitive\index{monoid!primitive ---} monoids with numerical\\
dimension functions}

Note first that $\mathbf{E}$ is, in fact, given by a
\emph{functor} from the category of all QO-systems%
\index{quasi-ordered (QO) -system} to the category of \cm s.
This functor obviously preserves direct limits. Let $\utr$ be
the partial ordering associated with an antisymmetric
QO-system $\vv<P,\tr>$ (as at the beginning of
Section~\ref{S:MonFromQO}), let
$P=P^{(0)}\cup P^{(1)}$ be the associated decomposition.

We shall from now on use the monoid
$\ZZb=\ZZ^+\cup\{\infty\}$%
\index{zzzeeplusbar@$\ZZb$|ii}, endowed with its natural
structure of \cm. The objective of what follows is to
represent any primitive monoid by a partially ordered, additive
monoid of functions to $\ZZb$.

For every antisymmetric QO-system $\vv<P,\tr>$ and every $p\in P$,
let $\dnw p$\index{dzzownrwp@$\dnw p$, $\dnw_Pp$|ii}
(or $\dnw_Pp$ if $P$ is not understood)
denote the principal lower set generated by $p$ for the natural
ordering $\utr$ of $P$, that is, the set
\(\{q\in P\mid q\utr p\}\) and then let $\mathbf{V}_P$ be
the set of all finite unions of principal lower sets of $P$.
Note that
$\mathbf{V}_P$ is a distributive join-subsemilattice of the
powerset algebra of $P$. Put
\(S_x=\{p\in P\mid x(p)\ne0\}\), for all $x\in(\ZZb)^P$.

Let $\Ft(P)$\index{FzzofPtilde@$\Ft(P)$|ii} be the set of all
mappings $x\colon P\to\ZZb$ satisfying the following four
conditions:

\begin{itemize}
\item $S_x\in\mathbf{V}_P$;
\item $x$ is \emph{antitone}, that is, $p\utr q$ implies
that $x(p)\geq x(q)$.
\item $x[P^{(0)}]\subseteq\{0,\infty\}$;
\item ${x}^{-1}[\NN]$ is an \emph{incomparable} subset
of $P$, that is, for any elements $p$ and $q$ of
${x}^{-1}[\NN]$ such that $p\ne q$, neither $p\utr q$
nor $q\utr p$ holds.
\end{itemize}

\begin{lemma}\label{L:strucFt(P)} The set
$\Ft(P)$ is an additive submonoid of
$(\ZZb)^P$, closed under $\vee$. Furthermore, if $P$ is
\emph{finite}, then $\Ft(P)$ is
closed under
$\wedge$.
\end{lemma}

\begin{proof} Let $x$, $y\in\Ft(P)$.
For all
$z\in\{x\wedge y,x\vee y,x+y\}$, it is clear that
$z[P^{(0)}]\subseteq\{0,\infty\}$ and that $z$ is antitone.
Furthermore, since $S_{x+y}=S_{x\vee y}=S_x\cup S_y$, $S_z$
belongs to $\mathbf{V}_P$ if $z=x+y$ or $z=x\vee y$. If $P$ is
finite, then $\mathbf{V}_P$ is the set of all lower subsets of
$P$, thus $S_{x\wedge y}\in\mathbf{V}_P$. It remains to check
that ${z}^{-1}[\NN]$ is incomparable. Suppose otherwise,
that is, there are $p\dtr q$ in $P$ such that both
$z(p)$ and $z(q)$ belong to $\NN$. For $z\in\{x\vee y,x+y\}$,
this implies that all the elements $x(p)$, $x(q)$, $y(p)$, and
$y(q)$ are finite, thus (since $x$, $y\in\Ft(P)$)
$x(q)=y(q)=0$, whence $z(q)=0$, a contradiction. Now suppose
that $z=x\wedge y$. Without loss of generality,
$x(p)=z(p)\in\NN$. If $x(q)\leq y(q)$, then $x(q)=z(q)\in\NN$,
which contradicts the fact that $x\in\Ft(P)$: therefore,
$x(q)>y(q)$. But $x(p)\in\NN$, thus $x(q)\notin\NN$, whence
$x(q)=\infty$, so that $x(p)<x(q)$, which contradicts the fact
that $x$ is antitone.\end{proof}

For every $x\in(\ZZb)^P$, let $\Gamma_x$ the set
of all \emph{maximal} elements of $S_x$.
Note that $\Gamma_x$ is finite and $S_x=\dnw_P\Gamma_x$, for all
$x\in\Ft(P)$. Then
put\index{FzzofP@$\mathbf{F}(P)$|ii}
\[
\mathbf{F}(P)=\bigl\{x\in\Ft(P)\mid
x[\Gamma_x\cap P^{(1)}]\subseteq\NN\bigr\}.
\]

\begin{lemma}\label{L:strucF(P)} The subset $\mathbf{F}(P)$
is closed under $+$ and $\vee$.
\end{lemma}

\begin{proof}
For $x$, $y\in\mathbf{F}(P)$, we prove that
$x+y\in\mathbf{F}(P)$ (the proof for $\vee$ is similar). Let
$p\in\Gamma_{x+y}\cap P^{(1)}$, it suffices to prove that
$(x+y)(p)<\infty$. Otherwise, we have, without loss of
generality, $x(p)=\infty$, thus (since $x\in\mathbf{F}(P)$ and
$p\in P^{(1)}$) $p\notin\Gamma_x$; it follows that there exists
$q$ such that $p\dtr q$ and $x(q)>0$;
therefore, $(x+y)(q)>0$, which contradicts
\(p\in\Gamma_{x+y}\).
\end{proof}

Now, for all $p\in P$, define $f_p\colon P\to\ZZb$
the following way: for all $q\in P$, put
\[
f_p(q)=\left\{
\begin{array}{ll}
\infty,&\text{if $q\dtr p$,}\\
1,&\text{if $q=p$ and $p\ntr p$,}\\
\infty,&\text{if $q=p$ and $p\tr p$,}\\
0,&\text{if $q\nutr p$}.
\end{array}
\right.
\]
As for the $e_p$'s, we write $f_p^P$ instead of $f_p$ in
the case where there may be any ambiguity on $P$.

\begin{lemma}\label{L:E(P)toF(P)} There exists a unique monoid
homomorphism
$\varphi_P\colon \mathbf{E}(P)\to\mathbf{F}(P)$ such that
$\varphi_P(e_p)=f_p$, for all $p\in P$.
\end{lemma}

\begin{proof}
It suffices to verify that all the
$f_p$'s belong to $\mathbf{F}(P)$ and that
$2f_p=f_p$ for all $p\in P^{(0)}$, while $f_p+f_q=f_q$ for all
$p\dtr q$ in $P$. This is straightforward.
\end{proof}

Now, define a map $\psi_P\colon
\mathbf{F}(P)\to\mathbf{E}(P)$ by the rule
\[\psi_P(x)=\sum_{p\in\Gamma_x\cap P^{(0)}}e_p+
\sum_{p\in\Gamma_x\cap P^{(1)}}x(p)e_p.\]

\begin{lemma}\label{L:F(P)toE(P)} Let $x\in\mathbf{F}(P)$ and
let $S$ and $T$ be finite subsets of $P$ such that
\(\Gamma_x\cap P^{(0)}\subseteq S\subseteq S_x\cap P^{(0)}\)
and \(\Gamma_x\cap P^{(1)}\subseteq T\subseteq
{x}^{-1}[\ZZ^+]\cap P^{(1)}\). Then we have
\[\psi_P(x)=\sum_{p\in S}e_p+
\sum_{p\in T}x(p)e_p.\]
\end{lemma}

\begin{proof}
It suffices to prove that
$e_p\ll\psi_P(x)$, for all
\(p\in S\setminus(\Gamma_x\cap P^{(0)})\),
while $x(p)e_p\ll\psi_P(x)$, for all
\(p\in({x}^{-1}[\ZZ^+]\cap P^{(1)})\setminus
(\Gamma_x\cap P^{(1)})\). Let us
begin with the first relation. Since
\(p\in S\setminus(\Gamma_x\cap P^{(0)})\) and
\(S\subseteq P^{(0)}\), we have $p\notin\Gamma_x$. Since
$S\subseteq S_x$, we have $p\in S_x$. Since
$S_x\in\mathbf{V}_P$, there exists $q\in\Gamma_x$
such that $p\dtr q$. Then
the $q$th summand defining $\psi_P(x)$ is either $e_q$ (if
$q\in P^{(0)}$) or
$x(q)e_q$ (if $q\in P^{(1)}$), thus it absorbs $e_p$. Now let
us check the second point; thus let $p$ be an element of
\(({x}^{-1}[\ZZ^+]\cap P^{(1)})\setminus (\Gamma_x\cap
P^{(1)})\). If $p\notin S_x$, then
$x(p)=0$ and the conclusion is trivial, so suppose that $p\in
S_x$. Since
$p\notin\Gamma_x$ and $S_x\in\mathbf{V}_P$, there exists
$q\in\Gamma_x$ such that $p\dtr q$. Then the $q$th summand
defining $\psi_P(x)$ is either $e_q$ (for $q\in P^{(0)}$) or
$x(q)e_q$ (for $q\in P^{(1)}$), and in both cases, it absorbs
$x(p)e_p$.
\end{proof}

\begin{proposition}\label{P:E(P)isF(P)}
The maps $\varphi_P$ and
$\psi_P$ are mutually inverse monoid isomorphisms between
$\mathbf{E}(P)$ and $\mathbf{F}(P)$.
\end{proposition}

\begin{proof}
It is immediate that
$\psi_P(0)=0$ and that $\psi_P(f_p)=e_p$, for all $p\in P$, so
that it suffices to verify that $\psi_P$ is a semigroup%
\index{semigroup} homomorphism. Thus let
$x$, $y\in\mathbf{F}(P)$. Put
\(S=(\Gamma_x\cup\Gamma_y)\cap P^{(0)}\),
\(K=(x+y)^{-1}[\ZZ^+]\) and
\(T=(\Gamma_x\cup\Gamma_y)\cap K\cap P^{(1)}\). Note that
\(\Gamma_{x+y}\subseteq\Gamma_x\cup\Gamma_y\). Thus, by using
the fact that all $e_p$ (for $p\in P^{(0)}$) are idem-multiple
and by Lemma~\ref{L:F(P)toE(P)}, we have
\begin{align*}
\psi_P(x+y)&=\sum_{p\in S}e_p+\sum_{p\in T}(x(p)+y(p))e_p\\
&=X+Y,
\end{align*}
where we put
\begin{align*}
X&=\sum_{p\in\Gamma_x\cap P^{(0)}}e_p
+\sum_{p\in\Gamma_x\cap K\cap P^{(1)}}x(p)e_p,\\
Y&=\sum_{p\in\Gamma_y\cap P^{(0)}}e_p
+\sum_{p\in\Gamma_y\cap K\cap P^{(1)}}y(p)e_p.
\end{align*}
Note that the caveat for $X$ (resp., $Y$) of being equal to
$\psi_P(x)$ (resp., $\psi_P(y)$) is the fact that the
second sum defining it ranges over all indices in
\(\Gamma_x\cap K\cap P^{(1)}\) (resp.,
\(\Gamma_y\cap K\cap P^{(1)}\)) rather than just
\(\Gamma_x\cap P^{(1)}\) (resp., \(\Gamma_y\cap P^{(1)}\)).
Thus it suffices to prove that
$x(p)e_p\ll\psi_P(x+y)$ for all
\(p\in\Gamma_x\cap P^{(1)}\setminus K\),
and, symmetrically, that
$y(p)e_p\ll\psi_P(x+y)$ for all
\(p\in\Gamma_y\cap P^{(1)}\setminus K\).
Let us, for example, prove the first
statement. Since $p\notin K$, we have, by definition,
$x(p)+y(p)=\infty$; but $p\in\Gamma_x$, thus $x(p)<\infty$; it
follows that $y(p)=\infty$. Since $S_y\in\mathbf{V}_P$, there
exists $q\in\Gamma_y$ such that $p\dtr q$.
Thus, if $q\in P^{(0)}$, then the
$q$th summand defining $Y$ is equal to $e_q$ and it absorbs
$x(p)e_p$. If $q\in P^{(1)}$, then, since $q\in\Gamma_y$, we
have $y(q)\in\NN$, thus, since $x(q)=0$, $q\in K$, so that the
$q$th summand defining $Y$ is equal to $y(q)e_q$, and again, it
absorbs $x(p)e_p$. Thus, again, by Lemma~\ref{L:F(P)toE(P)}, we
obtain that \(X+Y=\psi_P(x)+\psi_P(y)\).
\end{proof}

\section{Approximation by finite QO-systems%
\index{quasi-ordered (QO) -system}}

In this section, we shall set up a certain number of tools that will
be used to establish certain properties satisfied by \emph{all}
primitive monoids, \emph{via finitely generated} primitive
monoids.

If $P$ is any antisymmetric QO-system%
\index{quasi-ordered (QO) -system} and $Q\subseteq P$, let
$\ee$%
\index{ezzpsilonQP@$\ee$|ii} be the canonical monoid
homomorphism from
$\mathbf{E}(Q)$ to
$\mathbf{E}(P)$, so that $\ee(e_q^Q)=e_q^P$ for all
$q\in Q$. Then let $\fft$%
\index{ezztatildeQP@$\fft$|ii} be the map from
$(\ZZb)^Q$ to $(\ZZb)^P$ defined by the rule
\[
(\forall x\in(\ZZb)^Q)(\forall q\in Q)\left[
\fft(x)(q)=\left\{
\begin{array}{ll}
x(q),&\text{ if }q\in Q\\
\infty,&\text{ if }q\in\dnw_PS_x\setminus Q\\
0,&\text{ if }q\in P\setminus(\dnw_PS_x\cup Q)
\end{array}
\right.\right]
\]
This will be the context of the following
Lemma~\ref{L:etaQP} and Lemma~\ref{L:diagram}:

\begin{lemma}\label{L:etaQP} The map
$\fft$ is an embedding of
\(\vv<(\ZZb)^Q,+,0,\leq>\) into\linebreak
\(\vv<(\ZZb)^P,+,0,\leq>\) ($\leq$ denotes the
componentwise ordering on both
$(\ZZb)^Q$ and $(\ZZb)^P$). Furthermore,
\(\fft[\Ft(Q)]\subseteq\Ft(P)\) and
\(\fft[\mathbf{F}(Q)]\subseteq\mathbf{F}(P)\).
\end{lemma}

\begin{proof} It is trivial that $\fft(0)=0$. Next, let $x$
and $y$ be elements of $(\ZZb)^Q$; we check that
\(\fft(x+y)(p)=\fft(x)(p)+\fft(y)(p)\) for all
$p\in P$. If
$p\in Q$, this amounts to proving that \((x+y)(p)=x(p)+y(p)\),
which is just the definition, so suppose that $p\notin Q$. If
\(p\notin\dnw_PS_x\cup\dnw_PS_y\), this amounts to proving
that $0=0+0$. If \(p\in\dnw_PS_x\cap\dnw_PS_y\), then
this amounts to proving that $\infty=\infty+\infty$.
If \(p\in\dnw_PS_x\setminus\dnw_PS_y\), this amounts to
proving that $\infty=\infty+\fft(y)(p)$; the
symmetric case
\(p\in\dnw_PS_y\setminus\dnw_PS_x\) is similar. Hence,
$\fft$ is a monoid homomorphism. Since
$\leq$ is also the algebraic\index{algebraic!preordering}
ordering on both
$(\ZZb)^Q$ and $(\ZZb)^P$,
$\fft$ is also a homomorphism for $\leq$.
Moreover, since
$x=\fft(x)\res_Q$ for all $x\in(\ZZb)^Q$, $\fft$
is, in fact, an order-embedding.

Now let us prove the second part of the lemma; thus let first
$x$ be an element of $\Ft(Q)$. First, we have
\[
\fft(x)[P^{(0)}]\subseteq\fft(x)[Q^{(0)}]\cup\{0,\infty\}
=x[Q^{(0)}]\cup\{0,\infty\}\subseteq\{0,\infty\}.
\]
Next, we have
\[
S_{\fft(x)}=S_x\cup(\dnw_PS_x\setminus Q)
=\dnw_PS_x\in\mathbf{V}_P,
\]
thus, in particular,
\(\Gamma_{\fft(x)}=\Gamma_x\). Furthermore,
\(\fft(x)^{-1}[\NN]={x}^{-1}[\NN]\) is incomparable. Thus, to
verify that $\fft(x)\in\Ft(P)$, it suffices to verify that
$\fft(x)$ is antitone. Hence let $p\dtr q$ in $Q$,
and suppose that $\fft(x)(p)<\fft(x)(q)$. In particular,
$\fft(x)(p)<\infty$ and $\fft(x)(q)>0$. Since
$\fft(x)\res_Q=x$ is antitone, it cannot happen that
both $p$ and
$q$ belong to $Q$. Suppose first that $p\in Q$, so that
$x(p)<\infty$ and $q\notin Q$. It follows that
$\fft(x)(q)=\infty$, so that $q\in\dnw_PS_x$, that is,
there exists $p'\in S_x$ such that $q\utr p'$. Thus
$p\dtr p'\in S_x$ in $Q$, thus $x(p)=\infty$, a contradiction.
Thus $p\notin Q$. Since $\fft(x)(p)<\infty$, we have
\(p\notin\dnw_PS_x\cup Q\). If $q\in Q$, then (since $p\utr q$
and $p\notin\dnw_PS_x$) $q\notin S_x$, thus
\(\fft(x)(q)=x(q)=0\), a contradiction; thus $q\notin Q$.
Since
$\fft(x)(q)>0$, we have \(q\in\dnw_PS_x\setminus Q\), thus,
since $p\utr q$,
\(p\in\dnw_PS_x\setminus Q\), so that $\fft(x)(p)=\infty$, a
contradiction again. This completes the proof that
$\fft(x)$ is antitone, thus that
\(\fft[\Ft(Q)]\subseteq\Ft(P)\). To obtain that
\(\fft[\mathbf{F}(Q)]\subseteq\mathbf{F}(P)\), it suffices to
prove that
\(\fft(x)[\Gamma_{\fft(x)}\cap P^{(1)}]\subseteq\NN\) for all
\(x\in\mathbf{F}(Q)\).
However, as we have seen,
\(\Gamma_{\fft(x)}=\Gamma_x\subseteq Q\) so that the conclusion
follows at once from the fact that \(x\in\mathbf{F}(Q)\).
\end{proof}

According to Lemma~\ref{L:etaQP}, we shall define $\ff$%
\index{ezztaQP@$\ff$|ii} as the restriction of $\fft$ from
$\mathbf{F}(Q)$ to
$\mathbf{F}(P)$.

\begin{lemma}\label{L:diagram} The following diagram is
commutative:
\[
\begin{CD}
\mathbf{E}(Q) @>{\ee}>> \mathbf{E}(Q)\\
@V{\varphi_Q}VV
@VV{\varphi_P}V\\
\mathbf{F}(Q) @>{\ff}>> \mathbf{F}(P)
\end{CD}
\]
\end{lemma}

\begin{proof}
Since the arrows of the diagram above are monoid
homomorphisms, it suffices to prove that
$\varphi_P\circ\ee$ and $\ff\circ\varphi_Q$ agree on all
elements of
$\mathbf{E}(Q)$ of the form
$e_q^Q$, where $q\in Q$. This, in turn, amounts to proving that
$\ff(f_q^Q)=f_q^P$, which is a simple verification.
\end{proof}

\begin{lemma}\label{L:invF(P)}
Let $P$ be an arbitrary antisymmetric QO-system%
\index{quasi-ordered (QO) -system}. Then the following holds:

\begin{itemize}
\item[\rm(i)] The equality
\[
x=\sum_{p\in\Gamma_x\cap P^{(0)}}f_p^P+
\sum_{p\in\Gamma_x\cap P^{(1)}}x(p)f_p^P
\]
holds for all $x\in\Ft(P)$.
\item[\rm(ii)] For all $x$ in $\Ft(P)$
(resp., $\mathbf{F}(P)$) and all $Q$ such that
$\Gamma_x\subseteq Q\subseteq P$, $x$ belongs to
\(\fft[\Ft(Q)]\) (resp., \(\ff[\mathbf{F}(Q)]\)).
\end{itemize}

\end{lemma}

\begin{proof}
(i) Let $y$ be the right hand side, let
$p\in P$; we verify that $x(p)=y(p)$. If $p\notin S_x$, then
$x(p)=y(p)=0$. If $p\in\Gamma_x$, then, since $\Gamma_x$ is
incomparable, $y(p)$ is equal either to $f_p(p)=\infty$, if
$p\in P^{(0)}$, or to $x(p)f_p(p)=x(p)$, if $p\in P^{(1)}$. In
both cases, $y(p)=x(p)$. If $p\in S_x\setminus\Gamma_x$, then
(since $S_x\in\mathbf{V}_P$) there exists $q\in\Gamma_x$ such
that $p\dtr q$. If
$x(p)<\infty$, then, since $x$ is antitone, $x(q)<\infty$, but
$q\in\Gamma_x$, thus $x(q)>0$, so that both $x(p)$ and $x(q)$
belong to $\NN$, contradicting $x\in\Ft(P)$; hence
$x(p)=\infty$. Moreover, $f_q(p)=\infty$, so that
$x(p)=\infty=y(p)$.

(ii) follows immediately from (i).
\end{proof}

Now, when $P$ is a \emph{finite} antisymmetric QO-system%
\index{quasi-ordered (QO) -system}, let us define, for all
$n\in\ZZ^+$, a map
$\rho_n\colon (\ZZb)^P\to\mathbf{F}(P)$ by the rule
\[
\rho_n(x)=\sum_{p\in\Gamma_x}(x(p)\wedge n)f_p,
\]
where $\wedge$ denotes the infimum operation in $\ZZ^+$.

\begin{lemma}\label{L:tech} Let $P$ be a \emph{finite}
antisymmetric QO-system%
\index{quasi-ordered (QO) -system},
let $x$, $y\in(\ZZb)^P$, and
let $n\in\ZZ^+$. Then the following statements hold:

\begin{itemize}
\item[\rm (i)] $\rho_n(x)=\sum_{p\in P}(x(p)\wedge n)f_p$.

\item[\rm (ii)] If $x\in\Ft(P)$, then
$\rho_n(x)\leq x$ (for the componentwise ordering).

\item[\rm (iii)] $x\in\mathbf{F}(P)$ if and only if
$\rho_n(x)=x$, for all large enough $n$.

\item[\rm (iv)] $\rho_n(x)+\rho_n(y)\leq\rho_{2n}(x+y)$.

\item[\rm (v)] $\rho_n(x+y)\leq\rho_n(x)+\rho_n(y)$.
\end{itemize}
\end{lemma}

\begin{proof}
(i) It suffices to prove that
$(x(p)\wedge n)f_p\ll\rho_n(x)$,
for all $p\in P\setminus\Gamma_x$.
If $p\notin S_x$, then $x(p)=0$ and the conclusion is trivial.
If $p\in S_x\setminus\Gamma_x$, then there exists $q\in\Gamma_x$
such that $p\dtr q$, so that
$(x(p)\wedge n)f_p\ll(x(q)\wedge n)f_q\leq\rho_n(x)$.

(ii) Immediate by Lemma~\ref{L:invF(P)} and the fact that all
$f_p$ (for $p\in P^{(0)}$) are idem-multiple.

(iii) It is trivial the the range of $\rho_n$ is contained in
$\mathbf{F}(P)$. Conversely, let
$x\in\mathbf{F}(P)$ and let
$n\in\NN$ such that
$n\geq x(p)$, for all $p\in\Gamma_x\cap P^{(1)}$. Then we have
\begin{align*}
\rho_n(x)&=\sum_{p\in\Gamma_x\cap P^{(0)}}(x(p)\wedge n)f_p+
\sum_{p\in\Gamma_x\cap P^{(1)}}(x(p)\wedge n)f_p\\
&=\sum_{p\in\Gamma_x\cap P^{(0)}}f_p+
\sum_{p\in\Gamma_x\cap P^{(1)}}x(p)f_p\\
&=x\qquad\hbox{(by Lemma~\ref{L:invF(P)})}.
\end{align*}

(iv) Follows immediately from (i) and the inequality
$x\wedge n+y\wedge n\leq(x+y)\wedge 2n$ (in
$\ZZ^+$).

(v) Follows immediately from (i) and the inequality
\((x+y)\wedge n\leq x\wedge n+y\wedge n\) (in
$\ZZ^+$).
\end{proof}

\section{Further properties of
primitive\index{monoid!primitive ---} monoids}

In this section, we shall apply the results of the previous
sections to obtain further properties of
primitive\index{monoid!primitive ---} monoids.
In \ref{L:algord}--\ref{P:pseudocan},
let $P$ be an antisymmetric QO-system%
\index{quasi-ordered (QO) -system}.

\begin{lemma}\label{L:algord} The
algebraic\index{algebraic!preordering} preordering of
$\mathbf{F}(P)$ is the componentwise ordering of
$\mathbf{F}(P)$.
\end{lemma}

\begin{proof}
Denote, temporarily, by $\leq_\mathrm{alg}$ the
algebraic\index{algebraic!preordering} preordering of
$\mathbf{F}(P)$. It is trivial that the
componentwise ordering on $\mathbf{F}(P)$ is coarser than the
algebraic\index{algebraic!preordering} ordering on
$\mathbf{F}(P)$. For the converse, let us first
consider the case where $P$ is
\emph{finite}. Let
$x$, $y\in\mathbf{F}(P)$ such that $x\leq y$. Let the mapping
$\bar z\colon P\to\ZZb$ be defined by the rule
\[
\bar z(p)=\text{largest element }t\in\ZZb
\text{ such that }x(p)+t=y(p).
\]
Then $x+\bar z=y$. There exists $n_0\in\NN$
such that $\rho_n(x)=x$ and
$\rho_n(y)=y$ for all $n\geq n_0$. But by Lemma~\ref{L:tech}, we have
\[
\rho_n(x)+\rho_n(\bar z)\leq\rho_{2n}(y)\leq
\rho_{2n}(x)+\rho_{2n}(\bar z),
\]
which can be written as
\[
x+\rho_n(\bar z)\leq y\leq x+\rho_{2n}(\bar z);
\]
but we also have $x+\rho_{2n}(\bar z)\leq y$
(by the above inequality applied to $2n$), whence
$x+\rho_{2n}(\bar z)=y$. But $\rho_{2n}(\bar z)$ belongs to
$\mathbf{F}(P)$, whence $x\leq_\mathrm{alg}y$.

In the general case, let again
$x$, $y\in\mathbf{F}(P)$ such that $x\leq y$. Put
$Q=\Gamma_x\cup\Gamma_y$. By Lemma~\ref{L:invF(P)}, there
exist elements $x'$ and $y'$ of $\mathbf{F}(Q)$ such that
$x=\ff(x')$ and $y=\ff(y')$. Since, by Lemma~\ref{L:etaQP},
the map $\ff$ is an order-embedding, we have
$x'\leq y'$, thus, by the finite case,
$x'\leq_\mathrm{alg}y'$, thus, applying
$\ff$, \(x\leq_\mathrm{alg}y\).
\end{proof}

\begin{corollary}\label{C:subdirect} Every
primitive\index{monoid!primitive ---} monoid admits an
embedding (for the monoid structure as well as for the
ordering) into a power of $\ZZb$.\qed
\end{corollary}

Note that we already know that (the
algebraic\index{algebraic!preordering} preordering on)
$\mathbf{E}(P)$ is antisymmetric.
Nevertheless, Lemma~\ref{L:algord} makes it possible, for example,
to obtain immediately the following result:

\begin{corollary}\label{C:unperf} The \cm\
$\mathbf{E}(P)$, equipped with its
algebraic\index{algebraic!preordering} preordering, is
\emph{unperforated}\index{unperforation}, that is, it
satisfies the axiom
 \begin{equation}
 (\forall x,y)(mx\leq my\Rightarrow x\leq y)
 \end{equation}
for all $m\in\NN$.\qed
\end{corollary}

In particular, $\mathbf{E}(P)$ is
\emph{separative\index{separativity} as a positively
preordered commutative monoid}, see \cite{Wehr94}, which means,
since we are discussing the
\emph{algebraic}\index{algebraic!preordering} preordering, that
it satisfies both axioms
 \begin{gather*}
 (\forall x,y)(2x=x+y=2y\Longrightarrow x=y),\\
 (\forall x,y)(x+y=2y\Longrightarrow x\leq y).
 \end{gather*}

We shall now investigate further axioms satisfied by all
$\mathbf{E}(P)$'s. The first one, named the
\emph{interval axiom}\index{interval axiom|ii} in
\cite{Wehr96a}, is the following:
 \[
 (\forall x,y_0,y_1,z)\bigl(z\leq x+y_0,x+y_1
 \Rightarrow(\exists y\leq y_0,y_1)(z\leq x+y)\bigr).
 \]

\begin{proposition}\label{P:intax} The \cm\
$\mathbf{E}(P)$, endowed with its
algebraic\index{algebraic!preordering} preordering, satisfies
the interval axiom.
\end{proposition}

\begin{proof}
It suffices to consider the case where
$P$ is finite and to prove that $\mathbf{F}(P)$ satisfies
the interval axiom. Let $\bar y=y_0\wedge y_1$
(componentwise infimum, calculated in $(\ZZb)^P$). Then
$z\leq x+\bar y$. Let $n\in\NN$ be large enough such that
$\rho_n(x)=x$ and $\rho_n(z)=z$. Put $y=\rho_n(\bar y)$. Then,
again by using Lemma~\ref{L:tech},
\[
z=\rho_n(z)\leq\rho_n(x)+y.
\]
But by Lemma~\ref{L:strucFt(P)},
$\bar y\in\Ft(P)$, thus, by
Lemma~\ref{L:tech},
$y\leq\bar y$, so that $y\leq y_0$, $y_1$.
\end{proof}

The second axiom, the
\emph{pseudo-cancellation\index{pseudo-cancellation|ii}
property}, has been considered, for example, in
\cite{ShWe94,Wehr92a,Wehr92b}. It is the following:
\[
(\forall x,y,z)\bigl(x+z\leq y+z
\Rightarrow(\exists t\ll z)(x\leq y+t)\bigr).
\]

\begin{proposition}\label{P:pseudocan} The \cm\
$\mathbf{E}(P)$, endowed with its
algebraic\index{algebraic!preordering} preordering, satisfies
the pseudo-cancellation property.
\end{proposition}

\begin{proof}
Again, it suffices to consider the case where
$P$ is finite and to prove, in that case, that
$\mathbf{F}(P)$ satisfies the
pseudo-cancellation property. Let
$\bar z=z/\infty$ be defined by the rule
 \[
 \bar z(p)=\left\{
   \begin{array}{ll}
   \infty,&\text{ if }z(p)=\infty\\
   0,&\text{ if }z(p)<\infty.
   \end{array}
 \right.
 \]
Then $\bar z\in\Ft(P)$ and
\(x+\bar z=y+\bar z\). Therefore, the following relations
 \begin{align*}
 x\leq x+\rho_n(\bar z)=\rho_n(x)+\rho_n(\bar z)
 &\leq\rho_{2n}(x+\bar z)\\
 &=\rho_{2n}(y+\bar z)\\
 &\leq\rho_{2n}(y)+\rho_{2n}(\bar z)\\
 &=y+\rho_{2n}(\bar z)
 \end{align*}
hold for all large enough \(n\in\NN\)
with \(\rho_{2n}(\bar z)\in\mathbf{F}(P)\) and
\(\rho_{2n}(\bar z)\leq\bar z\ll z\), so that
\(\rho_{2n}(\bar z)\ll z\).
\end{proof}

This gives us an alternative (and indirect) proof of the
refinement\index{refinement!property} property of
$\mathbf{E}(P)$: indeed, by
\cite[Theorem 3.2]{ShWe94} (see Claim 2 of the proof), every
\cm\ which, endowed with its
algebraic\index{algebraic!preordering} preordering, is
antisymmetric and satisfies both the interval
axiom\index{interval axiom} and the
pseudo-cancellation\index{pseudo-cancellation} property,
satisfies also the refinement\index{refinement!property}
property (it is called in several papers a
\emph{strong refinement\index{monoid!strong refinement ---|ii}
monoid}).

Now, the following result gives us a summary of the main
results of this chapter:

\begin{theorem}\label{T:list}
Let $M$ be a primitive monoid\index{monoid!primitive ---},
endowed with its algebraic ordering.
Then $\vv<M,+,0,\leq>$ satisfies the following properties:

\begin{itemize}
\item[\rm (i)] Antisymmetry:
\((\forall x,y)((x\leq y\text{ and }y\leq x)
\Longrightarrow x=y)\).

\item[\rm (ii)] Unperforation:
\((\forall x,y)(mx\leq my\Rightarrow x\leq y)\),
for all $m\in\NN$.

\item[\rm (iii)] Interval axiom\index{interval axiom}
(see Proposition~\ref{P:intax}).

\item[\rm (iv)] Pseudo-cancellation property (see
Proposition~\ref{P:pseudocan}).

\item[\rm (v)] Refinement\index{refinement!property} property,

\item[\rm (vi)] $\vv<M,+,0,\leq>$ embeds
into a power of $\ZZb$.\qed
\end{itemize}
\end{theorem}

\chapter{Dimension monoids of BCF\index{lattice!BCF ---}
lat\-tices}
\label{ShortLatt}

Random checking on small finite lattices leads to the natural
conjecture that the dimension monoid of a finite lattice is always
a refinement monoid, and, in fact, a primitive monoid. In this
chapter, we shall prove this conjecture as a particular case of a
more general result, see Corollary~\ref{C:dimprim}. The class of
lattices to which this result applies will be called the class of
BCF lattices; it contains all finite lattices.

\section{An alternative presentation of the dimension monoid}

In this section, we shall define BCF lattices. Furthermore, for every BCF
lattice $L$, we shall find an alternative presentation
(in terms of generators and relations) of the
dimension monoid $\DD L$, that will imply, with the help of
Lemma~\ref{L:genrel}, that $\DD L$ is a primitive monoid. The proof
is reminiscent of the classical proof of the Jordan-H\"older Chain
Condition in a finite semimodular lattice. Hence our
result can be considered as a version of the Jordan-H\"older Chain
Condition for \emph{any} finite lattice. The $\omega$-valued height
function on $L$ is replaced by the $\DD L$-valued dimension
function on $L$.

\begin{definition}\label{D:short}
A partially ordered set is
\emph{BCF}\index{BCF partially ordered set|ii}, if it has
no infinite bounded chain.
\end{definition}

Here, BCF stands for ``Bounded Chain Finite''.

In particular, every finite partially ordered set is
BCF.

\begin{lemma}\label{L:equivshort}
A partially ordered set is BCF if and only if the set
of its closed intervals, partially ordered under inclusion, is
well-founded.
\end{lemma}

\begin{proof} Let $\vv<P,\leq>$ be a partially ordered set.
Suppose first that $P$ is BCF and let
\(\vv<[a_n,\,b_n]\mid n\in\omega>\) be a decreasing
sequence of closed intervals of $P$. Thus $\vv<a_n\mid
n\in\omega>$ is increasing while $\vv<b_n\mid
n\in\omega>$ is decreasing. Since $P$ has no infinite bounded
chain, both sequences
\(\vv<a_n\mid n\in\omega>\) and
\(\vv<b_n\mid n\in\omega>\) are eventually constant; thus
so is the sequence \(\vv<[a_n,\,b_n]\mid n\in\omega>\).

Conversely, let $C$ be an infinite bounded chain of $P$; we
show that the set of closed intervals of $P$ is not well
founded under inclusion. Let $a\leq b$ in $P$ such that
$C\subseteq[a,\,b]$. If $C$ is well-founded (resp., dually
well-founded), then there is a strictly increasing (resp.,
strictly decreasing) sequence
$\vv<c_n\mid n\in\omega>$ of elements of $C$. Define a
strictly decreasing sequence
$\vv<I_n\mid n\in\omega>$ of closed intervals of $P$ by
putting
$I_n=[c_n,\,b]$ in the first case and \(I_n=[a,\,c_n]\) in the
second case.
\end{proof}

In \ref{D:caustic}--\ref{L:rect}, let $L$ be a BCF%
\index{lattice!BCF ---} lattice.

\begin{definition}\label{D:caustic} A pair $\vv<a,b>$ of
elements of $L$ is \emph{caustic}\index{caustic!pair|ii}, if
$a\parallel b$\index{Pzzarallel@$a\parallel b$|ii}
(that is, $a$ and $b$ are incomparable) and
the following statements hold:
\begin{align*} (\forall x\in\open{a\wedge b,\,a})&(x\vee
b=a\vee b);\\
(\forall y\in\open{a\wedge b,\,b})&(a\vee
y=a\vee b);\\
(\forall x'\in\open{a,\,a\vee b})&(x'\wedge
b=a\wedge b);\\
(\forall y'\in\open{b,\,a\vee b})&(a\wedge
y'=a\wedge b).
\end{align*}
\end{definition}

\begin{definition}\label{D:path} A \emph{path} (resp.,
\emph{antipath}) of $L$ is a finite sequence
$\alpha\colon n+1\to L$ such that
\((\forall i<n)(\alpha(i)\prec\alpha(i+1))\) (resp.,
\((\forall i<n)(\alpha(i+1)\prec\alpha(i))\)). If $\alpha$ is
either a path or an antipath, let $\alpha^*$ be its
\emph{opposite}, defined by the rule
\(\alpha^*(i)=\alpha(n-i)\) (where
\(n=\mathrm{lh}(\alpha)\), all $i<n+1$). We say that
\(n=\mathrm{lh}(\alpha)\) is the \emph{length} of $\alpha$;
furthermore, if $a=\alpha(0)$ and $b=\alpha(n)$, we write
$a\to^\alpha b$, and we say that $\alpha$ is a path from $a$ to
$b$. Note that by the assumption on $L$, there exists a path
from $a$ to~$b$, for all $a\leq b$ in $L$. If
\(\alpha=\vv<a,a_1,\ldots,a_{m-1},b>\) and
\(\beta=\vv<b,b_1,\ldots,b_{n-1},c>\) are paths, we shall
put \(\alpha\beta=
\vv<a,a_1,\ldots,a_{m-1},b,b_1,\ldots,b_{n-1},c>\). A
\emph{caustic path}%
\index{caustic!path|ii} is a quadruple
\(\gamma=\vv<\alpha_0,\beta_0,\alpha_1,\beta_1>\) of paths for
which there exists a (necessarily unique) caustic pair $\vv<a,b>$
(the \emph{base} of $\gamma$) such that
\(a\wedge b\to^{\alpha_0}a\to^{\alpha_1}a\vee b\) and
\(a\wedge b\to^{\beta_0}b\to^{\beta_1}a\vee b\).
\end{definition}

Now let $\DD' L$ be the \cm\ defined by generators $\|b-a\|$
(for $a\prec b$ in $L$) and, for every caustic
path $\vv<\alpha_0,\beta_0,\alpha_1,\beta_1>$ with base
$\vv<a,b>$, the relations
\begin{align}
\|\beta_0(1)-\beta_0(0)\|&=\|\alpha_1^*(0)-\alpha_1^*(1)\|,
\label{Eq:ShortRel1}\\
\|\beta_0(j+1)-\beta_0(j)\|&\ll\|\beta_0(1)-\beta_0(0)\|
\quad(\mathrm{all}\ j\in\{1,\ldots,\mathrm{lh}(\beta_0)-1\}),
\label{Eq:ShortRel2}\\
\|\beta_1^*(j)-\beta_1^*(j+1)\|&\ll
\|\beta_1^*(0)-\beta_1^*(1)\|
\quad(\mathrm{all}\ j\in\{1,\ldots,\mathrm{lh}(\beta_1)-1\}).
\label{Eq:ShortRel3}
\end{align}
Note that since
$\vv<\beta_0,\alpha_0,\beta_1,\alpha_1>$ is also a caustic
path with base $\vv<b,a>$, the symmetric set
of relations
\begin{align*}
\|\alpha_0(1)-\alpha_0(0)\|&=\|\beta_1^*(0)-\beta_1^*(1)\|,\\
\|\alpha_0(i+1)-\alpha_0(i)\|&\ll\|\alpha_0(1)-\alpha_0(0)\|
\quad(\mathrm{all}\ i\in\{1,\ldots,\mathrm{lh}(\alpha_0)-1\})\\
\|\alpha_1^*(i)-\alpha_1^*(i+1)\|&\ll
\|\alpha_1^*(0)-\alpha_1^*(1)\|
\quad(\mathrm{all}\ i\in\{1,\ldots,\mathrm{lh}(\alpha_1)-1\})
\end{align*} also holds in $\DD' L$. Now, by
Lemma~\ref{L:genrel}, we have immediately the following result:

\begin{proposition}\label{P:V'(L)prim} The \cm\ $\DD' L$ is
primitive\index{monoid!primitive ---}.\qed
\end{proposition}

In particular, $\DD' L$ is a strong
refinement\index{monoid!strong refinement ---} monoid, it is
antisymmetric, unperforated, \emph{etc.} (see
Theorem~\ref{T:list}). We shall now prove that $\DD' L$ is
isomorphic to
$\DD L$, \emph{via}
$\|b-a\|\leftrightarrow\Dim(a,b)$.

If $\alpha$ is a path of $L$, let the \emph{arclength}%
\index{arclength|ii} of $\alpha$ be the following finite sum
in $M$:
\[
\|\alpha\|=\sum_{i<\mathrm{lh}(\alpha)}
\|\alpha(i+1)-\alpha(i)\|.
\]
We say that a closed interval $[a,\,b]$ of $L$ is
\emph{rectifiable}\index{rectifiable|ii}, if any two paths
from $a$ to $b$ have the same arclength, and
then denote by $\|b-a\|$ the common value of this arclength.
Note that if $a=b$ or $a\prec b$, then $[a,\,b]$ is trivially
rectifiable.

\begin{lemma}\label{L:rect} Let $c\leq d$ in $L$. Then the
following holds:

\begin{itemize}
\item[\rm (i)] Every closed subinterval of $[c,\,d]$ is
rectifiable.

\item[\rm (ii)] The equality
$\|a\vee b-a\|=\|b-a\wedge b\|$ holds for all $a$, $b\in[c,\,d]$.
\end{itemize}
\end{lemma}

\begin{proof}
By Lemma~\ref{L:equivshort}, the set of closed intervals of $L$,
endowed with containment, is well-founded. Thus, we can argue by
$\subseteq$-induction on $[c,\,d]$. So suppose that the property holds
for all strict closed subintervals of $[c,\,d]$. In
particular, for every $[x,\,y]\subset[c,\,d]$ (where
$\subset$ denotes
\emph{strict} inclusion), the notation
$\|y-x\|$ is well-defined. Most of our effort will be devoted
to prove the following claim:

\setcounter{claim}{0}
\begin{claim}
For all elements $a$, $b$ of
$[c,\,d]$ such that $a\ne c$ and $b\ne d$, $[a,\,a\vee b]$ and
$[a\wedge b,\,b]$ are rectifiable\index{rectifiable} and
$\|a\vee b-a\|=\|b-a\wedge b\|$.
\end{claim}

\begin{cproof}
The conclusion is trivial if $a$ and $b$
are comparable, so suppose that $a\parallel b$. Since
\(c<a\leq a\vee b\leq d\) and \(c\leq a\wedge b\leq b<d\), the
rectifiability statement is obvious. If
\([a\wedge b,\,a\vee b]\subset[c,\,d]\), then the conclusion
follows from the induction hypothesis; thus suppose that
$a\wedge b=c$ and $a\vee b=d$.

Suppose first that $\vv<a,b>$ is a caustic pair%
\index{caustic!pair}. Let
\(\vv<\alpha_0,\beta_0,\alpha_1,\beta_1>\)
be a caustic path%
\index{caustic!path} with base
$\vv<a,b>$. Then, by the definition of the relations
(\ref{Eq:ShortRel1}), (\ref{Eq:ShortRel2}),
(\ref{Eq:ShortRel3}) in $\DD' L$, we have
\begin{align*}
\|d-a\|&=\|\alpha_1^*(0)-\alpha_1^*(1)\|+
\sum_{1\leq i<\mathrm{lh}(\alpha_1)}
\|\alpha_1^*(i)-\alpha_1^*(i+1)\|\\
&=\|\alpha_1^*(0)-\alpha_1^*(1)\|\\
&=\|\beta_0(1)-\beta_0(0)\|\\
&=\|\beta_0(1)-\beta_0(0)\|
+\sum_{1\leq j<\mathrm{lh}(\beta_0)}
\|\beta_0(j+1)-\beta_0(j)\|\\
&=\|b-c\|,
\end{align*}
and we are done.

Now suppose that $\vv<a,b>$ is not a caustic%
\index{caustic!pair} pair. Then four cases can occur:

\begin{itemize}
\item[\textbf{\textit{Case 1.}}] There exists
\(x\in\open{c,\,a}\) such that $x\vee b<d$.

Put \(\bar x=(x\vee b)\wedge a\) (so that $\bar x\geq x$). Then
\([\bar x,\,x\vee b]\nearrow[a,\,d]\) and
\([\bar x,\,d]\subset[c,\,d]\); thus, by the induction
hypothesis, \(\|d-a\|=\|x\vee b-\bar x\|\). Moreover,
\([c,\,b]\nearrow[\bar x,\,x\vee b]\) and
\([c,\,x\vee b]\subset[c,\,d]\), thus, by the induction
hypothesis, \(\|x\vee b-\bar x\|=\|b-c\|\). Therefore,
\(\|d-a\|=\|b-c\|\).\smallskip

\item[\textbf{\textit{Case 2.}}]
There exists \(x\in\open{b,\,d}\) such that $c<x\wedge a$.

Put \(\bar x=(x\wedge a)\vee b\) (so that $\bar x\leq x$). Then
\([x\wedge a,\,\bar x]\nearrow[a,\,d]\) and
\([x\wedge a,\,d]\subset[c,d]\); thus, by the induction
hypothesis, \(\|d-a\|=\|\bar x-x\wedge a\|\). Moreover,
\([c,\,b]\nearrow[x\wedge a,\,\bar x]\) and
\([c,\,\bar x]\subset[c,\,d]\), thus, by the induction
hypothesis,
\(\|\bar x-x\wedge a\|=\|b-c\|\). Therefore,
\(\|d-a\|=\|b-c\|\).\smallskip

\item[\textbf{\textit{Case 3.}}] There exists
\(y\in\open{c,\,b}\) such that $y\vee a<d$.

Put \(\bar y=(y\vee a)\wedge b\) (so that
$\bar y\geq y$). Then \([c,\,\bar y]\nearrow[a,\,y\vee a]\) and
\([c,\,a\vee\bar y]\subset[c,\,d]\), thus, by the induction
hypothesis,
\(\|\bar y-c\|=\|y\vee a-a\|\). Moreover,
\([\bar y,\,b]\nearrow[y\vee a,\,d]\) and
\([\bar y,\,d]\subset[c,\,d]\), thus, by the induction
hypothesis,
\(\|d-y\vee a\|=\|b-\bar y\|\). Therefore,
\begin{align*}
\|d-a\|&=\|d-y\vee a\|+\|y\vee a-a\|\\
&=\|b-\bar y\|+\|\bar y-c\|\\
&=\|b-c\|.
\end{align*}

\item[\textbf{\textit{Case 4.}}] There exists
\(y\in\open{a,\,d}\) such that $c<y\wedge b$.

Let \(\bar y=a\vee(y\wedge b)\) (so that $\bar y\leq y$). Then
\([c,\,y\wedge b]\nearrow[a,\,\bar y]\) and \([c,\,\bar
y]\subset[c,\,d]\), thus, by the induction hypothesis,
\(\|\bar y-a\|=\|y\wedge b-c\|\). Moreover,
\([y\wedge b,\,b]\nearrow[\bar y,\,d]\) and
\([y\wedge b, d]\subset[c,\,d]\), thus, by the induction
hypothesis,
\(\|d-\bar y\|=\|b-y\wedge b\|\). Therefore,
\begin{align*}
\|d-a\|&=\|d-\bar y\|+\|\bar y-a\|\\
&=\|b-y\wedge b\|+\|y\wedge b-c\|\\
&=\|b-c\|.\tag*{\qedc}
\end{align*}

\end{itemize}
\renewcommand{\qedc}{}\end{cproof}

Now, we prove that the interval $[c,\,d]$ itself is
rectifiable\index{rectifiable}. Thus let $\alpha$ and $\beta$
be paths from
$c$ to $d$, we prove that \(\|\alpha\|=\|\beta\|\). As observed
earlier, this is trivial if $c=d$ or $c\prec d$. Thus suppose
that both $\alpha$ and $\beta$ have length at least $2$. Put
$a=\alpha(1)$, $b=\beta(1)$ and let $\alpha'$, $\beta'$ be the
paths of respective length $\mathrm{lh}(\alpha)-1$,
$\mathrm{lh}(\beta)-1$ defined by the rules
\(\alpha'(i)=\alpha(i+1)\) and \(\beta'(j)=\beta(j+1)\). If
$a=b$, then, since \([a,\,d]\subset[c,\,d]\), we have
\(\|\alpha'\|=\|\beta'\|\), whence
\(\|\alpha\|=\|a-c\|+\|\alpha'\|=\|b-c\|+\|\beta'\|=\|\beta\|\).
Thus suppose that $a\ne b$. It follows that $a\wedge b=c$. Put
$e=a\vee b$ and let $\alpha''$, $\beta''$, and $\gamma$ be
paths respectively from $a$ to $e$, from $b$ to $e$ and from
$e$ to $d$. Since $c<a$ and $c<b$, both intervals $[a,\,d]$
and $[b,\,d]$ are rectifiable\index{rectifiable}, whence
\(\|\alpha'\|=\|\alpha''\|+\|\gamma\|\) and
\(\|\beta'\|=\|\beta''\|+\|\gamma\|\). Moreover, by
Claim~1, $\|\alpha''\|=\|b-c\|$ and
$\|\beta''\|=\|a-c\|$. Therefore, the following holds:
\begin{align*}
\|\alpha\|&=\|a-c\|+\|\alpha'\|\\
&=\|a-c\|+\|b-c\|+\|\gamma\|\\
&=\|b-c\|+\|\beta'\|\\
&=\|\beta\|,
\end{align*}
so that we have proved rectifiability of
$[c,\,d]$. Thus condition (i) of Lemma~\ref{L:rect} is
satisfied. As to condition (ii), it is trivial if $a$ and $b$
are comparable. If $a$ and $b$ are not comparable, then they are both
distinct from $c$ and
$d$ and Claim~1 applies. Thus (ii) is satisfied as
well and Lemma~\ref{L:rect} is proved.
\end{proof}

We can now prove the following result:

\begin{theorem}\label{T:VisV'}
Let $L$ be a BCF\index{lattice!BCF ---} lattice. Then there exists
an isomorphism from
$\DD L$ onto $\DD' L$ that
sends every
$\Dim(a,b)$ (for $a\leq b$ in $L$) to
$\|b-a\|$.
\end{theorem}

\begin{proof} By Lemma~\ref{L:rect}, every closed interval of
$L$ is rectifiable\index{rectifiable}, thus $\|b-a\|$ is
well-defined for all
$a\leq b$ in $L$. Then, it is easy to see that the function
$\vv<x,y>\mapsto\|y-x\|$ satisfies (D0) and
(D1), while (D2) follows from
Lemma~\ref{L:rect}.(ii). Therefore, there exists a unique
monoid homomorphism $\varphi\colon\DD L\to\DD' L$
such that for all $a\leq b$ in $L$,
$\varphi(\Dim(a,b))=\|b-a\|$.

Conversely, we prove that the relations (\ref{Eq:ShortRel1}),
(\ref{Eq:ShortRel2}), and (\ref{Eq:ShortRel3}) are satisfied by
the corresponding elements of $\DD L$. Indeed, let
\(\vv<\alpha_0,\beta_0,\alpha_1,\beta_1>\) be a caustic path%
\index{caustic!path} of
$L$, with base $\vv<a,b>$. Let
\(j\in\{1,\ldots,\mathrm{lh}(\beta_0)\}\). Then
\(a\wedge\beta_0(j)=a\wedge b\) and
\(a\vee\beta_0(j)=a\vee b\), thus
\(\Dim(a,a\vee b)=\Dim(a\wedge b,\beta_0(j))\). Therefore,
taking $j=1$, we obtain
\(\Dim(a,a\vee b)=\Dim(\beta_0(0),\beta_0(1))\).
Thus the relation
\(\Dim(\beta_0(j),\beta_0(j+1))\ll\Dim(\beta_0(0),\beta_0(1))\)
holds for all $j\geq 1$ in
\(\{1,\ldots,\mathrm{lh}(\beta_0)-1\}\).
Similarly with $\alpha_1$ instead of $\beta_0$, we obtain that
\(\Dim(a\wedge b,b)=\Dim(\alpha_1^*(i),a\vee b)\) for all
\(i\in\{1,\ldots,\mathrm{lh}(\alpha_1)\}\), so that
\(\Dim(a,a\vee b)=\Dim(a\wedge b,b)
=\Dim(\alpha_1^*(1),\alpha_1^*(0))\), thus
\(\Dim(\beta_0(0),\beta_0(1))=
\Dim(\alpha_1^*(1),\alpha_1^*(0))\). Hence the relations
(\ref{Eq:ShortRel1}), (\ref{Eq:ShortRel2}) and
(\ref{Eq:ShortRel3}) are satisfied by the
$\Dim(a,b)$, $a\prec b$ in $L$. So there exists a unique
monoid homomorphism
$\psi\colon \DD' L\to\DD L$ such
that \(\psi(\|b-a\|)=\Dim(a,b)\), for all $a\prec b$ in $L$.
Hence it is obvious that $\varphi$ and $\psi$ are mutually
inverse.
\end{proof}

\section{Further consequences of the
primitivity of BCF\index{lattice!BCF ---} lat\-tices}

This section will be mainly devoted to derive some consequences
of Theorem~\ref{T:VisV'}. We start it by pointing the connection
with Chapter~\ref{PrimMon}, which follows immediately from
Lemma~\ref{L:genrel}:

\begin{corollary}\label{C:dimprim}
The dimension monoid of a BCF\index{lattice!BCF ---} lattice is a
primitive\index{monoid!primitive ---} monoid.\qed
\end{corollary}

By using the results of Chapter~\ref{PrimMon}, in particular,
Theorem~\ref{T:list}, we obtain the following corollary:

\begin{corollary}\label{C:pptiesV(L)} Let $L$ be a BCF%
\index{lattice!BCF ---} lattice. Then
$\DD L$ is an antisymmetric
refinement monoid\index{monoid!refinement ---} satisfying the
interval axiom\index{interval axiom} and the
pseudo-cancellation\index{pseudo-cancellation} property;
furthermore, it admits an embedding (for the monoid
structure as well as for the ordering) into a power of
$\ZZb$.\qed
\end{corollary}

By composing the embedding of $\DD L$ into a power of $\ZZb$
and the projections, one obtains a family
\(\vv<D_p\mid p\in P>\) (where $P$ is the QO-system%
\index{quasi-ordered (QO) -system} of all
pseudo-indecomposable%
\index{pseudo-indecomposable (PI)} elements of
$\DD L$) of mappings (``\emph{dimension functions}'' of
$L$)\linebreak
\(D_p\colon\S(L)\to\ZZb\) that satisfy axioms (D0), (D1) and (D2)
introduced in Chapter~\ref{DimMon} and such that if
\(\vv<[a_i,\,b_i]\mid i<m>\) and
\(\vv<[c_j,\,d_j]\mid j<n>\) are finite sequences of
elements of $\S(L)$, then
\(\sum_i\Dim(a_i,b_i)\leq\sum_j\Dim(c_j,d_j)\) if and only if
\(\sum_iD_p(a_i,b_i)\leq\sum_jD_p(c_j,d_j)\) holds for all
$p\in P$ (and similarly
for the \emph{equality} of dimension words).

\begin{corollary}\label{C:shortsimple}
Let $L$ be a nontrivial simple%
\index{lattice!simple ---} BCF\index{lattice!BCF ---}
lattice. Then exactly one of the following two possibilities
can occur:

\begin{itemize}
\item[\rm (a)] $L$ is modular%
\index{lattice!modular (not necessarily complemented) ---} and
$\DD L\cong\ZZ^+$.

\item[\rm (b)] $L$ is not modular%
\index{lattice!modular (not necessarily complemented) ---} and
$\DD L\cong\two$.
\end{itemize}
\end{corollary}

\begin{proof} By definition, the fact that $L$ is
simple\index{lattice!simple ---} nontrivial means exactly that
\(\con L\cong\two\), whence also
\(\ccon L\cong\two\). Therefore, by
Corollary~\ref{C:congquotV}, $\DD L$ is a nontrivial simple%
\index{monoid!simple ---} conical\index{monoid!conical ---}
\cm,
that is, it is nonzero and any two nonzero elements are
equivalent modulo $\asymp$. On the other hand, there exists an
antisymmetric QO-system $\vv<P,\tr>$%
\index{quasi-ordered (QO) -system} such that
\(\DD L\cong\mathbf{E}(P)\).
If $p$ and $q$ are elements of $P$, then
\(e_p\asymp e_q\) in \(\mathbf{E}(P)\), thus $p=q$, so that
$|P|\leq1$; since $\DD L\ne\{0\}$, we obtain that $P$ is a
singleton, say
$P=\{p\}$. If \(p\in P^{(0)}\), that is, $2e_p=e_p$, then
$\DD L\cong\two$; by Proposition~\ref{P:shortmod}, $L$ cannot
be modular\index{lattice!modular (not necessarily
complemented) ---}. If
$p\in P^{(1)}$, then $\DD L\cong\ZZ^+$; by
Proposition~\ref{P:dirfin},
$L$ is modular\index{lattice!modular (not necessarily
complemented) ---}.
\end{proof}

\begin{corollary}\label{C:SimpleVmodShort} Let $L$ be a
simple\index{lattice!simple ---} lattice. Then the following
holds:
\begin{itemize}
\item[\rm (a)] $L$ has a prime interval and it
is modular
\index{lattice!modular (not necessarily complemented) ---}
if and only if \(\DD L\cong\ZZ^+\).

\item[\rm (b)] Suppose that $L$ is
\Vmod\index{lattice!V-modular ---} and that $L$ has a non
modular\index{lattice!modular (not necessarily complemented)
---} BCF%
\index{lattice!BCF ---} subinterval. Then
\(\DD L\cong\two\).
\end{itemize}
\end{corollary}

\begin{proof} Part (a) is exactly Corollary~\ref{C:succVhom}.
Now let us prove (b); thus let $K$ be a
BCF\index{lattice!BCF ---} non
modular
\index{lattice!modular (not necessarily complemented) ---}
subinterval of $L$. Then \(\DD K\) is not cancellative, thus the
set of relations defining \(\DD K\) contains at least one of the
relations (\ref{Eq:ShortRel2}) or (\ref{Eq:ShortRel3});
suppose that it is (\ref{Eq:ShortRel2}). Let
$\vv<\alpha_0,\beta_0,\alpha_1,\beta_1>$ the corresponding
caustic\index{caustic!path} path. Put \(u=\beta_0(0)\),
\(v=\beta_0(1)\), \(w=\beta_0(2)\),
\(\varepsilon=\Dim(u,v)\), \(\eta=\Dim(v,w)\). Then, by
(\ref{Eq:ShortRel2}), we have
\(\eta\ll\varepsilon\). On the other hand, it follows from
Corollary~\ref{C:simpleVL} that
\(\DD L=\ZZ^+\cdot\varepsilon=\ZZ^+\cdot\eta\). Thus there
exists \(k\in\NN\) with \(\varepsilon=k\eta\). Since
\(\eta\ll\varepsilon\), we have \(k\eta=(k+1)\eta\), thus, by
an easy induction, \(2k\eta=k\eta\). Since
\(k\eta=\varepsilon\), it follows that
\(2\varepsilon=\varepsilon\). Since
\(\DD L=\ZZ^+\cdot\varepsilon\), it follows that
\(\DD L\cong\two\).
\end{proof}

Now recall that a lattice $L$ is \emph{geometric}%
\index{lattice!geometric ---|ii}, if it is
algebraic\index{lattice!algebraic ---}, the compact%
\index{compact (element of a lattice)} elements are exactly
the finite suprema of atoms, and $L$ is (upper)
semimodular\index{lattice!semimodular ---} (see
\cite{Birk93,Grat} for more details).

\begin{corollary}\label{C:SimpleGeom}
Let $L$ be a
nontrivial simple\index{lattice!simple ---} geometric lattice.
Then exactly one of the following two possibilities can occur:
\begin{itemize}
\item[\rm (a)] $L$ is modular and $\DD L\cong\ZZ^+$.

\item[\rm (b)] $L$ is not modular and
$\DD L\cong\two$.
\end{itemize}
\end{corollary}

\begin{proof} Note first that $L$ is relatively complemented,
thus it is
\Vmod\index{lattice!V-modular ---}
(see Proposition~\ref{P:ExVmod}). Furthermore, let $K$ be the
sublattice of compact elements of $L$. Since
$L$ is not modular%
\index{lattice!modular (not necessarily complemented) ---},
neither is $K$. But $K$ is a BCF\index{lattice!BCF ---}
subinterval of
$L$. Therefore, by Corollary~\ref{C:SimpleVmodShort}.(b),
\(\DD L\cong\two\).
\end{proof}

In particular, we can immediately compute the dimension monoids
of all partition lattices: if, for every set $X$,
$\Pi_X$\index{PzzartX@$\Pi_X$|ii} denotes the lattice of all
partitions of $X$, then
\(\DD(\Pi_X)=\{0\}\), if $X=\varnothing$,
\(\DD(\Pi_X)\cong\ZZ^+\), if \(|X|\in\{2,3\}\), and
\(\DD(\Pi_X)\cong\two\), if $|X|\geq4$.

The presentation given by the relations
(\ref{Eq:ShortRel1}), (\ref{Eq:ShortRel2}), and
(\ref{Eq:ShortRel3}) also makes it possible to compute easily
dimension monoids of frequently encountered finite lattices:

\begin{example}\label{E:finiteex}
Let us first chose to
represent any antisymmetric QO-system%
\index{quasi-ordered (QO) -system} $P$ by the graph of the
ordering of $P$, where the elements of $P^{(0)}$ are marked
with filled circles (while those of $P^{(1)}$ are marked with
hollow circles). Now, for every positive integer $n$, let
$\mathfrak{C}_n$\index{Czzhainn@$\mathfrak{C}_n$|ii} denote the
$n$-element chain
$\{0,1,\ldots,n-1\}$. Using the presentation of
$\DD L\cong\DD' L$ given at the
beginning of this chapter, one obtains easily that
$\DD(\mathfrak{C}_n\amalg\mathfrak{C}_1)\cong\mathbf{E}(P_n)$
for all $n\in\{2,3\}$,
where $P_2$ and $P_3$ are the
(finite) antisymmetric QO-systems%
\index{quasi-ordered (QO) -system} whose graphs are
represented on Figure~\ref{Fig:QOsys}.

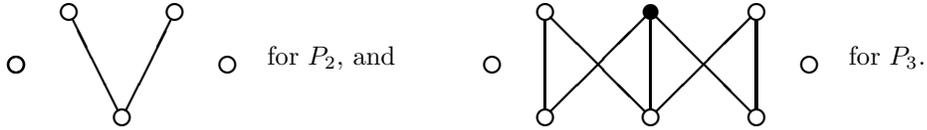
\begin{figure}[hbt]
\begin{picture}(250,80)(30,-20)
\thicklines
\put(38.5,2.6){\line(-1,2){17}}
\put(41.5,2.6){\line(1,2){17}}
\put(0,20){\circle{6}}
\put(80,20){\circle{6}}
\put(20,40){\circle{6}}
\put(60,40){\circle{6}}
\put(0,20){\circle{6}}
\put(40,0){\circle{6}}
\put(95,20){for $P_2$, and}
\put(180,20){\circle{6}}
\put(200,0){\circle{6}}
\put(200,40){\circle{6}}
\put(240,0){\circle{6}}
\put(240,40){\circle*{6}}
\put(280,0){\circle{6}}
\put(280,40){\circle{6}}
\put(300,20){\circle{6}}
\put(200,3){\line(0,1){34}}
\put(202.12,2.12){\line(1,1){35.76}}
\put(202.12,37.88){\line(1,-1){35.76}}
\put(240,3){\line(0,1){37}}
\put(240,40){\line(1,-1){37.88}}
\put(242.12,2.12){\line(1,1){35.76}}
\put(280,3){\line(0,1){34}}
\put(315,20){for $P_3$.}
\end{picture}
\caption{Some QO-systems}\label{Fig:QOsys}
\end{figure}

Recall that the pictures of
$\mathfrak{C}_2\amalg\mathfrak{C}_1$ and
$\mathfrak{C}_3\amalg\mathfrak{C}_1$ are respectively, from
the left to the right, represented on
Figure~\ref{Fig:FreeChain} (see \cite[Chapter VI]{Grat}).

\begin{figure}[hbt]
\begin{picture}(100,260)(80,-130)
\thicklines
\put(30,-80){\circle{6}}
\put(0,-50){\circle{6}}
\put(60,-50){\circle{6}}
\put(30,-20){\circle{6}}
\put(30,20){\circle{6}}
\put(110,0){\circle{6}}
\put(0,50){\circle{6}}
\put(60,50){\circle{6}}
\put(30,80){\circle{6}}

\put(27.88,-77.88){\line(-1,1){25.76}}
\put(2.12,-47.88){\line(1,1){25.76}}
\put(2.12,52.12){\line(1,1){25.76}}
\put(32.12,77.88){\line(1,-1){25.76}}
\put(57.88,47.88){\line(-1,-1){25.76}}
\put(62.12,47.88){\line(1,-1){45.76}}
\put(30,17){\line(0,-1){34}}
\put(27.88,22.12){\line(-1,1){25.76}}
\put(32.12,-22.12){\line(1,-1){25.76}}
\put(32.12,-77.88){\line(1,1){25.76}}
\put(62.12,-47.88){\line(1,1){45.76}}

\put(200,0){\circle{6}}
\put(200,100){\circle{6}}
\put(200,-100){\circle{6}}
\put(220,120){\circle{6}}
\put(220,-120){\circle{6}}
\put(230,70){\circle{6}}
\put(230,-70){\circle{6}}
\put(235,35){\circle{6}}
\put(235,-35){\circle{6}}
\put(250,20){\circle{6}}
\put(250,-20){\circle{6}}
\put(250,50){\circle{6}}
\put(250,-50){\circle{6}}
\put(250,90){\circle{6}}
\put(250,-90){\circle{6}}
\put(265,35){\circle{6}}
\put(265,-35){\circle{6}}
\put(285,55){\circle{6}}
\put(285,-55){\circle{6}}
\put(340,0){\circle{6}}

\put(202.12,2.12){\line(1,1){30.76}}
\put(202.12,-2.12){\line(1,-1){30.76}}
\put(202.12,102.12){\line(1,1){15.76}}
\put(202.12,-102.12){\line(1,-1){15.76}}
\put(202.12,97.88){\line(1,-1){25.76}}
\put(202.12,-97.88){\line(1,1){25.76}}
\put(222.12,117.88){\line(1,-1){25.76}}
\put(222.12,-117.88){\line(1,1){25.76}}
\put(232.12,67.88){\line(1,-1){15.76}}
\put(232.12,-67.88){\line(1,1){15.76}}
\put(232.12,72.12){\line(1,1){15.76}}
\put(232.12,-72.12){\line(1,-1){15.76}}
\put(237.12,37.12){\line(1,1){10.76}}
\put(237.12,32.88){\line(1,-1){10.76}}
\put(237.12,-37.12){\line(1,-1){10.76}}
\put(237.12,-32.88){\line(1,1){10.76}}
\put(250,-17){\line(0,1){34}}
\put(252.12,22.12){\line(1,1){10.76}}
\put(252.12,-22.12){\line(1,-1){10.76}}
\put(252.12,47.88){\line(1,-1){10.76}}
\put(252.12,-47.88){\line(1,1){10.76}}
\put(252.12,87.88){\line(1,-1){30.76}}
\put(252.12,-87.88){\line(1,1){30.76}}
\put(267.12,37.12){\line(1,1){15.76}}
\put(267.12,-37.12){\line(1,-1){15.76}}
\put(287.12,52.88){\line(1,-1){50.76}}
\put(287.12,-52.88){\line(1,1){50.76}}
\end{picture}
\caption{Finite lattices freely generated by chains}
\label{Fig:FreeChain}
\end{figure}
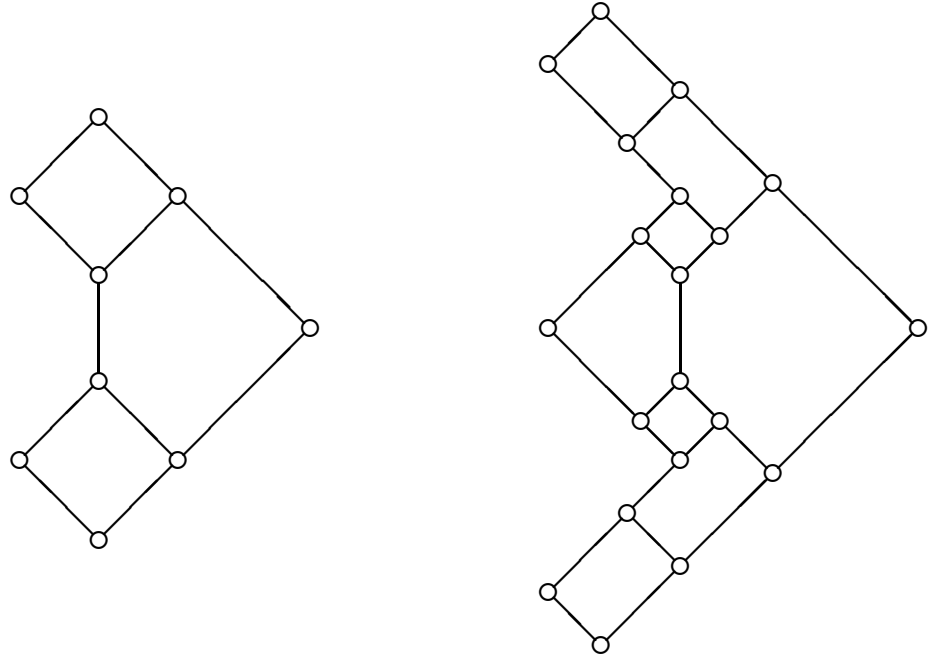

In particular, $\DD(\mathfrak{C}_1)=\{0\}$ and
$\DD(\mathfrak{C}_2)\cong\ZZ^+$, thus one sees immediately on
the example of $\mathfrak{C}_2\amalg\mathfrak{C}_1$ that the
$\DD$ functor does \emph{not} preserve arbitrary coproducts
(compare this with Corollary~\ref{C:presdiag}).
\end{example}

\section{The rectangular extension of a finite lattice}

As in \cite{GrSc}, we say that the \emph{rectangular extension}%
\index{rectangular extension|ii}%
\index{rzzectL@$\rect L$|ii}
of a finite lattice $L$ is the following direct product:
\[
\rect L=\prod_{p\in\MM(\con L)}L/p,
\]
where \(\MM(\con L)\) denotes the set of all (non-coarse)
meet-irreducible congruences of $L$. Hence, the natural
lattice homomorphism \(j_L\colon L\hookrightarrow\rect L\)
(defined by \(j_L(x)=\vv<[x]_p\mid p\in\MM(\con L)>\)) is
exactly the canonical subdirect decomposition of $L$ (into
subdirectly irreducible factors)---this is because $L$ is finite.
It is proved in \cite{GrSc} that $j_L$ has the
\emph{congruence extension property}%
\index{congruence!extension property|ii} (CEP), that is,
for every congruence $\alpha$ of $L$, there exists a (not
necessarily unique) congruence $\beta$ of $\rect L$ whose
inverse image under $j_L$ is equal to $\alpha$. This means
exactly that the semilattice homomorphism from \(\con L\) to
\(\con(\rect L)\) induced by $j_L$ is an embedding. This
suggests the following definition:

\begin{definition}\label{D:DEP} Let \(f\colon K\to L\) be a
lattice homomorphism. We say that $f$ has the
\emph{Dimension Extension Property (DEP)}
\index{dimension!extension property|ii}, if the monoid
homomorphism \(\DD f\colon \DD K\to\DD L\) is an embedding of
preordered monoids (for the algebraic preorderings of $\DD K$ and
$\DD L$).
\end{definition}

In particular, the DEP implies the CEP. We can now formulate
the following strengthening of the previously mentioned result
of G. Gr\"atzer and E.~T. Schmidt in \cite{GrSc}:

\begin{theorem}\label{T:SDDEP} Let $n$ be a non-negative
integer and let
\(f\colon L\hookrightarrow\prod_{i<n}L_i\) be a subdirect
decomposition of a BCF lattice $L$. Then $f$ has the DEP.
\end{theorem}

\begin{proof}
Put \(\widetilde{L}=\prod_{i<n}L_i\), and let
$\theta_i$ be the kernel of the natural surjective homomorphism
from $L$ onto $L_i$, for
all \(i<n\). Note first that the map \(\DD f\) is
given, modulo the natural identifications of
\(\DD(\widetilde{L})\) with \(\prod_{i<n}\DD(L_i)\)
(see Proposition~\ref{P:Vfunctor}) and of \(\DD(L_i)\) with
\(\DD L/\DD\theta_i\) (see Proposition~\ref{P:latt<>mon}) by the
rule
\[
\DD(f)(\alpha)=
\vv<\alpha\ \mathrm{mod}\,\DD\theta_i\mid i<n>.
\]
Since \(\DD L\) is antisymmetric
(see Corollary~\ref{C:pptiesV(L)}), it suffices to prove the
following statement:
\(\DD(f)(\alpha)\leq\DD(f)(\beta)\) implies that
\(\alpha\leq\beta\), for
all elements $\alpha$ and $\beta$ of \(\DD L\). The assumption
means that, for all
\(i<n\), there exists \(\gamma_i\in\DD\theta_i\) such that
\(\alpha\leq\beta+\gamma_i\). Since \(\DD L\) satisfies the
interval axiom\index{interval axiom}
(see Corollary~\ref{C:pptiesV(L)}), there exists \(\gamma\in\DD L\)
such that \(\alpha\leq\beta+\gamma\) and \(\gamma\leq\gamma_i\)
for all $i$. The latter statement implies that
\(\gamma\in\DD\theta_i\) holds for all \(i<n\), thus, if
\(\rho\colon \DD L\twoheadrightarrow\ccon L\) is the natural
homomorphism, \(\rho(\gamma)\leq\theta_i\) (use
Corollary~\ref{C:congquotV}). This holds for all
\(i<n\), thus, since $f$ is an embedding (so that
\(\bigcap_{i<n}\theta_i=\mathrm{id}_L\)),
\(\rho(\gamma)\) is the trivial congruence of $L$,
that is, \(\gamma=0\). It follows immediately that
\(\alpha\leq\beta\).
\end{proof}

\begin{corollary}\label{C:RectDEP} Let $L$ be a finite
lattice. Then the natural rectangular extension map $j_L$ has
the DEP\index{dimension!extension property}.\qed
\end{corollary}

\chapter[Sectionally complemented modular lattices]%
{Basic properties of sectionally complemented modular
lat\-tices}\label{BasicSCL}

We shall recall in this chapter some basic properties of sectionally
complemented modular lattices, emphasizing the relation of
\emph{perspectivity}\index{perspectivity}. Most of the results
presented here are well-known and come mainly from
the references \cite{FMae58,Neum60}, we restate them
here mainly for convenience. We shall omit the proofs in
most cases and refer the reader to the monographs
\cite{Birk93,Grat,FMae58,Neum60}.

\section{Independence and continuity}

In any lattice with zero, define, as usual, a partial addition
$\oplus$%
\index{ozzplus@$\oplus$!on elements of a lattice|ii} by
\(a\oplus b=c\) if and only if
\(a\wedge b=0\) and \(a\vee b=c\). Although $\oplus$ is always
commutative, it is not always associative, see
Chapter~\ref{MaxCommQuot} for the definition of associativity of a
partial operation. This is
shown very easily by the following folklore result:

\begin{proposition}\label{P:oplus} Let $L$ be a modular%
\index{lattice!modular (not necessarily complemented) ---}
lattice with zero. Then the partial addition $\oplus$ of
$L$ is associative. If $L$ is sectionally complemented, then
the converse holds.
\end{proposition}

\begin{proof} Suppose first that $L$ is modular%
\index{lattice!modular (not necessarily complemented) ---}. Let
$a$, $b$, and $c$ be elements of $L$ such that
\((a\oplus b)\oplus c\) is defined. Thus
\((a\oplus b)\wedge c=0\), thus \emph{a fortiori}
\(b\wedge c=0\),
that is,
\(b\oplus c\) is defined. Furthermore,
\begin{align*}
a\wedge(b\oplus c)&=a\wedge(b\vee c)\\
&=a\wedge(a\vee b)\wedge(b\vee c)\\
&=a\wedge\Bigl[b\vee\bigl((a\vee b)\wedge c\bigr)\Bigr]
\qquad\text{(by the modularity of $L$)}\\
&=a\wedge b\\
&=0,
\end{align*}
whence $a\oplus(b\oplus c)$ is defined. Of
course, both values $(a\oplus b)\oplus c$ and
$a\oplus(b\oplus c)$ equal $a\vee b\vee c$.

Conversely, suppose that $L$ is sectionally complemented and
that $\oplus$ is associative. Let $a$, $b$, $c\in L$ such that
$a\geq c$ and put $u=a\wedge(b\vee c)$ and
$v=(a\wedge b)\vee c$. Since $u\geq v$, there exists $w$ such
that $v\oplus w=u$. Let $b'$ such that
\((a\wedge b)\oplus b'=b\). It is easy to verify that
\([a\wedge b,\,b]\nearrow [u,\,b\vee c]\) and
\([a\wedge b,\,b]\nearrow[v,\,b\vee c]\), thus
\(b\vee c=u\oplus b'=v\oplus b'\). Therefore,
\(b\vee c=(w\oplus v)\oplus b'\). By assumption, we also have
\(b\vee c=w\oplus(v\oplus b')\), that is,
\(b\vee c=w\oplus(b\vee c)\). Hence $w=0$, that is, $u=v$.
Thus $L$ is modular%
\index{lattice!modular (not necessarily complemented)
---}.\end{proof}

Proposition~\ref{P:oplus} is the basis of the theory of
\emph{independence}\index{independent} in modular%
\index{lattice!modular (not necessarily complemented) ---}
lattices. Since the coming chapters deal with lattices
satisfying some completeness assumptions, it is convenient to
formulate this in a somewhat more general context.

\begin{definition}\label{D:MJContLatt} Let $L$ be a lattice
and let $\kappa$ be a cardinal number.

\begin{itemize}

\item $L$ is \emph{$\kappa$-\mcont}%
\index{lattice!meet-continuous ---!$\kappa$-\mcont\ ---|ii},
if for every
$a\in L$ and every increasing $\kappa$-sequence
$\vv<b_\xi\mid \xi<\kappa>$ of elements of $L$, the
supremum
$\bigvee_{\xi<\kappa}b_\xi$ exists and the following equality
\[
a\wedge\bigvee_{\xi<\kappa}b_\xi=
\bigvee_{\xi<\kappa}(a\wedge b_\xi).
\]
holds.

\item $L$ is \emph{$<\kappa$-\mcont}%
\index{lattice!meet-continuous ---!$<\kappa$-\mcont\ ---|ii},
if for every $\alpha<\kappa$, $L$ is
$\alpha$-\mcont.

\item $L$ is \emph{\mcont}%
\index{lattice!meet-continuous ---|ii}, if it is
$\kappa$-\mcont\
for every cardinal number $\kappa$.

\item One defines dually the concepts of
\emph{$\kappa$-join-con\-ti\-nu\-ity} and
\emph{$<\kappa$-join-con\-ti\-nu\-ity}%
\index{lattice!join-continuous ---!$<\kappa$-\jcont\ ---|ii}%
\index{lattice!join-continuous ---!$\kappa$-\jcont\ ---|ii}%
\index{lattice!join-continuous ---|ii}.

\item If {\rm (P)} is a lattice-the\-o\-ret\-i\-cal property,
we say that
$L$ is \emph{conditionally\index{conditionally (P)|ii}
{\rm (P)}}, if every bounded closed interval of $L$ has {\rm (P)}.

\end{itemize}
\end{definition}

\begin{note}
Although `continuity' of a given lattice has been defined, in early
times, as the conjunction of meet-continuity and join-continuity,
the reader has to be warned that in a somehow unfortunate way, this
is no longer the case (see \cite{GHKL80} for more details about
continuous%
\index{lattice!continuous ---} lattices). Note also that every
lattice is (vacuously)
$<\aleph_0$-\mcont\ as well as $<\aleph_0$-\jcont.
\end{note}

\begin{definition}\label{D:indep}
Let $L$ be a lattice with zero. Then a family
$\vv<a_i\mid i\in I>$ of elements of $L$ is
\emph{independent}\index{independent|ii}, if the map $\varphi$
from the generalized Boolean algebra
$[I]^{<\omega}$ of all finite subsets of $I$ to $L$
defined by \(\varphi(J)=\bigvee_{i\in J}a_i\) is a lattice
homomorphism.
\end{definition}

\begin{lemma}\label{L:indep}
Let $L$ be a lattice with zero. Then the following statements hold:

\begin{itemize}

\item[\rm (i)]
A family of elements of $L$ is independent if
and only if all its finite subfamilies are
independent.

\item[\rm (ii)] Suppose that $L$ is modular%
\index{lattice!modular (not necessarily complemented) ---}
and let $n$ be a positive integer. Then a finite sequence
$\vv<a_i\mid i\leq n>$ of elements of $L$ is
independent if and only if
$\vv<a_i\mid i<n>$ is independent and
\(a_n\wedge\bigvee_{i<n}a_i=0\).\qed

\end{itemize}
\end{lemma}

In particular, for modular%
\index{lattice!modular (not necessarily complemented) ---}
lattices with zero, the notation
$\oplus_{i<n}a_i$ is nonambiguous and the term that it
represents is defined if and only if the sequence
$\vv<a_i\mid i<n>$ is independent.

\begin{lemma}\label{L:GenAsso}
Let $\kappa$ be a cardinal number and let $L$ be a
$\kappa$-\mcont\
modular\index{lattice!modular (not necessarily complemented)
---} lattice. Let $I$ be a set of cardinality $<\kappa$ and let
$I=\bigcup_{j\in J}I_j$ be a partition of~$I$. Let
$\vv<a_i\mid i\in I>$ be a family of elements of $L$. Then
$\vv<a_i\mid i\in I>$ is independent
if and only if all subfamilies $\vv<a_i\mid i\in I_j>$
(for all $j\in J$) are independent and the
family
$\vv<b_j\mid j\in J>$ is independent,
where \(b_j=\bigvee_{i\in I_j}a_i\).\qed
\end{lemma}

\begin{lemma}\label{L:Cpleteindep}

Let $L$ be a modular%
\index{lattice!modular (not necessarily complemented) ---}
lattice with zero, let $\kappa$ be a cardinal number, let $n$
be a positive integer, and let $I$ be a set of cardinality at
most $\kappa$. Furthermore, let
\(\vv<a_i\mid i\in I>\) be an
independent family of elements of $L$. Then
the map $\psi$ from the product
$\prod_{i\in I}[0,\,a_i]$ to $L$ defined by the rule
\[
\psi(\vv<x_i\mid i\in I>)=\bigvee_{i\in I}x_i
\]
is a lattice embedding.\qed
\end{lemma}

\section{Perspectivity and projectivity}

\begin{definition}\label{D:persp}
Let $L$ be a lattice, let
$a$ and $b$ be elements of $L$.
\begin{itemize}
\item[(a)] $a$ is \emph{perspective}\index{perspectivity|ii} to
$b$, in notation
\(a\sim b\)%
\index{PzzerspElt@$\sim$ (for elements)|ii}, if there exists
\(x\in L\) such that
\(a\wedge x=b\wedge x\) and
\(a\vee x=b\vee x\).

\item[(b)] $a$ is
\emph{subperspective}\index{subperspectivity|ii} to $b$, in
notation
\(a\lesssim b\)%
\index{szzubpers@$\lesssim$|ii}, if there exists \(x\in L\)
such that
\(a\wedge x\leq b\wedge x\) and
\(a\vee x\leq b\vee x\).
\end{itemize}

As usual, the transitive closure of the relation of
perspectivity is called
\emph{projectivity}\index{projectivity!of elements|ii}, in
symbol $\approx$%
\index{PzzrojElt@$\approx$ (for elements)|ii}.

Furthermore, for every non-negative integer $n$, we shall
denote by $\sim_n$%
\index{Pzzersp2k@$\sim_k$|ii} (resp., $\lesssim_n$)%
\index{szzubpersk@$\lesssim_k$|ii} the $n$th power of $\sim$
(resp., $\lesssim$) as a binary relation (so that
\(\sim_0=\lesssim_0=\mathrm{id}_L\)).
\end{definition}

The following lemma originates in
\cite[Theorem I.3.1]{Neum60}.

\begin{lemma}\label{L:persp}
Let $L$ be a relatively complemented modular lattice, let $a$
and $b$ be elements of $L$, let
$c$ and $d$ be elements of $L$ such that $c\leq a\wedge b$
and $d\geq a\vee b$. Then $a$ and $b$ are perspective if and only it
there exists $x\in L$ such that
\begin{equation}
a\wedge x=b\wedge x=c\ \text{ and }\
a\vee x=b\vee x=d.\tag*{\qed}
\end{equation}
\end{lemma}

The following result is usually stated for complemented modular
lattices (see, for example, \cite[Theorem I.6.1]{Neum60}). However,
as we shall now show, it also remains true in the
general context of relatively complemented
(not necessarily modular) lattices
(see \cite[Lemma III.1.3 and Exercise III.1.3]{Grat}).

\begin{proposition}\label{P:lesssim}
Let $L$ be a relatively
complemented lattice, let $a$ and $b$ be elements of $L$. Then
the following are equivalent:
\begin{itemize}
\item[\rm (a)] \(a\lesssim b\);

\item[\rm (b)] there exists \(y\leq b\) such that \(a\sim y\);

\item[\rm (c)] there exists \(x\geq a\) such that \(x\sim b\).
\end{itemize}
\end{proposition}

\begin{proof}
To be relatively complemented is self-dual, thus
it suffices to prove that (a) and (b) are equivalent. It is
clear that (b) implies (a), so let us prove that (a) implies
(b). Let $a$ and $b$ be elements of
$L$ such that \(a\lesssim b\). By definition, there exists $c$
such that \(a\wedge c\leq b\wedge c\) and
\(a\vee c\leq b\vee c\). Let $x$ be a relative complement of
\(a\vee c\) in the interval \([c,\,b\vee c]\), let
$y$ be a relative complement of \(x\wedge b\) in the interval
\([a\wedge c,\,b]\). Then \(y\leq b\) and
\begin{align*}
y\wedge x&=y\wedge b\wedge x=a\wedge c=
a\wedge(a\vee c)\wedge x=a\wedge x,\\
y\vee x&=y\vee(x\wedge b)\vee x=
b\vee x=b\vee c\vee x=b\vee c=a\vee c\vee x=a\vee x.
\end{align*}
The computations can be followed on Figure~\ref{Fig:RelComp}.
Note that, unlike most figures
displayed in this work, Figure~\ref{Fig:RelComp} draws a
\emph{partial sublattice} of $L$, as opposed to a sublattice.
For example, there is no reason \emph{a priori}
for the elements $c\vee(x\wedge b)$ and $x$ to be equal.
On the other hand, it is the case that \((a\vee c)\wedge x=c\).
We mark a meet, $w=u\wedge v$, or a join, $w=u\vee v$, to
be correct by drawing an arc at $w$ in the triangle $u$, $v$,
$w$. This convention may be helpful, although it is sometimes
ambiguous.

\begin{figure}[hbt]
\begin{picture}(200,220)(-60,-80)
\thicklines
\put(0,0){\circle{6}}
\put(60,-60){\circle{6}}
\put(60,0){\circle{6}}
\put(120,0){\circle{6}}
\put(0,60){\circle{6}}
\put(30,30){\circle{6}}
\put(90,30){\circle{6}}
\put(120,60){\circle{6}}
\put(60,60){\circle{6}}
\put(60,120){\circle{6}}

\put(0,3){\line(0,1){54}}
\put(120,3){\line(0,1){54}}
\put(60,63){\line(0,1){54}}
\put(60,-57){\line(0,1){54}}
\put(2.12,62.12){\line(1,1){55.76}}
\put(62.12,-57.88){\line(1,1){55.76}}
\put(62.12,117.88){\line(1,-1){55.76}}
\put(2.12,-2.12){\line(1,-1){55.76}}
\put(62.12,2.12){\line(1,1){25.76}}
\put(32.12,32.12){\line(1,1){25.76}}
\put(92.12,32.12){\line(1,1){25.76}}
\put(32.12,27.88){\line(1,-1){25.76}}
\put(62.12,57.88){\line(1,-1){25.76}}
\put(2.12,57.88){\line(1,-1){25.76}}

\put(-5,0){\makebox(0,0)[r]{$a$}}
\put(60,-70){\makebox(0,0){$a\wedge c$}}
\put(60,130){\makebox(0,0){$b\vee c$}}
\put(-5,60){\makebox(0,0)[br]{$a\vee c$}}
\put(25,25){\makebox(0,0)[tr]{$c$}}
\put(55,0){\makebox(0,0)[r]{$b\wedge c$}}
\put(65,65){\makebox(0,0)[l]{$x$}}
\put(93,27){\makebox(0,0)[tl]{$x\wedge b$}}
\put(125,0){\makebox(0,0)[l]{$y$}}
\put(125,60){\makebox(0,0)[l]{$b$}}

\qbezier(45.86,105.86)(60,100)(74.14,105.86)
\qbezier(45.86,-45.86)(60,-40)(74.14,-45.86)
\qbezier(49.39,10.61)(60,15)(70.61,10.61)
\qbezier(19.39,40.61)(30,45)(40.61,40.61)
\qbezier(79.39,40.61)(90,45)(100.61,40.61)
\qbezier(0,45)(5.74,46.14)(10.61,49.39)
\qbezier(109.39,49.39)(114.26,46.14)(120,45)

\end{picture}
\caption{A partial configuration in $L$}\label{Fig:RelComp}
\end{figure}
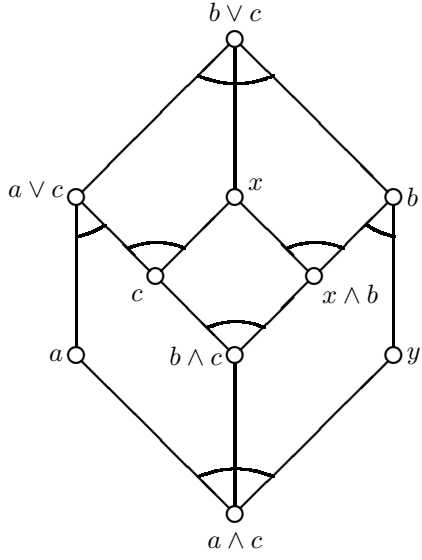

We obtain that \(a\sim y\) and \(y\leq b\).
\end{proof}

\begin{definition}\label{D:axis} Let $a$ and $b$ be
elements of a sectionally complemented modular lattice $L$.
Then an
\emph{axis of perspectivity} between $a$
and $b$ is any element $s$ such that \(a\oplus s=b\oplus s\);
in addition, we say that $s$ is \emph{proper}%
\index{proper!axis of perspectivity|ii}, if
\(a\oplus s=b\oplus s=a\vee b\). Moreover, the
\emph{perspective map}\index{perspective!isomorphism|ii} from
$[0,\,a]$ to $[0,\,b]$ with axis
$s$ is the map $\tau_{ab}^s$ from $[0,\,a]$ to $[0,\,b]$
defined by the rule
\[
(\forall x\leq a)
\bigl[\tau_{ab}^s(x)=(x\oplus s)\wedge b\bigr].
\]
\end{definition}

\begin{note}
If $s$ is an axis of
perspectivity\index{perspectivity} between $a$ and $b$, then
$s\wedge(a\vee b)$ is a proper axis of
perspectivity\index{perspectivity} between $a$ and $b$,
defining the same perspective
isomorphism.
\end{note}

\begin{definition}\label{D:ProjIso} Let $L$ be a sectionally
complemented modular lattice, let
\(n\in\omega\), and let \(\vv<a_i\mid i\leq n>\) be a
finite sequence of elements of $L$. A \emph{projective
isomorphism}%
\index{projective!isomorphism|ii} from \([0,\,a_0]\) onto
\([0,\,a_n]\) is any map of the form
\(T_{n-1}\circ\cdots\circ T_0\), where \(T_i\)
is a perspective isomorphism
from \([0,\,a_i]\) onto
\([0,\,a_{i+1}]\), for all \(i<n\).
\end{definition}

The following lemma summarizes the basic properties of
perspective maps:

\begin{lemma}\label{L:PerspMap} Let $L$ be a sectionally
complemented modular lattice, let $s$ be an axis of
perspectivity\index{perspectivity} between $a$ and $b$ in $L$.
Then
$\tau=\tau_{ab}^s$ is a lattice isomorphism from $[0,\,a]$ onto
$[0,\,b]$, with inverse $\tau_{ba}^s$. Moreover, for every
$x\leq a$, $s$ is an axis of
perspectivity\index{perspectivity} between $x$ and
$\tau(x)$ and the restriction of $\tau$ from $[0,\,x]$ to
$[0,\,\tau(x)]$ is equal to $\tau_{x,\tau(x)}^s$.\qed
\end{lemma}

Although there are many examples where
perspectivity\index{perspectivity} is not transitive,
transitivity appears in the following important special case,
see \cite[Theorem I.3.4]{Neum60}:

\begin{lemma}\label{L:TransOrthPersp}
Let $L$ be a sectionally
complemented modular lattice and let
$\vv<a,b,c>$ be a triple of elements of $L$ such that
\((a\vee b)\wedge(b\vee c)=b\). If $a\sim b$ and $b\sim c$,
then $a\sim c$. Moreover, if $S$ (resp., $T$) is a
perspective isomorphism from $[0,\,a]$ to
$[0,\,b]$ (resp., from $[0,\,b]$ to $[0,\,c]$), then $T\circ S$
is a perspective isomorphism from $[0,\,a]$ to $[0,\,c]$.
\end{lemma}

This lemma is usually stated in the case where
$\vv<a,b,c>$ is \emph{independent} but it
is sometimes used in the more general context above; see, for
example, the remark following Lemma~\ref{L:Jonsson}.

\begin{proof} By Lemma~\ref{L:persp}, there are elements $u$
and $v$ of $L$ such that \(a\oplus u=b\oplus u=a\vee b\) and
\(b\oplus v=c\oplus v=b\vee c\). Then put
\(w=(u\vee v)\wedge(a\vee c)\). One verifies, as usual, that
\(a\wedge w=c\wedge w=0\) and \(a\vee w=c\vee w=a\vee c\);
this means that
\begin{equation}\label{Eq:PerspAxis}
w=(u\vee v)\wedge(a\vee c)
\text{ is an axis of perspectivity between }
a\text{ and }c.
\end{equation}
Furthermore, let \(x\leq a\) and put \(y=S(x)\),
and \(z=T(y)\). We first prove that
\begin{equation}\label{Eq:QuasiOrth}
(x\vee y)\wedge(y\vee z)=y.
\end{equation}
Indeed, we have
\[
(x\vee y)\wedge(y\vee z)\leq(a\vee b)\wedge(b\vee c)=b,
\]
and moreover, by modularity,
\begin{align*}
x\vee y&=x\vee((x\oplus u)\wedge b)=
(x\oplus u)\wedge(x\vee b),\\
y\vee z&=y\vee((y\oplus v)\wedge c)=
(y\oplus v)\wedge(y\vee c),
\end{align*}
whence, using the fact that \(y=(x\oplus u)\wedge b\),
\begin{align*}
(x\vee y)\wedge(y\vee z)&=(x\vee y)\wedge(y\vee
z)\wedge b\\
&=(x\oplus u)\wedge(y\oplus v)\wedge
b\wedge(y\vee c)\\
&=y\wedge(y\oplus v)\wedge(y\vee c)\\
&=y,
\end{align*}
which proves (\ref{Eq:QuasiOrth}). In addition,
\(u_x=u\wedge(x\vee y)\) (resp.,
\(v_x=v\wedge(y\vee z)\)) is a proper axis of
perspectivity\index{perspectivity} between $x$ to $y$ (resp.,
between $y$ and $z$), thus, by (\ref{Eq:QuasiOrth}) (applied
to (\ref{Eq:PerspAxis})),
\(w_x=(u_x\vee v_x)\wedge(x\vee y)\) is an axis of
perspectivity\index{perspectivity} between $x$ and $z$. Since
\(w_x\leq w\) and \(x\wedge w=z\wedge w=0\), we also have
\(x\oplus w=z\oplus w\); whence \(T\circ S\) is the unique
perspective\index{perspective!isomorphism} isomorphism from
\([0,\,a]\) onto \([0,\,c]\) with axis $w$.
\end{proof}

Similarly, finite disjoint joins do not necessarily preserve
perspectivity\index{perspectivity}, but they do in the
following important particular case:

\begin{lemma}\label{L:BasicAddPersp} Let $\kappa$ be a
cardinal number and let $L$ be a
$\kappa$-\mcont%
\index{lattice!meet-continuous ---!$\kappa$-\mcont\ ---}
sectionally complemented modular lattice. Let $I$ be a set of
cardinality at most $\kappa$ and let
\(\vv<a_i\mid i\in I>\) and
\(\vv<b_i\mid i\in I>\) be families of elements of
$L$. If \(a_i\sim b_i\), for all \(i\in I\), and if the family
\(\vv<a_i\vee b_i\mid i\in I>\) is
independent, then
\(\oplus_{i\in I}a_i\sim\oplus_{i\in I}b_i\).\qed
\end{lemma}

\begin{corollary}\label{C:BasicAddPersp} Let \(\vv<a,b,c>\)
be an independent triple in a sectionally
complemented modular lattice. Then
\(a\sim b\) if and only if
\(a\oplus c\sim b\oplus c\).\qed
\end{corollary}

\chapter[Dimension for relatively complemented modular lattices]%
{Dimension monoids of relatively complemented modular
lat\-tices}\label{RelCompl}

This chapter starts with the simple observation that
the dimension function on a relatively complemented lattice with
zero is, really, a \emph{unary} function, because
if $a$ and $b$ are elements of any lattice $L$ with
zero, then, for any sectional complement $x$ of
$a\wedge b$ in $a\vee b$, the equality $\Dim(a,b)=\Dim(x)$ holds.
To obtain relations among those elements that define the dimension
monoid of $L$ is not that easy; the proof given in
Proposition~\ref{P:defrel} requires $L$ be relatively
complemented (but not necessarily modular). However, the further
consequences of this result will require, in general, modularity of
$L$. Among these consequences will be the fact that if $L$ is a
sectionally complemented modular lattice, then the range of the
unary dimension map $x\mapsto\Dim(x)=\Dim(0,x)$ is an \emph{ideal}
of $\DD L$, see Corollary~\ref{C:dimVmeas}, generalized to
relatively complemented modular lattices in
Corollary~\ref{P:splitV}.

Furthermore, it will turn out that for sectionally complemented
modular lattices, the equality of two dimension words,
$\sum_{i<m}\Dim(a_i)=\sum_{j<n}\Dim(b_j)$, can be conveniently
expressed in terms of the relation of \emph{projectivity by
decomposition}, see Definition~\ref{D:projdec}. A refining of
projectivity by decomposition, the relation of \emph{perspectivity
by decomposition}, will be used in Section~\ref{S:PerspDec} to
obtain further results about the action of the $\DD$ functor on
\emph{maps} (as opposed to lattices).

\section[Projectivity by decomposition]%
{Equality of dimension words; projectivity%
\index{projectivity!by decomposition} by decomposition}

The following simple result gives a more convenient
presentation of \(\DD L\) for a
relatively complemented (not necessarily modular) lattice with
$0$.

\begin{proposition}\label{P:defrel} Let $L$ be a relatively
complemented lattice with $0$. Then the following relations

\begin{itemize}
\item[\rm (D$'$0)]
\index{D'012@(D$'$0), (D$'$1), (D$'$2)|ii}
\(\Dim(0)=0\);

\item[\rm (D$'$1)]
\(\Dim(a\oplus b)=\Dim(a)+\Dim(b)\) (for all disjoint $a$,
$b\in L$);

\item[\rm (D$'$2)]
\(\Dim(a)=\Dim(b)\) (for all $a$, $b\in L$ such that \(a\sim b\))
\end{itemize}

\noindent are defining relations of the monoid%
\index{monoid}
\(\DD L\), and the dimension
range\index{dimension range} of $L$ equals
\(\{\Dim(x)\mid x\in L\}\).
\end{proposition}

\begin{proof}
Let $M$ be the monoid\index{monoid}
defined by generators $\|a\|$ (for $a\in L$) and relations
(D$'$0), (D$'$1) and (D$'$2).
We begin by proving that $M$ is commutative.
Let $a$, $b\in L$, let us prove that $\|a\|$ and
$\|b\|$ commute. By (D$'$1), this is true
if \(a\wedge b=0\). Now, if \(a\leq b\), then there exists $c$
such that
\(a\oplus c=b\); since \(a\wedge c=0\), \(\|a\|\) and
\(\|c\|\) commute, but
\(\|a\|\) commutes with itself, whence \(\|a\|\) commutes with
\(\|a\|+\|c\|\), that is, with \(\|b\|\). In the
general case, there exists $c$ such that
\((a\wedge b)\oplus c=b\); since \(a\wedge b\leq a\),
\(\|a\wedge b\|\) and
\(\|a\|\) commute; since \(c\wedge a=0\), \(\|c\|\)
commutes with
\(\|a\|\) as well; whence \(\|b\|=\|a\wedge b\|+\|c\|\)
commutes with \(\|a\|\).

Next, the unary function $\Dim$ satisfies (D$'$0), (D$'$1) and
(D$'$2). Indeed, observe that if \(a\wedge b=0\), then
\([0,\,a]\nearrow[b,\,a\oplus b]\), whence
\(\Dim(b,a\oplus b)=\Dim(a)\), so that
\(\Dim(a\oplus b)=\Dim(a)+\Dim(b)\). Moreover, if $a\sim b$,
that is, there exists $c$ such that
\(a\oplus c=b\oplus c\), then
\(\Dim(a)=\Dim(c,a\oplus c)=
\Dim(c,b\oplus c)=\Dim(b)\). This makes it possible to define a
monoid homomorphism
\(\psi\colon M\to\DD L\) such
that \(\psi(\|x\|)=\Dim(x)\) holds for all \(x\in L\).

Conversely, let $a\leq b$ in $L$. Then any two relative
complements $x$ and $y$ of $a$ in $b$ are
perspective\index{perspectivity} (with axis $a$), thus, by
(D$'$2), $\|x\|=\|y\|$. Note that we also
have $\Dim(a,b)=\Dim(x)$, whence the statement about the
dimension range. This makes it possible to define a map $\mu$ on
$\diag L$, by \(\mu(\vv<a,b>)=\|x\|\), for any $x$ such that
$a\oplus x=b$. We verify that $\mu$ satisfies
(D0), (D1) and (D2). Condition (D0) is trivially satisfied.
For all $a$, $b\in L$, let $x\in L$ such that
\((a\wedge b)\oplus x=b\). Then \(a\oplus x=a\vee b\), so that
\(\mu(\vv<a\wedge b,b>)=\|x\|=\mu(\vv<a,a\vee b>)\), so that
(D2) holds. Finally, let \(a\leq b\leq c\) in
$L$. Let $x$ be a relative complement of $a$ in \([0,\,b]\),
then let $z$ be a relative complement of $b$ in \([x,\,c]\) and
then, let $y$ be a relative complement of $x$ in \([0,\,z]\).
Then
\(a\oplus x=b\), thus \(\mu(\vv<a,b>)=\|x\|\). Further, we have
\begin{align*}
a\vee z&=a\vee x\vee z=b\vee z=c,\\
a\wedge z&=a\wedge b\wedge z=a\wedge x=0,
\end{align*}
thus \(a\oplus z=c\), whence
\(\mu(\vv<a,c>)=\|z\|\). Furthermore, we have
\begin{align*}
b\wedge y&=b\wedge z\wedge y=x\wedge y=0,\\
b\vee y&=b\vee x\vee y=b\vee z=c,
\end{align*}
thus \(b\oplus y=c\), whence
\(\mu(\vv<b,c>)=\|y\|\). But
\(z=x\oplus y\), whence we obtain that
\(\mu(\vv<a,c>)=\mu(\vv<a,b>)+\mu(\vv<b,c>)\). The computations
can be followed on Figure~\ref{Fig:DimRelComp}.

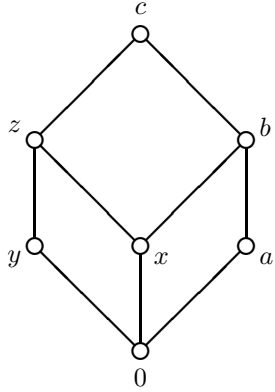
\begin{figure}[hbt]
\begin{picture}(100,140)(-20,-10)
\thicklines
\put(40,0){\circle{6}}
\put(40,40){\circle{6}}
\put(0,40){\circle{6}}
\put(0,80){\circle{6}}
\put(80,40){\circle{6}}
\put(80,80){\circle{6}}
\put(40,120){\circle{6}}

\put(42.12,2.12){\line(1,1){35.76}}
\put(37.88,2.12){\line(-1,1){35.76}}
\put(40,3){\line(0,1){34}}
\put(0,43){\line(0,1){34}}
\put(80,43){\line(0,1){34}}
\put(37.88,42.12){\line(-1,1){35.76}}
\put(42.12,42.12){\line(1,1){35.76}}
\put(2.12,82.12){\line(1,1){35.76}}
\put(77.88,82.12){\line(-1,1){35.76}}

\put(40,-10){\makebox(0,0){$0$}}
\put(-5,35){\makebox(0,0)[r]{$y$}}
\put(85,35){\makebox(0,0)[l]{$a$}}
\put(85,85){\makebox(0,0)[l]{$b$}}
\put(-5,85){\makebox(0,0)[r]{$z$}}
\put(45,35){\makebox(0,0)[l]{$x$}}
\put(40,130){\makebox(0,0){$c$}}

\end{picture}
\caption{Relative complements in $L$}\label{Fig:DimRelComp}
\end{figure}

So $\mu$ also satisfies (D1). Hence there exists
a unique monoid homomorphism
\(\varphi\colon \DD L\to M\)
such that
\(\varphi(\Dim(a,b))=\mu(\vv<a,b>)\) holds for all
\(a\leq b\) in $L$. In particular, \(\varphi(\Dim(x))=\|x\|\)
for all \(x\in L\). Since
\(\{\Dim(x)\mid x\in L\}\) generates \(\DD L\),
\(\varphi\) and \(\psi\) are mutually inverse.
\end{proof}

It follows, using the notation of Definition~\ref{D:dimsgrp},
that if $L$ is any sectionally complemented modular lattice,
then
$\DD L$ is naturally isomorphic
to $\DD(L,\oplus,\sim)$, where
$\oplus$, $\sim$ are as above. The
following lemma is basically well-known. For the
refinement\index{refinement!property} property, we refer for
an example to \cite[Theorem 4.1]{CrJo64} or
\cite[Proposition 1.1]{AGPO} (the context is slightly different
but adapts easily here) or to \cite[Theorem 3.1]{Halp39a} for
the much stronger infinite
refinement\index{refinement!property} property in
continuous\index{continuous!geometry} geometries.

In \ref{L:simref}--\ref{C:dimVmeas}, let
$L$ be a sectionally complemented modular lattice.

\begin{lemma}\label{L:simref}
The perspectivity\index{perspectivity} relation
$\sim$ is refining on $L$ and
$\vv<L,\oplus,\sim>$ satisfies the
refinement\index{refinement!property} property.
\end{lemma}

\begin{proof}
The fact that $\sim$ is refining is well-known
and follows, for example, immediately from
Lemma~\ref{L:PerspMap}. Now let $a$, $b$, $c$, and $d$ be
elements of $L$ such that $a\oplus b=c\oplus d$. Thus
$a\oplus b=a\vee c\vee d$. But there exist $c'\leq c$ and
$d'\leq d$ such that \(a\vee c=a\oplus c'\) and
\((a\vee c)\vee d=(a\vee c)\oplus d'\), whence
\(a\oplus b=a\oplus c'\oplus d'\). In particular,
\(b\sim c'\oplus d'\) thus there are \(\bar
c'\sim c'\) and
\(\bar d'\sim d'\) such that \(b=\bar c'\oplus\bar d'\). Let
$c''$ and $d''$ such that
\(c=c'\oplus c''\) and \(d=d'\oplus d''\). Then
\(a\oplus c'\oplus d'=a\oplus b=c\oplus d= c''\oplus d''\oplus
c'\oplus d'\), so that
\(a\sim c''\oplus d''\) and thus there are
$\bar c''\sim c''$ and $\bar d''\sim d''$ such that \(a=\bar
c''\oplus \bar d''\). Hence we have the following
$\sim$-refinement\index{refinement!matrix} matrix:
\[
\begin{tabular}{|c|c|c|}
\cline{2-3}
\multicolumn{1}{l|}{} & $c$ & $d$\tvi\\
\hline
$a$\tvi & $\bar c''/c''$ & $\bar d''/d''$\\
\hline
$b$\tvi & $\bar c'/c'$ & $\bar d'/d'$\\
\hline
\end{tabular}
\]
We conclude the proof by Lemma~\ref{L:refppty}.
\end{proof}

In particular, we can infer, using Theorem~\ref{T:refmon}
and Lemma~\ref{L:refwoO}, that $\DD L$ satisfies the
refinement\index{refinement!property} property, but this
already follows from Theorem~\ref{T:V(L)ref}. In fact, we are
interested in the additional information carried by the
$\DD(L,\oplus,\sim)$-representation.

\begin{definition}\label{D:projdec}
Let $\approxeq$ be the
binary relation (\emph{projectivity by decomposition})
\index{projectivity!by decomposition|ii}%
\index{PzzrojDec@$\approxeq$|ii}
on $L$ defined by \(a\approxeq b\) if and only if there are
decompositions \(a=\oplus_{i<n}a_i\), \(b=\oplus_{i<n}b_i\)
such that \(a_i\approx b_i\) for all \(i<n\).
\end{definition}

\begin{corollary}\label{C:eqwords} Let $\vv<a_i\mid i<m>$
and $\vv<b_j\mid j<n>$ be finite sequences of
elements of
$L$. Then $\sum_{i<m}\Dim(a_i)=\sum_{j<n}\Dim(b_j)$ if and only
if there exists a
$\approxeq$-refinement\index{refinement!matrix} matrix
of the following form:
\[
\begin{tabular}{|c|c|c|c|c|}
\cline{2-5}
\multicolumn{1}{l|}{} & $b_0$ & $b_1$ & $\ldots$ &
$b_{n-1}$\tvi\\
\hline
$a_0$\tvi & $c_{00}/d_{00}$ & $c_{01}/d_{01}$ & $\ldots$ &
$c_{0,n-1}/d_{0,n-1}$\\
\hline
$a_1$\tvi & $c_{10}/d_{10}$ & $c_{11}/d_{11}$ & $\ldots$ &
$c_{1,n-1}/d_{1,n-1}$\\
\hline
$\vdots$ & $\vdots$ & $\vdots$ & $\ddots$ & $\vdots$\\
\hline
$a_{m-1}$\tvi & $c_{m-1,0}/d_{m-1,0}$ & $c_{m-1,1}/d_{m-1,1}$ &
$\ldots$ & $c_{m-1,n-1}/d_{m-1,n-1}$\\
\hline
\end{tabular}
\]
\end{corollary}

\begin{proof} It is obvious that the given condition is
necessary. Conversely, suppose that
\(\sum_{i<m}\Dim(a_i)=\sum_{j<n}\Dim(b_j)\). To avoid
trivialities, suppose that $m>0$ and $n>0$. Putting
\(p=\vv<a_0,\ldots,a_{m-1}>\) and
\(q=\vv<b_0,\ldots,b_{n-1}>\) yields, if $\to$ and $\equiv$
are defined as in Lemma~\ref{L:basicto}, that
$p\equiv q$, that is, there exist a finite sequence
$r$ of elements of $L$ and $\alpha$, $\beta$ such that
$p\to^\alpha r$ and $q\to^\beta r$. Thus there are
$u_k$, $v_k\approx r(k)$ (for all $k\in\mathrm{dom}(r)$) such that
\(a_i=\oplus_{k\in{\alpha}^{-1}\{i\}}u_k\) and
\(b_j=\oplus_{k\in{\beta}^{-1}\{j\}}v_k\),
for all $i<m$ and all $j<n$. Now, for all
\(\vv<i,j>\in m\times n\), put
\(c_{ij}=\oplus_{k\in{\alpha}^{-1}\{i\}\cap
{\beta}^{-1}\{j\}}u_k\) and
\(d_{ij}=\oplus_{k\in{\alpha}^{-1}\{i\}\cap
{\beta}^{-1}\{j\}}v_k\), then the conclusion holds for
those values of $c_{ij}$ and~$d_{ij}$.
\end{proof}

\begin{corollary}\label{C:eqdim}
$\Dim(a)=\Dim(b)$ if and only if $a\approxeq b$, for all $a$ and $b$ in
$L$.\qed
\end{corollary}

Note that the following corollary is a special feature of
sectionally complemented modular lattices; indeed, it is
not valid in general modular lattices (even for finite chains
with at least three elements). We recall that
\emph{\Vmeas s} have been defined by
Dobbertin, see \cite{Dobb83}:

\begin{corollary}\label{C:dimVmeas}
The dimension map
$\Dim\colon L\to\DD L$ is a
\Vmeas\index{Vmeas@\Vmeas|ii},
that is, if $c\in L$ and
$\alpha$, $\beta\in\DD L$ such that
$\Dim(c)=\alpha+\beta$, then there exist $a$, $b\in L$ such
that $\alpha=\Dim(a)$,
$\beta=\Dim(b)$, and $c=a\oplus b$. In particular, the dimension
range of $L$, that is, $\{\Dim(x)\mid x\in L\}$, is a
lower subset of $\DD L$.
\end{corollary}

\begin{proof}
There are finite sequences $\vv<a_i\mid i<m>$ and
$\brack{b_j\mid j<n}$ of elements of $L$ such that
$\alpha=\sum_{i<m}\Dim(a_i)$ and
$\beta=\sum_{j<n}\Dim(b_j)$, therefore, by
Lemma~\ref{C:eqwords}, there exists a
$\approxeq$-refinement\index{refinement!matrix} matrix as
follows:
\[
\begin{tabular}{|c|c|}
\cline{2-2}
\multicolumn{1}{l|}{} & $c$\tvi\\
\hline
$a_i\ (i<m)$ & $a_i/x_i$\tvi\\
\hline
$b_j\ (j<n)$ & $b_j/y_j$\tvi\\
\hline
\end{tabular}
\]
Thus \(c=\oplus_{i<m}x_i\oplus\oplus_{j<n}y_j\).
In particular,
$a=\oplus_{i<m}x_i$ and $b=\oplus_{j<n}y_j$ are defined.
Furthermore, the following equalities hold:
\begin{gather*}
\alpha=\sum_{i<m}\Dim(a_i)=\sum_{i<m}\Dim(x_i)=\Dim(a),\\
\beta=\sum_{j<n}\Dim(b_j)=\sum_{j<n}\Dim(y_j)=\Dim(b),\\
c=a\oplus b.\tag*{\qed}
\end{gather*}
\renewcommand{\qed}{}
\end{proof}

This result can be generalized right away to relatively
complemented modular lattices:

\begin{proposition}\label{P:splitV}
Let $L$ be a relatively
complemented modular lattice, let
$u\leq v$ in $L$, let $\alpha$, $\beta$ be elements of
$\DD L$ such that
$\Dim(u,v)=\alpha+\beta$. Then there exists
$x\in[u,\,v]$ such that $\alpha=\Dim(u,x)$ and
$\beta=\Dim(x,v)$. In particular, the dimension
range\index{dimension range} of $L$ is a lower subset of
$\DD L$.
\end{proposition}

\begin{proof}
Note that $L$ is the direct union of its closed
intervals (viewed as convex sublattices), thus, since the
$\DD$ functor preserves direct limits
(see Proposition~\ref{P:Vfunctor}), it suffices to prove the result
for a \emph{bounded} lattice $L$, say with bounds $0$
and $1$ (note indeed that every convex sublattice of $L$ is
still relatively complemented modular). Since $L$ is
sectionally complemented, there exists
$c$ such that $u\oplus c=v$, so that $\Dim(u,v)=\Dim(c)$. By
Corollary~\ref{C:dimVmeas}, there are $a$ and $b$ in $L$ such
that $c=a\oplus b$, $\Dim(a)=\alpha$, and $\Dim(b)=\beta$;
whence
\(\Dim(u,u\oplus a)=\alpha\) and \(\Dim(u\oplus a,v)=\beta\).
Hence the conclusion holds for $x=u\oplus a$.
\end{proof}

\section[Perspectivity by decomposition]%
{Lattice ideals\index{ideal!of a semilattice};
perspectivity by decomposition}
\label{S:PerspDec}

Let us recall first that if $L$ is a sectionally complemented
modular lattice, a
\emph{neutral ideal}%
\index{ideal!neutral --- (in lattices)|ii} of
$L$ is an ideal\index{ideal!of a semilattice} $I$ of $L$ such
that $x\sim y$ and $y\in I$ implies that
$x\in I$. Then it is well-known that
congruences\index{congruence!lattice ---} of $L$ are in
one-to-one correspondence with neutral
ideals of $L$, in the
following way: with a congruence\index{congruence!lattice ---}
$\theta$ of $L$, associate the neutral
ideal \(\{x\in L\mid\vv<0,x>\in\theta\}\); conversely, with every
neutral ideal $I$ of
$L$, associate the congruence\index{congruence!lattice ---}
$\theta$ defined by
\[
\vv<x,y>\in\theta\Leftrightarrow (\exists u\in I)
(x\vee u=y\vee u).
\]
The details can be found in \cite{Birk93} or
\cite{Grat}. If $x$ is an element of $L$, we shall denote by
$\Theta(x)$%
\index{Tzzhetax@$\Theta(x)$|ii} the neutral
ideal generated by
$x$. It can easily be seen that the natural map from
$\DD L$ onto
$\ccon L$ (see
Corollary~\ref{C:congquotV}) sends
$\Dim(x)$ to $\Theta(x)$.

\begin{corollary}\label{C:V(I)inV(L)}
Let $L$ be a sectionally
complemented modular lattice and let $I$ be a neutral ideal of $L$.
Then the natural embedding from $\DD I$ into $\DD L$ is a
\Vemb\index{Vemb@\Vemb}. Furthermore, if we identify
$\DD I$ with its image in $\DD L$ under this isomorphism,
we have \(\DD(L/I)\cong\DD L/\DD I\).
\end{corollary}

\begin{proof} Let $f$ be the natural homomorphism from
$\DD I$ to $\DD L$. By
Proposition~\ref{P:exV-hom}, $f$ is a \Vhom\index{Vhom@\Vhom}.
Therefore, in order to prove that $f$ is a
\Vemb\index{Vemb@\Vemb}, it suffices, by Lemma~\ref{L:cancref},
to prove that the restriction of $f$ to the dimension range of
$I$ is one-to-one. Thus let $a$ and $b$ be elements of $I$; we
prove that
\(\Dim_L(a)=\Dim_L(b)\) implies that \(\Dim_I(a)=\Dim_I(b)\).
Indeed, by Corollary~\ref{C:eqdim}, the hypothesis means that
\(a\approxeq b\) in $L$. However, since $I$ is a
neutral ideal of $L$, it
is easy to see that this is equivalent to the fact that
\(a\approxeq b\) in $I$. Thus, $f$
is a \Vemb\index{Vemb@\Vemb}.

Furthermore, by Proposition~\ref{P:latt<>mon},
\(\DD(L/I)\cong\DD L/\DD\theta\), where
$\theta$ is the congruence\index{congruence!lattice ---}
associated with $I$. But
\(\Dim(x)=\Dim(0,x)\in\DD\theta\) for all $x\in I$, so that
\(\DD I\subseteq\DD\theta\);
conversely, for all elements
$x$ and $y$ of $L$ such that \(x\equiv y\pmod{I}\), there
exists \(u\in L\) such that
\((x\wedge y)\oplus u=x\vee y\), thus \(u\in I\), so that
\(\Dim(x,y)=\Dim(0,u)\in\DD I\): hence we have proved that
\(\DD\theta=\DD I\). It follows that
\(\DD(L/I)\cong\DD L/\DD I\).
\end{proof}

However, there are examples (see
Corollary~\ref{C:V(K)notembV(L)}) of sectionally complemented
modular lattices $L$ and convex sublattices $K$ of $L$ such
the natural map $\DD K\to\DD L$
is not one-to-one. To conclude this chapter, we shall give a
class of lattices in which this strange behaviour does not
happen.

\begin{definition}\label{D:persdec} Let $L$ be a modular%
\index{lattice!modular (not necessarily complemented) ---}
lattice with $0$. Then elements $a$ and $b$ of $L$ are said
to be \emph{perspective by decomposition}%
\index{perspectivity by decomposition|ii}, which we shall
write $a\simeq b$%
\index{PzzersDec@$\simeq$|ii}, if there are
$n\in\omega$ and decompositions \(a=\oplus_{i<n}a_i\),
\(b=\oplus_{i<n}b_i\) such that $a_i\sim b_i$ for all $i<n$.
\end{definition}

Thus perspectivity by decomposition $\simeq$
trivially implies projectivity by decomposition%
\index{projectivity!by decomposition} $\approxeq$; furthermore, the
converse holds if and only if $\simeq$ is transitive. This
is, in fact, easily seen to be equivalent to the fact that
$x\sim_2y$ implies that $x\simeq y$, for all
$x$, $y$. We shall meet large classes of lattices, in
Chapter~\ref{CtbleMeetCo}, where this
happens.

Proposition~\ref{P:Vemb} yields a whole class of \Vemb s%
\index{Vemb@\Vemb} between dimension monoids of sectionally
complemented modular lattices:

\begin{proposition}\label{P:Vemb}
Let $L$ be a sectionally
complemented modular lattice, let $K$ be an
ideal\index{ideal!of a semilattice} of $L$. Suppose that
perspectivity by decomposition is transitive in
$L$. Then the natural map \(f\colon\DD K\to\DD L\) is a
\Vemb\index{Vemb@\Vemb}.
\end{proposition}

\begin{proof} By Proposition~\ref{P:exV-hom}, $f$ is a
\Vhom\index{Vhom@\Vhom}. Thus, by Corollary~\ref{C:dimVmeas}
and Lemma~\ref{L:onetoone}, it suffices to prove that the
restriction of $f$ to the dimension range of $K$ is one-to-one.
Thus let $x$ and $y$ be elements of $K$ such that
\(f(\Dim_K(x))=f(\Dim_K(y))\). This means, of course, that
$\Dim_L(x)=\Dim_L(y)$. By Corollary~\ref{C:eqdim}, we have
$x\approxeq y$ in $L$, thus, by
assumption, $x\simeq y$, so that there are $n\in\omega$ and
decompositions \(x=\oplus_{i<n}x_i\) and
\(y=\oplus_{i<n}y_i\) such that $x_i\sim y_i$ in $L$, for all
$i<n$. Since $K$ is an ideal\index{ideal!of a semilattice} of
$L$, both $x_i$ and $y_i$ belong to $K$ and thus they are also
perspective\index{perspectivity} in $K$. Therefore,
$x\simeq y$ in $K$, so that \(\Dim_K(x)=\Dim_K(y)\).
\end{proof}

\chapter[Normal equivalences]%
{Normal equivalences; dimension
monoids of regular\index{ring!(von Neumann) regular ---}
rings}\label{RegRings}

If $R$ is any regular ring, then the lattice $\mathcal{L}(R_R)$ of
all principal right ideals of $R$ is a complemented modular
lattice. Thus one can define on $\mathcal{L}(R_R)$ the two
relations of perspectivity ($\sim$) and projectivity ($\approx$),
that may or may not be the same. Another important equivalence
relation on $\mathcal{L}(R_R)$ is the \emph{isomorphy} relation
($\cong$), implied by the relations $\sim$ and $\approx$ above,
but, in general, strictly larger.

Now comes an important point: the isomorphy relation on
$\mathcal{L}(R_R)$ is defined in module-theoretical terms,
but it cannot, \emph{a priori}, be defined in lattice-theoretical
terms. This may suggest that there is no nontrivial dimension
theory of non-coordinatizable complemented modular lattices.
However, one moment's reflection shows that this guess is too
pessimistic: indeed, non-Arguesian projective planes have a perfectly
understood (and trivial) dimension theory, that does not rely at all
on the possibility of coordinatization!

This is the main motivation for introducing a substitute to the
isomorphy relation for lattices that are not necessarily
coordinatizable. This will be achieved by extracting a
lattice-theoretical property of the relation of isomorphy, which we
will call \emph{normality}. As it will turn out, there may be more
than one normal equivalence on a given complemented modular
lattice, see Corollary~\ref{C:RegRingCounterex}; there may also
be no such relation, see Section~\ref{S:NonNorm}. Normal lattices
will be defined as the lattices with at least one normal
equivalence. The theory of normal equivalences will culminate in
Theorem~\ref{T:MeetContNor}, which implies, in particular, that
every countably meet-continuous complemented modular lattice is
normal.

\section{Normal equivalences; normal lat\-tices}

We start with the definition of a normal equivalence.

\begin{definition}\label{D:normal}
Let $L$ be a modular lattice with zero. A
\emph{normal equivalence}\index{normal!equivalence|ii} on
$L$ is an equivalence relation $\equiv$ on $L$ satisfying the
following properties:

\begin{itemize}
\item[\rm (E1)]\index{E1@(E1), (E2), (E3)|ii}
$\equiv$ is additive and refining.

\item[\rm (E2)]
$x\sim y$ implies that $x\equiv y$, for all $x$, $y\in L$.

\item[\rm (E3)]
$x\perp y$ and $x\equiv y$ implies that $x\sim y$, for all $x$,
$y\in L$.
\end{itemize}

We say that a modular lattice $L$ with zero is
\emph{normal}\index{normal!lattice|ii}, if there exists at
least one normal equivalence on $L$.
\end{definition}

The reference \cite[page 98]{HeHu74} contains another,
completely unrelated, definition of normality for modular (not
necessarily complemented) lattices, but the contexts are
different enough to avoid any risk of confusion.

A fundamental class of normal modular lattices with
zero is given by the following result, contained in
\cite[Proposition 4.22]{Good91}.

\begin{lemma}\label{L:SubmNorm}
Let $M$ be a right module over a ring $R$. Then the relation of
module isomorphy on the lattice
$\mathcal{L}(M)$ of all submodules of $M$ is a normal equivalence
on $\mathcal{L}(M)$. Thus, $\mathcal{L}(M)$ is normal.
\end{lemma}

\begin{proof}
It is obvious that the relation of module isomorphy, $\cong$,
satisfies (E1) and (E2). Now we prove (E3). Let $A$ and $B$ be
submodules of
$M$, such that $A\cap B=\{0\}$ and $A\cong B$, we prove that
$A\sim B$. Let $f\colon A\to B$ be an isomorphism. Put
 \[
 C=\{x+f(x)\mid x\in A\}.
 \]
Then $C\in\mathcal{L}(M)$ and $A\oplus C=B\oplus C$.
\end{proof}

\begin{example}
Note that, for submodules $A\subseteq B$ and
$A'\subseteq B'$ of a module $M$, $\Dim(A,B)=\Dim(A',B')$ does
\emph{not} necessarily imply that $B/A\cong B'/A'$. For example,
for $R=\ZZ$ and $M=\ZZ\oplus\ZZ$, take $A=\{0\}\oplus\{0\}$,
$B=\ZZ\oplus\{0\}$, and $A'=2\ZZ\oplus\{0\}$, $B'=\ZZ\oplus\ZZ$.
Note that
$\{0\}\oplus\ZZ$ and $2\ZZ\oplus\{0\}$ are isomorphic and
independent, thus, by Lemma~\ref{L:SubmNorm}, they are
perspective; whence
 \[
 \Dim(\{0\}\oplus\{0\},\{0\}\oplus\ZZ)=
 \Dim(\{0\}\oplus\{0\},2\ZZ\oplus\{0\}).
 \]
It follows that
 \begin{align*}
  \Dim(A',B')&=
  \Dim(2\ZZ\oplus\{0\},\ZZ\oplus\{0\})+
  \Dim(\ZZ\oplus\{0\},\ZZ\oplus\ZZ)\\
  &=\Dim(2\ZZ\oplus\{0\},\ZZ\oplus\{0\})+
  \Dim(\{0\}\oplus\{0\},\{0\}\oplus\ZZ)\\
  &=\Dim(2\ZZ\oplus\{0\},\ZZ\oplus\{0\})+
  \Dim(\{0\}\oplus\{0\},2\ZZ\oplus\{0\})\\
  &=\Dim(\{0\}\oplus\{0\},\ZZ\oplus\{0\})\\
  &=\Dim(A,B),
 \end{align*}
while $B/A\cong\ZZ$ and $B'/A'\cong(\ZZ/2\ZZ)\oplus\ZZ$ are not
isomorphic.
\end{example}

In any sectionally complemented modular lattice,
projectivity by decomposition%
\index{projectivity!by decomposition} $\approxeq$ is an
equivalence relation satisfying both (E1) and (E2)
but not necessarily (E3), as proved by the
counterexample of Section~\ref{S:NonNorm}.

\begin{lemma}\label{L:charnorm}
Let $L$ be a sectionally complemented modular lattice. Then the
following are equivalent:

\begin{itemize}
\item[\rm (i)] There exists a normal\index{normal!equivalence}
equivalence on $L$.

\item[\rm (ii)] $x\approx y$ and
$x\perp y$ implies that $x\sim y$, for all $x$, $y\in L$.

\item[\rm (iii)] The relation $\approxeq$ of projectivity%
\index{projectivity!by decomposition} by decomposition is a
normal\index{normal!equivalence} equivalence on~$L$.
\end{itemize}
\end{lemma}

\begin{proof} (i)$\Rightarrow$(ii). By (E2), it
is trivial that every normal\index{normal!equivalence}
equivalence on $L$ contains $\approx$. The conclusion
follows then from (E3).

(ii)$\Rightarrow$(iii). The only nontrivial point is to verify
that $\approxeq$ satisfies
(E3). Thus let $x$, $y\in L$ such that
$x\approxeq y$ and $x\perp y$. By
definition, there are decompositions $x=\oplus_{i<n}x_i$ and
$y=\oplus_{i<n}y_i$ such that
$x_i\approx y_i$ for all $i<n$. Since
$x_i\perp y_i$, it follows from the hypothesis that
$x_i\sim y_i$. Since $x\perp y$, the
finite sequence
$\vv<x_i\vee y_i\mid i<n>$ is
independent. Therefore, by
Lemma~\ref{L:BasicAddPersp},
\(\oplus_{i<n}x_i\sim\oplus_{i<n}y_i\), that is, $x\sim y$.

(iii)$\Rightarrow$(i) is trivial.\end{proof}

Thus a sectionally complemented modular lattice is
normal\index{normal!lattice} if and only if any two
independent projective%
\index{projectivity!of elements} elements are
perspective\index{perspectivity}.

A convenient sufficient (but not necessary) condition for
normality\index{normal!lattice} is the following one:

\begin{proposition}\label{P:SuffNorm}
Let $L$ be a sectionally
complemented modular lattice. If perspectivity by
decomposition\index{perspectivity by decomposition} is
transitive in $L$, then $L$ is normal.
\end{proposition}

\begin{proof} Let $a$, $b\in L$ such that $a\approx b$ and
$a\perp b$. By assumption, there are decompositions
\(a=\oplus_{i<n}a_i\) and
\(b=\oplus_{i<n}b_i\) such that
$a_i\sim b_i$, for all $i<n$. It follows from
Lemma~\ref{L:BasicAddPersp} that $a\sim b$.
\end{proof}

Definition~\ref{D:normal} of normal
lattices can be extended to lattices without zero, by noting
the following result:

\begin{lemma}\label{L:sublNorm} Let $L$ be a
normal sectionally complemented modular
lattice and let $K$ be a convex sublattice of $L$ with a least
element. Then $K$ is normal.
\end{lemma}

\begin{proof}
Denote by $0_K$, $\perp_K$, and \(\oplus^K\)
($0_L$, $\perp_L$, and \(\oplus^L\), resp.)
the zero element, independence
relation and independent sum in $K$ (resp., $L$). Since $K$ is
a convex sublattice of $L$, the relation of
perspectivity\index{perspectivity} in $K$ is the restriction to
$K$ of the relation of perspectivity\index{perspectivity} in
$L$. Let $x$, $y$ be elements of $K$ such that
\(\Dim_K(x)=\Dim_K(y)\) and
\(x\perp_K y\). There are elements $x'$ and $y'$ of $L$ such
that \(x=0_K\oplus x'\) and
\(y=0_K\oplus^Ly'\). Thus we have \(\Dim_L(x')=\Dim_L(y')\) and
\(x'\perp_Ly'\); therefore, since $L$ is
normal, \(x'\sim y'\). By
Corollary~\ref{C:BasicAddPersp},
it follows that \(x\sim y\). So we have
proved that $K$ is normal.
\end{proof}

This shows that the following definition is an extension of
Definition~\ref{D:normal}:

\begin{definition}\label{D:GenNorm}
A relatively complemented modular lattice $L$ is
\emph{normal}\index{normal!lattice|ii}, if every bounded
interval of $L$ is normal.
\end{definition}

\begin{corollary}\label{C:GenNorm} Every convex sublattice of
a normal relatively complemented modular
lattice is normal.\qed
\end{corollary}

The following proposition lists some basic preservation
results of normality:

\begin{proposition}\label{P:PresNorm} The class of all
normal relatively complemented modular
lattices is closed under the following operations:
\begin{itemize}
\item[\rm (a)] direct limits;

\item[\rm (b)] reduced products;

\item[\rm (c)] homomorphic images.
\end{itemize}
\end{proposition}

\begin{proof} (a) and (b) are obvious. Now, let us prove (c).
Let $\theta$ be a congruence\index{congruence!lattice ---} of a
normal relatively complemented modular
lattice $L$, we prove that \(L/\theta\) is
normal. Since
\(L/\theta\) is the direct limit of all lattices
\(M/\theta_M\), where $M$ ranges over all closed intervals of
$L$ and where we put \(\theta_M=\theta\cap(M\times M)\), it
suffices to consider the case where $L$ is complemented
modular. Then $\theta$ is the
congruence\index{congruence!lattice ---} associated with some
neutral ideal\index{ideal!neutral --- (in lattices)} $I$ of
$L$. Let \(\pi\colon L\twoheadrightarrow L/I\) be the canonical
projection.

\setcounter{claim}{0}
\begin{claim}
Let $a$ and $b$ be elements of
$L$. If \(\pi(a)\sim\pi(b)\) in \(L/I\), then there
exist elements $u$ and $v$ of $L$ such that
\begin{gather*}
u\leq a,\qquad v\leq b;\\
\pi(u)=\pi(a),\qquad\pi(v)=\pi(b);\\
u\sim v.
\end{gather*}
\end{claim}

\begin{cproof}
By definition, there exists
\(c\in L\) such that both \(a\wedge c\) and \(b\wedge c\)
belong to $I$ and \(a\vee c\equiv b\vee c\pmod{I}\). The latter
means that there exists \(d\in I\) such that
\(a\vee c\vee d=b\vee c\vee d\). But
\(a\wedge(c\vee d)\equiv a\wedge c\pmod{I}\), thus
\(a\wedge(c\vee d)\in I\); similarly,
\(b\wedge(c\vee d)\in I\), so that one can replace $c$ by
\(c\vee d\) and thus assume that both \(a\wedge c\) and
\(b\wedge c\) belong to $I$ and that \(a\vee c=b\vee c\). Now
let $u$ (resp., $v$) be any sectional complement of
\(a\wedge c\) (resp., \(b\wedge c\)) in $a$ (resp., $b$). Then it
is immediate that \(u\oplus c=v\oplus c=a\vee c=b\vee c\),
whence $u$ and $v$ satisfy the required conditions.
\end{cproof}

Then one deduces easily by induction the following claim:

\begin{claim}
Let $a$ and $b$ be elements of
$L$. If \(\pi(a)\approx\pi(b)\) in \(L/I\), then there exist
elements $u$ and $v$ of $L$ such that
\begin{gather*}
u\leq a,\qquad v\leq b;\\
\pi(u)=\pi(a),\qquad\pi(v)=\pi(b);\\
u\approx v.\tag*{\qedc}
\end{gather*}
\end{claim}

Now let \(\boldsymbol{a}\) and \(\boldsymbol{b}\) be elements
of \(L/I\) such that \(\boldsymbol{a}\perp\boldsymbol{b}\) and
\(\boldsymbol{a}\approx\boldsymbol{b}\) in \(L/I\). By
Claim~2 above, there are elements $a$ of
$\boldsymbol{a}$ and $b$ of $\boldsymbol{b}$ such that
\(a\approx b\). Let $T$ be a
projective\index{projective!isomorphism} isomorphism from
\([0,\,a]\) onto \([0,\,b]\). Since $\boldsymbol{a}$ and
$\boldsymbol{b}$ are independent, we have \(a\wedge b\in I\).
Now let $u$ be a sectional complement of \(a\wedge b\) in $a$,
and put \(v=T(u)\); thus \(u\approx v\). Furthermore,
\(u\perp b\), thus \(u\perp v\). Since $L$ is
normal\index{normal!lattice}, it follows that
\(u\sim v\). But \(a\equiv
u\pmod{I}\) and
\(b\equiv v\pmod{I}\), thus
\(\boldsymbol{a}\sim\boldsymbol{b}\).
\end{proof}

There are quite large classes of normal\index{normal!lattice}
sectionally complemented modular lattices. For example, every
principal ideal\index{ideal!of a ring} lattice of a
von~Neumann  regular\index{ring!(von Neumann) regular ---}
ring is normal\index{normal!lattice}
(\cite[Proposition 4.22]{Good91}; see Lemma~\ref{L:SubmNorm}).
So is every countably complete sectionally complemented modular
lattice that is $\aleph_0$-\mcont%
\index{lattice!meet-continuous ---!$\aleph_0$-\mcont\ ---}
as well as $\aleph_0$-\jcont%
\index{lattice!join-continuous ---!$\aleph_0$-\jcont\ ---}%
---in fact, in this context,
perspectivity\index{perspectivity} is transitive, see
\cite{Halp38}, but the verification of
normality\index{normal!equivalence} of
$\approxeq$ is much easier than the
proof of the full transitivity of
perspectivity\index{perspectivity}, while it provides a
stepping stone towards a proof of the latter result. Finally,
any complete and \mcont%
\index{lattice!meet-continuous ---} sectionally complemented
modular lattice is normal\index{normal!lattice} (this is a
result of von~Neumann and Halperin, see \cite{HaNe40}; see also
\cite[Satz II.3.7]{FMae58}). Wider classes of
normal\index{normal!lattice} lattices will be encountered
throughout this work.

\section[Poset of normal equivalences]%
{The partially ordered
set of normal\index{normal!equivalence} equivalences}
\label{S:PosetNorm}

For every sectionally complemented modular lattice $L$, denote
by \(\NEq(L)\)\index{NzzEqL@$\NEq(L)$|ii} the set of all
normal\index{normal!equivalence} equivalences on $L$,
partially ordered under inclusion. The following trivial
statement summarizes what we know at this point about
\(\NEq(L)\):

\begin{proposition}\label{P:BasicNEq} Let $L$ be a sectionally
complemented modular lattice. Then the following properties
hold:
\begin{itemize}
\item[\rm (i)] If \(\NEq(L)\) is
nonempty, then it has a minimal element (namely,
$\approxeq$).

\item[\rm (ii)] The union of any nonempty directed subset of
\(\NEq(L)\) belongs to
\(\NEq(L)\).\qed
\end{itemize}
\end{proposition}

We shall now study the effect of $\NEq$
on the \emph{dual} lattice of a lattice.
A weaker form of Theorem~\ref{T:NormEqDual}, stating that the dual
of a normal\index{normal!lattice} complemented modular lattice is
normal, was communicated to the author by Christian Herrmann.

Let us first state a useful (and straightforward) lemma:

\begin{lemma}\label{L:ModCube}
Let $L$ be a bounded modular
lattice, let $a$, $b$, and $c$ be elements of $L$ such that
\((a\wedge b)\oplus c=1\) and \(a\vee b=1\). Put
\(a'=a\wedge c\) and \(b'=b\wedge c\). Then
\(\vv<a\wedge b,a',b'>\) is independent and the following equalities hold:
\begin{equation*}
a=(a\wedge b)\oplus a';\qquad
b=(a\wedge b)\oplus b';\qquad
c=a'\oplus b'.\tag*{\qed}
\end{equation*}
\end{lemma}

\begin{figure}[hbt]
\begin{picture}(200,200)(-40,-100)
\thicklines
\put(0,25){\circle{6}}
\put(50,25){\circle{6}}
\put(50,75){\circle{6}}
\put(100,25){\circle{6}}

\put(0,-25){\circle{6}}
\put(50,-25){\circle{6}}
\put(50,-75){\circle{6}}
\put(100,-25){\circle{6}}

\put(2.12,27.12){\line(1,1){45.76}}
\put(2.12,22.88){\line(1,-1){45.76}}
\put(52.12,22.88){\line(1,-1){45.76}}
\put(0,-22){\line(0,1){44}}
\put(2.12,-22.88){\line(1,1){45.76}}
\put(2.12,-27.12){\line(1,-1){45.76}}
\put(50,-72){\line(0,1){44}}
\put(52.12,-72.88){\line(1,1){45.76}}
\put(52.12,22.88){\line(1,-1){45.76}}
\put(52.12,-22.88){\line(1,1){45.76}}
\put(52.12,72.88){\line(1,-1){45.76}}
\put(50,28){\line(0,1){44}}
\put(100,-22){\line(0,1){44}}

\put(-5,25){\makebox(0,0)[r]{$a$}}
\put(-5,-25){\makebox(0,0)[r]{$a'$}}
\put(55,-25){\makebox(0,0)[l]{$a\wedge b$}}
\put(105,-25){\makebox(0,0)[l]{$b'$}}
\put(105,25){\makebox(0,0)[l]{$b$}}
\put(55,25){\makebox(0,0)[l]{$c$}}
\put(50,85){\makebox(0,0){$a\vee b=1$}}
\put(50,-85){\makebox(0,0){$a'\wedge b'=0$}}
\end{picture}
\caption{A cube configuration}\label{Fig:Cube}
\end{figure}
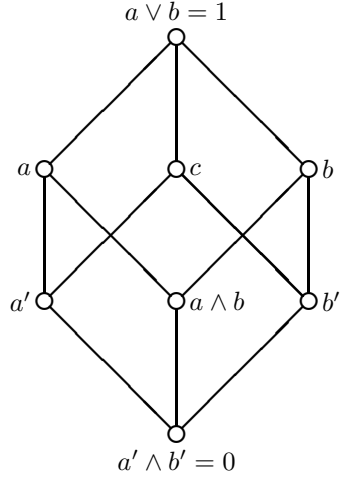

The configuration of the parameters in Lemma~\ref{L:ModCube}
is shown on Figure~\ref{Fig:Cube}.

Now let $L$ be a complemented modular lattice. Let
$\alpha$ be an equivalence relation on $L$. For all
$x$, $y\in L$, we shall write \(x\equiv_\alpha y\) for
\(\vv<x,y>\in\alpha\). If $\alpha$ is a
normal\index{normal!equivalence} equivalence on $L$, let us
define a binary relation $\alpha^*$ on $L$ by putting
\begin{equation}\label{Eq:Defadual}
\vv<x,y>\in\alpha^*\Leftrightarrow(\exists x',y')
(x\oplus x'=y\oplus y'=1\text{ and }x'\equiv_\alpha y').
\end{equation}

Note that if $x''$ (resp., $y''$) is another complement of $x$
(resp., $y$), then \(x'\sim x''\) and \(y'\sim y''\), thus we
also have \(x''\equiv_\alpha y''\). Hence the definition
(\ref{Eq:Defadual}) of $\alpha^*$ can be replaced by the
following:

\begin{equation}\label{Eq:Defadual2}
\vv<x,y>\in\alpha^*\Leftrightarrow(\forall x',y') (x\oplus
x'=y\oplus y'=1\Rightarrow x'\equiv_\alpha y').
\end{equation}

It follows easily from this that $\alpha^*$ is an equivalence
relation on $L$.

\begin{lemma}\label{L:AlphastNorm}
The equivalence $\alpha^*$
is a normal\index{normal!equivalence} equivalence on the dual
lattice \(L^\mathrm{op}\).
\end{lemma}

\begin{proof} Before starting the proof, note that both
perspectivity\index{perspectivity} and
projectivity\index{projectivity!of elements} are self-dual;
furthermore, elements $a$ and $b$ of $L$ are independent in
$L^\mathrm{op}$ if and only if \(a\vee b=1\) (and then their
``independent sum'' in $L^\mathrm{op}$ is
\(a\wedge b\)). Now we check the defining properties
successively, by repeated use of Lemma~\ref{L:ModCube}.
\smallskip

(E1) Let first $a_i$, $b_i$ (for \(i<2\)) be
elements of $L$ such that \(a_0\equiv_{\alpha^*}a_1\),
\(b_0\equiv_{\alpha^*}b_1\), and  \(a_i\vee b_i=1\)
(for all \(i<2\)). Let $c_i$ (for \(i<2\)) such that
\((a_i\wedge b_i)\oplus c_i=1\) in $L$ and put
\(a'_i=a_i\wedge c_i\) and \(b'_i=b_i\wedge c_i\). By
Lemma~\ref{L:ModCube}, we have \(a_i\oplus b'_i=1\), thus,
since \(a_0\equiv_{\alpha^*}a_1\), we have
\(b'_0\equiv_{\alpha}b'_1\). One can prove similarly that
\(a'_0\equiv_{\alpha}a'_1\). Since $\alpha$ is additive, one
obtains that
\[
c_0=a'_0\oplus b'_0\equiv_\alpha a'_1\oplus b'_1=c_1.
\]
Since \((a_i\wedge b_i)\oplus c_i=1\), we obtain that
\(a_0\wedge b_0\equiv_{\alpha^*}a_1\wedge b_1\). Thus
$\alpha^*$ is additive.

Next, let $a_0$, $a_1$, and $b$ in $L$ such that
\(a_0\vee a_1=1\) and \(a_0\wedge a_1=b\). Let $a'$ such that
\((a_0\wedge a_1)\oplus a'=1\) and put \(a'_i=a_i\wedge a'\)
(for all \(i<2\)). Since \(a_0\wedge a_1\equiv_{\alpha^*}b\), we
have, by definition, \(a'\equiv_\alpha b'\). Since
\(a'=a'_0\oplus a'_1\) and $\alpha$ is refining, there are
$b'_0$ and $b'_1$ such that \(b'=b'_0\oplus b'_1\) and
\(a'_i\equiv_\alpha b'_i\) (for all \(i<2\)). Note that
\(\vv<b,b'_0,b'_1>\) is independent with join $1$. Put
\(b_i=b\oplus b'_i\) (for all \(i<2\)). Then
\(b_i\oplus b'_{1-i}=a_i\oplus a'_{1-i}=1\) holds for all \(i<2\),
thus, since \(a'_{1-i}\equiv_\alpha b'_{1-i}\), we obtain that
\(a_i\equiv_{\alpha^*}b_i\). Furthermore, \(b_0\wedge b_1=b\);
this proves that $\alpha^*$ is refining.
\smallskip

(E2) Let \(a\sim b\) in $L^\mathrm{op}$ (thus \(a\sim b\) in
$L$). There exists $c$ such that \(a\oplus c=b\oplus c=1\).
Since
\(c\equiv_\alpha c\), we also have \(a\equiv_{\alpha^*}b\).
\smallskip

(E3) Let \(a\equiv_{\alpha^*}b\) be such that
\(a\vee b=1\). Let
$c$ such that \((a\wedge b)\oplus c=1\) and put
\(a'=a\wedge c\) and \(b'=b\wedge c\). Then
\(a\oplus b'=b\oplus a'=1\), thus \(a'\equiv_\alpha b'\).
Since \(a'\perp b'\) and $\alpha$ is
normal\index{normal!equivalence}, we obtain that
\(a'\sim b'\). Thus, since \(\vv<a\wedge b,a',b'>\) is
independent and by Corollary~\ref{C:BasicAddPersp},
\((a\wedge b)\oplus a'\sim(a\wedge b)\oplus b'\), that is,
\(a\sim b\).
\end{proof}

Once this is proved, it is easy to verify that
\(\alpha\mapsto\alpha^*\) is order-preserving (for the
inclusion), and that \(\alpha^{**}=\alpha\). Therefore, we
obtain the following result:

\begin{theorem}\label{T:NormEqDual}
Let $L$ be a complemented
modular lattice. Then the partially ordered sets
\(\NEq(L)\) and
\(\NEq(L^\mathrm{op})\) are
isomorphic.\qed
\end{theorem}

\begin{corollary}\label{C:NormEqDual1} Let $L$ be a relatively
complemented modular lattice. If $L$ is
normal\index{normal!lattice}, then so is the dual lattice
\(L^\mathrm{op}\).\qed
\end{corollary}

\begin{corollary}\label{C:NormEqDual2} Let $L$ be a
complemented modular lattice. If there exists a unique
normal\index{normal!equivalence} equivalence on $L$, then
there exists a unique normal\index{normal!equivalence}
equivalence on~\(L^\mathrm{op}\).\qed
\end{corollary}

\section{Normality\index{normal!lattice} and dimension
embedding for ideals\index{ideal!of a semilattice}}

As we shall see in this section, normality---more precisely,
uniqueness of the normal equivalence---can be used to prove
``dimension embedding'' theorems, such as Proposition~\ref{P:conv}.
In the same spirit, it will follow from Corollary~\ref{C:uniqnor}
that in many coordinatizable lattices, the relation of isomorphy on
ideals \emph{can} be defined in lattice-theoretical terms, more
precisely, in terms of the square $\sim_2$ of the relation of
perspectivity.

We start with a variant of Proposition~\ref{P:Vemb} for the case
where $K$ is not necessarily an ideal\index{ideal!of a lattice} of
$L$, but merely a convex sublattice:

\begin{proposition}\label{P:conv} Let $L$ be a relatively
complemented modular lattice and let
$K$ be a convex sublattice of $L$ with a smallest element. We
make the following additional assumptions:

\begin{itemize}
\item[\rm (i)] $L$ is normal\index{normal!lattice}.

\item[\rm (ii)] There exists at most one
normal\index{normal!equivalence} equivalence on $K$.
\end{itemize}
\noindent Then there exists exactly one
normal\index{normal!equivalence} equivalence on $K$ and the
natural map\linebreak
\(f\colon \DD K\to\DD L\) is a \Vemb\index{Vemb@\Vemb}.
\end{proposition}

\begin{proof}
By Proposition~\ref{P:exV-hom}, $f$ is a
\Vhom\index{Vhom@\Vhom}. Let
$\equiv$ be the equivalence relation on
$K$ defined by the rule
\[
x\equiv y\Longleftrightarrow\Dim_L(0_K,x)=\Dim_L(0_K,y).
\]

\setcounter{claim}{0}
\begin{claim}
The relation $\equiv$ is a
normal\index{normal!equivalence} equivalence on $K$.
\end{claim}

\begin{cproof}
The equality
\[
\Dim_K(0_K,x)+\Dim_K(0_K,y)=
\Dim_K(0_K,x\wedge y)+\Dim_K(0_K,x\vee y)
\]
holds for all $x$, $y\in K$,
thus, applying $f$, we obtain that
\[
\Dim_L(0_K,x)+\Dim_L(0_K,y)=
\Dim_L(0_K,x\wedge y)+\Dim_L(0_K,x\vee y).
\]
In particular, if $x\perp_Ky$, then
\(\Dim_L(0_K,x)+\Dim_L(0_K,y)=\Dim_L(0_K,x\vee y)\). This
implies immediately that $\equiv$ is additive.

Let us now prove that $\equiv$ is refining. Thus let
$a_0$, $a_1$, $b\in K$ be such that \(a_0\oplus^Ka_1\equiv b\),
that is,
\(\Dim_L(0_K,a_0)+\Dim_L(0_K,a_1)=f(\Dim_K(0_K,b))\). Since $f$
is a \Vhom\index{Vhom@\Vhom}, there are elements $\alpha_i$
(for $i<2$) of $\DD K$ such that
\(\Dim_K(0_K,b)=\alpha_0+\alpha_1\) and
\(f(\alpha_i)=\Dim_L(0_K,a_i)\) (for all $i<2$). Since $\Dim_K$ is
a \Vmeas\index{Vmeas@\Vmeas}\ (see Corollary~\ref{C:dimVmeas}),
there are elements
$b_0$ and $b_1$ of $K$ such that \(\alpha_i=\Dim_K(0_K,b_i)\)
(for all $i<2$) and
\(b=b_0\oplus^Kb_1\). Thus,
\(\Dim_L(0_K,b_i)=f(\alpha_i)=\Dim_L(0_K,a_i)\) for all $i<2$,
that is,
\(a_i\equiv b_i\). This completes the proof that $\equiv$ is
refining.

Let $x$, $y\in K$ such that \(x\sim_Ky\). This implies that
\(\Dim_K(0_K,x)=\Dim_K(0_K,y)\), thus, applying~$f$,
\(\Dim_L(0_K,x)=\Dim_L(0_K,y)\), that is, \(x\equiv y\).
Thus $\equiv$ contains $\sim_K$.

Finally, let $x$, $y\in K$ be such that \(x\perp_K y\) and
\(x\equiv y\). Since $L$ is the direct union of all its closed
intervals, there exists an interval
\(L'=[a,\,b]\) of $L$ such that \(a\leq 0_K\leq x\vee y\leq b\)
and \(\Dim_{L'}(0_K,x)=\Dim_{L'}(0_K,y)\). Let $x'$ (resp.,
$y'$) be a relative complement of $0_K$ in the interval
\([a,\,x]\) (resp., \([a,\,y]\)); we thus obtain that
\(x'\perp_{L'}y'\) and
\(\Dim_{L'}(a,x')=\Dim_{L'}(a,y')\). But $L'$ is
normal\index{normal!lattice}, thus
\(x'\sim y'\). Therefore, by Corollary~\ref{C:BasicAddPersp},
\(x\sim y\).
\end{cproof}

It follows from Claim~1 above and the assumption on
$K$ that there exists exactly one
normal\index{normal!equivalence} equivalence on $K$;
furthermore, $\equiv$ is exactly the relation of projectivity
by decomposition%
\index{projectivity!by decomposition} in $K$,
that is,
\(\Dim_L(0_K,x)=\Dim_L(0_K,y)\) if and only if
\(\Dim_K(0_K,x)=\Dim_K(0_K,y)\), for all $x$, $y\in K$.
Thus, the restriction of $f$
to the dimension range of $K$ is one-to-one. By
Corollary~\ref{C:dimVmeas} and Lemma~\ref{L:onetoone}, it
follows that $f$ is one-to-one.
\end{proof}

\begin{lemma}\label{L:decomp}
Let $L$ be a sectionally
complemented modular lattice. Let $n\in\NN$ and let $a$,
$u_i$ (for $i<n$) be elements of $L$ such that
\(a\leq\oplus_{i<n}u_i\). Put
\(\bar u_i=\oplus_{j<i}u_j\), for all $i\leq n$; furthermore,
let $a_i$ satisfy
\((a\wedge\bar u_i)\oplus a_i=a\wedge\bar u_{i+1}\),
for all $i<n$.
Then \(a=\oplus_{i<n}a_i\) and
\(a_i\perp\bar u_i\) and \(a_i\leq\bar u_{i+1}\), for all $i<n$.
\end{lemma}

\begin{proof} The second part of the statement is immediate.
For the first part, it is easy to prove by induction on
$k\leq n$ that
\(a\wedge\bar u_k=\oplus_{i<k}a_i\); the conclusion follows by
taking $k=n$.\end{proof}

\begin{proposition}\label{P:uniqnor}
Let $L$ be a
complemented modular lattice, let $\equiv$ be a
normal equivalence on $L$. Let $n\geq 2$ in $\NN$
and let
\(\vv<u_i\mid i<n>\), \(\vv<u'_i\mid i<n>\)
be finite sequences of
elements of $L$ satisfying the following properties:

\begin{itemize}
\item[\rm (i)] \(\oplus_{i<n}u_i=1\);

\item[\rm (ii)] $u_i\equiv u'_i$
and $u'_i\leq\oplus_{j\ne i}u_j$, for all $i<n$.
\end{itemize}

Then $a\equiv b$ if and
only if there are decompositions
\(a=\oplus_{i<n^2}a_i\), \(b=\oplus_{i<n^2}b_i\)
such that for all $i<n^2$, \(a_i\sim_2b_i\), for all $a$, $b$ in $L$.
\end{proposition}

\begin{proof}
First, if
$a$ and $b$ are decomposed as above, then
$a\approxeq b$, thus $a\equiv b$
since $\equiv$ is an additive equivalence relation
containing $\sim$.

Conversely, put, as in Lemma~\ref{L:decomp},
\(\bar u_i=\oplus_{j<i}u_j\) (for all $i\leq n$).
>Furthermore, define a subset $S$ of $L$ by
\[
S=\{x\in L\mid (\exists i<n) (x\leq\bar u_{i+1}
\text{ and }x\perp\bar u_i)\}.
\]
Now let
$a$, $b\in L$ such that $a\equiv b$. By using twice
Lemma~\ref{L:decomp} and the fact that $\equiv$ is
refining, it is easy to verify that there are
decompositions
\(a=\oplus_{i<n^2}a_i\) and \(b=\oplus_{i<n^2}b_i\)
such that all $a_i$, $b_i$'s lie in $S$ and
$a_i\equiv b_i$, for all $i$. Therefore, to
prove the conclusion,
it is sufficient to prove that if $a$ and $b$ are
elements of $S$, then $a\equiv b$ implies that
$a\sim_2b$. Let $i$, $j<n$ such that
\(a\leq\bar u_{i+1}\), \(a\perp\bar u_i\),
\(b\leq\bar u_{j+1}\), and \(b\perp\bar u_j\); we
prove that $a\sim_2b$. Without loss of generality,
$i\leq j$. If
$i<j$, then $a\leq\bar u_j$, so that
$b\perp\bar u_j$ implies that $a\perp b$, whence,
since $\equiv$ is normal\index{normal!equivalence},
$a\sim b$. Now suppose
that $i=j$. Since
\(a\oplus\bar u_i\leq\bar u_{i+1}
=u_i\oplus\bar u_i\), $a$ is
perspective\index{perspectivity} to a part of $u_i$.
Thus, by assumption, there is
$x\leq u'_i$ such that $a\equiv x$. Since
$a\equiv b$, we also have $b\equiv x$. In addition,
\(a\leq\oplus_{k\leq i}u_k\) and
\(x\leq\oplus_{k\ne i}u_k\) and the $u_k$'s are
independent, whence
\(a\wedge x\leq\oplus_{k<i}u_k=\bar u_i\). But
\(a\perp\bar u_i\), whence $a\wedge x=0$. Since
$\equiv$ is normal\index{normal!equivalence},
$a\sim x$. Similarly,
$b\sim x$. It follows
that $a\sim_2b$.
\end{proof}

By applying Proposition~\ref{P:uniqnor} to the
relation of projectivity by decomposition%
\index{projectivity!by decomposition}
$\approxeq$, one obtains immediately the following
consequence:

\begin{corollary}\label{C:uniqnor} Let $L$ be a
complemented modular lattice. Let $n\geq 2$ in $\NN$
and let
\(\vv<u_i\mid i<n>\) be a finite sequence of
elements of $L$ satisfying the following properties:

\begin{itemize}
\item[\rm (i)] \(\oplus_{i<n}u_i=1\);

\item[\rm (ii)]
\(\Dim(u_i)\leq\sum_{j\ne i}\Dim(u_j)\), for all $i<n$.
\end{itemize}

Then there exists at most one normal equivalence on
$L$. Furthermore, if such an equivalence
$\equiv$ exists, then
$a\equiv b$ if and only if there are decompositions
\(a=\oplus_{i<n^2}a_i\), \(b=\oplus_{i<n^2}b_i\)
such that for all $i<n^2$, \(a_i\sim_2b_i\), for all
$a$, $b$ in $L$.\qed
\end{corollary}

\section{A non-normal\index{normal!non- --- lattice}
modular ortholattice\index{ortholattice}}
\label{S:NonNorm}

In this section, we shall present an example
of a non-normal%
\index{normal!non- --- lattice} modular
ortholattice\index{ortholattice} of closed subspaces of the
Hilbert space of infinite, countable dimension. The basic
underlying idea of this example is contained in the papers
\cite{DaHW72,Herr81}; the example itself is the lattice studied in
\cite{BrRo92}.

First some notation and terminology. We say that an
\emph{$M_3$-diamond}\index{diamond} of a lattice $L$ is a
quintuple
\(\vv<u,v,a,b,c>\) of elements of $L$ satisfying the two
conditions
\[
a\wedge b=a\wedge c=b\wedge c=u,\qquad
a\vee b=a\vee c=b\vee c=v
\]
(this notation for diamonds\index{diamond} is slighly
different from the one used in
Section~\ref{S:LatMonInd}---this is to make it look more
symmetric in $a$, $b$, and $c$).

We start with a classical (and straightforward) lemma.

\begin{lemma}\label{L:Subdiam} Let \(\vv<u,v,a,b,c>\) be an
$M_3$-diamond\index{diamond} of a modular lattice $L$. Let
$\alpha$ (resp., $\beta$) be the
perspective\index{perspective!isomorphism} isomorphism from
\([u,\,c]\) onto \([u,\,a]\) (resp., \([u,\,b]\)) with axis $b$
(resp., $a$). Let
\(c'\in[u,\,c]\) and put \(a'=\alpha(c')\), \(b'=\beta(c')\),
\(v'=a'\vee b'\). Then the following properties hold:
\begin{itemize}
\item[\rm (i)] \(\vv<u,v',a',b',c'>\) is an
$M_3$-diamond\index{diamond}.

\item[\rm (ii)] \([a',\,a]=\alpha[c',\,c]\) and
\([b',\,b]=\beta[c',\,c]\); in particular, the three intervals
\([a',\,a]\), \([b',\,b]\), and \([c',\,c]\) are pairwise
projective\index{projectivity!of intervals} (thus
isomorphic).\qed
\end{itemize}
\end{lemma}

We next outline G. Bruns and M. Roddy's construction as it is
presented in \cite{BrRo92}. Let $H$ be an
infinite-dimensional, separable real Hilbert space, with an
orthonormal basis consisting of vectors $e_i$ and $f_i$, for
$i\in\omega$.

Denote by \(\mathcal{L}_{\mathrm{cl}}(H)\) the (non
modular) ortholattice\index{ortholattice} of all closed
subspaces of $H$. For every subset $X$ of $H$, denote by
\([X]\) the closed subspace of $H$ generated by $X$. Then let
$a$, $b$, and $c$ (denoted in \cite{BrRo92} by
$A$, $C^\bot$, and $D$, respectively) be the following elements of
\(\mathcal{L}_{\mathrm{cl}}(H)\):
\begin{itemize}
\item[{}] \(a=[\{f_i\mid i\in\omega\}]\).

\item[{}] \(b=[\{f_0\}\cup\{e_i-2f_{i+1}\mid i\in\omega\}]\).

\item[{}] \(c=[\{e_i+f_i\mid i\in\omega\}]\).
\end{itemize}

As in \cite{BrRo92}, denote by \(x+y\) (resp., \(x\vee y\)) the
sum, as vector subspaces, of any vector subspaces $x$ and
$y$ of $H$ (resp., the closure of \(x+y\)). Next,
put \(u=[f_0]\), \(v=u^\bot\), \(\zero=\{0\}\) and \(\one=H\).
Then the relevant information about $a$, $b$, and $c$ is
contained in the following two results, see
\cite[Lemma 1]{BrRo92}:

\begin{lemma}
\begin{multline*}
a+b=a+b^\bot=a+c=a+c^\bot=a^\bot+b=a^\bot+c\\
=a^\bot+c^\bot=b+c=b+c^\bot=b^\bot+c=b^\bot+c^\bot=\one,
\end{multline*}
while
\begin{equation*}
a^\bot+b^\bot=v\prec\one.\tag*{\qed}
\end{equation*}
\end{lemma}

In particular, we see that \(a\wedge b(=a\cap b)=u\). Next, put
\[
M=\{\zero,a,a^\bot,b,b^\bot,c,c^\bot,\one\}
\]
and let $L$ be the sublattice of
\(\mathcal{L}_{\mathrm{cl}}(H)\) of those closed subspaces of
$H$ that are congruent to some element of $M$ modulo the ideal
of finite-dimensional subspaces of $H$. Then one proves the
following result, see \cite[Lemma 2]{BrRo92}:

\begin{lemma}\label{L:SumisJoin}
The lattice $L$ is a
sub-ortholattice\index{ortholattice} of
\(\mathcal{L}_{\mathrm{cl}}(H)\). Furthermore, \(x+y\) is a closed subspace
of $H$ (thus
\(x\vee y=x+y\)), for all $x$, $y\in L$.\qed
\end{lemma}

It follows, in particular, that $L$ is a complemented modular
(and even \emph{Arguesian}%
\index{lattice!Arguesian ---}) lattice.\smallskip

Now, we define inductively elements $a_n$, $a'_n$, $b_n$,
$b'_n$ (for \(n\in\omega\)) of $L$ as follows:
\begin{itemize}
\item[{}] \(a_0=a\) and \(b_0=b\).

\item[{}] \(a'_n=(a_n\vee c)\wedge a^\bot\)
and \(b'_n=(b_n\vee c)\wedge b^\bot\), for all \(n\in\omega\).

\item[{}]
\(a_{n+1}=(a'_n\vee b^\bot)\wedge a\) and
\(b_{n+1}=(b'_n\vee a^\bot)\wedge b\), for all \(n\in\omega\).
\end{itemize}

\begin{lemma}\label{L:Relataa'bb'}
The following properties hold, for all \(n\in\omega\):
\begin{itemize}
\item[\rm (a)] \(a_{n+1}\prec a_n\) and \(b_{n+1}\prec b_n\).

\item[\rm (b)]
\(a'_{n+1}\prec a'_n\) and \(b'_{n+1}\prec b'_n\).
\end{itemize}

Furthermore, all elements $a_m$, $a'_m$, and $b_n$, $b'_n$
(for $m$, $n\in\omega$) are mutually projective%
\index{projectivity!by decomposition}.
\end{lemma}

\begin{proof}
For \(n=0\), we have \(a_0=a\), \(a'_0=a^\bot\), and
\(a_1=a\wedge v\). If \(a_0=a_1\), then \(a\leq v\), thus,
since \(a^\bot\leq v\), we obtain that \(v=\one\), a contradiction;
hence, since \(v\prec\one\) and $L$ is modular, we obtain that
\(a_1\prec a_0\). Furthermore, one deduces from modularity of
$L$ that we have
\[
a_1\oplus b^\bot=(a\wedge v)\oplus b^\bot =(a\oplus
b^\bot)\wedge v=v=a^\bot\oplus b^\bot,
\]
thus \(a_1\sim a'_0\). One can prove
similarly that \(b_1\prec b_0\) and \(b_1\sim b'_0\).
Furthermore,
\(a'_0\oplus c=b'_0\oplus c=\one\), thus \(a'_0\sim b'_0\).
\smallskip

Now suppose that \(a_{n+1}\prec a_n\leq a\). Since
\(a\oplus c=a^\bot\oplus c=\one\), \(a'_n\) (resp.,
\(a'_{n+1}\)) is the image of \(a_n\) (resp., \(a_{n+1}\))
under the perspective\index{perspective!isomorphism}
isomorphism from
\([0,\,a]\) onto \([0,\,a^\bot]\) with axis $c$; it follows
that
\(a'_{n+1}\prec a'_n\). Similarly, \(a_1\oplus
b^\bot=a^\bot\oplus b^\bot=v\), thus
\(a_{n+1}\) (resp., \(a_{n+2}\)) is the image of \(a'_n\)
(resp., \(a'_{n+1}\)) under the
perspective\index{perspective!isomorphism} isomorphism from
\([0,\,a^\bot]\) onto \([0,\,a_1]\) with axis \(b^\bot\); thus
\(a_{n+2}\prec a_{n+1}\). Furthermore, the proof above shows
that \(a_n\sim a'_n\) and
\(a'_n\sim a_{n+1}\)
(with respective axes $c$ and $b^\bot$). The proof for
$b_n$ and $b'_n$ is similar.
\end{proof}

Now we can state the concluding theorem of this section:

\begin{theorem}\label{T:NonNorm}
The lattice $L$ is not normal%
\index{normal!non- --- lattice}.
\end{theorem}

\begin{proof}
It follows from Lemma~\ref{L:Relataa'bb'} that
all $a'_m$ and
$b'_n$ are mutually projective%
\index{projectivity!of elements}. Since \(a'_m\perp b'_n\), if
suffices to prove that \(a'_m\sim b'_n\) implies that \(m=n\)
(so that, for example, \(a'_0\approx b'_1\) and
\(a'_0\perp b'_1\) but \(a'_0\not\sim b'_1\)).
\smallskip

Denote by $\equiv$ the lattice congruence of $L$ generated by
the ideal of finite dimensional subspaces, that is,
\begin{multline*}
x\equiv y\text{ if and only if there exists a}\\
\text{finite-dimensional subspace }z\text{ such that }x+z=y+z.
\end{multline*}

Let $m$, $n\in\omega$ and suppose that \(a'_m\sim b'_n\). Put
\(u=a'_m\) and \(v=b'_n\). There exists \(w\in L\) such that
\(u\oplus w=v\oplus v=u\oplus v\). Let $\bar w$ be the unique
element of $M$ which is congruent to $w$ modulo $\equiv$ and
put \(w'=w\wedge\bar w\). Since
\(u\oplus w=v\oplus w=u\oplus v\) and \(u\equiv a^\bot\),
\(v\equiv b^\bot\), we cannot have \(w\equiv a^\bot\), or
\(w\equiv b^\bot\), or \(w\equiv\zero\), or \(w\equiv\one\).
Since \(\bar w\in M\), it follows that
\(\bar w\in\{a,b,c,c^\bot\}\). In particular, it follows that
\begin{equation}\label{Eq:whereisw}
w\text{ is a complement of both }a^\bot\text{ and }b^\bot.
\end{equation}

Furthermore, let $k$ (resp., $l$) denote the (finite) height of
the interval \([w',\,w]\) (resp.,
\([w',\,\bar w]\)). By Lemma~\ref{L:Subdiam}, there are
elements \(u'\leq u\) and \(v'\leq v\) such that
\(\vv<0,u'\vee v',u',v',w'>\) is an
$M_3$-diamond\index{diamond} and all three intervals
\([u',\,u]\), \([v',\,v]\), and \([w',\,w]\) are mutually
projective\index{projectivity!of intervals}; in particular,
they all have the same (finite) height~$k$. Then it follows
from (\ref{Eq:whereisw}) that we have
\begin{align*}
\mathsf{height}[u'\oplus w',\,\one]&=
\mathsf{height}[u'\oplus w',\,u'\oplus\bar w]+
\mathsf{height}[u'\oplus\bar w,\,u\oplus\bar w]+\\
&\phantom{=}+
\mathsf{height}[u\oplus\bar w,\,a^\bot\oplus\bar w]\\
&=l+k+m,
\end{align*} and, similarly,
\(\mathsf{height}[v'\oplus w',\,\one]=l+k+n\). However,
\(u'\oplus w'=v'\oplus w'\), thus \(m=n\) and we are done.
\end{proof}

The dimension monoid of $L$ is easy to compute: it is
isomorphic to \(\ZZ^+\cup\{nc\mid n\in\NN\}\), where $c$ is an
element satisfying \(1+c=c<2c<3c<\cdots\). \emph{Via} this
isomorphism, for all \(x\in L\), \(\Dim(x)=\dim(x)\), if $x$ is
finite dimensional, while \(\Dim(x)=c\), otherwise.

\section{Dimension monoids of regular%
\index{ring!(von Neumann) regular ---} rings}\label{S:RegRing}

It is not always the case that for a regular ring $R$, the
dimension monoid $V(R)$ of the ring $R$ is isomorphic to the
dimension monoid $\DD\mathcal{L}(R)$ of the lattice
$\mathcal{L}(R)$, see Corollary~\ref{C:RegRingCounterex}. However,
we shall see in this section that there is always a canonical
homomorphism from the latter monoid onto the former monoid.

We start with the observation that if $R$ is a regular%
\index{ring!(von Neumann) regular ---}
ring, then, with the notations of
Definition~\ref{D:dimsgrp}, we have
\(V(R)=\DD(\mathcal{L}(R_R),\oplus,\cong)\). Any element $M$
of $\FP(R)$ determines a natural monoid
homomorphism \(\pi_M\colon\DD\mathcal{L}(M)\to V(R)\)%
\index{pzziM@$\pi_M$, $\pi_R$, $\pi_{nR}$|ii}
sending every $\Dim(X)$ (for $X\in\mathcal{L}(M)$) to $[X]$.
To ease the notation, we shall write $\pi_R$
\index{pzziM@$\pi_M$, $\pi_R$, $\pi_{nR}$|ii} instead of
$\pi_{R_R}$.

\begin{proposition}\label{P:piMVhom}
The map $\pi_M$ is a \Vhom\index{Vhom@\Vhom}.
\end{proposition}

\begin{proof}
It is immediate that $\pi_R$ is a \Vhom\index{Vhom@\Vhom} at
every element of
$\DD\mathcal{L}(M)$ of the form
$\Dim(X)$, $X\in\mathcal{L}(M)$. Since
$V(R)$ is a
refinement\index{monoid!refinement ---} monoid, we obtain the
conclusion by Lemma~\ref{L:V-homat}.
\end{proof}

In particular, if $[M]$ is an order-unit of
$V(R)$, then $\pi_M$ is surjective.
For example, this is the case when
$M=nR_R$, for some $n\in\NN$. On the other hand, even
$\pi_R$ may not
necessarily be one-to-one (see
Corollary~\ref{C:RegRingCounterex}).

\begin{proposition}\label{P:pinRiso} Let $R$ be a regular%
\index{ring!(von Neumann) regular ---} ring. Then for every
integer \(n\geq 2\),
\(\pi_{nR}\) is an
isomorphism from \(\DD\mathcal{L}(nR_R)\) onto \(V(R)\).
\end{proposition}

\begin{proof} First, $nR$ endowed with the relation of
isomorphism $\cong$ satisfies the assumption of
Proposition~\ref{P:uniqnor} (the $u_i$'s are the elements of
the canonical homogeneous\index{homogeneous!basis} basis of
$nR$ and $v_i=\oplus_{j\ne i}u_j$), thus $\cong$ and
$\approxeq$ coincide, so that, by
Corollary~\ref{C:eqdim}, $[X]=[Y]$
if and only if $\Dim(X)=\Dim(Y)$, for all $X$ and $Y$ in
$\mathcal{L}(nR_R)$. This proves that the
restriction of $\pi_{nR}$ to the dimension
range of $\mathcal{L}(nR_R)$ is
one-to-one. We obtain the conclusion by Lemma~\ref{L:onetoone}.
\end{proof}

\begin{remark}
The congruence\index{congruence!semilattice}
semilattice analogue of this result (that is, one takes
the maximal semilattice quotient of each side) holds even for
$n=1$, see \cite[Corollary 2.4]{WehrC}.
\end{remark}

We recall the following definition from \cite{TuWe}:

\begin{definition}\label{D:CongSplit}
A lattice $L$ has \emph{commuting congruences}%
\index{lattice!with commuting congruences|ii}, if for any
congruences $\alpha$ and $\beta$ of $L$,
$\alpha\circ\beta=\beta\circ\alpha$.
\end{definition}

Note that the class of lattices with commuting congruences is
fairly large:

\begin{proposition}\label{P:SuffCS}
The following statements hold:
\begin{itemize}

\item[\rm (a)] Every lattice which is either relatively
complemented or sectionally complemented has commuting congruences.

\item[\rm (b)] Every atomistic lattice has commuting congruences.

\item[\rm (c)] The class of lattices with commuting congruences is
closed under direct limits.\qed

\end{itemize}
\end{proposition}

The origin of the following result can be traced back to
\cite{PTWe}, where it was proved for the so-called
\emph{congruence-splitting}%
\index{lattice!congruence splitting ---}
lattices. It was later strengthened
to lattices with commuting congruences, see \cite{TuWe}.

\begin{corollary}\label{C:NoReprLatt}
Let $\mathcal{V}$ be
any non-distributive lattice variety. For any set $X$, let
\(F_{\mathcal{V}}(X)\)%
\index{FzzreeLattX@$F_{\mathcal{V}}(X)$|ii} be the free
$\{0,1\}$-lattice in $\mathcal{V}$ on $X$. If
\(|X|\geq\aleph_2\), then there exists no lattice $L$
with commuting congruences such that
\(\DD L\cong\DD(F_{\mathcal{V}}(X))\).
\end{corollary}

\begin{proof}
If \(\DD L\cong\DD(F_{\mathcal{V}}(X))\), then,
by Corollary~\ref{C:congquotV}, we also have
\(\ccon L\cong\ccon(F_{\mathcal{V}}(X))\), a contradiction by
the results of \cite{TuWe}.
\end{proof}

\begin{corollary}\label{C:NoReprRing}
Let $\mathcal{V}$ and $X$ satisfy the hypothesis of
Corollary~\ref{C:NoReprLatt}. Then there exists no regular%
\index{ring!(von Neumann) regular ---} ring $R$ such that
\(V(R)\cong\DD(F_{\mathcal{V}}(X))\).\qed
\end{corollary}

\begin{remark}
These last two results yields many new
refinement\index{monoid!refinement ---} monoids (of size at
least $\aleph_2$) with order-unit which cannot be represented
as \(V(R)\) for $R$ regular (the
first examples had been found in \cite{WehrB}). These
refinement\index{monoid!refinement ---} monoids, the
\(F_{\mathcal{V}}(X)\), are not always cancellative (for
example, if $\mathcal{V}$ is the variety of all modular
\index{lattice!modular (not necessarily complemented) ---}
lattices). In fact, we do not know whether
\(F_{\mathcal{V}}(X)\) is always a
refinement\index{monoid!refinement ---} monoid: although it is
certainly the case if $\mathcal{V}$ is \emph{modular}%
\index{lattice!modular (not necessarily complemented) ---}, it
is also the case if $\mathcal{V}$ is \emph{generated by a
finite lattice}, because then, by well-known standard results
of universal algebra,
\(F_{\mathcal{V}}(X)\) is finite whenever $X$ is finite, so
that in that case, by Corollary~\ref{C:dimprim},
\(\DD(F_{\mathcal{V}}(X))\) is a
primitive\index{monoid!primitive ---}
refinement\index{monoid!refinement ---} monoid.
\end{remark}

The argument of the proof of Proposition~\ref{P:pinRiso} shows
that for every \(M\in\FP(R)\) such that
$\mathcal{L}(M)$ satisfies the
hypothesis of Proposition~\ref{P:uniqnor} with respect to the
normal\index{normal!equivalence} equivalence relation $\cong$,
$\pi_M$ is a
\Vemb\index{Vemb@\Vemb}, so that if $M=nR_R$ and $n\geq1$, then
$\pi_M$ is in that case
an isomorphism. Our next result will make it possible to find a
large class of regular\index{ring!(von Neumann) regular ---} rings
such that $M=R_R$ satisfies this condition.

\begin{proposition}\label{P:piRiso} Let $R$ be a regular%
\index{ring!(von Neumann) regular ---} ring. Suppose that
there are principal right ideals\index{ideal!of a ring}
$I$ and $J$ of $R$ such that
\(I\cap J=\{0\}\) and both $[I]$ and $[J]$ are order-units of
$V(R)$. Then $\mathcal{L}(R_R)$,
endowed with the relation of isomorphy $\cong$, satisfies the
hypotheses of Proposition~\ref{P:uniqnor}. Therefore,
$\pi_R$ is an
isomorphism from $\DD\mathcal{L}(R_R)$ onto $V(R)$.
\end{proposition}

\begin{proof} Without loss of generality,
$I\oplus J=R_R$. Since $V(R)$ is a
refinement\index{monoid!refinement ---} monoid, there are, by
Corollary~\ref{C:charasymp}, decompositions in
$\mathcal{L}(R_R)$ of the form
\[
I=\bigoplus_{i\in X}I_i\ \text{ and }\
J=\bigoplus_{j\in Y}J_j,
\]
with finite sets $X$ and $Y$ (\emph{which we may suppose
disjoint}), along with a subset $\Gamma\subseteq X\times Y$ of
domain $X$ and image $Y$ such that
$I_i\cong J_j$, for all
$\vv<i,j>\in\Gamma$. This yields a decomposition
\[
R_R=\bigoplus_{k\in X\cup Y}U_k,
\]
where we put $U_k=I_k$, if $k\in X$ and $U_k=J_k$, if
$k\in Y$. Then put $V_k=\bigoplus_{l\ne k}U_l$. Since
$\Gamma$ has domain $X$, for all $k\in X$, there exists
$l\in Y$ such that $U_k\cong U_l$, whence $[U_k]\leq[V_k]$.
Symmetrically, the same conclusion holds for all $k\in Y$. By
Proposition~\ref{P:uniqnor}, the conclusion follows.
\end{proof}

\begin{corollary}\label{C:caseRsimple} Let $R$ be a simple%
\index{ring!simple ---} regular\index{ring!(von Neumann) regular ---} ring.
Then $\pi_R$ is an isomorphism.
\end{corollary}

\begin{proof} If there are nonzero elements $\alpha$ and
$\beta$ of $V(R)$ such that
$[R]=\alpha+\beta$, then the conclusion follows immediately
from Proposition~\ref{P:piRiso}. Otherwise, for every nonzero
principal right ideal\index{ideal!of a ring} $I$ of $R$,
$I\cong R_R$; if $I\ne R$, then there exists a nonzero
principal right ideal\index{ideal!of a ring} $J$ such that
$I\oplus J=R$, so that
$[I]+[J]=[R]$ with both $[I]$ and $[J]$ nonzero, a
contradiction: thus $I=R$. Hence, $R$ is a division ring and
the conclusion is trivial.\end{proof}

Here is another class of rings for which
$\pi_R$ is an isomorphism:

\begin{proposition}\label{P:piiso(unreg)} Let $R$ be a
unit-regular\index{ring!unit-regular ---} ring. Then
$\pi_R$ is an isomorphism.
\end{proposition}

\begin{proof}
By \cite[Proposition 6.8]{Good91}, two principal
ideals\index{ideal!of a ring} of $R$ are isomorphic if and
only if they are perspective\index{perspectivity}. The
conclusion follows by Lemma~\ref{L:onetoone} and
Proposition~\ref{P:piMVhom}.
\end{proof}

\chapter[Locally finitely distributive lattices]%
{Locally finitely distributive%
\index{distributive!locally finitely ---} relatively
complemented modular lat\-tices}
\label{FinDistr}

The basic idea of the proof of Theorem~\ref{T:MeetContNor} is the
following. For any sectionally complemented modular lattice $L$, we
shall find a large \emph{normal} ideal $I$ of $L$, the so-called
\emph{normal kernel} of $L$ (see Chapter~\ref{NormKer}), and we
shall prove that the quotient lattice $L/I$ is ``small'', in
particular, normal. The relevant concept of smallness is here the
concept of
\emph{$n$-distributivity}, see A.~P. Huhn \cite{Huhn72}. Thus,
this section will be devoted to the basic dimension
theory of $n$-distributive relatively complemented modular
lattices.

\section{Lattice index and monoid index}\label{S:LatMonInd}

The main result of this section is that for a given complemented
modular lattice $L$, the least $n$ such that $L$ is $n$-distributive
is equal to the index of $\Dim(1_L)$ in the monoid $\DD L$, see
Proposition~\ref{P:twoindex}. This will lead to a
monoid-theoretical treatment of certain questions related to
$n$-distributivity.

We shall first recall Huhn's definition of $n$-distributivity.

\begin{definition}\label{D:n-Distr} Let $n$ be a positive
integer. A modular%
\index{lattice!modular (not necessarily complemented) ---}
lattice $L$ is \emph{$n$-distributive}%
\index{distributive!$n$- ---}, if it satisfies the following
identity:
\[
x\vee\bigwedge_{i\leq n}y_i=
\bigwedge_{i\leq n}\left(
x\vee\bigwedge_{j\leq n;\ j\ne i}y_j
\right).
\]
We say that $L$ is \emph{finitely distributive}%
\index{distributive!finitely ---}, if it is
$n$-distributive for some $n\in\NN$. We say that $L$ is
\emph{locally finitely distributive}%
\index{distributive!locally finitely ---}, if every closed
interval of $L$ is finitely distributive.
\end{definition}

In particular, a lattice is distributive if and only if it is
$1$-distributive.

The concept of $n$-distributive lattice%
\index{distributive!$n$- ---} is very closely related to the
concept of an \emph{diamond}. Since the definition
of a diamond is not completely uniform in the
literature, we shall state our own definition here:

\begin{definition}\label{D:diam}
Let $n$ be a positive integer
and let $L$ be a modular%
\index{lattice!modular (not necessarily complemented) ---}
lattice. An \emph{$n$-diamond}\index{diamond|ii} of $L$ is a
finite sequence
\[
\delta=\vv<a_0,a_1,\ldots,a_{n-1},e>
\]
of $n+1$ elements of $L$ such that, if we put
$0_\delta=\bigwedge_{i<n}a_i$ and
$1_\delta=\bigvee_{i<n}a_i$, then
$a_i\wedge e=0_\delta$ and $a_i\vee e=1_\delta$ for all $i<n$,
and $\vv<a_0,\ldots,a_{n-1}>$ is
independent above
$0_\delta$. We say that the diamond $\delta$ is
\emph{trivial}, if $0_\delta=1_\delta$.
\end{definition}

The following theorem is proved in A.~P. Huhn
\cite[Satz 2.1]{Huhn72}.

\begin{theorem}\label{T:nDistrdiam}
Let $n$ be a positive integer.
Then a modular%
\index{lattice!modular (not necessarily complemented) ---}
lattice is $n$-distributive%
\index{distributive!$n$- ---} if and only if it has no non
trivial $n+1$-diamond\index{diamond}.\qed
\end{theorem}

As a corollary, Huhn observes in particular that the
concept of $n$-distributivity is \emph{self-dual}, see
\cite[Satz 3.1]{Huhn72}. Note that the concept of
$n$-distributivity can be extended to
\emph{non-modular} lattices, and this concept is not self-dual for
$n\geq2$.

For our purposes here, we shall focus
on sectionally complemented lattices. Furthermore, we shall
not really need the last entry of a diamond\index{diamond};
thus, we formulate the following definition:

\begin{definition}\label{D:HomSys} Let $n$ be a non-negative
integer and let $L$ be a modular%
\index{lattice!modular (not necessarily complemented) ---}
lattice with zero. Then a
\emph{homogeneous\index{homogeneous!sequence|ii} sequence of
length $n$} of
$L$ is a finite \emph{independent} sequence
\(\sigma=\vv<a_0,a_1,\ldots,a_{n-1}>\) such that
\((\forall i<n)(a_0\sim a_i)\); put then
\(1_\sigma=\oplus_{i<n}a_i\). We say that the
homogeneous sequence
$\sigma$ is \emph{trivial}, if $1_\sigma=0$.
\end{definition}

\begin{lemma}\label{L:HomDiam} Let $n$ be a positive integer,
let $L$ be a sectionally complemented modular lattice. Then
$L$ is $n$-distributive%
\index{distributive!$n$- ---} if and only if it has no non
trivial homogeneous sequence of
length $n+1$.
\end{lemma}

We refer, for example, to \cite[3.3]{JiRo92} for more details
about this.

\begin{proof}
If there is no nontrivial
$n+1$-diamond\index{diamond}, then there is no nontrivial
homogeneous sequence of length
$n+1$: indeed, if
\(\sigma=\vv<a_i\mid i\leq n>\) is a
homogeneous sequence, then let $c_i$ be such that
\(a_0\oplus c_i=a_i\oplus c_i=a_0\oplus a_i\), for all
$i\in\{1,\ldots,n\}$, and put
\(e=\oplus_{i=1}^nc_i\). Then it is easy to verify that
\(\delta=\vv<a_0,\ldots,a_n,e>\) is an
$n+1$-diamond\index{diamond}, which is by assumption trivial,
thus $1_\sigma=1_\delta=0_\delta=0$. Conversely, suppose that
there is no nontrivial
homogeneous sequence of length
$n+1$ and let
\(\delta=\vv<a_0,\ldots,a_n,e>\) be an
$n+1$-diamond\index{diamond}. Let $x_i\in L$ be
such that
\(0_\delta\oplus x_i=a_i\), for all $i\leq n$. Then \(\vv<x_i\mid
i\leq n>\) is a homogeneous sequence of
length
$n+1$, thus $x_i=0$, for all $i$, thus $\delta$ is trivial. We
obtain the conclusion by Theorem~\ref{T:nDistrdiam}.
\end{proof}

\begin{definition}\label{D:LattInd}
Let $L$ be a sectionally
complemented modular lattice. Then we define
the \emph{index} of $x$ in $L$ for all $x\in L$ and we denote it by
$\ind_L(x)$%
\index{index!of an element in a lattice|ii}%
\index{izzndl@$\ind_L(x)$|ii}, as the largest
$n\in\ZZ^+$ such that there exists a nontrivial
homogeneous sequence of length $n$ below
$x$, if it exists, and $\infty$, otherwise.
\end{definition}

\begin{note}
In particular, we allow
homogeneous sequences of length
$0$, in order to define $\ind_L(0)=0$.
\end{note}

\begin{remark}
At this point, our index terminology is
consistent with the one used for regular%
\index{ring!(von Neumann) regular ---} rings: in particular,
one sees, by \cite[Theorem 7.2]{Good91}, that if $R$ is a
regular\index{ring!(von Neumann) regular ---} ring, then the
index of $R$ in $\mathcal{L}(R_R)$ as defined
above is equal to the index of nilpotence%
\index{index!of nilpotence|ii} of $R$,
that is, the supremum of the indexes of nilpotence
of the elements of $R$ (the \emph{index of nilpotence} of an element
$x$ of $R$ is the least positive integer $n$, if it exists, such
that $x^n=0$).
\end{remark}

In order to relate the lattice-the\-o\-ret\-i\-cal concept of
index and the monoid-the\-o\-ret\-i\-cal one, we shall need the
following lemma, which is a detailed version of a
result of J\'onsson, see B. J\'onsson \cite[Lemma 1.4]{Jons60}:

\begin{lemma}\label{L:Jonsson}
Let $a$ and $b$ be elements of a sectionally complemented
modular lattice such that $a\lesssim_2b$. Then there are
decompositions
\begin{equation}\label{Eq:asim2b}
\begin{split} a&=u_0\oplus u_1\oplus(\oplus_{i<4}a_i)\\
b&=u\oplus(\oplus_{i<4}b_i)\oplus h
\end{split}
\end{equation}
such that $u_0\sim u$, $u_1\sim u$, and $a_i\sim b_i$,
for all $i<4$.
\end{lemma}

\begin{proof}
It suffices to prove the result for
\(a\sim_2b\), that is, there exists $c$ such that \(a\sim c\)
and \(c\sim b\). We follow the original
proof of J\'onsson's Lemma \cite[Lemma 1.4]{Jons60}. So let $S$
(resp., $T$) be a perspective\index{perspective!isomorphism}
isomorphism from
\([0,\,a]\) onto \([0,\,c]\) (resp., \([0,\,c]\) onto
\([0,\,b]\)). As in the original proof of J\'onsson's Lemma,
consider a decomposition \(a=(a\wedge c)\oplus a'\) and observe
that \(S(a\wedge c)=a\wedge c\), so that one can reduce the problem
to \(a\wedge c=0\) provided that we can find \emph{three} summands
\(a_i\) (and \(b_i\)) instead of four in (\ref{Eq:asim2b}).
Similarly, one can reduce the problem to
\(a\wedge c=b\wedge c=0\) provided that we make it with \emph{two}
summands \(a_i\) (and \(b_i\)) instead of four in
(\ref{Eq:asim2b}). Next, put
\(a'=a\wedge(b\oplus c)\) and \(b'=b\wedge(a\oplus c)\), so
that $a'$ and $b'$ are perspective\index{perspectivity}
with axis $c$, and let
$a''$ such that
\(a=a'\oplus a''\). Put \(b''=TS(a'')\); since
\(\vv<a'',b,c>\) is independent, so is
\(\vv<a'',S(a''),b''>\), thus, by
Lemma~\ref{L:TransOrthPersp}, \(a''\sim b''\). Put
\(u=b'\wedge b''\) and let $h$ such that
\(b=(b'\vee b'')\oplus h\). Let \(b_0\) and \(b_1\) such that
\[
u\oplus b_0=b'\ \text{ and }\ u\oplus b_1=b''.
\]
Since \(a'\sim b'\) and \(a''\sim b''\), there are
decompositions
\(a'=u_0\oplus a_0\) and \(a''=u_1\oplus a_1\) such that
\(u_i\sim u\) (for all \(i<2\)) and
\(a_i\sim b_i\) (for all \(i<2\)).  Thus the
decompositions
\begin{align*} a&=u_0\oplus u_1\oplus a_0\oplus a_1\\
b&=u\oplus b_0\oplus b_1\oplus h
\end{align*} are as required.\end{proof}

\begin{remark} In fact, in the proof of Lemma~\ref{L:Jonsson},
it is also possible to prove that \(u_0\sim u_1\) (by using
Lemma~\ref{L:TransOrthPersp}), but we shall not need this.
\end{remark}

The following corollary is proved in \cite[Lemma 1.5]{Jons60}:

\begin{corollary}\label{C:Jonsson}
Let $a$ and $x$ be elements of a
sectionally complemented modular lattice $L$. Then
\(\Theta(x)\subseteq\Theta(a)\) if and only if there exists a
decomposition \(x=\oplus_{i<n}x_i\) such that
\(x_i\lesssim a\) for all \(i<n\).\qed
\end{corollary}

\begin{corollary}\label{C:IncrCong}
Let $L$ be a sectionally
complemented modular lattice and let
$n$ be a positive integer. Let \(\vv<a_i\mid i\leq n>\)
be a finite sequence of elements of $L$ such that
\(\Theta(a_0)\subseteq\Theta(a_i)\) for all
\(i\leq n\). Then there
are a positive integer $N$ and decompositions
\begin{equation}\label{Eq:a0,ai}
\begin{split}
a_0&=\bigvee_{j<N}x_{0j},\\
a_i&\geq\bigvee_{j<N}x_{ij}\qquad(1\leq i\leq n)
\end{split}
\end{equation} such that
$x_{ij}\sim x_{i+1,j}$ for all $i<n$ and all $j<N$.
\end{corollary}

\begin{proof}
By induction on $n$. It is trivial for $n=0$. Let
$\vv<a_i\mid i\leq n+1>$ be a finite sequence of elements
of $L$ such that
\(\Theta(a_0)\subseteq\Theta(a_i)\) for all $i\leq n+1$.
By the induction hypothesis,
there are a positive integer $N$ and decompositions as in
(\ref{Eq:a0,ai}) (in the range \(0\leq i\leq n\)) such
that $x_{ij}\sim x_{i+1,j}$ for all $i<n$ and all $j<N$. Now,
\(\Theta(x_{nj})=\Theta(x_{0j})\subseteq\Theta(a_0)
\subseteq\Theta(a_{n+1})\) for
all $j<N$, thus, by
Corollary~\ref{C:Jonsson}, there are a positive integer $m$,
decompositions \(x_{nj}=\bigvee_{k<m}x_{njk}\), and elements
$x_{n+1,j,k}\leq a_{n+1}$ such that \(x_{njk}\sim x_{n+1,j,k}\)
for all $j<N$ and all $k<m$. Now, by using
Lemma~\ref{L:PerspMap}, it is easy to find by downward
induction on $i\leq n$ decompositions
\(x_{ij}=\bigvee_{k<m}x_{ijk}\) such that
$x_{ijk}\sim x_{i+1,j,k}$ for all $i<n$, all $j<N$, and all $k<m$.
Thus the decompositions
\(a_i\geq\bigvee_{j<N;\ k<m}x_{ijk}\) (for all
$i\leq n+1$, with equality, if $i=0$) are as required.
\end{proof}

\begin{proposition}\label{P:twoindex} Let $L$ be a sectionally
complemented modular lattice. Then
\(\ind_L(x)=\ind_{\DD L}(\Dim(x))\)%
\index{index!of an element in a monoid}%
\index{index!of an element in a lattice} for all $x\in L$.
\end{proposition}

\begin{proof}
First the easy direction. If $n\in\NN$ be such that
$n\leq\ind_L(x)$%
\index{index!of an element in a lattice}, then there exists,
by definition, a nontrivial
homogeneous\index{homogeneous!sequence} sequence
$\vv<y_i\mid i<n>$ below
$x$; thus \(n\cdot\Dim(y_0)=\Dim(\oplus_{i<n}y_i)\leq\Dim(x)\)
with $\Dim(y_0)>0$, so that
\(n\leq\ind_{\DD L}(\Dim(x))\)%
\index{index!of an element in a monoid}. Conversely, let
$n\in\NN$ be such that
\(n\leq\ind_{\DD L}(\Dim(x))\)%
\index{index!of an element in a monoid}. By
Corollary~\ref{C:dimVmeas}, there exists a nonzero $y\in L$
such that
\(n\cdot\Dim(y)\leq\Dim(x)\). Again by
Corollary~\ref{C:dimVmeas}, there exists an
independent finite sequence
$\vv<y_i\mid i<n>$ such that $y_i\leq x$
and $\Dim(y_i)=\Dim(y)$, for all $i<n$. It follows that
\emph{a fortiori},
\(\Theta(y_i)=\Theta(y)\), thus, by Corollary~\ref{C:IncrCong},
there are $N\in\NN$ and decompositions
\begin{align*}
y_0&=\bigvee_{j<N}y_{0j},\\
y_i&\geq\bigvee_{j<N}y_{ij}\qquad(1\leq i<n)
\end{align*} such that
\(y_{ij}\sim y_{i+1,j}\) for all $i<n-1$ and all $j<N$.
Since $y_0\ne0$, there
exists $j<N$ such that $y_{0j}\ne0$. Then it follows from
Lemma~\ref{L:TransOrthPersp} that the finite sequence
$\vv<y_{ij}\mid i<n>$ is
homogeneous\index{homogeneous!sequence}, thus, since it lies
below $x$, $n\leq\ind_L(x)$%
\index{index!of an element in a lattice}.
\end{proof}

\begin{corollary}\label{C:nDistrCanc} Let $L$ be a sectionally
complemented modular lattice. If $L$ is locally finitely
distributive%
\index{distributive!locally finitely ---}, then every element
of \(\DD L\) has finite index
and
\(\DD L\) is the positive cone
of an Archimedean\index{Archimedean} dimension group%
\index{group!dimension ---}.
\end{corollary}

\begin{proof} By Proposition~\ref{P:twoindex}, every element
of the dimension range of $L$ has finite index. By
Lemma~\ref{L:addInd}, every element of
$\DD L$ has finite index. The conclusion follows by
Proposition~\ref{P:finIndDim}.
\end{proof}

\section{Transitivity of perspectivity\index{perspectivity}}

The main purpose of this section is to prove that in any relatively
complemented modular lattice in which every closed interval is
$n$-distributive for some $n$, the relation of perspectivity is
transitive, see Corollary~\ref{C:nDistrTrans}. This result is
technically easier to prove than the corresponding results proved
first by von~Neumann for continuous geometries and then by Halperin
for countably continuous geometries.

We start with a weak version of the desired result:

\begin{lemma}\label{L:LocDisPP} Let $L$ be a locally finitely
distributive%
\index{distributive!locally finitely ---} sectionally
complemented modular lattice. Then
$a\approxeq b$ implies that $a\simeq b$, for all $a$, $b\in L$.
\end{lemma}

\noindent(We refer the reader to Definitions \ref{D:projdec}
and \ref{D:persdec} for the definitions of $\approxeq$ and
$\simeq$).

\begin{proof}
We prove the result by induction on $\ind_L(a)$%
\index{index!of an element in a lattice}. It is trivial for
$\ind_L(a)=0$%
\index{index!of an element in a lattice} (in which case
\(a=b=0\)). So suppose that the conclusion has been proved for
$\ind_L(a)<n$, where $n$ is a positive integer; let
$a\approxeq b$ in $L$, where $\ind_L(a)=n$.
By definition, there are decompositions
\(a=\oplus_{i<m}a_i\) and \(b=\oplus_{i<m}b_i\) such that
$a_i\approx b_i$, for all $i<m$. By
Lemma~\ref{L:addInd}, the inequality
$\Dim(a_i)\leq\Dim(a)$ holds for all $i$, thus it suffices to
prove the conclusion for $a\approx b$. Since the index
function is perspectivity-invariant\index{perspectivity} (this
follows immediately from Proposition~\ref{P:twoindex}), an easy
induction proof shows that it suffices to consider the case
$a\sim_2b$. In this case, by Lemma~\ref{L:Jonsson}, there are
decompositions
\begin{align*}
a&=u_0\oplus u_1\oplus(\oplus_{i<4}a_i),\\
b&=u\oplus(\oplus_{i<4}b_i)\oplus h,
\end{align*}
such that $u_0$, $u_1$, and $u$ are mutually
perspective\index{perspectivity} and
$a_i\sim b_i$, for all $i<4$; thus, in
order to prove the conclusion, it suffices to prove that
$u_0\simeq h$. Since
$\Dim(a)=\Dim(b)$ and, by Corollary~\ref{C:nDistrCanc},
$\DD L$ is cancellative,
\(\Dim(u_0)=\Dim(h)\), which means, by Corollary~\ref{C:eqdim},
that \(u_0\approxeq h\). However, since
$u_0\oplus u_1\leq a$, the inequality $2\cdot\Dim(u_0)\leq\Dim(a)$
holds, so that it follows from Lemma~\ref{L:addInd} that
\(\ind_{\DD L}(\Dim(a))\geq 2\cdot\ind_{\DD L}(\Dim(u_0))\)%
\index{index!of an element in a monoid},
that is, \(\ind_L(a)\geq 2\cdot\ind_L(u_0)\). Since
$\ind_L(a)=n>0$, \(\ind_L(u_0)<\ind_L(a)\)%
\index{index!of an element in a lattice}. Therefore, it
follows from \(u_0\approxeq h\) and the induction hypothesis
that $u_0\simeq h$, which concludes the proof.
\end{proof}

\begin{corollary}\label{C:nDistrNorm} Every locally finitely
distributive%
\index{distributive!locally finitely ---} sectionally
complemented modular lattice is normal\index{normal!lattice}.
\end{corollary}

(We refer the reader to Definition~\ref{D:normal} for the
definition of a normal\index{normal!lattice} lattice).

\begin{proof}
By Lemma~\ref{L:LocDisPP} and Proposition~\ref{P:SuffNorm}.
\end{proof}

Both Lemma~\ref{L:LocDisPP} and Corollary~\ref{C:nDistrNorm}
will be strengthened in Corollary~\ref{C:nDistrTrans}.

\begin{proposition}\label{P:transperspchar} Let $L$ be a
relatively complemented modular lattice. Then
perspectivity\index{perspectivity} is transitive in $L$ if and
only if $L$ is normal\index{normal!lattice} and
$\DD L$ is cancellative.
\end{proposition}

We shall see that \emph{the assumption of normality is
redundant}\index{normal!lattice} in
Proposition~\ref{P:transperspchar}, but this is much harder to
prove so that we postpone it until
Theorem~\ref{T:transperspchar}.

\begin{proof}
Suppose, first, that
perspectivity\index{perspectivity} is transitive in $L$; we
prove that \(\DD L\) is
cancellative. Without loss of generality, $L$ is a bounded
lattice (thus complemented and modular). Since
perspectivity\index{perspectivity} is transitive, $L$ is
\emph{a fortiori} normal\index{normal!lattice}. Then, using the
same proof as in
\cite[IV.2.2]{FMae58} (stated there for continuous
geometries\index{continuous!geometry} but valid with merely
transitive perspectivity\index{perspectivity}; see also
\cite[Part I, Chapter VI]{Neum60}), one can prove the
following statement:

\setcounter{claim}{0}
\begin{claim}
The following assertions hold:
\begin{itemize}

\item[\rm (i)]
\(a\oplus x=b\oplus y\) and $a\sim b$ implies that $x\sim y$,
for all $a$, $b$, $x$, and $y$ in $L$.

\item[\rm (ii)]
$a\oplus x\sim b\oplus y$ and $a\sim b$ implies that
$x\sim y$, for all $a$, $b$, $x$, and $y$ in $L$.

\item[\rm (iii)] If $\vv<a_i\mid i<n>$ and
$\vv<b_i\mid i<n>$ are independent\index{independent}
finite sequences of elements of $L$ such that
$a_i\sim b_i$ (for all $i<n$), then
\(\oplus_{i<n}a_i\sim\oplus_{i<n}b_i\).
\qedc
\end{itemize}
\end{claim}

Using Corollary~\ref{C:eqdim}, it follows immediately that
$\Dim(a)=\Dim(b)$ if and only if $a\sim b$, for
all $a$, $b\in L$.
Then, using (ii) of Claim~1, one deduces immediately
that the dimension range of $L$ is cancellative (in the sense
of Lemma~\ref{L:cancref}). By Lemma~\ref{L:cancref}, it
follows that $\DD L$ is cancellative.\smallskip

Conversely, suppose that $L$ is normal\index{normal!lattice}
and that $\DD L$ is cancellative. We
first prove that $a\wedge b=u$ and
$\Dim(u,a)=\Dim(u,b)$ imply $a\sim b$,
for all $a$, $b$, $u\in L$. So put
\(L_u=\{x\in L\mid x\geq u\}\). By Corollary~\ref{C:eqdim},
there are decompositions (in $L_u$)
\(a=\oplus_{i<n}a_i\) and
\(b=\oplus_{i<n}b_i\) such that
\(a_i\approx b_i\), for all \(i<n\). But \(a_i\wedge
b_i=u\), thus, since $L_u$ is normal\index{normal!lattice},
\(a_i\sim b_i\). But since the finite sequence
\(\vv<a_0,\ldots,a_{n-1},b_0,\ldots,b_{n-1}>\) is
independent\index{independent} in $L_u$, this implies that
\(\oplus_{i<n}a_i\sim\oplus_{i<n}b_i\),
that is, \(a\sim b\), as claimed.

Now in the general case, let $a$, $b$, $c\in L$ be such that
\(a\sim b\) and \(b\sim c\). Put \(u=a\wedge b\wedge c\) and
\(v=a\vee b\vee c\), and put \(K=[u,\,v]\). Note that
\(a\sim b\) and \(b\sim c\) in $K$. It follows that
\(\Dim(u,a)=\Dim(u,b)=\Dim(u,c)\). Now let $a'$ and $c'$ in
$K$ be such that $K$ satisfies \(a=a\wedge c\oplus a'\) and
\(c=a\wedge c\oplus c'\). Then
\(\Dim(u,a)=\Dim(u,c)\) can be written as
\(\Dim(u,a\wedge c)+\Dim(u,a')=\Dim(u,a\wedge c)+\Dim(u,c')\),
thus, since \(\DD L\) is
cancellative,
\(\Dim(u,a')=\Dim(u,c')\). Since \(a'\perp c'\) in $K$, it
follows from previous paragraph that \(a'\sim c'\). But
\(\vv<a',c',a\wedge c>\) is independent\index{independent}
in $K$, thus, by Corollary~\ref{C:BasicAddPersp}, we also have
\(a'\oplus(a\wedge c)\sim c'\oplus(a\wedge c)\) in $K$,
that is, \(a\sim c\).
\end{proof}

\begin{corollary}\label{C:nDistrTrans}
In every locally finitely distributive%
\index{distributive!locally finitely ---} relatively
complemented modular lattice,
perspectivity\index{perspectivity} is transitive.
\end{corollary}

\begin{proof}
By Corollary~\ref{C:nDistrCanc},
Corollary~\ref{C:nDistrNorm}, and
Proposition~\ref{P:transperspchar}.
\end{proof}

\begin{lemma}\label{L:nDistrVemb} Let $L$ be a relatively
complemented modular lattice, let $K$ be a convex sublattice
of $L$. If $L$ is locally finitely distributive%
\index{distributive!locally finitely ---}, then the natural map
\(f\colon \DD K\to\DD L\) is a \Vemb\index{Vemb@\Vemb}.
\end{lemma}

\begin{proof} We already know that $f$ is a
\Vhom\index{Vhom@\Vhom}\  (see Proposition~\ref{P:exV-hom}), so
that it suffices to prove that $f$ is one-to-one. By
expressing a lattice as the direct union of its closed
intervals, one sees easily that it suffices to consider the
case where both $K$ and $L$ are bounded lattices (not
necessarily with the same $0$ and $1$). Thus both $K$ and $L$
are complemented modular lattices. By Lemma~\ref{L:cancref},
it suffices to prove that the restriction of $f$ to the
dimension range of $K$ is one-to-one. But this follows
immediately from Lemma~\ref{L:persp} and
Corollary~\ref{C:nDistrTrans}.
\end{proof}

We deduce the following result:

\begin{theorem}\label{T:FinIndnDistr}
Let $L$ be a modular lattice%
\index{lattice!modular (not necessarily complemented) ---}.
Then the following assertions hold:
\begin{itemize}

\item[\rm (i)] If every element of $\DD L$ has finite index,
then $L$ is locally finitely distributive%
\index{distributive!locally finitely ---}.

\item[\rm (ii)] If $L$ is relatively complemented and locally
finitely distributive%
\index{distributive!locally finitely ---}, then every element
of $\DD L$ has finite index.

\end{itemize}
\end{theorem}

\begin{proof} (i) Let $a\leq b$ in $L$ and let $n$ be the index
of
$\Dim(a,b)$ in $\DD L$. Let
\(\delta=\vv<a_0,\ldots,a_n,e>\) be an
$n+1$-diamond\index{diamond} in the interval $[a,\,b]$. Then we
have
\begin{align}
\Dim(a,b)&=\Dim(a,0_\delta)+\Dim(0_\delta,a_0)+\cdots
+\Dim(a_0\vee\cdots\vee a_{n-1},1_\delta)+\Dim(1_\delta,b)
\notag\\
&\geq\sum_{k\leq n}\Dim(0_\delta,a_k)\quad
\tag*{\text{(because the $a_k$'s are independent above
$0_\delta$)}}\\
&=(n+1)\cdot\Dim(0_\delta,a_0),\notag
\end{align}
thus, by the definition of $n$,
\(\Dim(0_\delta,a_0)=0\), that is, $\delta$ is trivial. We
end the proof by using Theorem~\ref{T:nDistrdiam}.

(ii) Express $L$ as the direct union of its closed intervals,
viewed as sublattices, the transition maps being the inclusion
maps. Every closed interval $K$ of $L$ is a locally finitely
distributive%
\index{distributive!locally finitely ---} complemented modular
lattice, thus, by Corollary~\ref{C:nDistrCanc}, every element
of $\DD K$ has finite index. Moreover, if $K_0$ and
$K_1$ are closed intervals of $L$ such that
$K_0\subseteq K_1$, then the natural homomorphism from
$\DD K_0$ to $\DD K_1$ is by Lemma~\ref{L:nDistrVemb} a
\Vemb\index{Vemb@\Vemb}. Since the $\DD$ functor preserves
direct limits (see Proposition~\ref{P:Vfunctor}), every element of
$\DD L$ has finite index.
\end{proof}

\begin{remark}[Christian Herrmann]
One cannot generalize
Theorem~\ref{T:FinIndnDistr}.(ii) to arbitrary modular
lattices. For
example, take for $L$ the (Arguesian%
\index{lattice!Arguesian ---}) lattice of all additive
subgroups of \(\ZZ^2\) (to prove the $2$-distributivity%
\index{distributive!$n$- ---} of $L$, one can use
Theorem~\ref{T:nDistrdiam} and elementary properties of
abelian groups). If one puts
\(a=\ZZ\times\{0\}\), \(b=\{0\}\times 2\ZZ\), and
\(c=2\ZZ\times\{0\}\), then \(a\sim b\) and \(b\sim c\), thus
\(\Dim(\{0\},c)=\Dim(\{0\},a)=\Dim(\{0\},c)+\Dim(c,a)\), while
\(\Dim(c,a)>0\). In particular, $\DD L$ is not cancellative.

Another example is provided by the lattice considered in
\cite{DaHW72}.
\end{remark}

We can also recover from Proposition~\ref{P:transperspchar}
the ring-theoretical result that if $R$ is a
regular ring\index{ring!(von Neumann) regular ---}, then
$V(R)$ is cancellative if and only if
perspectivity\index{perspectivity} is transitive in
$\mathcal{L}(2R_R)$ \cite[Theorem 4.24]{Good91}.

\begin{corollary}\label{C:RegRingCounterex}
There exists a regular%
\index{ring!(von Neumann) regular ---} ring $R$ such that
\(\DD\mathcal{L}(R_R)\not\cong V(R)\).
\end{corollary}

\begin{proof} Let $R$ be G.~M. Bergman's example
\cite[Example 4.26]{Good91}. Then
perspectivity\index{perspectivity} is transitive in
$\mathcal{L}(R_R)$, thus
$\DD\mathcal{L}(R_R)$ is cancellative by
Proposition~\ref{P:transperspchar}. Nevertheless, $R$ is not
unit-regular%
\index{ring!unit-regular ---}, that is,
$V(R)$ is not cancellative.
\end{proof}

Note that there are at least two distinct normal equivalences
on Bergman's lattice $\mathcal{L}(R_R)$: namely, isomorphy and
projectivity by decomposition.

The same example yields also the following result:

\begin{corollary}\label{C:V(K)notembV(L)}
There exist a
complemented modular lattice $L$ and an
ideal\index{ideal!of a semilattice}
$K$ of $L$ such that $K$ and $L$
are normal\index{normal!lattice} and the natural map
$f\colon \DD K\to\DD L$ is
surjective but not one-to-one.
\end{corollary}

\begin{proof}
Again, let $R$ be G.~M. Bergman's example
\cite[Example 4.26]{Good91}. Put $L=\mathcal{L}(2R_R)$ and
$K=\mathcal{L}(R_R)$; identify $K$
with an ideal\index{ideal!of a semilattice} of $L$
\emph{via} $I\mapsto I\times\{0\}$. It is clear that $f$ is
surjective. Nevertheless, $\DD L$ is, by
Proposition~\ref{P:pinRiso}, isomorphic to
$V(R)$, while $\DD K$ is not (by the proof of
Corollary~\ref{C:RegRingCounterex}).
\end{proof}

\chapter[Normal kernel]%
{The normal
kernel of a sectionally complemented modular lattice}
\label{NormKer}

We shall introduce in this chapter the last technical tools that
are required to establish a nontrivial dimension theory of
non-coordinatizable lattices. For any sectionally complemented
modular lattice $L$, we define the \emph{normal kernel}, $\nor L$,
of $L$, in a fairly straightforward fashion and we proceed by
proving that
$\nor L$ is \emph{large} in the sense that every homogeneous
sequence of at least four elements of $L$ is contained in $\nor L$,
see Proposition~\ref{P:mLinNor(L)}. To prove this result, we use a
powerful generalization of the von~Neumann
Coordinatization Theorem, due to B.~J\'onsson. The largeness of
the normal kernel makes it possible to infer many previously
inaccessible normality results, such as
Corollary~\ref{C:SimpleNorm}, Theorem~\ref{T:transperspchar}, and
Theorem~\ref{T:MeetContNor}.

\section{The normal kernel and the
ideals\index{ideal!of a semilattice} $mL$}

In this section, we shall define the normal kernel, $\nor L$, of a
sectionally complemented modular lattice $L$. Furthermore, we
shall prove that all terms of any homogeneous sequence of at least
four elements of $L$ belongs to $\nor L$, by using the J\'onsson
form of the von~Neumann Coordinatization Theorem.

\begin{definition}\label{D:NorKer} Let $L$ be a sectionally
complemented modular lattice. Define the
\emph{normal\index{normal!kernel|ii} kernel}%
\index{nzzormkern@$\nor L$|ii}
of $L$ as a subset of $L$:
\[
\nor L=\{x\in L\mid \Theta(x)\text{ is normal}\}.
\]
Recall that $\Theta(x)$ is an
ideal\index{ideal!of a semilattice}, thus a sublattice, of $L$.
\end{definition}

\begin{lemma}\label{L:NorGen}
Let $L$ be a sectionally
complemented modular lattice, let \(a\in L\), let $S$ be a
lower subset of $L$ such that the following assertions hold:
\begin{itemize}

\item[\rm (i)] Every element of $L$ which is
subperspective\index{subperspectivity} to
$a$ is a finite join of elements of $S$.

\item[\rm (ii)]
If \(x\approx y\) and \(x\perp y\), then
\(x\sim y\), for all elements $x$ and $y$ of $S$.

\end{itemize} Then $a$ belongs to \(\nor L\).
\end{lemma}

\begin{proof}
First, every element of \(\Theta(a)\) is, by
Corollary~\ref{C:Jonsson}, a finite join of elements of $L$
which are subperspective\index{subperspectivity} to $a$, thus,
by assumption (i), a finite join of elements of $S$. Since $L$
is sectionally complemented modular, every element of
\(\Theta(a)\) can be written in the form \(\oplus_{i<n}x_i\),
where \(x_i\lesssim a\), for all \(i<n\).

Now let $x$ and $y$ be elements of \(\Theta(a)\) such that
\(x\approx y\) and \(x\perp y\). By
the previous paragraph, there are decompositions
\(x=\oplus_{i<m}x_i\) and
\(y=\oplus_{j<n}y_j\) such that all \(x_i\)'s and all \(y_j\)'s
belong to $S$. Since \(x\approx y\), there exists, by
Corollary~\ref{C:eqwords}, a
$\approxeq$-refinement\index{refinement!matrix} matrix of the
following form:
\[
\begin{tabular}{|c|c|c|c|c|}
\cline{2-5}
\multicolumn{1}{l|}{} & $y_0$ & $y_1$ & $\ldots$ &
$y_{n-1}$\tvi\\
\hline
$x_0$\tvi & $u_{00}/v_{00}$ & $u_{01}/v_{01}$ & $\ldots$ &
$u_{0,n-1}/v_{0,n-1}$\\
\hline
$x_1$\tvi & $u_{10}/v_{10}$ & $u_{11}/v_{11}$ & $\ldots$ &
$u_{1,n-1}/v_{1,n-1}$\\
\hline
$\vdots$ & $\vdots$ & $\vdots$ & $\ddots$ & $\vdots$\\
\hline
$x_{m-1}$\tvi & $u_{m-1,0}/v_{m-1,0}$ & $u_{m-1,1}/v_{m-1,1}$ &
$\ldots$ & $u_{m-1,n-1}/v_{m-1,n-1}$\\
\hline
\end{tabular}
\]
It follows that both \(u_{ij}\) and
\(v_{ij}\) belong to $S$ and \(u_{ij}\approxeq v_{ij}\), for all
\(i<m\) and all \(j<n\); thus, since $S$ is a lower subset of $L$,
by Lemma~\ref{L:charnorm} and by the assumption on $S$,
\(u_{ij}\sim v_{ij}\). Since
\(x=\oplus_{i<m;\ j<n}u_{ij}\),
\(y=\oplus_{i<m;\ j<n}v_{ij}\), and \(x\perp y\), we obtain, by
Lemma~\ref{L:BasicAddPersp}, that \(x\sim y\). Thus $a$ belongs to
\(\nor L\).
\end{proof}

\begin{remark}
By applying Lemma~\ref{L:NorGen} to
\(S=\{x\in L\mid x\lesssim a\}\), we see that \emph{$a$ belongs
to \(\nor L\) if and only if for all
elements $x$ and $y$ of $L$ such that $x\lesssim a$ and
$y\lesssim a$, if $x\approx y$ and $x\wedge y=0$, then $x\sim y$}.
\end{remark}

\begin{proposition}\label{P:NorId} Let $L$ be a sectionally
complemented modular lattice. Then $\nor L$ is a neutral
ideal\index{ideal!neutral --- (in lattices)} of $L$.
Furthermore, \(\nor L\) is
normal\index{normal!lattice}.
\end{proposition}

\begin{proof}
Put \(N=\nor L\). It is obvious that
\(0\in N\) and that $N$ is closed under
perspectivity\index{perspectivity}. Furthermore, let
\(a\leq b\) in $L$ and \(b\in N\). By definition, \(\Theta(b)\) is a
normal\index{normal!lattice} lattice. Furthermore,
\(\Theta(a)\) is a neutral ideal
\index{ideal!neutral --- (in lattices)} of \(\Theta(b)\), thus
two elements of \(\Theta(a)\) are projective%
\index{projectivity!of elements} (resp.,
perspective\index{perspectivity}) in
\(\Theta(b)\) if and only if they are projective%
\index{projectivity!of elements} (resp.,
perspective\index{perspectivity}) in \(\Theta(a)\): therefore,
\(\Theta(a)\) is also a normal\index{normal!lattice} lattice,
so that
\(a\in N\).

Finally, let $a$ and $b$ be elements of $N$; we prove that
\(a\vee b\in N\). By replacing $b$ by a sectional complement
of \(a\wedge b\) in $b$, one may suppose without loss of
generality that \(a\perp b\). Put
\(S=\Theta(a)\cup\Theta(b)\); note that $S$ is a lower subset
of $L$, closed under perspectivity\index{perspectivity}.
Moreover, an easy combination of Lemma~\ref{L:decomp} and
Corollary~\ref{C:Jonsson} yields that every element of $L$
which is subperspective\index{subperspectivity} to
\(a\oplus b\) is a finite join of elements of $S$. Thus, by
Lemma~\ref{L:NorGen}, it suffices to verify that
for all $x$, $y\in S$, \(x\approx y\) and
\(x\perp y\) imply that \(x\sim y\). Since $x$ and $y$
belong simultaneously either to \(\Theta(a)\) or to
\(\Theta(b)\), the conclusion follows from the
normality\index{normal!lattice} of both \(\Theta(a)\) and
\(\Theta(b)\). Hence \(a\vee b\in N\).

At this point, we know that $N$ is a neutral
ideal\index{ideal!neutral --- (in lattices)} of $L$. Finally,
for all elements $x$ and $y$ of $N$ such that
\(x\perp y\) and \(x\approx y\) (the latter means the same in $N$
as in $L$), \(x\vee y\in N\) and \(\Theta(x\vee y)\) also
satisfies \(x\perp y\) and \(x\approx y\); since
\(\Theta(x\vee y)\) is normal\index{normal!lattice}, we
obtain the conclusion that \(x\sim y\). Hence $N$ is
normal\index{normal!lattice}.
\end{proof}

For the followiong definition, we refer to
B. J\'onsson \cite[Definition 1.1]{Jons60}.

\begin{definition}\label{D:PartFrame}
Let $L$ be a complemented
modular lattice, let $n$ be a positive integer.
\begin{itemize}

\item[(1)] A \emph{partial $n$-frame} of $L$ is a triple
\(\alpha=\vv<\boldsymbol{a},\boldsymbol{b},c>\), where
\(\boldsymbol{a}=\vv<a_i\mid i<n>\) is a finite sequence
of elements of $L$,
\(\boldsymbol{b}=\vv<b_{ij}\mid\vv<i,j>\in n\times n>\)
is a matrix of elements of $L$, $c$ is an element of $L$, and
the following assertions hold:
\begin{itemize}

\item[\rm (i)] \(1=(\oplus_{i<n}a_i)\oplus c\) (in particular,
$\boldsymbol{a}$ is independent\index{independent}).

\item[\rm (ii)]
\(a_i\oplus b_{ij}=a_i\vee a_j\) and \(b_{ij}=b_{ji}\),
for all $i$, $j<n$.

\item[\rm (iii)]
\(b_{ij}=(b_{ik}\vee b_{jk})\wedge(a_i\vee a_j)\),
for all $i$, $j$, $k<n$.

\end{itemize}
Then we put \(\alpha_i=a_i\),
\(\alpha_{ij}=b_{ij}\) and \(\alpha^*=c\), for all $i$, $j<n$.

\item[(2)] A partial $n$-frame $\alpha$ is \emph{large}, if
\(\alpha^*\in\Theta(\alpha_0)\).

\end{itemize}
By Corollary~\ref{C:Jonsson}, $\alpha$
is large if and only if
$c$ is a finite join of elements which are
subperspective\index{subperspectivity}
(Definition~\ref{D:persp}) to $a_0$.
\end{definition}

\begin{note}
If $\alpha$ is a partial $n$-frame, then
we necessarily have \(b_{ii}=0\), for all \(i<n\). Once this is
set, conditions (ii) and (iii) need to be checked only for
pointwise distinct values of $i$, $j$, and $k$.
\end{note}

We shall make use of the following important improvement, due
to J\'onsson, of von~Neumann's Coordinatization%
\index{von Neumann Coordinatization Theorem!
J\'onsson improvement}
Theorem, see B. J\'onsson \cite[Corollary 8.4]{Jons60}:

\begin{theorem}\label{T:Jonsson}
Let $L$ be a complemented
modular lattice having a large partial $4$-frame.
Then $L$ is \emph{coordinatizable},
that is, there exists a regular%
\index{ring!(von Neumann) regular ---} ring $R$ such that
$L\cong\mathcal{L}(R_R)$.\qed
\end{theorem}

In fact, we shall only use the consequence that if $L$ has a
large partial $4$-frame, then $L$ is
normal\index{normal!lattice} (recall that every
$\mathcal{L}(R_R)$ is normal%
\index{normal!lattice}, see Lemma~\ref{L:SubmNorm}).

\begin{definition}\label{D:mL} Let $L$ be a modular%
\index{lattice!modular (not necessarily complemented) ---}
lattice, let $m$ be a positive integer.
\begin{itemize}
\item[(i)] Let \(\diam{m}{L}\)%
\index{DzziamRel@$\diam{m}{L}$|ii} denote the
congruence\index{congruence!lattice ---} of $L$ generated by
all pairs \(\vv<0_\delta,1_\delta>\), where $\delta$ is an
$m$-diamond\index{diamond} of $L$.

\item[(ii)] Suppose that $L$ is sectionally complemented. Then
let $mL$ be the neutral ideal
\index{ideal!neutral --- (in lattices)} of $L$ generated by the
elements $a$ of $L$ such that there exists a
homogeneous\index{homogeneous!sequence} sequence
$\vv<a_i\mid i<m>$ with $a_0=a$.
\end{itemize}
\end{definition}

\begin{remark}
It is easy to see that
\(\diam{m+1}{L}\subseteq\diam{m}{L}\) and
\((m+1)L\subseteq mL\). It can also be proved that
\(L/\diam{m}{L}\) is the maximal
$(m-1)$-distributive%
\index{distributive!$n$- ---} quotient of $L$ (this is an
immediate consequence of
Theorem~\ref{T:ProjFrame} and of Theorem~\ref{T:nDistrdiam}).
\end{remark}

\begin{proposition}\label{P:DiammL} Let $L$ be a sectionally
complemented modular lattice and let
$m$ be a positive integer. Then the equivalence
\[
x\in mL\Longleftrightarrow x\equiv 0\pmod{\diam{m}{L}}
\]
holds for all \(x\in L\).
In particular, \(L/mL=L/\diam{m}{L}\).
\end{proposition}

\begin{proof} Let \(\vv<a_i\mid i<m>\) be a
homogeneous\index{homogeneous!sequence} sequence of $L$. As in
the beginning of the proof of Lemma~\ref{L:HomDiam}, there
exists a diamond\index{diamond} of the form
\(\delta=\vv<a_0,\ldots,a_{m-1},e>\). Since \(0_\delta=0\)
and \(a_0\leq\oplus_{i<m}a_i=1_\delta\), the relation
\(a_0\equiv 0\pmod{\diam{m}{L}}\) holds. This proves the direction
$\Rightarrow$.

Conversely, let \(x\in L\) be such that
\(x\equiv 0\pmod{\diam{m}{L}}\). By
\cite[Corollary III.1.4]{Grat}, there exist a decomposition
\(0=x_0\leq x_1\leq\cdots\leq x_n=x\) and
$m$-diamonds\index{diamond}
\(\delta_j\) (for \(j<n\)) such that
\([x_j,\,x_{j+1}]\wpr[0_{\delta_j},\,1_{\delta_j}]\) for all
\(j<n\). Let $y_j$ (resp., $z_j$) be such that
\(x_j\oplus y_j=x_{j+1}\) (resp.,
\(0_{\delta_j}\oplus z_j=1_{\delta_j}\)). Then $z_j$ is
perspective\index{perspectivity} to the finite join of $m$
elements of \(mL\), thus
\(z_j\in mL\). Since
\([0,\,y_j]\wpr[0,\,z_j]\),
we also have
\(y_j\in mL\). Since \(x=\oplus_{j<n}y_j\), it follows that
\(x\in mL\).

The end of the statement of Proposition~\ref{P:DiammL} follows
immediately from the fact that $L$ is sectionally complemented.
\end{proof}

\begin{lemma}\label{L:eltsmL}
Let $L$ be a sectionally complemented modular lattice. Let $S_m$
denote the set of all elements $x$ of $L$ such that there exists a
homogeneous\index{homogeneous!sequence} sequence
\(\vv<x_i\mid i<m>\) satisfying that
\(x=x_0\). Then $S_m$ is a lower subset of $L$,
and $mL$ is the set of finite joins of elements of $S_m$.
\end{lemma}

\begin{proof}
It is trivial that $S_m$ is a lower subset of $L$. Let
$I$ be the set of all finite joins of elements of $S_m$.
It is trivial that $I$ is closed under finite joins.

By Corollary~\ref{C:Jonsson}, in order to prove
that $I$ is a neutral ideal
\index{ideal!neutral --- (in lattices)} of $L$, it suffices to
prove that $x\in I$ implies that $y\in I$, for all \(x\in S_m\)
and all \(y\in L\) such that
\(y\lesssim x\), $y$ belongs to $I$. Since \(S_m\) is a lower
subset of $L$, it is sufficient to consider the case
\(x\sim y\). By definition, there exists a homogeneous
sequence \(\vv<x_i\mid i<m>\) such that \(x_0=x\). Put
\(u=y\vee\oplus_{i<m}x_i\). Since $x$ and $y$ are
perspective\index{perspectivity} and both lie below $u$, there
exists \(y'\leq u\) such that \(x\oplus y'=y\oplus y'=u\).
Furthermore, there exists $x'$ such that
\(\oplus_{0<i<m}x_i\leq x'\leq u\) and
\(x\oplus x'=u\). In particular, \(x\oplus x'=x\oplus y'=u\),
thus \(x'\sim y'\). Since \(\vv<x_i\mid 0<i<m>\) is a
homogeneous sequence of length $m-1$ below $x'$, there exists
a homogeneous sequence \(\vv<y_i\mid 0<i<m>\) below $y'$
such that $x_i\sim y_i$, for all $i$ with $0<i<m$. It follows
that the $y_i$, $0\leq i<m$, are mutually projective. A
simple application of Corollary~\ref{C:IncrCong} yields then
that $y$ is a finite join of elements of $S_m$; thus it
belongs to $I$.
\end{proof}

\begin{proposition}\label{P:mLinNor(L)}
Let $L$ be a sectionally complemented modular lattice. Then
\(4L\subseteq\nor L\).
\end{proposition}

\begin{proof}
Put $N=\nor L$. By definition, $4L$ is
the neutral ideal\index{ideal!neutral --- (in lattices)}
generated by $S_4$, thus, by Proposition~\ref{P:NorId}, it
suffices to prove that
$S_4\subseteq N$. Thus, let
\(\boldsymbol{a}=\vv<a_i\mid i<4>\) be a
homogeneous\index{homogeneous!sequence} sequence of length $4$
in $L$; we prove that $\Theta(a_0)$ is
normal\index{normal!lattice}. Let
$x$, $y\in\Theta(a_0)$ be such that $x\perp y$ and
$x\approx y$. By definition, there are $n\in\omega$ and
elements $z_j$ (for \(0\leq j\leq n\)) such that $z_0=x$, $z_n=y$,
and \(z_j\sim z_{j+1}\), for all $j<n$. Define an element
$u$ of $L$ by
\[
u=(\oplus_{i<4}a_i)\vee\bigvee_{j\leq n}z_j,
\]
and let $c$ be an element of $L$ such that
\((\oplus_{i<4}a_i)\oplus c=u\). Since \(\vv<a_i\mid
i<4>\) is a homogeneous\index{homogeneous!sequence} sequence,
there are, by
\cite[Lemma II.5.3]{Neum60}, elements
$b_{ij}$ (for $i$, $j<4$) of $L$ such that
\(\boldsymbol{b}=\vv<b_{ij}\mid\vv<i,j>\in n\times n>\)
satisfies (i) and (ii) of the definition of a partial
$4$-frame (we put, of course, \(b_{ii}=0\), for all $i<4$).
Since \(\Theta(z_j)=\Theta(x)\subseteq\Theta(a_0)\), we have
\(\Theta(u)=\Theta(a_0)\), whence \(c\in\Theta(a_0)\).
Therefore, \(\vv<\boldsymbol{a},\boldsymbol{b},c>\) is a large
partial $4$-frame of $[0,\,u]$. By Theorem~\ref{T:Jonsson},
$[0,\,u]$ is coordinatizable, thus
normal\index{normal!lattice}. Since $x\perp y$ and
$x\approx y$ in
$[0,\,u]$, it follows that $x\sim y$. Hence $\Theta(a_0)$ is
normal\index{normal!lattice}, that is, $a_0\in N$.
\end{proof}

\section{Quotient of a lattice by its
normal\index{normal!lattice} kernel}

We shall now recall another very useful result, due to
R. Freese \cite{Free76} (see also the proof in A.~P. Huhn
\cite[Satz 2.1]{Huhn72}, and \cite[Corollary 3.17]{JiRo92}):

\begin{theorem}\label{T:ProjFrame}
Let $L$ and $M$ be modular lattices%
\index{lattice!modular (not necessarily complemented) ---},
let $n$ be a positive integer, and let
\(f\colon M\twoheadrightarrow L\) be a surjective homomorphism.
Let $\delta$ be an $n$-diamond\index{diamond} in
$L$. Then there exists an $n$-diamond\index{diamond}
$\hat\delta$ in $M$ such that
\(f(\hat\delta)=\delta\).\qed
\end{theorem}

\begin{corollary}\label{C:L/Nor(L)FinInd} The lattice
\(L/\nor L\) is $3$-distributive%
\index{distributive!$n$- ---}.
\end{corollary}

\begin{proof}
Let \(\pi\colon L\twoheadrightarrow L/\nor L\)
be the canonical homomorphism. If
$\delta$ is a $4$-diamond\index{diamond} in
\(L/\nor L\), then there exists a
$4$-diamond\index{diamond}
$\hat\delta$ of $L$ such that
\(\pi(\hat\delta)=\delta\). Put
\(\hat\delta=
\vv<\hat{a}_0,\hat{a}_1,\hat{a}_2,\hat{a}_3,\hat{e}>\) and
let $x_i$ (for $i<4$) be such that
\(0_{\hat\delta}\oplus x_i=\hat{a}_i\) (for all $i<4$). Then
\(\vv<x_i\mid i<4>\) is a homogeneous%
\index{homogeneous!sequence} sequence of length $4$ in $L$,
thus, by definition,
$x_i\in 4L$ (for all $i<4$). It follows then from
Proposition~\ref{P:mLinNor(L)} that all $x_i$'s lie in
$\nor L$. Applying $\pi$ yields that
$\delta$ is trivial. The conclusion follows by
Theorem~\ref{T:nDistrdiam}.
\end{proof}

By Corollary~\ref{C:nDistrTrans}, we deduce immediately
the following result:

\begin{corollary}\label{C:L/Nor(L)TransPers} Let $L$ be a
sectionally complemented modular lattice. Then
perspectivity\index{perspectivity} is transitive in
$L/\nor L$.\qed
\end{corollary}

In particular, both $\nor L$ and $L/\nor L$ are
normal\index{normal!lattice} lattices (although, by the
counterexample of
Section~\ref{S:NonNorm}, $L$ may not be normal%
\index{normal!non- --- lattice}). Nevertheless, our results will
be sufficient to obtain normality\index{normal!lattice} in
quite a number of cases.

About general lattices (without necessarily a zero), one can still
formulate the following corollary:

\begin{corollary}\label{C:Loverdiam}
Let $L$ be a relatively
complemented modular lattice and let
$m$ be a positive integer. Then the following properties hold:
\begin{itemize}
\item[\rm (i)] Perspectivity is transitive in
\(L/\diam{m}{L}\).

\item[\rm (ii)] Let $a$, $b\in L$ be such that
\(a\equiv b\pmod{\diam{4}{L}}\). If
\(\Dim(a\wedge b,a)=\Dim(a\wedge b,b)\), then \(a\sim b\). In
particular, every
\(\diam{4}{L}\)-equivalence class is
normal\index{normal!lattice}.
\end{itemize}
\end{corollary}

\begin{proof}
It follows immediately from
Theorem~\ref{T:ProjFrame} that
\(L/\diam{m}{L}\) has no non
trivial
$m$-diamond\index{diamond}. Therefore, by
Theorem~\ref{T:nDistrdiam}, \(L/\diam{m}{L}\) is
$(m-1)$-distributive%
\index{distributive!$n$- ---}. Therefore, (i) follows from
Corollary~\ref{C:nDistrTrans}.

Now let us prove (ii). Let $0$ (resp., $1$) denote the infimum
(resp., supremum) of all parameters in $L$ witnessing
\(a\equiv b\pmod{\diam{4}{L}}\) and
\(\Dim(a\wedge b,a)=\Dim(a\wedge b,b)\). Then those relations
also hold in \([0,\,1]\), so that one can suppose, without loss
of generality, that \(L=[0,\,1]\). Let $a'$ (resp., $b'$) be a
sectional complement of \(a\wedge b\) in $a$ (resp., $b$). Then
\begin{equation}\label{Eq:a'ortperpb'}
a'\perp b'\text{\quad and\quad}\Dim(a')=\Dim(b').
\end{equation}
Furthermore, it follows from
\(a\equiv b\pmod{\diam{4}{L}}\) that
\(a'\equiv b'\equiv 0\pmod{\diam{4}{L}}\), thus, by
Proposition~\ref{P:DiammL}, both $a'$ and $b'$ belong to
\(4L\). Then it follows from Proposition~\ref{P:mLinNor(L)}
that both $a'$ and $b'$ belong to \(\nor L\). Hence, it
follows from (\ref{Eq:a'ortperpb'}) that \(a'\sim b'\).
Therefore, by Corollary~\ref{C:BasicAddPersp}, \(a\sim b\).

The last statement of (ii) follows easily.
\end{proof}

\begin{corollary}\label{C:SimpleNorm} Every
simple\index{lattice!simple ---} relatively complemented
modular lattice is normal\index{normal!lattice}.
\end{corollary}

\begin{proof} Since $L$ is simple\index{lattice!simple ---},
\(\diam{4}{L}\) is equal either to
\(\mathrm{id}_L\) or to \(L\times L\). In the first case,  we
obtain the conclusion by Corollary~\ref{C:Loverdiam}.(i). In the
second case, every \(\diam{4}{L}\)-equivalence
class is equal to $L$, so that
the conclusion follows by Corollary~\ref{C:Loverdiam}.(ii).
\end{proof}

Proposition~\ref{P:NorId} and
Corollary~\ref{C:L/Nor(L)TransPers} also make it possible to improve
Proposition~\ref{P:transperspchar}:

\begin{theorem}\label{T:transperspchar}
Let $L$ be a relatively complemented modular lattice. Then
\(\DD L\) is cancellative if and only if
perspectivity\index{perspectivity} is transitive in $L$.
\end{theorem}

\begin{proof} By Proposition~\ref{P:transperspchar}, it
remains to prove that if \(\DD L\) is cancellative, then $L$ is
normal\index{normal!lattice}. Thus let \(p\leq q\) in $L$; we
prove that \(K=[p,\,q]\) is normal\index{normal!lattice}. So
let $a$ and $b$ be elements of
$K$ which are independent and projective%
\index{projectivity!of elements} in $K$; we shall prove that
\(a\sim b\). Denote by
\(\pi\colon K\twoheadrightarrow K/4K\) the canonical
projection. First, \(\pi(a)\approx\pi(b)\) in
\(K/4K\), thus, by  Proposition~\ref{P:DiammL} and
Corollary~\ref{C:Loverdiam}.(i), \(\pi(a)\sim\pi(b)\),
that is, there exists
\(c\in K\) such that \(a\wedge c,b\wedge c\in 4K\) and
\(a\vee c\equiv b\vee c\pmod{4K}\). Hence there exists
\(c'\in 4K\) such that \(a\vee c\vee c'=b\vee c\vee c'\),
thus, after having replaced $c$ by \(c\vee c'\), one can
suppose that \(a\vee c=b\vee c\). Let $a'$, $b'\in K$
be such that
$K$ satisfies \(a=a\wedge c\oplus a'\) and
\(b=b\wedge c\oplus b'\). Then $K$ satisfies that
\(a'\oplus c=b'\oplus c\), so that we obtain
\begin{equation}\label{Eq:a'simb'}
a'\sim b'.
\end{equation}
On the other hand, $a$ and $b$ are projective%
\index{projectivity!of elements} in $K$, thus
\(\Dim_K(p,a)=\Dim_K(p,b)\). Thus \(\Dim_L(p,a)=\Dim_L(p,b)\),
which can be written as
\begin{equation}\label{Eq:Decaprb}
\Dim_L(p,a\wedge c)+\Dim_L(a\wedge c,a)=
\Dim_L(p,b\wedge c)+\Dim_L(b\wedge c,b).
\end{equation}
However, by (\ref{Eq:a'simb'}), we also have
\[
\Dim_L(a\wedge c,a)=\Dim_L(p,a')=\Dim_L(p,b')=
\Dim_L(b\wedge c,b),
\]
thus, by applying to (\ref{Eq:Decaprb}) the hypothesis that
\(\DD L\) is cancellative,
\begin{equation}\label{Eq:paceqpbc}
\Dim_L(p,a\wedge c)=\Dim_L(p,b\wedge c).
\end{equation} However, \(a\wedge c\) and \(b\wedge c\) are
independent in $K$ and both belong to \(4K\), thus, by
Proposition~\ref{P:DiammL},
\(a\wedge c\equiv b\wedge c\pmod{\diam{4}{K}}\), thus
\emph{a fortiori}
\(a\wedge c\equiv b\wedge c\pmod{\diam{4}{L}}\). Therefore,
it follows from Corollary~\ref{C:Loverdiam}.(ii) that we have
\begin{equation}\label{Eq:acsimbc}
a\wedge c\sim b\wedge c.
\end{equation}
Finally, it follows from (\ref{Eq:a'simb'}),
(\ref{Eq:acsimbc}), and Lemma~\ref{L:BasicAddPersp} that
\(a\sim b\). So we have proved that
$L$ is normal\index{normal!lattice}.
\end{proof}

In order to state Corollary \ref{C:transperspchar},
we say that
an element $c$ of a \cm\ $M$ is \emph{cancellable}%
\index{cancellable|ii} in $M$, if
$M$ satisfies the statement
\[
(\forall x,y)(x+c=y+c\Longrightarrow x=y).
\]

\begin{corollary}\label{C:transperspchar}
Let $L$ be a
sectionally complemented modular lattice. Then every element
$c$ of $L$ such that \(\Dim_{\Theta(c)}(c)\) is cancellable
in \(\DD(\Theta(c))\), and thus belongs to \(\nor L\).
\end{corollary}

In particular, by Lemma~\ref{L:ApproxCanc} and
Proposition~\ref{P:twoindex}, every element of $L$ with finite
index satisfies the assumption above, thus belongs to
\(\nor L\).

\begin{proof} Put \(I=\Theta(c)\); then $I$ is a neutral
ideal\index{ideal!neutral --- (in lattices)} of $L$. Then
\(\Dim_I(c)\) is an order-unit of \(\DD I\) and
\(\DD I\) is a
refinement\index{monoid!refinement ---} monoid, thus, in order
to prove that \(\DD I\) is
cancellative, it suffices, by Lemma~\ref{L:cancref}, to prove
that the interval
\([0,\,\Dim_I(c)]\) is cancellative as a partial monoid.
Let $\xi$, $\eta$, $\zeta\in[0,\,\Dim_I(c)]$ be such that
\(\xi+\zeta=\eta+\zeta\leq\Dim_I(c)\). Since
\(\zeta\leq\Dim_I(c)\), this implies that
\(\xi+\Dim_I(c)=\eta+\Dim_I(c)\), thus, by assumption,
\(\xi=\eta\). Hence we have proved that
\(\DD I\) is cancellative.
Therefore, by Theorem~\ref{T:transperspchar},
\(\Theta(c)=I\) is normal\index{normal!lattice}, that is,
\(c\in N\).
\end{proof}

\section{The case of $\aleph_0$-meet-continuous%
\index{lattice!meet-continuous ---!$\aleph_0$-\mcont\ ---}
lat\-tices}

The following Theorem~\ref{T:MeetContNor} will yield us a new
class of normal\index{normal!lattice} lattices; its analogue
for
\mcont\index{lattice!meet-continuous ---} lattices was
proved in I. Halperin and J. von~Neumann
\cite[Theorem 6]{HaNe40}. We shall first prove a lemma.

\begin{lemma}\label{L:ctbleref}
Let $L$ be a sectionally
complemented modular conditionally $\aleph_0$-\mcont%
\index{lattice!meet-continuous ---!$\aleph_0$-\mcont\ ---}
lattice. Let
$\vv<a_n\mid n\in\omega>$ and
$\vv<b_n\mid n\in\omega>$ be independent\index{independent}
sequences of elements of $L$ such that
\(\oplus_ia_i=\oplus_jb_j\). Then there exists an infinite
$\sim_2$-refinement\index{refinement!matrix} matrix as
follows:
\[
\begin{tabular}{|c|c|}
\cline{2-2}
\multicolumn{1}{l|}{} & $b_n\,(n\in\omega)$\tvi\\
\hline
$a_m\,(m\in\omega)$\tvi & $x_{mn}/y_{mn}$\\
\hline
\end{tabular}
\]
\end{lemma}

\begin{proof}
We put
\(\bar a_n=\oplus_{i<n}a_i\) and \(\bar b_n=\oplus_{i<n}b_i\), for all
$n\leq\omega$.
Furthermore, put
\(c_{mn}=(\bar a_{m+1}\wedge\bar b_n)\vee
(\bar a_m\wedge\bar b_{n+1})\) and
\(d_{mn}=\bar a_{m+1}\wedge\bar b_{n+1}\), for all $m$, $n\in\omega$. Note that
\(c_{mn}\leq b_{mn}\), thus there exists $z_{mn}$ such that
\(c_{mn}\oplus z_{mn}=d_{mn}\).

\setcounter{claim}{0}
\begin{claim}
The equality
\[
(\bar a_m\wedge\bar b_n)\oplus(\oplus_{j<n}z_{mj})=
\bar a_{m+1}\wedge\bar b_n
\]
holds, for all $m<\omega$ and $n\leq\omega$
\end{claim}

\begin{cproof}
We prove the conclusion by induction on $n$. The case
$n=\omega$ follows immediately from the finite cases and the
continuity assumption on $L$, thus it suffices to deal with
the finite case. The conclusion is trivial for $n=0$. Suppose that
the conclusion has been proved for $n$. In particular,
\(\oplus_{j<n}z_{mj}\leq\bar a_{m+1}\wedge\bar b_n\), thus
\begin{align*}
(\bar a_m\wedge\bar b_{n+1})\wedge\oplus_{j<n}z_{mj}
&=(\bar a_m\wedge\bar b_{n+1})
\wedge (\bar a_{m+1}\wedge\bar b_n)\wedge\oplus_{j<n}z_{mj}\\
&=(\bar a_m\wedge\bar b_n)\wedge\oplus_{j<n}z_{mj}\\
&=0,
\end{align*}
so that
\(\vv<\bar a_m\wedge\bar b_{n+1},z_{m0},\ldots,z_{m,n-1}>\)
is independent\index{independent}. Furthermore,
\begin{align*}
(\bar a_m\wedge\bar b_{n+1})\oplus(\oplus_{j<n}z_{mj})
&=(\bar a_m\wedge\bar b_{n+1})\vee(\bar a_m\wedge\bar b_n)
\vee(\oplus_{j<n}z_{mj})\\
&=(\bar a_{m+1}\wedge\bar b_n)\vee
(\bar a_m\wedge\bar b_{n+1})\\
&=c_{mn},
\end{align*}
whence
\begin{align*} d_{mn}&=c_{mn}\oplus z_{mn}\\
&=\bigl((\bar
a_m\wedge\bar b_{n+1})\oplus (\oplus_{j<n}z_{mj})\bigr)\oplus
z_{mn}\\
&=(\bar a_m\wedge\bar b_{n+1})\oplus (\oplus_{j\leq
n}z_{mj}),
\end{align*}
which completes the proof of the Claim.
\end{cproof}

Similarly, one can prove the following:

\begin{claim}
The equality
\begin{equation*}
(\bar a_m\wedge\bar b_n)\oplus(\oplus_{i<m}z_{in})=
\bar a_m\wedge\bar b_{n+1}
\end{equation*}
holds for all $m\leq\omega$ and $n<\omega$.\qed
\end{claim}

In particular, using Claim~1 for $n=\omega$ yields
that \(\bar a_m\oplus(\oplus_{j<\omega}z_{mj})=
\bar a_{m+1}=\bar a_m\oplus a_m\), so that
\(\oplus_{j<\omega}z_{mj}\sim a_m\). Thus there exists a
decomposition \(a_m=\oplus_{j\in\omega}x_mj\) such that
$x_{mj}\sim z_{mj}$, for
all $j\in\omega$. Similarly,
Claim~2 yields that there exists a decomposition
\(b_n=\oplus_{i\in\omega}y_{in}\) such that \(z_{in}\sim y_{in}\) holds for
all $i\in\omega$. Therefore, the
$x_{ij}$'s and the $y_{ij}$'s satisfy the required
conditions.\end{proof}

\begin{theorem}\label{T:MeetContNor}
Every conditionally $\aleph_0$-\mcont%
\index{lattice!meet-continuous ---!$\aleph_0$-\mcont\ ---}
relatively complemented modular lattice is
normal\index{normal!lattice}.
\end{theorem}

\begin{proof} Let $L$ be a conditionally
$\aleph_0$-\mcont%
\index{lattice!meet-continuous ---!$\aleph_0$-\mcont\ ---}
relatively complemented modular lattice. By
Proposition~\ref{P:PresNorm}.(a) and the fact that $L$ is the
direct union of its closed intervals, one can reduce the
problem to the case where $L$ is a \emph{bounded} lattice. Put
\(N=\nor L\). We start with a
claim:

\setcounter{claim}{0}
\begin{claim}
The ideal\index{ideal!of a semilattice}
$N$ is closed under countable suprema.
\end{claim}

\begin{cproof}
Since $N$ is an
ideal\index{ideal!of a semilattice} of $L$, it suffices to
prove that $N$ is closed under suprema of
independent\index{independent} infinite sequences. Thus let
\(\vv<a_n\mid n\in\omega>\) be an
independent\index{independent} sequence of elements of $N$, let
\(a=\oplus_na_n\), we prove that
\(a\in N\). Put \(S=\bigcup_n\Theta(a_n)\). An easy
combination of Lemma~\ref{L:ctbleref} and
Corollary~\ref{C:Jonsson} yields that every element of
\(\Theta(a)\) can be written as the supremum of an
independent\index{independent} countable sequence of elements
of $S$.

Thus let $x$ and $y$ be elements of \(\Theta(a)\) be such
that \(x\perp y\) and \(x\approx y\). By previous
paragraph, there are decompositions
\(x=\oplus_{m\in\omega}x_m\) and \(y=\oplus_{n\in\omega}y_n\)
such that all $x_m$'s and all $y_n$'s belong to $S$. Hence
there exists a decomposition \(x=\oplus_{n\in\omega}y'_n\),
where \(y_n\approx y'_n\), for all $n$.
Therefore, by Lemma~\ref{L:ctbleref}, there exists an infinite
$\approx$-refinement\index{refinement!matrix} matrix of the
following form:
\[
\begin{tabular}{|c|c|}
\cline{2-2}
\multicolumn{1}{l|}{} & $y_n\,(n\in\omega)$\tvi\\
\hline
$x_m\,(m\in\omega)$\tvi & $u_{mn}/v_{mn}$\\
\hline
\end{tabular}
\]
Now, for all non-negative integers $m$ and $n$, both
\(u_{mn}\) and \(v_{mn}\) belong to $S$, thus to $N$, and
\(u_{mn}\approx v_{mn}\). Since $N$ is
normal\index{normal!lattice}, it follows that
\(u_{mn}\sim v_{mn}\). Since
\(x=\oplus_{m,n}u_{mn}\) and
\(y=\oplus_{m,n}v_{mn}\), it follows from
Lemma~\ref{L:BasicAddPersp} that \(x\sim y\). Hence
\(a\in N\).
\end{cproof}

Now, we shall prove that \(N=L\). Thus let
\(\pi\colon L\twoheadrightarrow L/N\) be the canonical
projection. Let $a$, $b\in L$ be such that \(a\perp b\) and
\(a\approx b\). Let $T$ be a
projective%
\index{projective!isomorphism} isomorphism from
\([0,\,a]\) onto
\([0,\,b]\) (\emph{cf.} Definition~\ref{D:ProjIso}). Then
\(\pi(a)\approx\pi(b)\), thus, by
Corollary~\ref{C:L/Nor(L)TransPers}, \(\pi(a)\sim\pi(b)\). As
in the proof of Theorem~\ref{T:transperspchar}, there exists
\(c\in L\) such that both \(a\wedge c\) and \(b\wedge c\)
belong to $N$ and \(a\vee c=b\vee c\). Define a sequence
\(\vv<c_n\mid n\in\omega>\) by putting \(c_0=c\), and
\begin{equation}\label{Eq:Defcn}
c_{n+1}=c_n\vee T(a\wedge
c_n)\vee T^{-1}(b\wedge c_n)
\end{equation}
for all \(n\in\omega\).
Then it is easy to see that the sequence
\(\vv<c_n\mid n\in\omega>\) is increasing and that
\begin{equation}\label{Eq:Pptiescn}
a\wedge c_n\in
N\ \text{ and }\ b\wedge c_n\in N
\end{equation}
holds for all \(n\in\omega\).
Furthermore, it follows immediately from
(\ref{Eq:Defcn}) that
\begin{equation}\label{Eq:cnalmostOK}
T(a\wedge c_n)\leq b\wedge c_{n+1}
\ \text{ and }\
T^{-1}(b\wedge c_n)\leq a\wedge c_{n+1}.
\end{equation}
Now let \(c_\omega=\bigvee_{n\in\omega}c_n\).
Since $L$ is conditionally $\aleph_0$-\mcont%
\index{lattice!meet-continuous ---!$\aleph_0$-\mcont\ ---}, we
have
\(a\wedge c_\omega=\bigvee_{n\in\omega}(a\wedge c_n)\) and
\(b\wedge c_\omega=\bigvee_{n\in\omega}(b\wedge c_n)\), thus,
by (\ref{Eq:Pptiescn}) and Claim~1, both
\(a\wedge c_\omega\) and \(b\wedge c_\omega\) belong to $N$.
Furthermore, since $T$ is an isomorphism, it follows from
(\ref{Eq:cnalmostOK}) that
\(T(a\wedge c_\omega)=b\wedge c_\omega\), whence
\(a\wedge c_\omega\approx b\wedge c_\omega\). Therefore,
since $N$ is normal\index{normal!lattice} and by
(\ref{Eq:Pptiescn}), we obtain that
\begin{equation}\label{Eq:comOK}
a\wedge c_\omega\sim b\wedge c_\omega.
\end{equation}
But if $u$ and $v$ are elements of $L$ such that
\(a=(a\wedge c_\omega)\oplus u\) and
\(b=(b\wedge c_\omega)\oplus v\), we have
\(u\oplus c_\omega=v\oplus c_\omega\), so that \(u\sim v\).
Therefore, by (\ref{Eq:comOK}) and by
Lemma~\ref{L:BasicAddPersp}, we obtain
\(a\sim b\).
\end{proof}

The following analogue of
Theorem~\ref{T:transperspchar} will show us that in a quite
large class of lattices, all lattices are
normal\index{normal!lattice} and, in fact,
perspectivity\index{perspectivity} by decomposition is
transitive. This will prove
\emph{a posteriori} to be a generalization of both
Theorem~\ref{T:MeetContNor} and Theorem~\ref{T:normeqform},
but we shall not know this until Proposition~\ref{P:V(L)GCA}.
Nevertheless, it will make it possible to verify
further results, for example about
$\aleph_0$-\emph{join}-continuous%
\index{lattice!join-continuous ---!$\aleph_0$-\jcont\ ---}
lattices.

\begin{theorem}\label{T:AxImplNorm}
Let $L$ be a relatively complemented modular lattice. Suppose
that \(\DD L\) satisfies the following axiom:
\begin{equation}\label{Eq:HalfCan}
(\forall x,y,z)\bigl[x+z=y+z\Rightarrow(\exists t)
(2t\leq z\text{ and }x+t=y+t)\bigr].
\end{equation}
Then $L$ is normal\index{normal!lattice}. In
addition, if $L$ has a zero and $a$ and $b$ are elements of
$L$, then \(a\approxeq b\) if and only if there are
decompositions
\[
a=a_0\oplus a_1,\quad b=b_0\oplus b_1
\]
such that \(a_0\sim b_0\) and \(a_1\sim b_1\). Furthermore, if
$a'$ (resp., $b'$) is any sectional complement of \(a\wedge b\)
in $a$ (resp., $b$), one can take \(a_0\perp b_0\),
\(a_0\geq a'\), \(b_0\geq b'\), and $a_1$, $b_1\leq a\wedge b$.
\end{theorem}

\begin{proof}
We first prove normality\index{normal!lattice}.
Thus let \(p\leq q\) in $L$; we prove that \(K=[p,\,q]\) is
normal\index{normal!lattice}. So let $a$ and $b$ be elements of
$K$ that are independent and projective%
\index{projectivity!of elements} in $K$; we shall prove that
\(a\sim b\). Denote by \(\pi\colon K\twoheadrightarrow K/4K\)
the canonical projection. As in the proof of
Theorem~\ref{T:transperspchar}, one proves that there exists
\(c\in K\) such that $a\wedge c$, $b\wedge c\in 4K$ and
$a\vee c=b\vee c$. As in the proof of
Theorem~\ref{T:transperspchar}, let $a'$ (resp., $b'$) be such
that \(a=a\wedge c\oplus a'\) and \(b=b\wedge c\oplus b'\).
Then $K$ satisfies \(a'\oplus c=b'\oplus c\), thus
\(a'\sim b'\). Put
\(\gamma=\Dim_L(p,a')=\Dim_L(p,b')\). Then
\[
\Dim_L(a\wedge c,a)=\gamma=\Dim_L(b\wedge c,b),
\]
thus
\[
\Dim_L(p,a\wedge c)+\gamma=\Dim_L(p,a)=\Dim_L(p,b)
=\Dim_L(p,b\wedge c)+\gamma.
\]
By (\ref{Eq:HalfCan}), it follows that there exists
\(\delta\in\DD L\) such that
\begin{gather}
\Dim_L(p,a\wedge c)+\delta=\Dim_L(p,b\wedge c)+\delta,
\label{Eq:acbcDel}\\
2\delta\leq\gamma.\label{Eq:DelSmall}
\end{gather} Then, by using (\ref{Eq:DelSmall}) and
Proposition~\ref{P:splitV} (applied to $L$) and then by taking
sectional complements in $K$, one obtains a decomposition
\(a=a^*\oplus u_0\oplus u_1\) in $K$ such that
\(\Dim_L(p,u_0)=\Dim_L(p,u_1)=\delta\).

Now let $T$ be the perspective\index{perspective!isomorphism}
isomorphism from \([p,\,a']\) onto \([p,\,b']\) with axis~$c$.
Put
\(b^*=T(a^*)\) and
\(v_i=T(u_i)\) (for all \(i<2\)). Then
\(b'=b^*\oplus v_0\oplus v_1\) and
\(\Dim_L(p,v_0)=\Dim_L(p,v_1)=\delta\). Since
\(\vv<u_0,u_1,v_0,v_1>\) is independent\index{independent}
in $K$, it follows from Corollary~\ref{C:IncrCong} that $u_0$
is a finite join of elements each of which is the first entry
of a homogeneous\index{homogeneous!sequence} sequence of length
$4$, that is, \(u_0\in 4K\). Put
\(u=u_0\), \(a''=a^*\oplus u_1\) and, similarly,
\(v=v_0\), \(b''=b^*\oplus v_1\). Then we obtain that
\begin{gather*}
a'=a''\oplus u,\quad b'=b''\oplus v,\\
b'=T(a'),\quad b''=T(a''),\quad\text{and}\quad v=T(u),\\
\Dim_L(p,u)=\delta,\quad\text{and}\quad u\in 4K.
\end{gather*}
It follows from (\ref{Eq:acbcDel}) that we have
\begin{multline}\label{Eq:pabcuv}
\Dim_L(p,a\wedge c\oplus u)=\Dim_L(p,a\wedge c)+\delta\\
=\Dim_L(p,b\wedge c)+\delta=\Dim_L(p,b\wedge c\oplus v).
\end{multline}
Since \(a\wedge c\), \(b\wedge c\), $u$, and
$v$ belong to $4K$, it follows from
Proposition~\ref{P:DiammL} that
\[
a\wedge c\oplus u\equiv b\wedge c\oplus v\pmod{\diam{4}{K}},
\]
thus \emph{a fortiori}
\[
a\wedge c\oplus u\equiv b\wedge c\oplus v\pmod{\diam{4}{L}}.
\]
Since \(a\wedge c\oplus u\) and \(b\wedge c\oplus v\) are
independent in $K$ (the first element is \(\leq a'\) while the
second is \(\leq b'\)), it follows from
Corollary~\ref{C:Loverdiam}.(ii) and from (\ref{Eq:pabcuv})
that \(a\wedge c\oplus u\sim b\wedge c\oplus v\). Finally,
since
\(a''\sim T(a'')=b''\) and $a$ and $b$ are independent in
$K$, it follows from Lemma~\ref{L:BasicAddPersp} that
\(a\sim b\), thus completing the proof of
normality\index{normal!lattice} of $L$.
\smallskip

Now suppose that $L$ has a zero element and let $a$ and $b$ be
elements of $L$ such that
\(a\approxeq b\). Put \(c=a\wedge b\) and let $a'$
and $b'$ be such that \(a=c\oplus a'\) and \(b=c\oplus b'\). Put
\(\gamma=\Dim(c)\). Then we have
\(a'\oplus c\approxeq b'\oplus c\), thus
\(\Dim(a')+\gamma=\Dim(b')+\gamma\). By (\ref{Eq:HalfCan}),
there exists \(\delta\in\DD L\)
such that
\(2\delta\leq\gamma\) and \(\Dim(a')+\delta=\Dim(b')+\delta\).
Since \(2\delta\leq\Dim(c)\), there exists, by
Corollary~\ref{C:dimVmeas}, a decomposition
\begin{equation}\label{Eq:decompc}
c=d_0\oplus d_1\oplus e,
\end{equation}
where
\(\Dim(d_0)=\Dim(d_1)=\delta\). Define the \(a_i\), \(b_i\)
(for \(i<2\)) by
\begin{align}
a_0&=a'\oplus d_0&\text{and}&&a_1&=d_1\oplus e,
\label{Eq:Defa0a1}\\
b_0&=b'\oplus d_1&\text{and}&&b_1&=d_0\oplus e.
\label{Eq:Defb0b1}
\end{align}
Then \(a_0\geq a'\), \(b_0\geq b'\), and
$a_1$, $b_1\leq c$. In addition, since \(\vv<a',b',c>\) is
independent\index{independent}, so is
\(\vv<a',b',d_0,d_1>\), thus \(a_0\perp b_0\). But
\[
\Dim(a_0)=\Dim(a')+\delta=\Dim(b')+\delta=\Dim(b_0),
\]
thus, since $L$ is normal\index{normal!lattice},
\(a_0\sim b_0\).
Moreover,
\(d_0\perp d_1\) and \(\Dim(d_0)=\Dim(d_1)=\delta\), thus,
again by normality\index{normal!lattice} of $L$,
\(d_0\sim d_1\). Therefore, since
\(\vv<d_0,d_1,e>\) is independent\index{independent} and by
Corollary~\ref{C:BasicAddPersp}, \(a_1\sim b_1\). Finally, it
follows from (\ref{Eq:decompc}), (\ref{Eq:Defa0a1}), and
(\ref{Eq:Defb0b1}) that \(a=a_0\oplus a_1\) and
\(b=b_0\oplus b_1\).
\end{proof}

We refer the reader to Appendix~\ref{App:Norm} for a summary
of known classes of normal\index{normal!lattice} relatively
complemented modular lattices.

\chapter[Dimension for $\aleph_0$-continuous lattices]%
{Dimension monoids
of $\aleph_0$-meet- and $\aleph_0$-join-con\-tin\-u\-ous%
\index{lattice!meet-continuous ---!$\aleph_0$-\mcont\ ---}%
\index{lattice!join-continuous ---!$\aleph_0$-\jcont\ ---}
lat\-tices}
\label{CtbleMeetCo}

Once the normality result of Theorem~\ref{T:MeetContNor} is
obtained, the work of figuring out the dimension theory of
countably meet-continuous complemented modular lattices becomes
much easier. As we shall see in Section~\ref{S:UniqNorm}, on
any such lattice, there exists, in fact, a \emph{unique} normal
equivalence. Of course, for coordinatizable lattices, this normal
equivalence is just isomorphy of submodules.

\section{Uniqueness of the normal\index{normal!equivalence}
equivalence}\label{S:UniqNorm}

The main technical result about
$\aleph_0$-\mcont%
\index{lattice!meet-continuous ---!$\aleph_0$-\mcont\ ---}
sectionally complemented modular lattices is that they admit
\emph{exactly one} normal\index{normal!equivalence}
equivalence. To establish this result, we start out by the
following lemma, which is proved
exactly the same way as in \cite[Lemma 1.3]{Good82}:

\begin{lemma}\label{L:lattcan}
Let $L$ be a sectionally
complemented modular lattice and let $\equiv$ be an additive
and refining equivalence relation on $L$ containing
$\sim$, with associated
preordering $\preceq$ (that is, \(a\preceq b\) if and only
if there exists
\(x\leq b\) such that \(a\equiv x\)). Then for all elements
$a$, $b$, and $c$ of $L$ such that
$a\perp b$ and \(a\oplus c\equiv b\oplus c\), there are
decompositions
\begin{align*}
a=a'\oplus a'',\quad b=b'\oplus b'',\quad&\text{ and }\quad
c=c'\oplus c''\\
\intertext{such that}
a'\equiv b',\quad a''\oplus c''\equiv b''\oplus c'',\quad
&\text{ and }\quad a''\oplus b''\preceq c'.\tag*{\qed}
\end{align*}
\end{lemma}

\begin{theorem}\label{T:normeqform}
Let $L$ be a sectionally
complemented, modular, conditionally $\aleph_0$-\mcont%
\index{lattice!meet-continuous ---!$\aleph_0$-\mcont\ ---}
lattice. Then there exists exactly one
normal\index{normal!equivalence} equivalence
$\equiv$ on $L$, defined by the following rule: \(a\equiv b\)
if and only if there are \(a_i\) and \(b_i\) (for \(i<2\)) such
that \(a=a_0\oplus a_1\), \(b=b_0\oplus b_1\), and
\(a_i\sim b_i\) (for all \(i<2\)).
\end{theorem}

\begin{proof}
The existence of a
normal\index{normal!equivalence} equivalence follows from
Theorem~\ref{T:MeetContNor}. Conversely, let $\equiv$ be a
normal\index{normal!equivalence} equivalence on $L$. We start
with the argument of
\cite[Theorem 1.4]{Good91}. Let $a$ and
$b$ be elements of $L$ such that \(a\equiv b\). Put
\(c_0=a\wedge b\) and let \(a_0\) and \(b_0\) be elements
of $L$ such that
\(a_0\oplus c_0=a\) and \(b_0\oplus c_0=b\). Note that
\(\vv<a_0,b_0,c_0>\) is independent\index{independent}.
Using Lemma~\ref{L:lattcan} inductively, we obtain sequences
\(\vv<a_n\mid n\in\omega>\),
\(\vv<a'_n\mid n\in\omega>\), \(\vv<b_n\mid
n\in\omega>\),
\(\vv<b'_n\mid n\in\omega>\), \(\vv<c_n\mid
n\in\omega>\),
\(\vv<c'_n\mid n\in\omega>\) of elements of $L$ satisfying the
property (for \(\preceq\) is the natural preordering
associated with \(\equiv\)) that the following relations
\begin{gather}
a_n=a'_n\oplus a_{n+1},\qquad b_n=b'_n\oplus
b_{n+1},\qquad c_n=c'_n\oplus c_{n+1},\label{Eq:Indabc1}\\
a'_n\equiv b'_n,\qquad a_n\oplus
c_n\equiv b_n\oplus c_n,\qquad a_{n+1}\oplus b_{n+1}
\preceq c'_n
\label{Eq:Indabc2}
\end{gather}
hold for all \(n\in\omega\).

Put \(a'=\oplus_{n\in\omega}a'_n\) and
\(b'=\oplus_{n\in\omega}b'_n\).
For all \(n\in\omega\), the  relations
\(a'_n\equiv b'_n\), \(a'_n\leq a_0\), \(b'_n\leq b_0\), and
\(a_0\perp b_0\), hold, so that \(a'_n\perp b'_n\); since $\equiv$
is normal\index{normal!equivalence}, we obtain that
\(a'_n\sim b'_n\). By
Lemma~\ref{L:BasicAddPersp}, it follows from this that
 \begin{equation}\label{Eq:a'b'}
 a'\sim b'
 \end{equation}
holds. Let $a''$, $b''$ be such that
\(a_0=a'\oplus a''\) and \(b_0=b'\oplus b''\). The equality
\[
(\oplus_{i\leq n}a'_i)\oplus a_{n+1}=a_0=a'\oplus a''
\geq(\oplus_{i\leq n}a'_i)\oplus a''
\]
holds for all \(n\in\omega\), so that
\(a''\lesssim a_{n+1}\), thus \emph{a fortiori}
\(a''\preceq a_{n+1}\). Similarly, \(b''\preceq b_{n+1}\),
whence, by (\ref{Eq:Indabc2}),
\[
a''\oplus b''\preceq a_{n+1}\oplus b_{n+1}\preceq c'_n.
\]
Thus there exists \(w_n\) such that
\begin{equation}\label{Eq:Defwn}
w_n\equiv a''\oplus b''\ \text{ and }\
w_n\leq c'_n.
\end{equation}
Since \(\equiv\) is refining, there are \(u_n\)
and \(v_n\) such that
\begin{gather}
u_n\equiv a''\ \text{ and }\ v_n\equiv b'',
\label{Eq:Inequn,vn}\\
w_n=u_n\oplus v_n.\label{Eq:Defun,vn}
\end{gather}
Moreover, by (\ref{Eq:Indabc1}) and
(\ref{Eq:Defwn}), there exists $w$ such that
 \begin{equation}\label{Eq:Defw}
 c_0=w\oplus(\oplus_nw_n).
 \end{equation}
Since \(a_0\perp c_0\), the sequence
\(\vv<a_0,c'_0,c'_1,\ldots>\) is
independent\index{independent}, thus, by (\ref{Eq:Defwn}) and
(\ref{Eq:Defun,vn}), so is
\emph{a fortiori} \(\vv<a'',u_0,u_1,\ldots>\); since the
entries of the latter are pairwise $\equiv$-equivalent and
$\equiv$ is normal\index{normal!equivalence}, its entries are
pairwise perspective\index{perspectivity}. Similarly,
\(b''\), \(v_0\), \(v_1\), \emph{etc.}, are mutually
perspective\index{perspectivity}. Now put
\begin{align}
r_1&=a''\oplus(\oplus_nu_{2n}),&r_2&=\oplus_nu_{2n+1}
\label{Eq:Defr1r2}\\
r_3&=\oplus_nv_{2n},&r_4&=\oplus_nv_{2n+1}.\label{Eq:Defr3r4}
\end{align} Thus, by (\ref{Eq:Defun,vn}), (\ref{Eq:Defw}),
(\ref{Eq:Defr1r2}), and (\ref{Eq:Defr3r4}),
\begin{equation}\label{Eq:aExp}
a=a'\oplus a''\oplus c_0
=a'\oplus r_1\oplus r_2\oplus r_3\oplus r_4\oplus w.
\end{equation}
Similarly, put
\begin{align}
s_1&=\oplus_nu_{2n+1},&s_2&=\oplus_nu_{2n}
\label{Eq:Defs1s2}\\
s_3&=\oplus_nv_{2n+1},&s_4&=b''\oplus(\oplus_nv_{2n}).
\label{Eq:Defs3s4}
\end{align}
Then, by (\ref{Eq:Defun,vn}), (\ref{Eq:Defw}),
(\ref{Eq:Defs1s2}), and (\ref{Eq:Defs3s4}),
 \begin{equation}\label{Eq:bExp}
 b=b'\oplus b''\oplus c_0
 =b'\oplus s_1\oplus s_2\oplus s_3\oplus s_4\oplus w.
 \end{equation}
Furthermore, by using Lemma~\ref{L:BasicAddPersp}, we
can see that the relation
 \begin{equation}\label{Eq:risi}
 r_i\sim s_i
 \end{equation}
holds for all \(i\in\{1,2,3,4\}\).
however, this is still not as strong as the desired conclusion.

Now, put
 \begin{align*}
 x_0&=a'\oplus r_1\oplus r_4,&y_0&=b'\oplus
 s_1\oplus s_4\\
 x_1&=r_2\oplus r_3\oplus w,&y_1&=s_2\oplus s_3\oplus w.
 \end{align*}
We obviously have \(a=x_0\oplus x_1\) and
\(b=y_0\oplus y_1\). Now, by (\ref{Eq:Indabc1}),
(\ref{Eq:Defwn}), (\ref{Eq:Defun,vn}) and (\ref{Eq:Defw}),
\(\vv<a_0,b_0,\oplus_nu_n,\oplus_nv_n>\) is
independent\index{independent}, thus so is
\[
\vv<a',b',a'',b'',\oplus_nu_{2n},\oplus_nu_{2n+1},
\oplus_nv_{2n},\oplus_nv_{2n+1}>;
\]
thus \(\vv<a',b',r_1,r_4,s_1,s_4>\) is
independent\index{independent}. By (\ref{Eq:a'b'}),
(\ref{Eq:risi}), and Lemma~\ref{L:BasicAddPersp}, it follows that
 \begin{equation}
 a'\oplus r_1\oplus r_4\sim b'\oplus s_1\oplus s_4,
 \end{equation}
that is, \(x_0\sim y_0\). Similarly, \(\vv<r_2,r_3,s_2,s_3,w>\) is
independent\index{independent}, by
(\ref{Eq:Defw}), (\ref{Eq:Defr1r2}), (\ref{Eq:Defr3r4}),
(\ref{Eq:Defs1s2}), and (\ref{Eq:Defs3s4}),
 thus, by (\ref{Eq:risi}),
 \begin{equation}
 r_2\oplus r_3\oplus w\sim s_2\oplus s_3\oplus w,
 \end{equation}
that is, \(x_1\sim y_1\). The conclusion follows from
(\ref{Eq:aExp}) and (\ref{Eq:bExp}).
\end{proof}

By using Corollary~\ref{C:NormEqDual2} (and reducing the
problem to the bounded case), one obtains immediately the
following result:

\begin{corollary}\label{C:UniNormJoin}
Let $L$ be a
sectionally complemented, modular, conditionally
$\aleph_0$-\jcont%
\index{lattice!join-continuous ---!$\aleph_0$-\jcont\ ---}
lattice. Then there exists a unique normal equivalence%
\index{normal!equivalence} on $L$.\qed
\end{corollary}

This unique normal\index{normal!equivalence} equivalence on $L$
is again perspectivity by decomposition%
\index{perspectivity by decomposition}, see
the second remark following Theorem~\ref{T:DimComplGCA}.

\begin{corollary}\label{C:ctblepiiso}
Let $R$ be a
$\aleph_0$-right-con\-tin\-u\-ous%
\index{ring!$\aleph_0$-right-continuous ---}
regular\index{ring!(von Neumann) regular ---} ring. Then
\(\pi_R\) is an
isomorphism from \(\DD\mathcal{L}(R_R)\) onto
\(V(R)\). Moreover, for all
principal ideals\index{ideal!of a ring} $I$ and $J$ of $R$,
\(I\cong J\) if and only if there are decompositions
\(I=I_0\oplus I_1\) and \(J=J_0\oplus J_1\) such that
\(I_0\sim J_0\) and \(I_1\sim J_1\).
\end{corollary}

\begin{proof}
Since the relation of isomorphy is a normal equivalence on
\(\mathcal{L}(R_R)\) (see Lemma~\ref{L:SubmNorm}),
the second statement is an
immediate consequence of Theorem~\ref{T:normeqform}. It
trivially follows that \(I\cong J\)
if and only if
\(\Dim(I)=\Dim(J)\), for all elements $I$ and $J$ of
\(\mathcal{L}(R_R)\). The conclusion follows by
Lemma~\ref{L:onetoone}.
\end{proof}

Corollary~\ref{C:ctblepiiso} is strengthened in
Corollary~\ref{C:DirFinReg}.

\begin{corollary}\label{C:ctbleVemb} Let $L$ be a relatively
complemented modular lattice and let
$K$ be a convex sublattice of $L$. Suppose that $L$ is either
conditionally $\aleph_0$-\mcont%
\index{lattice!meet-continuous ---!$\aleph_0$-\mcont\ ---} or
conditionally $\aleph_0$-\jcont%
\index{lattice!join-continuous ---!$\aleph_0$-\jcont\ ---}.
Then the natural map \(f\colon\DD K\to\DD L\) is a
\Vemb\index{Vemb@\Vemb}.
\end{corollary}

\begin{proof}
We first settle the conditionally
$\aleph_0$-\mcont%
\index{lattice!meet-continuous ---!$\aleph_0$-\mcont\ ---}
case. As every lattice is the direct union of its closed
intervals, we see easily that it is sufficient to consider the
case of bounded $K$ and $L$.
Therefore, $K$ is complemented, modular,
$\aleph_0$-\mcont%
\index{lattice!meet-continuous ---!$\aleph_0$-\mcont\ ---}.
By Theorem~\ref{T:normeqform}, there exists at most one normal
equivalence on $K$. Furthermore, by assumption on $L$ and by
Theorem~\ref{T:normeqform}, $L$ is
normal\index{normal!lattice}. The conclusion follows by
Proposition~\ref{P:conv}.

In the conditionally $\aleph_0$-\jcont%
\index{lattice!join-continuous ---!$\aleph_0$-\jcont\ ---}
case, the result follows immediately from
Proposition~\ref{P:V(Lop)}.
\end{proof}

\begin{remark}
The number of two pieces used in the
decompositions \(a=a_0\oplus a_1\), \(b=b_0\oplus b_1\), of
Theorem~\ref{T:normeqform}, is
\emph{optimal}: indeed, in case $\sim$ is
transitive, \(\DD L\) is cancellative (see
\cite[pp. 91--92]{FMae58} or
Proposition~\ref{P:transperspchar}), but there are easy
examples of $\aleph_0$-\mcont%
\index{lattice!meet-continuous ---!$\aleph_0$-\mcont\ ---}
complemented modular lattices whose dimension monoid is not
cancellative (take the lattice of all vector subspaces of any
infinite-dimensional vector space). A key point about the
failure of cancellation is the fact that the \(r_i\)'s and the
\(s_i\)'s in the proof of Theorem~\ref{T:normeqform} may not be
zero. Thus, cancellativity of \(\DD L\) would be implied, for
example, by the fact that any infinite
independent\index{independent} sequence of pairwise
perspective\index{perspectivity} elements is the zero
sequence. This is the case, for example, in sectionally
complemented modular lattices that are $\aleph_0$-\mcont%
\index{lattice!meet-continuous ---!$\aleph_0$-\mcont\ ---}
as well as $\aleph_0$-\jcont%
\index{lattice!join-continuous ---!$\aleph_0$-\jcont\ ---}
(\cite[Theorem I.3.8]{Neum60} or \cite[Satz IV.2.1]{FMae58}).
\end{remark}

\begin{definition}\label{D:WeakFin}
A sectionally complemented modular lattice $L$ is
\emph{weakly finite}, if every independent\index{independent}
sequence of pairwise perspective\index{perspectivity} elements
of $L$ is trivial (that is, the zero sequence).
\end{definition}

The following proposition is a
generalization of Halperin's result of transitivity of
perspectivity\index{perspectivity} in $\aleph_0$-continuous
geometries%
\index{czzontj@$\aleph_0$-continuous geometry}.

\begin{proposition}\label{P:persptran}
Let $L$ be a sectionally complemented, modular,
$\aleph_0$-\mcont%
\index{lattice!meet-continuous ---!$\aleph_0$-\mcont\ ---}
lattice. If $L$ is weakly finite, then
perspectivity\index{perspectivity} is transitive in $L$.
\end{proposition}

\begin{proof}
It suffices to prove that \(a\approxeq b\) implies that
\(a\sim b\), for all elements $a$
and $b$ in $L$. Let us consider the
decompositions obtained in the proof of
Theorem~\ref{T:normeqform}. As mentioned above, our hypothesis
yields that \(a''\), \(b''\), the \(r_i\)'s, and the
\(s_i\)'s are all equal to zero, thus \(a_0=a'\) and
\(b_0=b'\), whence
\(a_0\sim b_0\). Since
\(\vv<a_0,b_0,c_0>\) is independent\index{independent}, it
follows from Lemma~\ref{L:BasicAddPersp} that
\(a_0\oplus c_0\sim b_0\oplus c_0\),
that is, \(a\sim b\).\end{proof}

Another use of Theorem~\ref{T:normeqform} is an easy
computation of the dimension monoid of any \emph{modular}
geometric\index{lattice!geometric ---} lattice. This could
have been carried out directly, but it would have required much
more computations.

We recall first some well-known results about geometric
lattices (see \cite{Grat} for more details). If $L$ is
a  geometric lattice, then the relation of perspectivity
\index{perspectivity}
on the set $P$ of points of $L$ is \emph{transitive}. Let $I$
denote the set of $\sim$-blocks of $P$.
For every $i\in I$, denote by $s_i$ the supremum
of all elements of $i$. Then $L$ is canonically
isomorphic to the direct product of all its closed intervals
$[0,\,s_i]$. For every element $x$ of $L$, denote by
$\dim(x)$ the cardinality of any maximal independent
\index{independent}
set of points of $[0,\,s_i]$. For every cardinal number
$\kappa$, denote by $M_\kappa$ the commutative monoid
defined by
\begin{equation*}
M_\kappa=
\begin{cases}
\{0\},&\text{if $\kappa=0$};\\
\ZZ^+,&\text{if $0<\kappa<\aleph_0$};\\
\ZZ^+\cup\{\aleph_\xi\mid\xi\leq\alpha\},&
\text{if $\kappa=\aleph_\alpha$}.
\end{cases}
\end{equation*}
In particular, note that $\kappa$ is an element
(and even an order-unit) of $M_\kappa$.

Then define a submonoid $D(L)$ of
\(\prod_{i\in I}M_{\dim(s_i)}\) as follows:
\[
D(L)=\left\{
\vv<x_i\mid i\in I>\in\prod_{i\in I}M_{\dim(s_i)}\mid
(\exists n\in\NN)(\forall i\in I)
(x_i\leq n\cdot\dim(s_i))
\right\}.
\]

Define also a
function $D$ from $L$ to \(D(L)\), by the
following rule:
\[
D(x)=\vv<\dim(x\wedge s_i)\mid i\in I>.
\]

Then the structure of the dimension monoid of a geometric
modular lattice is summarized by the following result:

\begin{proposition}\label{P:DimGeom}
Let $L$ be a geometric\index{lattice!geometric ---} modular
lattice. Then there exists a unique monoid isomorphism
\(\varphi\colon\DD L\to D(L)\) such that
\(\varphi(\Dim(x))=D(x)\) holds for all \(x\in L\).
\end{proposition}

\begin{proof}
It is obvious that the map $D$ satisfies the axioms
(D$'$0), (D$'$1) and (D$'$2),
thus there exists a unique monoid homomorphism
\(\varphi\colon\DD L\to D(L)\) such that
\(\varphi(\Dim(x))=D(x)\) holds for all \(x\in L\).

In order to prove that $\varphi$ is a \Vhom%
\index{Vhom@\Vhom}, it
suffices to prove that if $K$ is an indecomposable geometric
modular lattice, \(c\in K\), and $\alpha$,
\(\beta\in M_{\dim(K)}\) such that \(\dim(c)=\alpha+\beta\),
then there exists a decomposition \(c=a\oplus b\) in $K$ such
that \(\dim(a)=\alpha\) and \(\dim(b)=\beta\). Let $A$ and
$B$ be disjoint sets of respective cardinalities $\alpha$
and $\beta$ and let \(\vv<c_\xi\mid\xi\in A\cup B>\)
be a maximal independent
\index{independent} subset of points of $c$. Then
\(a=\oplus_{\xi\in A}c_\xi\) and \(b=\oplus_{\xi\in B}c_\xi\)
satisfy the required conditions.

It follows easily that $\varphi$ is surjective.

Conversely, let $\equiv$ denote the equivalence relation on
$L$ defined by the following rule:
\[
x\equiv y\text{ if and only if }D(x)=D(y).
\]
We claim that $\equiv$ is normal%
\index{normal!equivalence}. All individual items of the
definition of normality are obvious, except perhaps (E3).
Now, to verify (E3), it suffices to prove that if $K$ is an
indecomposable modular geometric
\index{lattice!geometric ---}
lattice, if $x$ and $y$ are elements of $K$ such that
\(\dim(x)=\dim(y)\), and
$x\perp y$, then $x\sim y$. If $K$ is either trivial or a
projective line, then this is trivial. If $K$ is Arguesian
and has at least three independent
\index{independent} points, then there exists
a vector space $V$ (over some skew field) such that
$K$ is isomorphic to the lattice of all vector subspaces of
$V$ and the proof of our claim is then very easy and similar
to the proof of Lemma~\ref{L:SubmNorm}. If $K$ is
non-Arguesian, then $K$ is necessarily a projective plane, so
that either \(x=y=0\), in which case our statement is
trivial, or both $x$ and $y$ are points, in which case, since
$K$ is indecomposable, $x$ and $y$ are perspective. So, in
every case, our claim is verified.

It follows from Theorem~\ref{T:normeqform} that $\equiv$ is
projectivity by decomposition%
\index{projectivity!by decomposition}. In particular,
\(D(x)=D(y)\) implies that \(\Dim(x)=\Dim(y)\), that is, the
restriction of $\varphi$ to the dimension range of $L$ is
one-to-one. It is obvious that $D(L)$ is a refinement monoid%
\index{monoid!refinement ---}.
Therefore, by Lemma~\ref{L:onetoone}, $\varphi$ is also
one-to-one. Hence, $\varphi$ is an isomorphism.
\end{proof}

\section[Dimension $\aleph_0$-meet-continuous]%
{Dimension
monoids of $\aleph_0$-meet-con\-tin\-u\-ous%
\index{lattice!meet-continuous ---!$\aleph_0$-\mcont\ ---}
lat\-tices}

The main result of this section, Proposition~\ref{P:V(L)GCA},
states that if $L$ is any conditionally
$\aleph_0$-\mcont, sectionally complemented modular lattice, then
$\DD L$ is a GCA and the dimension function is countably additive.
To prepare for the proof of this result, we first prove
Proposition~\ref{P:adddim}.

\begin{proposition}\label{P:adddim}
Let $L$ be a sectionally complemented modular
$\aleph_0$-\mcont%
\index{lattice!meet-continuous ---!$\aleph_0$-\mcont\ ---}
lattice. Then for any independent\index{independent}
sequences \(\vv<a_n\mid n\in\omega>\) and
\(\vv<b_n\mid n\in\omega>\) of elements of $L$ such that
\((\forall n\in\omega)(\Dim(a_n)=\Dim(b_n))\), we have
\(\Dim(\oplus_na_n)=\Dim(\oplus_nb_n)\).
\end{proposition}

\begin{proof} Let \(\approxeq_\omega\)%
\index{PzzrojDecInf@$\approxeq_\omega$|ii} be the binary
relation defined on $L$ by \(a\approxeq_\omega b\) if and
only if there are decompositions \(a=\oplus_{n\in\omega}a_n\),
\(b=\oplus_{n\in\omega}b_n\) such that
\(a_n\approx b_n\) holds for all \(n\in\omega\).
Then it is clear that
\(\approxeq_\omega\) is additive (it is even
countably additive)%
\index{additive (binary relation)!countably ---}. Furthermore,
by using Lemma~\ref{L:ctbleref}, one can see easily that
\(\approxeq_\omega\) is refining (and even
countably refining) and transitive. It trivially contains
$\sim$. Finally, if
\(a\approxeq_\omega b\) and \(a\perp b\), then,
if both decompositions
\(a=\oplus_{n\in\omega}a_n\) and
\(b=\oplus_{n\in\omega}b_n\) witness
\(a\approxeq_\omega b\), the relations
\(a_n\perp b_n\) and
\(a_n\approx b_n\) hold for all $n$, thus, since $L$ is
normal\index{normal!lattice}, \(a_n\sim b_n\). Therefore, by
Lemma~\ref{L:BasicAddPersp}, \(a\sim b\). This proves that
\(\approxeq_\omega\) is a
normal\index{normal!equivalence} equivalence on $L$. Therefore,
by Theorem~\ref{T:normeqform},
\(a\approxeq_\omega b\) if and only if
\(\Dim(a)=\Dim(b)\). The conclusion follows
immediately.\end{proof}

\begin{note} In the case of $\aleph_0$-right-con\-tin\-u\-ous%
\index{ring!$\aleph_0$-right-continuous ---}
regular\index{ring!(von Neumann) regular ---} rings, this
result has been proved first by Goodearl in the case of
directly%
\index{directly finite!element in a monoid} finite
regular\index{ring!(von Neumann) regular ---} rings in
\cite[Theorem 2.2]{Good82}, then in the general case by Ara
in \cite[Theorem 2.12]{Ara87}. The proof presented here is
purely lattice-the\-o\-ret\-i\-cal.
\end{note}

The remaining part of this section will be devoted to the
investigation of the structure of \(\DD L\) for $L$
sectionally complemented, modular,
$\aleph_0$-\mcont%
\index{lattice!meet-continuous ---!$\aleph_0$-\mcont\ ---}. We
shall use the results of
Section~\ref{S:GCAs} about generalized cardinal%
\index{generalized cardinal algebra (GCA)} algebras (GCA's).

\begin{proposition}\label{P:V(L)GCA}
Let $L$ be a sectionally complemented modular
$\aleph_0$-\mcont%
\index{lattice!meet-continuous ---!$\aleph_0$-\mcont\ ---}
lattice. Then
\(\DD L\) is a GCA%
\index{generalized cardinal algebra (GCA)}. Moreover, if
\(\vv<a_n\mid n\in\omega>\) is any increasing sequence of
elements of $L$ and
\(a=\bigvee_{n\in\omega}a_n\), then
\(\Dim(a)=\bigvee_{n\in\omega}\Dim(a_n)\).
\end{proposition}

\begin{proof} Let $U$ be the dimension range of $L$. By
Proposition~\ref{P:adddim}, one can define a partial infinite
addition on $U$ by putting \(\alpha=\sum_n\alpha_n\) if and
only there are elements $a$ and \(a_n\) (for \(n\in\omega\)) in
$L$ such that
\(a=\oplus_na_n\), \(\Dim(a)=\alpha\), and
\(\Dim(a_n)=\alpha_n\), for all $n$. Thus,
\(\Dim(a)=\sum_n\Dim(a_n)\) if and only if there are
\(x_n\in L\) (for all \(n\in\omega\)) such that
\(a=\oplus_nx_n\) and \(\Dim(x_n)=\Dim(a_n)\) for all $n$,
for all \(a\in L\) and
\(\vv<a_n\mid n\in\omega>\in L^\omega\).
Then one proves that $U$ equipped with this infinite
addition (and the restriction of the addition of
\(\DD L\)) is indeed a GCA%
\index{generalized cardinal algebra (GCA)}. It is easy to
verify directly (GCA1,2,3). The postulate (GCA4) follows
immediately from Lemma~\ref{L:ctbleref}. Let us finally check
the remainder postulate. Thus let
\(\vv<\alpha_n\mid n\in\omega>\) and
\(\vv<\beta_n\mid n\in\omega>\) be sequences of
elements of $U$ such that
\(\alpha_n=\beta_n+\alpha_{n+1}\) holds for all $n$. Then it is
easy to construct inductively sequences
\(\vv<a_n\mid n\in\omega>\) and
\(\vv<b_n\mid n\in\omega>\) of elements of $L$ such that
\(\alpha_n=\Dim(a_n)\), \(\beta_n=\Dim(b_n)\), and
\(a_n=b_n\oplus a_{n+1}\), for all $n$. Thus,
\(\oplus_{k<n}b_k\leq a_0\) holds for all $n$, so that
\(\oplus_{k\in\omega}b_k\leq a_0\); whence there exists
\(a\in L\) such that \(a\oplus(\oplus_{k\in\omega}b_k)=a_0\). It
follows that
\begin{align*}
a\oplus(\oplus_{k\geq n}b_k)\oplus(\oplus_{k<n}b_k)&=a_0\\
&=a_n\oplus(\oplus_{k<n}b_k)
\end{align*}
holds for all \(n\in\omega\).

so that \(a\oplus(\oplus_{k\geq n}b_k)\sim a_n\), whence
\begin{align*}
\Dim(a)+\sum_k\beta_{n+k}&=\Dim(a)+\sum_k\Dim(b_{n+k})\\
&=\Dim(a)+\Dim(\oplus_kb_{n+k})\\
&=\Dim(a\oplus(\oplus_kb_{n+k}))\\
&=\Dim(a_n)\\
&=\alpha_n,
\end{align*}
which completes the proof of the remainder postulate in \(U\).

Furthermore, by using the fact that \(\DD L\) is a
refinement\index{monoid!refinement ---} monoid generated (as a
monoid) by its lower subset
$U$, one can extend the infinite addition of \(U\) to
\(\DD L\) by defining
\(\sum_n\alpha_n=\alpha\) if and only if there are
\(m\in\omega\) and decompositions
\(\alpha=\sum_{i<m}\alpha^i\),
\(\alpha_n=\sum_{i<m}\alpha^i_n\) (for all \(n\in\omega\)) with all
the elements \(\alpha^i\), \(\alpha^i_n\) in $U$, and
\(\alpha^i=\sum_{n\in\omega}\alpha^i_n\) (for all \(i<m\)) in
\(U\). The verifications are tedious but quite
straightforward. Finally, the last continuity statement
follows immediately from Proposition~\ref{P:adddim}
and Proposition~\ref{P:basicGCA1}.(i).
\end{proof}

\section{Further extensions}

In this section, we shall extend the previous results to
lattices that are not necessarily
$\aleph_0$-\mcont%
\index{lattice!meet-continuous ---!$\aleph_0$-\mcont\ ---}.

\begin{theorem}\label{T:DimComplGCA}
Let $L$ be a relatively
complemented modular lattice. If $L$ is either conditionally
$\aleph_0$-\mcont%
\index{lattice!meet-continuous ---!$\aleph_0$-\mcont\ ---} or
conditionally
$\aleph_0$-\jcont%
\index{lattice!join-continuous ---!$\aleph_0$-\jcont\ ---},
then
\(\DD L\) is a GCA%
\index{generalized cardinal algebra (GCA)}.
\end{theorem}

\begin{proof}
Let us first see the case where $L$ is conditionally
$\aleph_0$-\mcont%
\index{lattice!meet-continuous ---!$\aleph_0$-\mcont\ ---}.
Express $L$ as the direct union of its closed intervals, viewed
as ($\aleph_0$-\mcont%
\index{lattice!meet-continuous ---!$\aleph_0$-\mcont\ ---},
complemented, modular) lattices. By
Proposition~\ref{P:V(L)GCA}, the dimension monoids of each of
these is a GCA%
\index{generalized cardinal algebra (GCA)}. Furthermore, by
Corollary~\ref{C:ctbleVemb}, if
\(K_0\) and \(K_1\) are closed intervals of $L$ such that
\(K_0\subseteq K_1\), then the natural map from \(\DD K_0\) to
\(\DD K_1\) is a
\Vemb\index{Vemb@\Vemb}. Therefore, by
Proposition~\ref{P:Vfunctor},
\(\DD L\) is a GCA%
\index{generalized cardinal algebra (GCA)}.

The conclusion in the conditionally
$\aleph_0$-\jcont%
\index{lattice!join-continuous ---!$\aleph_0$-\jcont\ ---}
case follows immediately from Proposition~\ref{P:V(Lop)}.
\end{proof}

\begin{remark}
Even in easy cases, \(\DD L\) may not be a CA%
\index{cardinal algebra (CA)}, because infinite sums may not
all be defined. In fact, an infinite sum
\(\sum_k\alpha_k\) is defined if and only if the set of all
partial sums \(\sum_{k<n}\alpha_k\) (for \(n\in\omega\)) is
bounded.
\end{remark}

\begin{remark}
Let $L$ be a lattice satisfying the hypothesis of
Theorem~\ref{T:DimComplGCA}. Then $\DD L$ is a GCA, thus, by
Proposition~\ref{P:basicCA}.(a) and Theorem~\ref{T:AxImplNorm},
there exists a unique normal equivalence $\equiv$ on $L$, defined
by $a\equiv b$ if and only if there are decompositions
$a=a_0\oplus a_1$, $b=b_0\oplus b_1$, such that $a_i\sim b_i$, for
all $i<2$.
\end{remark}

\begin{corollary}\label{C:GenEmbVs}
Let $K$ be a quotient
lattice of any relatively complemented modular lattice that is
either conditionally $\aleph_0$-\mcont%
\index{lattice!meet-continuous ---!$\aleph_0$-\mcont\ ---} or
conditionally $\aleph_0$-\jcont%
\index{lattice!join-continuous ---!$\aleph_0$-\jcont\ ---}.
Then \(\DD K\) satisfies
{\rm \ref{P:basicCA}.(a)}, thus {\rm (\ref{Eq:HalfCan})} and
$K$ is normal\index{normal!lattice}.
\end{corollary}

\begin{proof}
Put \(K=L/\theta\), where \(\theta\) is a
congruence\index{congruence!lattice ---} of $L$ and
$L$ is either conditionally
$\aleph_0$-\mcont%
\index{lattice!meet-continuous ---!$\aleph_0$-\mcont\ ---} or
conditionally
$\aleph_0$-\jcont%
\index{lattice!join-continuous ---!$\aleph_0$-\jcont\ ---}.
Then $K$ is the direct limit of all \(M/\theta_M\), where
$M$ ranges over all bounded intervals of $L$ and where we put
\(\theta_M=\theta\cap(M\times M)\); it follows that it is
sufficient to solve the problem in case where $L$ is
\emph{bounded}. However, in this case, it follows from
Theorem~\ref{T:DimComplGCA} that \(\DD L\) is a
GCA\index{generalized cardinal algebra (GCA)}; therefore, by
Proposition~\ref{P:basicCA},
\(\DD L\) satisfies
\ref{P:basicCA}.(a). Since \(\DD\theta\) as defined
before Proposition~\ref{P:congsonV} is an ideal\index{ideal!of
a monoid} of \(\DD L\) and, by
Proposition~\ref{P:latt<>mon}, \(\DD K\cong\DD L/\DD\theta\),
it follows easily that \(\DD K\) also satisfies
\ref{P:basicCA}.(a), thus, trivially, (\ref{Eq:HalfCan}).
Normality\index{normal!lattice} of $K$ follows then from
Theorem~\ref{T:AxImplNorm}.
\end{proof}

\begin{corollary}\label{C:CtbleDirFin} Let $L$ be a quotient
lattice of any relatively complemented modular lattice which is
either conditionally $\aleph_0$-\mcont%
\index{lattice!meet-continuous ---!$\aleph_0$-\mcont\ ---} or
conditionally $\aleph_0$-\jcont%
\index{lattice!meet-continuous ---!$\aleph_0$-\mcont\ ---}. If
all elements of the dimension range of $L$ are directly%
\index{directly finite!element in a monoid} finite, then
\(\DD L\) is cancellative.\qed
\end{corollary}

The following corollary extends both \cite[Theorem 2.7]{Ara87}
and Corollary~\ref{C:ctblepiiso}:

\begin{corollary}\label{C:DirFinReg} Let $R$ be a
regular\index{ring!(von Neumann) regular ---} ring which is
either $\aleph_0$-left-con\-tin\-u\-ous%
\index{ring!$\aleph_0$-left-continuous ---} or
$\aleph_0$-right-con\-tin\-u\-ous%
\index{ring!$\aleph_0$-right-continuous ---}, let $I$ be a
two-sided ideal\index{ideal!of a ring} of $R$; put
\(S=R/I\).
\begin{itemize}

\item[\rm(a)] The monoid \(V(S)\)
satisfies \ref{P:basicCA}.(a).

\item[\rm(b)] \(I\cong J\) if and only if there are decompositions
\(I=I_0\oplus I_1\) and \(J=J_0\oplus J_1\) such that
\(I_0\sim J_0\) and \(I_1\sim J_1\), for all principal
ideals\index{ideal!of a ring}
$I$ and $J$ of $S$.

\item[\rm(c)] The map \(\pi_S\) is an isomorphism; thus,
\(V(S)\cong\DD\mathcal{L}(S_S)\).

\item[\rm(d)] If $S$ is directly%
\index{directly finite!element in a monoid} finite, then it is
unit-regular\index{ring!unit-regular ---}.

\end{itemize}
\end{corollary}

\begin{proof} By \cite[Proposition 1.4]{AGPO}, \(V(I)\) is
an ideal\index{ideal!of a monoid} of
\(V(R)\) and \(V(S)\) is isomorphic to
\(V(R)/V(I)\), so that, to prove
(a), it is sufficient to prove that \(V(R)\) satisfies
\ref{P:basicCA}.(a). If $R$ is
$\aleph_0$-right-con\-tin\-u\-ous%
\index{ring!$\aleph_0$-right-continuous ---}, this follows
from Corollary~\ref{C:ctblepiiso}. Now, suppose that $R$ is
$\aleph_0$-left-con\-tin\-u\-ous%
\index{ring!$\aleph_0$-left-continuous ---}; put
\(L=\mathcal{L}(R_R)\). Then
$L$ is a $\aleph_0$-\jcont%
\index{lattice!join-continuous ---!$\aleph_0$-\jcont\ ---}
complemented modular lattice, thus, by
Corollary~\ref{C:GenEmbVs},
\(\DD L\) satisfies
\ref{P:basicCA}.(a). Then it follows from
Theorem~\ref{T:AxImplNorm} that any isomorphic principal
right ideals\index{ideal!of a ring} $I$ and $J$ of $R$ are
perspective by decomposition%
\index{perspectivity by decomposition}, thus \emph{a fortiori}
projective by decomposition%
\index{projectivity!by decomposition}; this means that the
restriction of
\(\pi_R\) to the
dimension range of $L$ is one-to-one, thus, by
Lemma~\ref{L:onetoone}, \(\pi_R\) is an isomorphism; in
particular,
\(V(R)\cong\DD L\) satisfies \ref{P:basicCA}.(a).
This proves (a).

Then (b) follows from Theorem~\ref{T:AxImplNorm} and (c)
follows from the argument used above to prove that
\(\pi_R\) is an
isomorphism. Finally, (d) follows immediately from (a).
\end{proof}

\begin{remark}

A word of warning about possible generalizations of all this to
$\kappa$-\mcont%
\index{lattice!meet-continuous ---!$\kappa$-\mcont\ ---} (or
even [totally] \mcont%
\index{lattice!meet-continuous ---}) lattices, where
$\kappa$ is any infinite cardinal. Even in the case where
\(\DD L\) is a complete lattice,
\emph{the dimension function may no longer be order-continuous}
(although, nevertheless, the
$\kappa$-analogue of Proposition~\ref{P:adddim} is easily seen
to hold, with a proof which is
\emph{mutatis mutandis} the same as for the countable case).
For example, let $K$ be any division ring and let
$E$ be a vector space of uncountable dimension over \(K\). Let
\(L=\PP(E)\)\index{PzzofE@$\PP(E)$|ii} be the lattice of all
vector subspaces of $E$. Then
\(\DD L\) is isomorphic to the
additive monoid of all cardinals (including the finite ones)
below the cardinality of
$I$, thus it is, in particular, a complete chain. Then
\(\Dim(X)\) is just the dimension of $X$, for all
\(X\in L\). Moreover, $L$ is a
complemented modular \mcont%
\index{lattice!meet-continuous ---} lattice. Nevertheless, if
\(\vv<X_\xi\mid \xi<\omega_1>\) is any strictly increasing
sequence of countably dimensional elements of $L$ and
\(X=\bigcup_{\xi<\omega_1}X_\xi\), then \(\Dim(X)=\aleph_1\),
while \(\Dim(X_\xi)=\aleph_0\) for all \(\xi<\omega_1\).

\end{remark}

\section{The case of $\vv<\kappa,\lambda>$-geometries}
\label{S:ContGeom}

We have seen in Proposition~\ref{P:persptran} that a
sectionally complemented, modular,
$\aleph_0$-\mcont%
\index{lattice!meet-continuous ---!$\aleph_0$-\mcont\ ---}
lattice $L$ has transitive
perspectivity\index{perspectivity}
(that is, \(\DD L\) is cancellative, see
Theorem~\ref{T:transperspchar}) if and only if $L$ is weakly
finite (Definition~\ref{D:WeakFin}). Accordingly, we
shall introduce the following definition, inspired by the
classical definition of a continuous%
\index{continuous!geometry} geometry:

\begin{definition}\label{D:kContGeom}
Let $\kappa$ and
$\lambda$ be (possibly finite) cardinal numbers. A
\emph{$\vv<\kappa,\lambda>$-geometry}%
\index{kzzlgeom@$\vv<\kappa,\lambda>$-geometry|ii} is a
relatively complemented modular lattice which is conditionally
$\kappa$-\mcont%
\index{lattice!meet-continuous ---!$\kappa$-\mcont\ ---}
as well as conditionally $\lambda$-\jcont%
\index{lattice!join-continuous ---!$\kappa$-\jcont\ ---} and
such that perspectivity\index{perspectivity} is transitive in
$L$ (\emph{cf.} Definition~\ref{D:MJContLatt}). This
is then generalized in the obvious way for case
$\kappa=\infty$ or for $\lambda=\infty$.
\end{definition}

In particular, by the results of \cite{Halp38}, every
continuous\index{continuous!geometry} geometry (or even every
relatively complemented modular lattice which is
conditionally $\aleph_0$-\mcont%
\index{lattice!meet-continuous ---!$\aleph_0$-\mcont\ ---}
as well as conditionally $\aleph_0$-\jcont%
\index{lattice!join-continuous ---!$\aleph_0$-\jcont\ ---}) is
a \(\vv<\aleph_0,\aleph_0>\)-geometry. A
\(\vv<\aleph_0,1>\)-geometry (resp., a
\(\vv<1,\aleph_0>\)-geometry) is a conditionally
$\aleph_0$-\mcont%
\index{lattice!meet-continuous ---!$\aleph_0$-\mcont\ ---}
(resp., conditionally $\aleph_0$-\jcont%
\index{lattice!meet-continuous ---!$\aleph_0$-\mcont\ ---})
relatively complemented modular lattice $L$ such that
perspectivity\index{perspectivity} is transitive in
$L$ (or, equivalently if $L$ has a zero, $L$ has no nontrivial
homogeneous\index{homogeneous!sequence} sequence). A
\(\vv<1,1>\)-geometry is just a relatively complemented modular
lattice with transitive perspectivity\index{perspectivity}.

\begin{lemma}\label{L:GenAdddim}
Let $\kappa$ be an infinite
cardinal, let $L$ be a sectionally complemented modular
$\kappa$-\mcont%
\index{lattice!meet-continuous ---!$\kappa$-\mcont\ ---}
lattice. Then for any independent%
\index{independent} sequences
\(\vv<a_\xi\mid \xi<\kappa>\) and
\(\vv<b_\xi\mid \xi<\kappa>\) of elements of $L$ such that
\(\Dim(a_\xi)=\Dim(b_\xi)\) (for all \(\xi<\kappa\)), we have
\(\Dim(\oplus_\xi a_\xi)=\Dim(\oplus_\xi b_\xi)\).
\end{lemma}

\begin{proof}
The proof is similar to the proof of Proposition
\ref{P:adddim} (the corresponding version of
Lemma~\ref{L:ctbleref} is proved in the same way).
\end{proof}

\begin{lemma}\label{L:CGtrim}
Let $L$ be a sectionally
complemented modular lattice with transitive
perspectivity\index{perspectivity}. Let \(a\leq b\) in $L$ and
\(\gamma\in\DD L\) be such that
\(\Dim(a)+\gamma\leq\Dim(b)\). Then there exists an element
\(c\leq b\) independent with $a$ such that \(\Dim(c)=\gamma\).
\end{lemma}

\begin{proof}
Let \(d\in L\) be such that \(a\oplus d=b\). Since,
by Proposition~\ref{P:transperspchar},
\(\DD L\) is cancellative, we
have
\(\gamma\leq\Dim(d)\), thus, by Corollary~\ref{C:dimVmeas},
there exists \(c\leq d\) such that
\(\Dim(c)=\gamma\).
\end{proof}

\begin{theorem}\label{T:K0(L)MsC}
Let $\kappa$ and $\lambda$ be nonzero cardinal numbers,
one of which is infinite, and let $L$ be a
\(\vv<\kappa,\lambda>\)-geometry%
\index{kzzlgeom@$\vv<\kappa,\lambda>$-geometry}. Then
\(K_0(L)\) is a monotone%
\index{monotone@monotone $\kappa$-completeness}
\(\max\{\kappa,\lambda\}\)-complete dimension group%
\index{group!dimension ---}. Furthermore, the dimension
function
$\Dim$ satisfies both of the following properties:

\begin{itemize}
\item[\rm(i)] $\Dim$ is $\kappa$-upper continuous%
\index{continuous!$\kappa$-upper- ---|ii},
that is, for every \(a\in L\) and every increasing
$\kappa$-sequence \(\vv<b_\xi\mid \xi<\kappa>\) of
elements of $L$, if \(b=\bigvee_{\xi<\kappa}b_\xi\), then
\begin{equation}\label{Eq:kcont}
\Dim(a,b)=\bigvee_{\xi<\kappa}\Dim(a,b_\xi)
\end{equation}
holds in $K_0(L)$.

\item[\rm(ii)] $\Dim$ is $\lambda$-lower continuous%
\index{continuous!$\lambda$-lower- ---|ii}, that is, for
every
\(b\in L\) and every decreasing
$\lambda$-sequence \(\vv<a_\eta\mid \eta<\lambda>\) of
elements of $L$, if \(a=\bigwedge_{\eta<\lambda}a_\eta\), then
\begin{equation}\label{Eq:lcont}
\Dim(a,b)=\bigvee_{\eta<\lambda}\Dim(a_\eta,b)
\end{equation}
holds in $K_0(L)$.
\end{itemize}
\end{theorem}

\noindent In particular, this generalizes the continuity
results of \cite{Halp39b}.

\begin{proof}
For any bounded closed intervals $I$ and $J$ of
$L$ such that \(I\subseteq J\), the natural map from
\(\DD I\) to \(\DD J\) is, by
Corollary~\ref{C:ctbleVemb}, a \Vemb\index{Vemb@\Vemb}. It
follows easily that it suffices to prove the theorem in
case $L$ is \emph{complemented} and \emph{modular}.

Let \(U=\{\Dim(x)\mid x\in L\}\) be the dimension range of
$L$. Note that
\(K_0(L)\) is an
interpolation\index{group!interpolation ---} group
(see Theorem~\ref{T:V(L)ref} and
Proposition~\ref{P:transperspchar}) and that $U$ generates
\(\DD L=K_0(L)^+\) as a monoid.

Next, we shall first prove that $U$ satisfies both (i) and
(ii) (we shall deduce monotone%
\index{monotone@monotone $\kappa$-completeness}
\(\max\{\kappa,\lambda\}\)-completeness later). Observe that
the dual lattice of $L$ is
$\lambda$-\mcont%
\index{lattice!meet-continuous ---!$\kappa$-\mcont\ ---},
thus, using the fact that
\(\DD L\) is cancellative, it
is easy to see that it suffices to prove (\ref{Eq:kcont}) for
\(a=0\). In order to do this, it is clearly sufficient to prove
the following claim:

\setcounter{claim}{0}
\begin{claim}
Let \(\theta\leq\kappa\) be an infinite cardinal.
\begin{itemize}
\item[\rm(i)] Let \(\vv<\alpha_\xi\mid \xi<\theta>\)
be an increasing $\theta$-sequence of elements of $U$. Then
there exists an increasing $\theta$-sequence
\(\vv<a_\xi\mid \xi<\theta>\) of elements of $L$ such that
\(\Dim(a_\xi)=\alpha_\xi\) holds for all \(\xi<\theta\).

\item[\rm(ii)] Let \(\vv<a_\xi\mid \xi<\theta>\) be
an increasing $\theta$-sequence of elements of $L$, let
\(a=\bigvee_{\xi<\theta}a_\xi\). Then
\begin{equation*}
\Dim(a)=\bigvee_{\xi<\theta}\Dim(a_\xi)
\end{equation*}
holds in $U$.
\end{itemize}
\end{claim}

\begin{cproof}
We argue by induction on $\theta$.
For (i), construct inductively the
\(a_\xi\)'s as follows. If all \(a_\eta\) (for \(\eta<\xi\)) have
been constructed in such a way that they form an increasing
$\xi$-sequence and \(\Dim(a_\eta)=\alpha_\eta\), for all
\(\eta<\xi\), put first \(b_\xi=\bigvee_{\eta<\xi}a_\eta\).
Then, by the induction hypothesis ((ii) at the ordinal
\(\mathrm{cf}(\xi)\)), we have
\(\Dim(b_\xi)=\bigvee_{\eta<\xi}\Dim(a_\eta)\) in $U$, so that in
particular, \(\Dim(b_\xi)\leq\alpha_\xi\). Then, it follows, by
Lemma~\ref{L:CGtrim}, that there exists
\(a_\xi\geq b_\xi\) such that \(\Dim(a_\xi)=\alpha_\xi\).
\smallskip

Now let us see (ii). It is easy to see that there exists an
independent\index{independent} $\theta$-sequence
\(\vv<a'_\xi\mid \xi<\theta>\) such that the equality
\(a_\xi=\oplus_{\eta\leq\xi}a'_\eta\) holds for all
\(\xi<\theta\). Let \(b\in L\) be such that
\(\Dim(a_\xi)\leq\Dim(b)\) holds for all \(\xi<\theta\).
We construct inductively
\(b'_\xi\), \(\xi<\theta\) as follows. Let
\(\xi<\theta\) and suppose having constructed an
independent\index{independent} sequence
\(\vv<b'_\eta\mid \eta<\xi>\) such that
\(\Dim(b'_\eta)=\Dim(a'_\eta)\) and \(b'_\eta\leq b\) hold for all
\(\eta<\xi\). Then, by
Lemma~\ref{L:GenAdddim},
\[
\Dim\left(\oplus_{\eta<\xi}b'_\eta\right)=
\Dim\left(\oplus_{\eta<\xi}a'_\eta\right)\leq
\Dim\left(\oplus_{\eta<\xi}a'_\eta\right)+\Dim(a'_\xi)=
\Dim(a_\xi)\leq\Dim(b),
\]
thus, by Lemma~\ref{L:CGtrim}, there exists \(b'_\xi\leq b\)
independent with \(\oplus_{\eta<\xi}b'_\eta\) such that
\(\Dim(b'_\xi)=\Dim(a'_\xi)\). This completes the inductive
construction of the \(b'_\xi\)'s. Then, it follows again by
Lemma~\ref{L:GenAdddim} that we have
\begin{equation}
\Dim(a)=
\Dim\left(\oplus_{\xi<\theta}a'_\xi\right)=
\Dim\left(\oplus_{\xi<\theta}b'_\xi\right)\leq
\Dim(b).\tag*{\qedc}
\end{equation}
\renewcommand{\qedc}{}
\end{cproof}

It follows immediately that $U$ is monotone%
\index{monotone@monotone $\kappa$-completeness}
$\kappa$-complete. Therefore, by Lemma~\ref{L:MsC,UtoG},
\(K_0(L)\) is itself monotone%
\index{monotone@monotone $\kappa$-completeness}
$\kappa$-complete.
\end{proof}

Since an interpolation\index{group!interpolation ---} group is
Dedekind-complete if and only if it is monotone%
\index{monotone@monotone $\kappa$-completeness}
$\kappa$-complete, for all $\kappa$, we obtain the following
result:

\begin{corollary}\label{C:K0(L)(cpletecase)}
Let $L$ be a relatively complemented modular lattice with
transitive perspectivity\index{perspectivity} that is either
\mcont%
\index{lattice!meet-continuous ---} or
\jcont%
\index{lattice!join-continuous ---}. Then
\(K_0(L)\) is a Dedekind-complete
lattice-ordered group\index{group!lattice-ordered ---}.\qed
\end{corollary}

It can be shown that the classical analysis of the center of
$L$ in the case where $L$ is a
continuous\index{continuous!geometry} geometry (that works in
fact if $L$ is merely a \mcont%
\index{lattice!meet-continuous ---} or \jcont%
\index{lattice!join-continuous ---} complemented modular
lattice with transitive perspectivity\index{perspectivity})
has a precise counterpart in the
world of lattice-ordered groups%
\index{group!lattice-ordered ---}. This can be outlined
as follows. We put, as usual, \(x\triangledown y\), if
\(\Theta(x)\cap\Theta(y)=\{0\}\), for all elements
$x$ and $y$ of $L$. For every subset $X$ of $L$,
put \(X^{\triangledown}=
\{s\in L\mid (\forall x\in X)(x\triangledown s)\}\) and say
that an ideal\index{ideal!of a semilattice} $I$ of $L$ is
\emph{closed}, if
\(I=I^{\triangledown\triangledown}\). Then the closed
ideals\index{ideal!of a semilattice} of
$L$ correspond exactly, \emph{via} the dimension map, to polar
subsets of the Dedekind-complete lattice-ordered group%
\index{group!lattice-ordered ---}
\(K_0(L)\). This leads to methods which
can be applied successfully to a complete analysis of the
dimension function on $L$. Furthermore, these methods can be
extended to arbitrary sectionally complemented, modular,
conditionally meet-continuous lattices. However, to show this in
detail would lead us too far away from the methods used in this
work.

\appendix

\chapter[Dimension for lattice embeddings]%
{A review on the $\DD$ functor and
lattice embeddings}\label{App:DimEmb}

Let $L$ be a lattice and let $K$ be a sublattice of $L$. Let
\(f\colon\DD K\to\DD L\) be the natural homomorphism. Then
the following assertions hold:

\begin{itemize}
\item[(1)] (see Corollary~\ref{C:VKidVL}) If $L$ is
\Vmod\ and $K$ is a convex sublattice of $L$, then the range
of the canonical map from \(\DD K\) to \(\DD L\) is an ideal
of \(\DD L\).

\item[(2)] (see Corollary~\ref{C:RectDEP}) If $K$ is BCF
and if the inclusion map from $K$ into
\(L=\prod_{i<n}L_i\) (for $n$ a positive integer) is a subdirect
decomposition of $K$, then
$f$ is an order-embedding of commutative preordered monoids.

\item[(3)] (see Corollary~\ref{C:V(I)inV(L)}) If $L$ is sectionally
complemented modular and $K$ is a neutral ideal of $L$, then
$f$ is a \Vemb. Furthermore, if we identify $\DD K$ with its
image in $\DD L$ under this isomorphism, then
\(\DD(L/K)\cong\DD L/\DD K\).

\item[(4)] (see Proposition~\ref{P:Vemb})
If $L$ is sectionally complemented modular, $K$
is an ideal of $L$ and perspectivity by decomposition is
transitive in $L$, then the natural map
$f\colon \DD K\to\DD L$ is a \Vemb.

\item[(5)] (see Proposition~\ref{P:conv}) Suppose that $L$ is
relatively complemented modular and that
$K$ is a convex sublattice of $L$ with a smallest element. If
the following additional assumptions hold:
\begin{itemize}
\item[\rm (i)] $L$ is normal;

\item[\rm (ii)] there exists at most one normal equivalence on
$K$;
\end{itemize}
then $f$ is a \Vemb.

\item[(6)] (see Lemma~\ref{L:nDistrVemb}) Suppose that $L$ is a
relatively complemented modular lattice and that $K$ is a
convex sublattice of $L$. If $L$ is locally finitely
distributive, then $f$ is a \Vemb.

\item[(7)] (see Corollary~\ref{C:ctbleVemb}) Suppose that $L$ is a
relatively complemented modular lattice and that
$K$ is a convex sublattice of $L$. Suppose that $L$ is either
conditionally $\aleph_0$-\mcont\ or
conditionally $\aleph_0$-\jcont.
Then $f$ is a \Vemb.

\end{itemize}

See also Proposition~\ref{P:exV-hom} and
Corollary~\ref{C:V(K)notembV(L)}.

\chapter[Lattice and ring dimension]%
{When is the lattice dimension
equivalent to the ring dimension?}\label{App:RLDim}

Let $R$ be a regular ring. We review here sufficient
conditions for the natural \Vhom\ from
\(\DD\mathcal{L}(R_R)\) onto \(V(R)\) to be an isomorphism.
Note that it is not always an isomorphism, by
Corollary~\ref{C:RegRingCounterex}.

\begin{itemize}
\item[(1)] (see Proposition~\ref{P:piRiso}) There are
principal right ideals $I$ and $J$ of $R$ such that
\(I\cap J=\{0\}\) and both $[I]$ and $[J]$ are order-units of
$V(R)$.

\item[(2)] (see Proposition~\ref{P:piiso(unreg)})
$R$ is unit-regular. Furthermore, in this case, if $I$ and
$J$ are principal right ideals of $R$, then $I\cong J$ if and
only if $I\sim J$.

\item[(3)] (see Corollary~\ref{C:caseRsimple}) $R$ is simple.

\item[(4)] (see Corollary~\ref{C:DirFinReg})
$R$ is a quotient ring of a regular ring that is
either $\aleph_0$-right-con\-tin\-u\-ous or
$\aleph_0$-left-con\-tin\-u\-ous.
Furthermore, in that case, if $I$ and
$J$ are principal right ideals of $R$, then $I\cong J$ if and
only if there are decompositions \(I=I_0\oplus I_1\) and
\(J=J_0\oplus J_1\) such that \(I_0\sim J_0\) and
\(I_1\sim J_1\).

\end{itemize}

\chapter[Review of normality]%
{A review of normality for relatively complemented modular
lattices}\label{App:Norm}

Recall that a sectionally complemented modular lattice is
\emph{normal}, if it satisfies the following sentence:
\[
(\forall x,y)\bigl((x\approx y\text{ and }x\wedge y=0)
\Rightarrow x\sim y\bigr),
\]
and a relatively complemented modular lattice is normal if
all its closed intervals, viewed as sublattices, are normal.

If $L$ is a relatively complemented modular lattice, then
it was already known that each of the following conditions is
sufficient for $L$ to be normal:

\begin{itemize}
\item[(1)] $L$ is \emph{coordinatizable}: that is, there exists
a regular ring $R$ such that \(L\cong\mathcal{L}(R_R)\), see
Lemma~\ref{L:SubmNorm}. This is the case, for example,
if $R$ admits a large partial $4$-frame or if $L$ is
Arguesian, and admits a large partial $3$-frame, see
\cite{Jons60}.

\item[(2)] $L$ is conditionally \mcont\ (J. von~Neumann and
I. Halperin \cite{HaNe40}).

\item[(3)] $L$ is conditionally $\aleph_0$-\mcont\ and
conditionally $\aleph_0$-\jcont\ (I. Halperin \cite{Halp38}).
\end{itemize}

The present work improves some of these classical results,
with the following list of sufficient conditions:

\begin{itemize}
\item[(4)] $L$ is simple (see Corollary~\ref{C:SimpleNorm}).
There is a similar statement about coordinatization, in
\cite[Corollary 8.5]{Jons60}, that requires dimension at
least~$4$.

\item[(5)] $L$ is locally finitely distributive
(see Corollary~\ref{C:nDistrNorm}). In fact, in this case,
perspectivity is transitive in $L$
(see Corollary \ref{C:nDistrTrans}).

\item[(6)] $L$ is conditionally $\aleph_0$-\mcont ---without any
additional join-continuity assumption (Theorem~\ref{T:MeetContNor}).
This improves conditions (2) and (3) above.

\item[(7)] $L$ is conditionally $\aleph_0$-\jcont\
(see Corollary~\ref{C:GenEmbVs}); see also
Corollary~\ref{C:NormEqDual1}).

\item[(8)] The dimension monoid \(\DD L\) satisfies the
following sentence:
\[
(\forall x,y,z)\bigl[x+z=y+z\Rightarrow(\exists t)
(2t\leq z\text{ and }x+t=y+t)\bigr].
\]
This result is contained in Theorem~\ref{T:AxImplNorm}. In
particular, if \(\DD L\) is cancellative, then $L$ is normal,
but, of course, there is in that case a stronger statement:
namely, \(\DD L\) is cancellative if and only if
perspectivity is transitive in $L$
(see Theorem~\ref{T:transperspchar}).
\end{itemize}

The class of all normal
relatively complemented modular lattices is
self-dual (see Corollary~\ref{C:NormEqDual1}) and closed under
convex sublattices (see Corollary~\ref{C:GenNorm}), direct limits,
reduced products and homomorphic images
(see Corollary~\ref{P:PresNorm}).

Both Corollary~\ref{C:nDistrNorm} and
Theorem~\ref{T:MeetContNor} admit a common strengthening as
follows: let $L$ be a sectionally complemented modular
lattice; we say that a sequence\linebreak
\(\vv<a_n\mid n\in\omega>\) of elements of $L$ is
\emph{tame}, if there exists a bounded
independent sequence
\(\vv<b_n\mid n\in\omega>\) of elements of $L$ such that
\(a_n\leq b_n\) and \(b_{n+k}\lesssim b_n\) hold for all $n$,
$k\in\omega$. Now, suppose that \(\bigvee_{n\in\omega}a_n\) exists,
for any tame sequence \(\vv<a_n\mid n\in\omega>\),
and that the equality \(x\wedge\bigvee_{n\in\omega}a_n=
\bigvee_{n\in\omega}\left(x\wedge\bigvee_{k<n}a_k\right)\)
holds for all \(x\in L\).
Then $L$ is normal. The proof of this is essentially
the same as that of Theorem~\ref{T:MeetContNor}, one just
has to be careful to verify that the relevant sequences are
tame. Note that if $L$ is locally finitely distributive, then
every tame sequence is eventually constant with value zero.

Finally, there exists a non-normal complemented modular
lattice, see Section~\ref{S:NonNorm}.
Moreover, this example is a modular ortholattice.

\chapter{Problems and comments}\label{App:Probl}

\begin{problem}
A \emph{dimension word} is an expression of the form
\[
\sum_{i<k}
\Dim(P_i(x_0,\ldots,x_{n-1}),Q_i(x_0,\ldots,x_{n-1})),
\]
where the $P_i$ and $Q_i$ are lattice polynomials. If $s$ and
$t$ are dimension words, is the problem whether $s=t$, resp.,
$s\leq t$, in every lattice (or, equivalently, in every free
lattice) \emph{decidable}?
\end{problem}

Note that the word problem in free lattices (which is known to
be decidable, see \cite{Whit41}) is a very particular case of
this ``dimension word problem''. Note also that the
corresponding problem for congruences is easily seen to be
equivalent to the word problem for finitely presented
lattices, and the latter is also known to be decidable, see,
for example, \cite{FJNa95,Mcki43}.

\begin{problem}
Let $K$ be a sublattice of a lattice $L$. We say that $L$ is a
\emph{dimension-preserving extension} of $K$, if the natural
map from \(\DD K\) to \(\DD L\) is an isomorphism. Does every
nontrivial lattice admit a proper dimension-preserving
extension?
\end{problem}

The corresponding problem for \emph{congruences} has an
affirmative answer, see \cite{GrWe}. However, the
extension considered is not dimension-preserving in general.

\begin{problem}
Is the dimension monoid \(\DD L\) of a lattice $L$ always a
refinement\index{monoid!refinement ---} monoid?
\end{problem}

By the results of Chapters~\ref{ModLatt} and \ref{ShortLatt},
any counterexample would have to be non-modular%
\index{lattice!modular (not necessarily complemented) ---}
with at least one infinite bounded chain.

\begin{problem}\label{Pb:V(mod)anyth}
Let $M$ be any
\crm\index{monoid!conical refinement ---}. Does there exist a
modular%
\index{lattice!modular (not necessarily complemented) ---}
lattice $L$ such that \(M\cong\DD L\)?
\end{problem}

Note that this problem is probably very difficult, because a
positive answer to it would, by Corollary~\ref{C:congquotV},
imply a positive solution to the Congruence Lattice
Problem\index{congruence!Congruence Lattice Problem}. Note
also that by the ordered vector space counterexample of
\cite{WehrB}, one cannot hope to solve this problem positively
by using, for example,
\emph{sectionally complemented lattices}: indeed, the map
$\Dim$ satisfies the identity
\(\Dim(x\wedge y)+\Dim(x\vee y)=\Dim(x)+\Dim(y)\); thus the
negative conclusion follows by \cite[Corollary 2.9]{WehrB}.
On the other hand, the answer to Problem \ref{Pb:V(mod)anyth}
is not even known for \emph{countable} \crm s%
\index{monoid!conical refinement ---} $M$.

\begin{problem}\label{Pb:finiteVLP}
Let $P$ be a finite antisymmetric QO-system%
\index{quasi-ordered (QO) -system}. Does there exist a finite
lattice $L$ such that
\(\DD L\cong\mathbf{E}(P)\)?
(We conjecture an affirmative answer.)
\end{problem}

Note that the analogue of this problem for
\emph{congruences}\index{congruence!lattice ---} has been
solved positively (see
\cite[Theorem II.3.17]{Grat}, where it is proved that every
finite distributive lattice is isomorphic to the
congruence\index{congruence!lattice} lattice of some finite
lattice).

\begin{problem}\label{Pb:trperdec}
Let $L$ be a sectionally complemented modular lattice. Does the
transitivity of the relation of perspectivity by decomposition%
\index{perspectivity by decomposition} in $L$ imply any
additional axiom (other than conicality or
refinement\index{refinement!property}) in
$\DD L$? (We conjecture a negative answer.)
\end{problem}

\begin{problem}\label{Pb:perdec2}
Let $L$ be a sectionally complemented modular lattice and let $a$
and $b$ be elements of $L$. If $a\simeq b$, do there exist elements
$x_0$, $x_1$, $y_0$, and $y_1$ of $L$ such that $a=x_0\oplus x_1$
and $b=y_0\oplus y_1$?
\end{problem}

\begin{problem}\label{Pb:NormLatt}
Let $L$ be a complemented modular lattice that either is
\emph{complete} (one does not assume meet-continuity) or has a
homogeneous basis with $3$ elements. Is $L$
normal\index{normal!lattice}? In the first case, is \(\DD L\) a GCA%
\index{generalized cardinal algebra (GCA)}?
\end{problem}

In case $L$ is complete, one can say more.
Indeed, by the results of Amemiya and Halperin, see
\cite[Corollary of Theorem 8.1]{AmHa59}, $L$ is the direct
product of three lattices $L_1$, $L_2$, and $L_3$, where $L_1$
is $\aleph_0$-\mcont%
\index{lattice!meet-continuous ---!$\aleph_0$-\mcont\ ---},
$L_2$ is $\aleph_0$-\jcont%
\index{lattice!join-continuous ---!$\aleph_0$-\jcont\ ---}
and $L_3$ has no nontrivial
homogeneous\index{homogeneous!sequence} sequence (these have
to be defined properly, see \cite[Definition 5.1]{AmHa59}).
Thus it suffices to solve the problem in case $L$
has no nontrivial homogeneous\index{homogeneous!sequence}
sequence. Of course, by Theorem~\ref{T:AxImplNorm} and
Proposition~\ref{P:basicCA}.(a), if \(\DD L\) is a GCA%
\index{generalized cardinal algebra (GCA)}, then
$L$ is normal\index{normal!lattice}.

We refer to Appendix~\ref{App:Norm} for a review of
sufficient conditions for normality.

\begin{problem}
If $L$ is a complemented modular lattice, investigate the set
\(\NEq(L)\) of all normal equivalences on $L$, ordered under
inclusion.
\end{problem}

\begin{problem}
Let $K$ be a complemented modular lattice and
let $\equiv$ be a normal equivalence on $K$. Does there exist
a complemented modular lattice $L$ such that $K$ is an ideal of
$L$ and the equivalence
\[
x\equiv y\Leftrightarrow x\approxeq y\text{ in }L
\]
holds for all $x$ and $y$ in $K$?
\end{problem}

\begin{problem}\label{Pb:kcpleteV(L)}
Let $\kappa$ be an
infinite cardinal and let $L$ be a sectionally complemented,
modular, $\kappa$-\mcont%
\index{lattice!meet-continuous ---!$\kappa$-\mcont\ ---}
lattice. Does every increasing bounded $\kappa$-sequence of
elements of \(\DD L\) admit a supremum?
\end{problem}

This is proved for \(\kappa=\aleph_0\) in
Theorem~\ref{T:DimComplGCA}. Moreover, it can be proved with
rather different methods ``for $\kappa=\infty$'',
that is, if $L$ is any conditionally \mcont%
\index{lattice!meet-continuous ---} relatively complemented
modular lattice, then every subset of
\(\DD L\) admits a supremum---in
this case, one can, in fact, elucidate completely the structure
of \(\DD L\).

\begin{problem}\label{Pb:Repres}
Let $M$ be a conical\index{monoid!conical ---} \cm\ which is, in
addition, a GCA%
\index{generalized cardinal algebra (GCA)}. Does there exist a
sectionally complemented, modular, conditionally
$\aleph_0$-\mcont%
\index{lattice!meet-continuous ---!$\aleph_0$-\mcont\ ---}
lattice $L$ such that
\(M\cong\DD L\)?
(We conjecture a negative answer.)
\end{problem}

By the results of \cite{WehrB}, there are \crm s%
\index{monoid!conical refinement ---} (and even positive cones
of partially ordered vector spaces with
interpolation\index{interpolation!property}) that do not
appear as \(\DD L\) for any relatively
complemented modular lattice $L$. The construction of these
counterexamples are ``free constructions'', so that one may
think that it would be sufficient to generalize these free
constructions to the countably infinite case. For example,
such a free construction has been used by the author to solve
Tarski's simple%
\index{cardinal algebra!simple ---} cardinal algebra problem
in \cite{Wehr96b}.

\begin{problem}\label{Pb:VGeomLatt} Study the dimension monoid
of \emph{geometric}%
\index{lattice!geometric ---} lattices.
\end{problem}

Although it is fairly easy to give a complete description of
the dimension monoid in the \emph{modular} case (also in the
\emph{finite} case), see Proposition~\ref{P:DimGeom}, the non
modular case seems to be considerably harder. One partial
result is, of course, that if $L$ is a simple%
\index{lattice!simple ---} non-modular
geometric\index{lattice!geometric ---} lattice, then
\(\DD L\cong\two\)
(see Corollary~\ref{C:SimpleGeom}), but we do not know,
for example,
whether the dimension monoid of an indecomposable non-modular
geometric\index{lattice!geometric ---} lattice $L$ is always a
\emph{semilattice} (thus isomorphic to the
congruence\index{congruence!semilattice} semilattice of $L$).

\begin{problem}\label{Pb:4Lcoord} Let $L$ be any sectionally
complemented modular lattice. Is
\(4L\) coordinatizable, that is, does there exist a
regular\index{ring!(von Neumann) regular ---} ring (not
necessarily unital) $R$ such that
\(L\cong\mathcal{L}(R_R)\)? If $L$ is
Arguesian\index{lattice!Arguesian ---}, is \(3L\)
coordinatizable? (We conjecture an affirmative answer.)
\end{problem}

\backmatter


\begin{theindex}

  \item additive (binary relation), \ii{26}
    \subitem countably ---, 105
  \item algebraic
    \subitem preordering, \ii{xiii}, 4, 16--18, 39, 43, 45--47
  \item Archimedean, \ii{17}, 19, 20, 86
  \item arclength, \ii{50}
  \item associative
    \subitem partial operation, \ii{25}
  \item $\leq_{\mathrm{alg}}$, \ii{xiii}
  \item $\asymp$, \ii{xiii}

  \indexspace

  \item BCF partially ordered set, \ii{49}
  \item $\mathbf{B}(L)$, \ii{12}

  \indexspace

  \item cancellable, \ii{96}
  \item cardinal algebra
    \subitem simple ---, 119
  \item cardinal algebra (CA), \ii{21}, 22, 23, 107
  \item caustic
    \subitem pair, viii, \ii{50}, 51
    \subitem path, \ii{50}, 51--53
  \item commutative
    \subitem partial operation, \ii{25}
  \item compact (element of a lattice), \ii{xiv}, 9, 54
  \item conditionally (P), \ii{58}
  \item congruence
    \subitem Congruence Lattice Problem, \ii{ix}, 117
    \subitem extension property, \ii{55}
    \subitem lattice, vii, ix, xiv, 9, 118
    \subitem lattice ---, 7--12, 67, 71, 93, 107, 118
    \subitem monoid ---, xiv, 1, 9
    \subitem semilattice, xiv, 14, 79, 119
  \item continuous
    \subitem $\aleph_0$-join- ---, \see {lattice}{xi}
    \subitem $\aleph_0$-right- ---, \see {ring}{xi}
    \subitem $\kappa$-upper- ---, \ii{109}
    \subitem $\lambda$-lower- ---, \ii{109}
    \subitem geometry, \ii{viii}, ix, x, xii, 65, 87, 108, 110
      \subsubitem indecomposable ---, viii, ix
    \subitem right- ---, \see {ring}{xi}
  \item $\ccon L$, \ii{xiv}, 1, 7
  \item $\mathfrak{C}_n$, \ii{54}
  \item $\C(I)$, \ii{7}
  \item $\con L$, \ii{xiv}
  \item $\aleph_0$-continuous geometry, 104

  \indexspace

  \item (D$'$0), (D$'$1), (D$'$2), \ii{63}
  \item (D0), (D1), (D2), \ii{1}
  \item diamond, 76--78, \ii{83}, 84, 88, 93--95
  \item dimension
    \subitem extension property, \ii{55}, 56
  \item dimension monoid
    \subitem of a lattice, vii, \ii{1}
    \subitem of a ring, ix, xv
  \item dimension range, x, \ii{1}, 37, 63, 66
  \item dimension semigroup of a \rps, \ii{26}, 29
  \item directly finite
    \subitem element in a monoid, 18, \ii{37}, 106--108
    \subitem ring, xi
  \item distributive
    \subitem $n$- ---, x, xi, 83, 84, 89, 93--95
    \subitem finitely ---, 83
    \subitem locally finitely ---, 83, 86--88
  \item $\diag L$, \ii{xiv}
  \item $\diam{m}{L}$, \ii{93}
  \item $\Dim(a,b)$, $\Dim_L(a,b)$, \ii{1}
  \item $\Dim^+(a,b)$, \ii{3}
  \item $\Dim(x)$, \ii{1}
  \item $\mathrm{Distr}(a,b,c)$, \ii{3}
  \item $\mathrm{Distr}^*(a,b,c)$, \ii{3}
  \item $\DD L$, $L$ lattice, \ii{1}
  \item $s\dnw=a$, \ii{25}
  \item $\dnw p$, $\dnw_Pp$, \ii{41}
  \item $\DD(S,\oplus,\sim)$ (for \rps s), \ii{26}
  \item $\DD(\theta)$ ($\theta$ congruence), \ii{7}

  \indexspace

  \item (E1), (E2), (E3), \ii{69}
  \item $\mathbf{E}(P)$, \ii{39}
  \item $\ee$, \ii{43}
  \item $\equiv_I$, \ii{xiv}
  \item $\ff$, \ii{44}
  \item $\fft$, \ii{43}

  \indexspace

  \item $\fin(M)$, \ii{19}
  \item $\mathbf{F}(P)$, \ii{41}
  \item $\Ft(P)$, \ii{41}
  \item $\FP(R)$, \ii{xv}
  \item $F_{\mathcal{V}}(X)$, \ii{80}

  \indexspace

  \item generalized cardinal algebra (GCA), xii, \ii{21}, 22, 106, 107, 
		118
    \subitem closure of a ---, 22
  \item geometric lattice, \see {lattice}{viii}
  \item group
    \subitem dimension ---, x, \ii{xiv}, 19--21, 86, 109
    \subitem Grothendieck ---, xv, 1
    \subitem interpolation ---, \ii{xiv}, 19, 23, 24, 109, 110
    \subitem lattice-ordered ---, 20, 21, 110
    \subitem simplicial ---, 21

  \indexspace

  \item homogeneous
    \subitem basis, x, 79
    \subitem sequence, xii, \ii{84}, 86, 93, 94, 99, 108, 118

  \indexspace

  \item ideal
    \subitem neutral --- (in lattices), \ii{67}, 71, 92--94, 96
    \subitem of a lattice, 74
    \subitem of a monoid, \ii{xiv}, 7, 9--12, 18, 19, 107, 108
    \subitem of a ring, x, 72, 80, 81, 103, 107, 108
    \subitem of a semilattice, \ii{xiv}, 67, 68, 74, 89, 91, 98, 110
  \item independent, x, 57, \ii{58}, 87, 88, 92, 97--105, 109, 110
  \item index
    \subitem of an element in a lattice, \ii{84}, 86, 87
    \subitem of an element in a monoid, 19, 21, 86, 87
    \subitem of nilpotence, \ii{84}
  \item interpolation
    \subitem $\kappa$- --- property, \ii{23}
    \subitem property, \ii{xiii}, xiv, 20, 118
  \item interval axiom, viii, \ii{46}, 47, 53, 56
  \item $\Id M$, \ii{xiv}
  \item $\Idc M$, \ii{xiv}
  \item $\ind_L(x)$, \ii{84}
  \item $\ind_M(x)$, \ii{19}
  \item $\infty a$, \ii{21}

  \indexspace

  \item $\vv<\kappa,\lambda>$-geometry, \ii{108}, 109
  \item $K_0(L)$, \ii{1}

  \indexspace

  \item lattice
    \subitem algebraic ---, \ii{xiv}, 54
    \subitem Arguesian ---, 77, 89, 119
    \subitem BCF ---, \see {BCF}{vii}, viii, 49, 52--54
    \subitem complemented ---, \ii{xiv}
    \subitem complete ---, viii, \ii{xiv}
    \subitem congruence splitting ---, 79
    \subitem continuous ---, 58
    \subitem geometric ---, viii, \ii{54}, 104, 105, 119
    \subitem join-continuous ---, viii, \ii{58}, 110
      \subsubitem $<\kappa$-\jcont\ ---, \ii{58}
      \subsubitem $\aleph_0$-\jcont\ ---, xi, 72, 99, 101, 103, 104, 
		107, 108, 118
      \subsubitem $\kappa$-\jcont\ ---, \ii{58}, 108
    \subitem meet-continuous ---, viii, \ii{58}, 72, 96, 108, 110, 118
      \subsubitem $<\kappa$-\mcont\ ---, \ii{58}
      \subsubitem $\aleph_0$-\mcont\ ---, xi, xii, 72, 96, 97, 98, 101, 
		103--108, 118
      \subsubitem $\kappa$-\mcont\ ---, \ii{58}, 61, 108, 109, 118
    \subitem modular (not necessarily complemented) ---, vii, viii, 3, 
		9, 10, 12, 26, 35--38, 53, 54, 57, 58, 67, 80, 83, 84, 
		88, 93, 94, 117
    \subitem normal ---, \see {normal}{xi}
    \subitem relatively complemented ---, \ii{xiv}
    \subitem sectionally complemented ---, \ii{xv}
    \subitem semimodular ---, 54
    \subitem simple ---, viii, 12, 38, 53, 54, 95, 119
    \subitem V-modular ---, \ii{11}, 12, 53, 54
    \subitem with commuting congruences, \ii{79}
  \item lower subset, \ii{xiii}
  \item $\mathcal{L}(R_R)$, $\mathcal{L}(M)$, \ii{x}, \ii{xv}
  \item $L^\mathrm{op}$, \ii{2}

  \indexspace

  \item monoid, \ii{xiii}, 1, 2, 63
    \subitem conical ---, \ii{xiii}, xiv, 1, 2, 11, 19, 53, 118
    \subitem conical refinement ---, ix, xv, 18--20, 36, 117, 118
    \subitem primitive ---, viii, \ii{39}, 40, 45--47, 50, 53, 80
    \subitem refinement ---, vii--ix, \ii{xiv}, 15--19, 26, 32, 40, 53, 
		79, 80, 96, 105, 106, 117
    \subitem simple ---, 18, 53
    \subitem stably finite ---, 18
    \subitem strong refinement ---, \ii{47}, 50
  \item monotone $\kappa$-completeness, \ii{23}, 24, 109, 110
  \item $\mathbf{M}(I;X,Y)$, \ii{40}
  \item $\mathrm{Mod}(a,b,c)$, \ii{3}
  \item $M/I$, \ii{xiv}

  \indexspace

  \item normal
    \subitem equivalence, x, \ii{69}, 70, 72--76, 80, 101--103, 105
    \subitem kernel, \ii{91}
    \subitem lattice, xi, xii, \ii{69}, 70, \ii{71}, 72, 74, 75, 87, 89, 
		92--96, 98--100, 103, 105, 107, 118
    \subitem non- --- lattice, 76, 78, 95
  \item $\NN$, \ii{xiv}
  \item $\NEq(L)$, \ii{72}
  \item $\nor L$, \ii{91}

  \indexspace

  \item ortholattice, 76, 77
  \item $\omega$, \ii{xiii}
  \item $[\omega]^{<\omega}$, \ii{27}
  \item $\oplus$
    \subitem on elements of a lattice, \ii{57}
    \subitem on intervals of a lattice, \ii{35}
    \subitem partial operation (general), \ii{25}

  \indexspace

  \item perspective
    \subitem isomorphism, \ii{60}, 61, 76, 78, 85, 99
  \item perspectivity, viii, x, xi, 57, \ii{59}, 60, 61, 64, 65, 68, 
		70--73, 76, 81, 85--89, 92--95, 99, 102, 104, 108--110
  \item perspectivity by decomposition, \ii{67}, 71, 103, 108, 118
  \item projective
    \subitem isomorphism, \ii{60}, 72, 98
  \item projectivity
    \subitem by decomposition, 63, \ii{65}, 67, 70, 75--77, 105, 108
    \subitem of elements, xi, \ii{59}, 70, 73, 78, 92, 95, 96, 99
    \subitem of intervals, \ii{1}, 36, 37, 77, 78
  \item proper
    \subitem axis of perspectivity, \ii{59}
  \item pseudo-cancellation, viii, 23, \ii{46}, 47, 53
  \item pseudo-indecomposable (PI), \ii{39}, 40, 53
  \item $a\parallel b$, \ii{50}
  \item $\Pi_X$, \ii{54}
  \item $\perp$ (for \rps s), \ii{27}
  \item $\simeq$, \ii{67}
  \item $\sim_2$, \ii{x}
  \item $\sim_k$, \ii{26}, \ii{59}
  \item $\sim$ (for elements), \ii{59}
  \item $\sim$ (for intervals), \ii{1}
  \item $\pi_M$, $\pi_R$, $\pi_{nR}$, \ii{79}
  \item $\PP(E)$, \ii{108}
  \item $\approxeq$, \ii{65}
  \item $\approxeq_\omega$, \ii{105}
  \item $\approx$ (for elements), \ii{59}
  \item $\approx$ (for intervals), \ii{1}
  \item $\propto$, \ii{xiii}

  \indexspace

  \item quasi-ordered (QO) -system, \ii{39}, 40, 43--45, 53, 54, 117

  \indexspace

  \item rectangular extension, \ii{55}
  \item rectifiable, \ii{50}, 51, 52
  \item refinement
    \subitem matrix, \ii{xiv}, 16--18, 21, 22, \ii{26}, 27, 28, 30--32, 
		35, 65, 66, 92, 97, 98
    \subitem postulate, \ii{21}
    \subitem property, \ii{xiv}, 16, 17, 22, 25, \ii{26}, 27, 28, 
		30--32, 35, 37, 47, 65, 118
    \subitem Schreier ---, vii, 35, 36
  \item refining (binary relation), \ii{26}
    \subitem left ---, \ii{26}
    \subitem right ---, \ii{26}
  \item ring
    \subitem $\aleph_0$-left-continuous ---, xii, 107, 108
    \subitem $\aleph_0$-right-continuous ---, xi, 103, 106--108
    \subitem (von Neumann) regular ---, ix--xii, xv, 1, 69, 72, 79, 80, 
		84, 89, 93, 103, 106, 107, 119
    \subitem right-continuous ---, xi
    \subitem simple ---, 80
    \subitem unit-regular ---, 81, 89, 108
  \item $\rect L$, \ii{55}
  \item $R{\mbox{\bf-Noeth}}$, \ii{32}

  \indexspace

  \item semigroup, vii, \ii{xiii}, xiv, 26, 29, 30, 32, 35, 42
    \subitem inverse ---, 32
    \subitem partial ---, vii, xiii, \ii{25}, 26, 27, 32, 35
    \subitem refined partial ---, 25, \ii{26}, 27, 30--32, 35
  \item separativity, viii, 18, 46
  \item subperspectivity, \ii{59}, 91--93
  \item $\mathbb{S}$, \ii{27}
  \item $\hat S$, \ii{29}
  \item $\Sigma p$, \ii{27}
  \item $a\sd b$, \ii{13}
  \item $\S(L)$, \ii{xiv}
  \item $\lesssim$, \ii{59}
  \item $\lesssim_k$, \ii{59}
  \item $S^\circ$, \ii{xiii}

  \indexspace

  \item $\Theta(x)$, \ii{67}
  \item $\to^\alpha$, \ii{27}
  \item $\searrow$, \ii{1}
  \item $\nearrow$, \ii{1}
  \item $\two$, \ii{viii}

  \indexspace

  \item unperforation, viii, \ii{xiv}, 46

  \indexspace

  \item \Vemb, 18, 67, 74, 80, 88, 103, 107, 109
  \item \Vhom, \ii{18}, 37, 67, 68, 75, 79, 88, 105
  \item \Vmeas, x, \ii{66}, 75
  \item von Neumann Coordinatization Theorem
    \subitem  J\'onsson improvement, xi, 93
    \subitem original form, ix
  \item $\DD L$, $L$ lattice, 7
  \item $V(R)$, $R$ ring, \ii{xv}

  \indexspace

  \item weak projectivity, \ii{8}
  \item $\wpr$, $\wprd$, $\wpru$, \ii{8}

  \indexspace

  \item $\ZZ^+$, \ii{xiv}
  \item $\ZZ^+[B]$, \ii{13}
  \item $\ZZb$, \ii{41}

\end{theindex}

\listoffigures


\begin{thebibliography}{99}

\bibitem{AmHa59}
I. Amemiya and I. Halperin,
\emph{Complemented modular lattices}, Canad. J. Math.
\textbf{11} (1959), 481--520.

\bibitem{Ara87}
P. Ara,
\emph{Aleph-nought-continuous regular rings},
J. Algebra \textbf{109} (1987), 115--126.

\bibitem{AGPO}
P. Ara, K. R. Goodearl, E. Pardo,
and K. C. O'Meara,
\emph{Separative cancellation for projective modules over exchange rings},
Israel J. Math. \textbf{105} (1998), 105--137.

\bibitem{Birk93}
G. Birkhoff,
\emph{Lattice Theory},
Corrected reprint of  the
1967 third edition. American Mathematical Society
Colloquium  Publications, Vol. \textbf{25}. American
Mathematical Society, Providence, R. I., 1979.
vi+418 pp.

\bibitem{Blac90}
B. Blackadar,
\emph{Rational C*-algebras and nonstable K-theory},
Rocky Mountain J. Math. \textbf{20}, no. 2 (1990),
285--316.

\bibitem{Broo}
G. Brookfield,
\emph{Direct sum cancellation of noetherian modules},
to appear in J. Algebra.

\bibitem{BrRo92}
G. Bruns and M. Roddy,
\emph{A finitely generated modular ortholattice},
Canad. Math. Bull.
\textbf{35}, no. 1 (1992), 29--33.

\bibitem{Busc90}
C. Busqu\'e,
\emph{Two-sided ideals in right
self-injective regular rings},
J. Pure Appl. Algebra
\textbf{67} (1990), 209-245.

\bibitem{CrDi73}
P. Crawley and R. P. Dilworth,
\emph{Algebraic Theory of Lattices},
Englewood Cliffs, New Jersey: Prentice-Hall, Inc. VI, 1973. 201 pp.

\bibitem{CrJo64}
P. Crawley and B. J\'onsson,
\emph{Refinements for infinite direct decompositions of algebraic
systems},
Pacific J. Math. \textbf{14} (1964), 797-855.

\bibitem{DaHW72}
A. Day, C. Herrmann, and R. Wille,
\emph{On modular lattices with four generators},
Algebra Universalis \textbf{2} (1972), 317--323.

\bibitem{Dobb82}
H. Dobbertin,
\emph{On Vaught's criterion for isomorphisms of
countable Boolean algebras},
Algebra Universalis \textbf{15} (1982), 95--114.

\bibitem{Dobb83}
\bysame,
\emph{Refinement monoids, Vaught
monoids, and Boolean algebras},
Math. Ann. \textbf{265} (1983), 475--487.

\bibitem{Dobb84}
\bysame,
\emph{Primely generated regular
refinement monoids},
J. Algebra \textbf{91} (1984), 166--175.

\bibitem{Fill65}
P. A. Fillmore,
\emph{The dimension theory of
certain cardinal algebras},
Trans. Amer. Math. Soc.
\textbf{117} (1965), 21--36.

\bibitem{Free76}
R. Freese,
\emph{Planar sublattices of $\mathrm{FM}(4)$},
Algebra Universalis \textbf{6} (1976), 69--72.

\bibitem{Free79}
\bysame,
\emph{Projective geometries as projective modular
lattices},
Trans. Amer. Math. Soc. \textbf{251} (1979), 329--342.

\bibitem{Free80}
\bysame,
\emph{Free modular lattices},
Trans. Amer. Math. Soc. \textbf{261} (1980), 81--91.

\bibitem{Free87}
\bysame,
\emph{A decomposition theorem for modular lattices
containing an $n$-diamond},
Acta Sci. Math. \textbf{51} (1987), 57--71.

\bibitem{FJNa95}
R. Freese, J. Je\v zek, and J. B. Nation,
\emph{Free Lattices},
Mathematical Surveys and Monographs, Vol. \textbf{42}.
American Mathematical Society, Providence, RI, 1995.
viii+293 pp.

\bibitem{GHKL80}
G. Gierz, K.~H. Hofmann, K. Keimel, J.~D.
Lawson, M. Mislove, and D. S. Scott,
\emph{A Compendium of Continuous Lattices},
Springer-Verlag, Berlin-New York, 1980. xx+371 pp.

\bibitem{Good82}
K.~R. Goodearl,
\emph{Directly finite
aleph-nought-continuous regular rings},
Pacific J. Math. \textbf{100} (1982), 105--122.

\bibitem{Good86}
\bysame,
\emph{Partially Ordered Abelian Groups with Interpolation}.
Mathematical Surveys and Monographs, Vol. \textbf{20}. American
Mathematical Society, Providence,  R.I., 1986. xxii+336 pp.

\bibitem{Good91}
\bysame,
\emph{Von Neumann Regular Rings},
Second edition.  Robert E. Krieger
Publishing Co., Inc., Malabar, FL, 1991.
xviii+412 pp.

\bibitem{Good94}
\bysame,
\emph{$K_0$ of regular rings with bounded index of nilpotence},
in: \emph{Abelian Group Theory and Related Topics}
(R. G\"obel, P. Hill and W. Liebert, eds.),
Contemporary Mathematics \textbf{171} (1994), 173--199.

\bibitem{Good95}
\bysame,
\emph{von~Neumann regular rings and direct sum decomposition problems},
in: \emph{Abelian Groups and Modules, Padova 1994}
(A. Facchini and C. Menini, eds.),
Dordrecht (1995) Kluwer, 249--255.

\bibitem{Grat}
G. Gr\"atzer,
\emph{General Lattice Theory. Second Edition},
Birkh\"auser Verlag, Basel. 1998. xix+663 pp.

\bibitem{GrSc0}
G. Gr\"atzer and E. T. Schmidt,
\emph{Congruence lattices of lattices}, Appendix C in G. Gr\"atzer:
\emph{General Lattice Theory. Second Edition}, Birkh\"auser Verlag,
Basel. 1998. xix+663 pp.

\bibitem{GrSc}
\bysame,
\emph{Congruence-preserving extensions of finite lattices to
sectionally complemented lattices},
to appear in Proc. Amer. Math. Soc.

\bibitem{GrWe}
G. Gr\"atzer and F. Wehrung,
\emph{Proper congruence-preserving extensions of lattices},
to appear in Acta Math. Hungar.; see
\emph{Abstracts Amer. Math. Soc.} no. 97T-06-189.

\bibitem{Gril70}
P.~A. Grillet,
\emph{Interpolation properties and tensor product of semigroups},
Semigroup Forum \textbf{1} (1970), 162--168.

\bibitem{Halp38}
I. Halperin,
\emph{On the transitivity of
perspectivity in continuous geometries},
Trans. Amer. Math. Soc. \textbf{44} (1938), 537--562.

\bibitem{Halp39a}
\bysame,
\emph{Dimensionality in reducible geometries},
Ann. of Math., II, \textbf{40} (1939), 581--599.

\bibitem{Halp39b}
\bysame,
\emph{Additivity and continuity of perspectivity},
Duke Math. J. \textbf{5}
(1939), 503--511.

\bibitem{Halp61}
\bysame,
\emph{Complemented modular lattices},
Proc. Sympos. Pure Math.
\textbf{2} (1961), 51--64.

\bibitem{HaNe40}
I. Halperin and J. von~Neumann,
\emph{On the transitivity of perspective mappings},
Ann. of Math., II, \textbf{41} (1940), 87--93.

\bibitem{Herr81}
C. Herrmann,
\emph{A finitely generated modular
ortholattice},
Canad. Math. Bull. \textbf{24}, no. 2
(1981), 241--243.

\bibitem{Herr84}
\bysame,
\emph{On the arithmetic of projective coordinate systems},
Trans. Amer. Math. Soc. \textbf{284},
(1984), 759--785.

\bibitem{HeHu74}
C. Herrmann and A.~P. Huhn,
\emph{Lattices
of normal subgroups which are generated by frames},
Colloq. Math. Soc. J\'anos Bolyai
\textbf{14}, Lattice Theory, Szeged (Hungary) (1974),
97--136.

\bibitem{Huhn72}
A.~P. Huhn,
\emph{Schwach distributive Verb\"ande. I},
Acta Sci. Math. \textbf{33} (1972), 297--305.

\bibitem{Huhn74}
\bysame
\emph{Two notes on $n$-distributive lattices},
Colloq. Math. Soc. J\'anos Bolyai
\textbf{14}, Lattice Theory, Szeged (Hungary) (1974), 137--147.

\bibitem{Huhn89a}
\bysame
\emph{On the representation of distributive
algebraic lattices. II}, Acta Sci. Math.
\textbf{53}, no. 1-2 (1989), 3--10.

\bibitem{Huhn89b}
\bysame
\emph{On the representation of distributive
algebraic lattices. III}, Acta Sci. Math.
\textbf{53}, no. 1-2 (1989), 11--18.

\bibitem{Iwam44}
T. Iwamura,
\emph{On continuous geometries I},
Japan. J. Math.
\textbf{19} (1944), 57--71.

\bibitem{JiRo92}
P. Jipsen and H. Rose,
\emph{Varieties of lattices}, Lecture Notes in Mathematics,
Vol. \textbf{1533}.  Springer-Verlag, Berlin, 1992. x+162 pp.

\bibitem{Jons60}
B. J\'onsson,
\emph{Representations of complemented modular lattices},
Trans. Amer. Math. Soc. \textbf{97} (1960), 64--97.

\bibitem{Kalm83}
G. Kalmbach,
\emph{Orthomodular Lattices},
London Mathematical Society Monographs, Vol. \textbf{18},
Academic Press, Inc., London-New York, 1983. viii+390 pp.

\bibitem{Loom55}
L. H. Loomis,
\emph{The lattice-theoretic background of the dimension theory of
operator algebras},
Mem. Amer. Math. Soc. \textbf{18} (1955), 35 pp.

\bibitem{FMae58}
F. Maeda,
\emph{Kontinuierliche Geometrien} (German), die Grundlehren der
mathematischen  Wissenschaften in Einzeldarstellungen mit
besonderer Ber\"ucksichtigung  der Anwendungsgebiete, Bd.
\textbf{95}. Springer-Verlag,  Berlin-G\"ottingen-Heidelberg, 1958.
x+244 pp. (translated from Japanese by S. Crampe, G. Pickert and R.
Schauffler).

\bibitem{SMae55}
S. Maeda,
\emph{Dimension Functions on Certain
General Lattices},
J. Sci. Hiroshima Univ., Ser. A \textbf{19} (1955),
211--237.

\bibitem{Mcki43}
J.~C.~C. McKinsey,
\emph{The decision problem for
some classes of sentences without quantifiers},
J. Symbolic Logic \textbf{8} (1943), 61--76.

\bibitem{Neum60}
J. von~Neumann,
\emph{Continuous geometry}, foreword by Israel
Halperin. Princeton Mathematical Series, Vol. \textbf{25}
Princeton University  Press, Princeton, N. J. 1960
xi+299 pp.

\bibitem{Pier89}
R.~S. Pierce,
\emph{Countable Boolean Algebras},
in \emph{Handbook of Boolean Algebras}, edited by J. D. Monk
with R. Bonnet, Elsevier, 1989, 775--876.

\bibitem{PTWe}
M. Plo\v s\v cica, J. T\r uma, and F. Wehrung,
\emph{Congruence lattices of free lattices
in non-distributive varieties},
Colloq. Math. \textbf{76},
no. 2 (1998), 269--278.

\bibitem{RaSh92}
K.~P.~S. Bhaskara Rao and R.~M. Shortt,
\emph{Weak cardinal algebras},
Ann. New York Acad. Sci.,
\textbf{659} (1992), 156--162.

\bibitem{ShWe94}
R.~M. Shortt and F. Wehrung,
\emph{Common extensions of semigroup-valued char\-ges},
J. Math. Anal. Appl.
\textbf{187}, no. 1 (1994), 235--258.

\bibitem{Tars49}
A. Tarski,
\emph{Cardinal Algebras},
with an Appendix: Cardinal Products
of Isomorphism  Types, by Bjarni J\'onsson and Alfred Tarski.
Oxford University Press, New  York, N. Y., 1949. xii+326 pp.

\bibitem{TuWe}
J. T\r uma and F. Wehrung,
\emph{Simultaneous representations of semilattices
by lattices with commuting congruences}, preprint.

\bibitem{Wago85}
S. Wagon,
\emph{The Banach-Tarski paradox}, with a foreword by Jan
Mycielski.  Encyclopedia of Mathematics and its Applications,
\textbf{24}. Cambridge  University Press, Cambridge-New York,
1985. xvi+251 pp.

\bibitem{Wehr90}
F. Wehrung,
\emph{Th\'eor\`eme de Hahn-Banach et
paradoxes continus et discrets},
C. R. Acad. Sci. Paris S\'er. I Math., t. \textbf{310} (1990),
303--306.

\bibitem{Wehr92a}
\bysame,
\emph{Injective positively ordered
monoids I},
J. Pure Appl. Algebra
\textbf{83} (1992), 43--82.

\bibitem{Wehr92b}
\bysame,
\emph{Injective positively ordered
monoids II},
J. Pure Appl. Algebra
\textbf{83} (1992), 83--100.

\bibitem{Wehr94}
\bysame,
\emph{Restricted injectivity,
transfer property and decompositions of separative positively
ordered monoids},
Comm. Algebra \textbf{22}, no. 5 (1994), 1747--1781.

\bibitem{Wehr96a}
\bysame,
\emph{Monoids of intervals of
ordered abelian groups},
J. Algebra \textbf{182}
(1996), 287--328.

\bibitem{Wehr96b}
\bysame,
\emph{Monotone $\sigma$-complete groups with
unbounded refinement},
Fund. Math. \textbf{151} (1996), 177--187.

\bibitem{WehrA}
\bysame,
\emph{Embedding simple commutative
monoids into simple refinement monoids},
Semigroup Forum
\textbf{56} (1998), 104--129.

\bibitem{WehrB}
\bysame,
\emph{Non-measurability  properties of
interpolation vector spaces},
Israel J. Math. \textbf{103} (1998),
177--206.

\bibitem{WehrC}
\bysame,
\emph{A uniform refinement property for
congruence lat\-tices},
Proc. Amer. Math. Soc. \textbf{127}, no. 2
(1999), 363--370.

\bibitem{Whit41}
P.~M. Whitman,
\emph{Free lattices},
Ann. of Math., II, \textbf{42} (1941),
325--330.
\end{thebibliography}
\end{document}